\documentclass[11pt]{book}
\usepackage{geometry}                		
\usepackage{graphicx}									
\usepackage{amssymb}\vspace{0.25cm}
\usepackage{mathtools}
\usepackage{setspace}
\usepackage{tikz-cd}
\usepackage[mathscr]{euscript}
\newcommand{\abs}[1]{\left\vert#1\right\vert}
\usepackage[normalem]{ulem}
\usepackage{manfnt}
\usepackage{pgfplots}
\usepackage{pifont}
\usepackage[safe]{tipa}
\usepackage{marvosym}
\usepackage{scalerel,stackengine}
\usepackage{titlesec}
\usepackage{makeidx}
\usepackage{centernot}

\usepackage{xspace}
\usepackage[noauto]{chappg}
\usepackage[OT2,T1]{fontenc}
\usepackage{hyperref}
\usepackage{amsmath}
\usepackage[amsmath,thmmarks]{ntheorem}

\usepackage{tabularx}
\newcolumntype{Y}{>{\raggedright\arraybackslash}X}

\newcommand\mA{%
$A$\xspace
}

\newcommand\mE{%
$E$\xspace
}
\newcommand\mF{%
$F$\xspace
}
\newcommand\mG{%
$G$\xspace
}
\newcommand\mH{%
$H$\xspace
}
\newcommand\mI{%
$I$\xspace
}

\newcommand\mS{%
$S$\xspace
}

\newcommand\mU{%
$U$\xspace
}
\newcommand\mV{%
$V$\xspace
}

\newcommand\mX{%
$X$\xspace
}
\newcommand\mY{%
$Y$\xspace
}

\newcommand\A{%
\mathbb{A}
}
\newcommand\C{%
\mathbb{C}
}
\newcommand\E{%
\mathbb{E}
}
\newcommand\F{%
\mathbb{F}
}

\newcommand\I{%
\mathbb{I}
}
\newcommand\K{%
\mathbb{K}
}
\newcommand\LL{%
\mathbb{L}
}
\newcommand\M{%
\mathbb{M}
}
\newcommand\N{%
\mathbb{N}
}

\newcommand\Q{%
\mathbb{Q}
}
\newcommand\R{%
\mathbb{R}
}
\newcommand\T{%
\mathbb{T}
}

\newcommand\Z{%
\mathbb{Z}
}

\newcommand\bINV{%
\textbf{INV}\xspace
}

\newcommand\Aut{%
\text{Aut}\hspace{0.05 cm}
}

\newcommand\sA{%
\mathcal{A}
}
\newcommand\sB{%
\mathcal{B}
}
\newcommand\sC{%
\mathcal{C}
}
\newcommand\sF{%
\mathcal{F}
}
\newcommand\sS{%
\mathcal{S}
}
\newcommand\sU{%
\mathcal{U}
}
\newcommand\cl{%
\text{c$\ell$}
}
\newcommand\ab{%
\text{ab}
}
\newcommand\ur{%
\text{ur}
}
\newcommand\sep{%
\text{sep}
}
\newcommand\Img{%
\text{Im}\hspace{0.05cm}
}
\newcommand\Inv{%
\text{Inv}\hspace{0.05cm}
}
\newcommand\Gal{%
\text{Gal}
}
\newcommand\rec{%
\text{rec}
}
\newcommand\ord{%
\text{ord}
}
\newcommand\pr{%
\text{pr}
}
\newcommand\res{%
\text{res}
}
\newcommand\ind{%
\text{ind}
}
\newcommand\mods{%
\text{mod}
}
\newcommand\id{%
\text{id}
}
\newcommand\itr{%
\text{int}\hspace{0.05cm}
}
\newcommand\nth{%
\text{th}
}
\newcommand\tJ{%
\text{J}
}
\newcommand\tN{%
\text{N}
}
\newcommand\tT{%
\text{T}
}

\newcommand\GL{%
\text{GL}
}
\newcommand\SL{%
\text{SL}
}
\newcommand\sym{%
\text{sym}
}
\newcommand\trdeg{%
\text{trdeg}
}
\newcommand\Tee{%
\mathsf{T}
}

\newcommand\ra{%
\rightarrow
}
\newcommand\lra{%
\longrightarrow
}
\newcommand\lla{%
\longleftarrow
}

\newcommand\Lra{%
\Leftrightarrow
}
\newcommand\ds{%
\displaystyle
}

\newcommand\card{%
\text{card}\hspace{0.05cm}
}
\newcommand\con{%
\text{con}\hspace{0.05cm}
}
\newcommand\End{%
\text{End}\hspace{0.05cm}
}
\newcommand\fin{%
\text{fin}\hspace{0.05cm}
}
\newcommand\modxs{%
\text{mod}
}
\newcommand\modx{%
\text{mod}\hspace{0.05cm}
}
\newcommand\vol{%
\text{vol}
}
\newcommand\sgn{%
\text{sgn}\hspace{0.05cm}
}
\newcommand\spp{%
\text{sp}
}
\newcommand\ev{%
\text{ev}
}
\newcommand\tr{%
\text{tr}\hspace{0.05cm}
}
\newcommand\trs{%
\text{tr}
}
\newcommand\ins{%
\text{in}
}

\newcommand\un[1]{%
\underline{#1}\xspace
}
\newcommand\ov[1]{%
\overline{#1}
}
\newcommand\acdot{%
\abs{\hspace{0.05cm} \cdot \hspace{0.05cm}}
}
\newcommand\ncdot{%
\norm{\hspace{0.05cm} \cdot \hspace{0.05cm}}
}
\newcommand\restr[2]{%
{#1}|{#2}
}
\newcommand\compl[1]{%
\text{$\ov{#1}$}
}

\newcommand\vsx{%
\vphantom{\int_\int^\int}
}
\newcommand\vsy{%
\vphantom{\int}
}
\newcommand\hsx{%
\hspace{0.05cm}
}

\newcommand\boxtimesdmc{%
\raisebox{-.075cm}{$\boxtimes$}
}

\newcommand{\norm}[1]{\left\lVert #1 \right\rVert}

\usepackage{scalerel,stackengine}
\stackMath
\newcommand\reallywidehat[1]{%
\savestack{\tmpbox}{\stretchto{%
  \scaleto{%
    \scalerel*[\widthof{\ensuremath{#1}}]{\kern-.6pt\bigwedge\kern-.6pt}%
    {\rule[-\textheight/2]{1ex}{\textheight}}
  }{\textheight}%
}{0.5ex}}%
\stackon[1pt]{#1}{\tmpbox}%
}
\makeatletter
\DeclareRobustCommand\widecheck[1]{{\mathpalette\@widecheck{#1}}}
\def\@widecheck#1#2{%
    \setbox\z@\hbox{\m@th$#1#2$}%
    \setbox\tw@\hbox{\m@th$#1%
       \widehat{%
          \vrule\@width\z@\@height\ht\z@
          \vrule\@height\z@\@width\wd\z@}$}%
    \dp\tw@-\ht\z@
    \@tempdima\ht\z@ \advance\@tempdima2\ht\tw@ \divide\@tempdima\thr@@
    \setbox\tw@\hbox{%
       \raise\@tempdima\hbox{\scalebox{1}[-1]{\lower\@tempdima\box
\tw@}}}%
    {\ooalign{\box\tw@ \cr \box\z@}}}
\makeatother

\newtheoremstyle{xx}
  {4pt}
  {0pt}
  {\upshape}
  {\bfseries}
  {}
  { }
  {}
  
\makeatletter 
 \newtheoremstyle{myu}%
  {\upshape\item[ \indent\indent\bf\underline{\theorem@headerfont ##2:}]}%
\makeatother
\makeatletter 
 \newtheoremstyle{myn}%
  {\item[\hskip\labelsep \ \bf ##1 \theorem@headerfont ##2.]}%
\makeatother
\theoremstyle{myn}
\newtheorem{theoremn}{Theorem} 
\theoremstyle{myu}
{\upshape}
\newtheorem{x}[theoremn]{}

\title{\textbf{Local and Global Analysis}}
\author{Garth Warner\\
Department of Mathematics\\
University of Washington}
\date{}									

\titleformat{\chapter}[display]
{\normalfont\filcenter\huge\bfseries}{}{0pt}{\large}

\titleformat{\chapter}[display]
{\normalfont\filcenter\huge\bfseries}{}{0pt}{\large}

\usepackage[OT2,OT1]{fontenc} 
\newcommand\cyr
{
\renewcommand\rmdefault{wncyr} 
\renewcommand\sfdefault{wncyss} 
\renewcommand\encodingdefault{OT2} 
\normalfont
\selectfont
}

\DeclareTextFontCommand{\textcyr}{\cyr}

\makeindex 
\begin{document}

\maketitle                              

\titlespacing*{\chapter}{0pt}{-50pt}{40pt}
\setlength{\parskip}{0.1em}
\[
\text{\large\bf ACKNOWLEDGEMENTS}
\]
\setlength\parindent{2em}

Many thanks to Judith Clare Salzer for typing the manuscript on an IBM Selectric.

Recently David Clark converted the typwritten manuscript to AMS-TeX.  This was a monumental task and in so doing he made a number of constructive and useful suggestions which serve to enhance the exposition.  His careful scrutiny of the manuscript has been invaluable.

\vspace{2 cm}

\[
\text{\large\bf DEDICATION}
\]

This article is dedicated to the memory of Paul Sally.

\begingroup
\fontsize{8pt}{8pt}\selectfont

\[
\text{\large\bf CONTENTS}
\]
 
 \vspace{0.3cm}
 
 \normalsize

\qquad $\S1.\ $ \qquad ABSOLUTE VALUES\\

\qquad $\S2.\ $ \qquad TOPOLOGICAL FIELDS\\

\qquad $\S3.\ $ \qquad COMPLETIONS\\

\qquad $\S4.\ $ \qquad p-ADIC STRUCTURE THEORY\\

\qquad $\S5.\ $ \qquad LOCAL FIELDS\\

\qquad $\S6.\ $ \qquad HAAR MEASURE\\

\qquad $\S7.\ $ \qquad HARMONIC ANALYSIS\\

\qquad $\S8.\ $ \qquad ADDITIVE p-ADIC CHARACTER THEORY\\

\qquad $\S9.\ $ \qquad MULTIPLICATIVE p-ADIC CHARACTER THEORY\\

\qquad $\S10.$ \qquad TEST FUNCTIONS\\

\qquad $\S11.$ \qquad LOCAL ZETA FUNCTIONS: $\R^\times$ OR $\C^\times$\\

\qquad $\S12.$ \qquad LOCAL ZETA FUNCTIONS: $\Q^\times_p$\\

\qquad $\S13.$ \qquad RESTRICTED PRODUCTS\\

\qquad $\S14.$ \qquad ADELES AND IDELES\\

\qquad $\S15.$ \qquad GLOBAL ANALYSIS\\

\qquad $\S16.$ \qquad FUNCTIONAL EQUATIONS\\

\qquad $\S17.$ \qquad GLOBAL ZETA FUNCTIONS\\

\qquad $\S18.$ \qquad LOCAL ZETA FUNCTIONS (BIS)\\

\qquad $\S19.$ \qquad L-FUNCTIONS\\

\qquad $\S20.$ \qquad FINITE CLASS FIELD THEORY\\

\qquad $\S21.$ \qquad LOCAL CLASS FIELD THEORY\\

\qquad $\S22.$ \qquad WEIL GROUPS: THE ARCHIMEDEAN CASE\\

\qquad $\S23.$ \qquad WEIL GROUPS: THE NON-ARCHIMEDEAN CASE\\

\qquad $\S24.$ \qquad THE WEIL-DELIGNE GROUP\\

\qquad \qquad APPENDIX A: \  TOPICS IN TOPOLOGY\\

\qquad \qquad APPENDIX B: \  TOPICS IN ALGEBRA\\

\qquad \qquad APPENDIX C: \  TOPICS IN GALOIS THEORY\\

\qquad \qquad REFERENCES\\

\endgroup 

\chapter*{PREFACE}
\setlength\parindent{2em}

\ \indent The objective of this article is to give an introduction to p-adic analysis along the lines of Tate's thesis, as well as incorporating material of a more recent vintage, for example Weil groups.

\pagenumbering{bychapter}
\setcounter{chapter}{0}
\chapter{
$\boldsymbol{\S}$\textbf{1}.\quad  ABSOLUTE VALUES}
\setlength\parindent{2em}
\setcounter{theoremn}{0}

\begin{x}{\small\bf DEFINITION} \ 
Let $\F$ be a field $-$then an 
\un{absolute value}
\index{absolute value} 
(a.k.a. a valuation of order 1) is a function
\[
\acdot : \F \ra \R_{\ge 0}
\]
satisfying the following conditions.\\

\indent \un{AV-1} \quad		$\abs{a}  = 0 \Lra a = 0$.\\
\indent \un{AV-2} \quad 		$\abs{ab} = \abs{a} \abs{b}$.\\
\indent \un{AV-3} \quad		$\exists \ M > 0$:\\
\[
\abs{a+b} \le M \sup(\abs{a}, \abs{b}).
\]
\end{x}
\vspace{0.1cm}

\begin{x}{\small\bf EXAMPLE} \ 
Let $\F = \R$ or $\C$ with the usual absolute value $\acdot_\infty$ $-$then one can take $M = 2$.
\end{x} 
\vspace{0.1cm}

\begin{x}{\small\bf DEFINITION} \ 
The 
\un{trivial absolute value}
\index{trivial absolute value} 
is defined by the rule
\[
\abs{a} = 1 \quad \text{$\forall$ a $\ne 0.$}
\]
\end{x}
\vspace{0.1cm}

\begin{x}{\small\bf LEMMA} \ 
If $\acdot$  is an absolute value, then
\[
\abs{1} = 1.
\]
\end{x}
\vspace{0.1cm}

\begin{x}{\small\bf APPLICATION} \ 
If $a^n$ = 1, then
\[
\abs{a^n} = \abs{a}^n = \abs{1} = 1
\]
\[
\implies \abs{a} = 1.
\]
\end{x}
\vspace{0.1cm}

\begin{x}{\small\bf RAPPEL} \ 
Let \mG be a cyclic group of order $r < \infty$ $-$then the order of any subgroup of \mG is a divisor of $r$ and if $n \mid r$, then \mG possesses one and only one 
subgroup of order $n$ (and this subgroup is cyclic).
\end{x}
\vspace{0.1cm}

\begin{x}{\small\bf RAPPEL} \ 
Let \mG be a cyclic group of order $r < \infty -$then the 
\un{order}
\index{order} 
of $x$ $\in$ \mG is, by definition, 
$\#\langle x \rangle$, the latter being the smallest positive integer $n$ such that $x^n$ = 1.
\end{x}
\vspace{0.1cm}
 
\begin{x}{\small\bf SCHOLIUM} \ 
Every absolute value on a finite field $\F_q$ is trivial.

\vspace{0.1cm}

[In fact, $\F_q^\times$  is cyclic of order $q-1$.]
\end{x}
\vspace{0.1cm}

\begin{x}{\small\bf DEFINITION} \ 
Two absolute values $\acdot_1$, and $\acdot_2$ on a field $\F$ are 
\un{equivalent} 
\index{absolute values equivalence } 
if $\exists$  $r > 0$:
\[
\acdot_2 = {\acdot_1}^r.
\]

Note: Equivalence is an equivalence relation.]
\end{x}
\vspace{0.1cm}

\begin{x}{\small\bf \un{N.B.}} \ 
If $\acdot$ is an absolute value, then so is ${\acdot}^r  \ (r > 0)$, 
the $M$ per $\acdot$ being $M^r$ per ${\acdot}^r$ .
\end{x}
\vspace{0.1cm}

\begin{x}{\small\bf LEMMA} \ 
Every absolute value is equivalent to one with $M \le 2$.

\vspace{0.1cm}

PROOF \  Assume from the beginning that $M > 2$, hence
\[
M^r \le 2 			\quad \text{$(r > 0)$}
\]
if 
\[
r \log M \le \log 2
\]
or still, if
\[
r  \le \frac{\log 2}{\log M}			\quad (< 1).
\]
\end{x}
\vspace{0.1cm}

\begin{x}{\small\bf DEFINITION} \ 
An absolute value $\acdot$ satisfies the 
\un{triangle inequality}
\index{triangle inequality} 
if
\[
\abs{a + b} \le \abs{a} + \abs{b}.
\]
\end{x}

\vspace{0.1cm}

\begin{x}{\small\bf LEMMA} \ 
Suppse given a function $\acdot : \F \ra \mathbb{R}_{\ge 0}$ satisfying AV-1 and AV-2, 
$-$then AV-3 holds with $M \le 2$  iff the triangle inequality obtains.

PROOF \   Obviously, if 
\[
\abs{a + b} \le \abs{a} + \abs{b},
\]
then
 \[
 \abs{a + b} \le 2 \sup( \abs{a},  \abs{b}).
\]
In the other direction, by induction on $m$, 
\[
\bigl|\sum_{k=1}^{2^m} a_{k} \bigr | \le 2^m \sup\limits_{1 \leq k \leq 2^m} \abs{a_k}.
\]
Next, given $n$ choose $m$: $2^m$ $\ge n > 2^{m-1}$, so upon inserting $2^m - n$ zero summands,
\[
\allowdisplaybreaks
\begin{aligned}
\allowdisplaybreaks
\bigl|\sum_{k=1}^{n} a_{k} \bigr | 
&\le M \sup \bigl(\bigl|\sum_{k=1}^{2^{m-1}} a_{k}\bigr | , \bigl|\sum_{k=2^{m-1}+1}^{2^m} a_{k} \bigr | \bigr )\\
& \le 2 \sup \bigl(\bigl|\sum_{k=1}^{2^{m-1}} a_{k} \bigr | ,\bigl|\sum_{k=2^{m+1} + 1}^{2^{m-1} + 2^{m-1}} a_{k} \bigr | \bigr)\\
& \le 2 \sup \bigl(2^{m-1} \sup_{k \le 2^{m-1}}\abs{a_k},    2^{m-1} \sup_{k > 2^{m-1}}\abs{a_k}\bigr) \\
& \le 2 \cdot 2^{m-1} \sup_{1 \le k \le n}\abs{a_k} \\
& \le 2 \cdot n \cdot \sup_{1 \le k \le n}\abs{a_k}.
\end{aligned}
\]

\allowdisplaybreaks
I.e.
\allowdisplaybreaks
\[
\begin{aligned}
\bigl|\sum_{k=1}^{n} a_{k} \bigr |  \ 
&\le \ 2n \sup_{1 \le k \le n}\abs{a_k} \\
&\le \ 2n \sum_{k=1}^{n} \abs{a_{k}}. 
\end{aligned}
\]
In particular, 
\[
\bigl|\sum_{k=1}^{n} 1 \bigr | \ = \abs{n} \  \le  \ 2n.
\]
Finally, 
\[
\begin{aligned}
\allowdisplaybreaks
\abs{a + b}^n 
&= \  \abs{(a + b)^n} \quad \text{(AV-2)}\\
&=  \ \bigl|\sum_{k=0}^{n} {n \choose k}  a^kb^{n-k}\bigr |\\
&\leq \  2(n+1) \sum_{k=0}^{n} \bigl|{n \choose k}a^kb^{n-k}\bigr | \\
&\leq \  2(n+1) \sum_{k=0}^{n} \bigl|{n \choose k}\bigr | \bigl|a^kb^{n-k}\bigr |    \quad \text{(AV-2)}\\
&\le  \ 2(n+1)2 \sum_{k=0}^{n} {n \choose k} \bigl|a^kb^{n-k}\bigr | \\
&=  \ 4(n+1) (\abs{a} + \abs{b})^n
\end{aligned}
\]
$\implies$\\
\[
\begin{aligned}
\allowdisplaybreaks
\abs{a + b} 
&\le 4^{1/n}(n + 1)^{1/n}(\abs{a} + \abs{b})\\
&\ra (\abs{a} + \abs{b})  \quad \text{ $(n \ra \infty)$}.
\end{aligned}
\]
\end{x}
\vspace{0.1cm}

\begin{x}{\small\bf SCHOLIUM} \ 
Every absolute value is equivalent to one that satisfies the triangle inequality.
\end{x}
\vspace{0.1cm}

\begin{x}{\small\bf DEFINITION} \ 
A 
\un{place}
\index{place} 
of $\F$ is an equivalence class of nontrivial absolute values.
\end{x}
\vspace{0.1cm}

Accordingly, every place admits a representative for which the triangle inequality is in force. \\

\begin{x}{\small\bf DEFINITION} \ 
An absolute value $\acdot$  is 
\un{non-archimedean}
\index{non-archimedean} 
if it satisfies the 
\index{ultrametric inequality}
\un{ultrametric inequality}:
\[
\abs{a + b} \le \sup(\abs{a},\abs{b})	\quad (\text{so } M = 1).
\]
\end{x}
\vspace{0.1cm}

\begin{x}{\small\bf \un{N.B.}} \ 
A non-archimedean absolute value satisfies the triangle inequality.
\end{x}
\vspace{0.1cm}

\begin{x}{\small\bf LEMMA} \ 
Suppose that $\acdot$  is non-archimedean and let $\abs{b} < \abs{a} -$then
\[
\abs{a + b} = \abs{a}.
\]
\quad PROOF \ 
\[
\begin{aligned}
\abs{a} = \abs{(a + b) - b} \ 
&\le \ \sup(\abs{a + b}, \abs{b}) \\
&=\  \abs{a + b}
\end{aligned}
\]
since $\abs{a} \le \abs{b}$ is untenable.  Meanwhile
\[
\abs{a + b} \ \le \  \sup(\abs{a}, \abs{b}) \ = \  \abs{a}.
\]
\end{x}
\vspace{0.1cm}

\begin{x}{\small\bf EXAMPLE} \ 
Fix a prime $p$ and take $\F$ = $\mathbb{Q}$.  Given a rational number $x$ $\ne$ 0, write
\[
x = p^k\frac{m}{n} 	\qquad (k \in \Z),
\]
where $p \nmid m$, $p \nmid n$, and then define the 
\un{$p$-adic absolute value} 
\index{$p$-adic absolute value} 
$\acdot_p$ by the prescription
\[
\abs{x}_p = p^{-k}	\qquad (\abs{0}_p = 0).
\]
\quad[AV-1 is obvious.
To check AV-2, write 
\[
x = p^k\frac{m}{n},  \ y = p^\ell\frac{u}{v},
\]
where $m, n, u, v$ are coprime to $p$ $-$then
\[
xy = p^{k + \ell} \frac{mu}{nv}
\]
\indent\indent\indent$\implies$
\[
\abs{xy}_p = p^{-(k+\ell)} = p^{-k} p^{-\ell} = \abs{x}_p\abs{y}_p.
\]
As for AV-3, $\acdot_p$ satisfies the ultrametric inequality.  
To establish this, assume without loss of generality that $k \le \ell$ and write 
\[
\begin{aligned}
\indent
x + y 
&= p^k\bigl(\frac{m}{n}+ p^{\ell-k}\frac{u}{v}\bigr )\\
& = p^k\frac{mv + p^{\ell-k}nu}{nv}.
\end{aligned}
\]
\indent\textbullet \quad $\abs{x}_p \ne \abs{y}_p$, so $\ell - k > 0$, hence 
\[
mv + p^{\ell-k}nu
\]
is coprime to $p$ (otherwise, 
\[
\begin{aligned}
mv 
&= p^r N - p^{\ell-k}nu	\quad (r \ge 1)\\
&= p(p^{r-1} N - p^{\ell - k - 1} nu) 
\end{aligned}
\]
\indent\indent\indent\indent\indent\indent $\implies p | mv$)\\ 
\indent\indent\indent $\implies$\\
\[
\begin{aligned}
\abs{x + y}_p 
&= p^{-k}\\
&= \abs{x}_p\\
&= \sup(\abs{x}_p, \abs{y}_p),
\end{aligned}
\]
since
\[
\begin{aligned}
\ell - k > 0 
&\implies p^{-\ell} <  p^{-k}\\ 
&\implies \abs{y}_p < \abs{x}_p.
\end{aligned}
\]
\indent\textbullet  \quad $\abs{x}_p = \abs{y}_p$, so, $\ell = k$, hence
\[
\begin{aligned}
&mv + nu = p^rN		\quad (r \ge 0) \ (p \nmid N) \\
\implies\\
&x + y = p^{k+r} \frac{N}{nv}\\
\implies\\
& \abs{x + y}_p = p^{-k-r}.
\end{aligned}
\]
And
\[p^{-k-r} \le\left\{
\begin{array}{l l}
p^{-k} = \abs{x}_p\\
p^{-k} = \abs{y}_p
\end{array}
\right.\]
\indent\indent$\implies$
\[
\abs{x + y}_p \le \sup(\abs{x}_p,\abs{y}_p).]
\]
\end{x} 
\vspace{0.1cm}

\begin{x}{\small\bf REMARK} \ 
It can be shown that every nontrivial absolute value on $\mathbb{Q}$ is equivalent to a $\acdot_p$ for some $p$ or to $\acdot_{\infty}.$
\end{x}
\vspace{0.1cm}

\begin{x}{\small\bf LEMMA} \ 
$\forall  \ x \in \mathbb{Q}^\times$,
\[
\prod_{p \le \infty} \abs{x}_p = 1,
\]
all but finitely many of the factors being equal to 1.\\

\quad PROOF \ Write
\[
x \ = \ \pm {p_1}^{k_1} \dotsb {p_n}^{k_n}		\quad (k_1, \dotsb , k_n \in \Z)
\]
for pairwise distinct primes $p_j$ $-$then $\abs{x}_p$ = 1 if $p$ is not equal to any of the $p_j$.  
In addition,
\[
\abs{x}_{p_j} \ = \ p^{-k_j},  \ \ \abs{x}_\infty \ = \ {p_1}^{k_1} \dotsb {p_n}^{k_n}	
\]
\indent\indent$\implies$
\[
\begin{aligned}
\prod_{p \le \infty} \abs{x}_p 
&= \bigl(\prod_{j=1}^n p_j^{-k_j}\bigr ) \cdot {p_1}^{k_1} \dotsb {p_n}^{k_n}\\
& = 1.	
\end{aligned}
\]
\\
\end{x}
\vspace{0.1cm}

\begin{x}{\small\bf REMARK} \ 
If $p_1$,  $p_2$, are distinct primes, then $\acdot_{p_1}$ is not equivalent to  $\acdot_{p_2}$.\\

\indent [Consider the sequence $\{{p_1}^n\}:$
\[
\abs{p_1}_{p_1} = p_1^{-1} \implies \abs{p_1^n}_{p_1} = p_1^{-n}  \ra 0.
\]
Meanwhile,
\[
\abs{p_1}_{p_2} = \abs{p_2^0 p_1}_{p_2} = p_2^{-0} = 1 
\]
\[
\implies \abs{p_1^n}_{p_2} \equiv 1.]
\]
\end{x}
\vspace{0.1cm}

\begin{x}{\small\bf CRITERION} \ 
Let $\acdot$ be an absolute value on $\F$ $-$then  $\acdot$ is non-archimedean iff $\{\abs{n}: n \in \N\}$ is bounded.
\\
\vspace{0.25cm}
\indent [Note: In either case, $\abs{n}$ is bounded by 1$:$
\[
\abs{n} = \abs{1 + 1 + \dotsb + 1} \le 1.]
\]
\end{x}

\chapter{
$\boldsymbol{\S}$\textbf{2}.\quad  TOPOLOGICAL FIELDS}
\setlength\parindent{2em}
\setcounter{theoremn}{0}


\ \indent Let $\abs{\hspace{0.05cm}\cdot \hspace{0.05cm}}$ be an absolute value on a field $\F$. Given a $\in \F, r > 0$, put
\[
N_r(a) = \{b: \abs{b - a} < r\}.
\]

\begin{x}{\small\bf LEMMA} \ 
There is a topology on $\F$ in which a basis for the neighborhoods of a are the $N_r(a)$.\\

PROOF \  The nontrivial point is to show that given  $V \in$ $\sB_a$
($\sB_a$ = the set of open balls centered at $a$)\xspace
, 
there is a $V_0 \in \sB_a$ such that if $a_0 \in V_0$, then there is a $W \in \sB_{a_0}$ such that $W \subset V$.  
So let V = $N_r(a)$, $V_0 = N_{r/2M}(a)$, $W = N_{r/2M}(a_0)$ $(a_0 \in V_0) -$ then $W \subset V:$
\[
\begin{aligned}
b \in W \implies \abs{b - a} \ 
&= \ \abs{(b - a_0) + (a_0 - a)}\\
&\le \ M \sup(\abs{b - a_0}, \abs{a_0 - a}) \\
&\le \ M \sup(r/2M, r/2M)\\
&= \  M(r/2M) \\
&= \ r/2 \\
&< \  r.
\end{aligned}
\]
\end{x}

\begin{x}{\small\bf EXAMPLE} \ 
The topology induced by $\acdot$ is the discrete topology iff  $\acdot$ is the trivial absolute value.
\end{x}

\begin{x}{\small\bf FACT} \ 
Absolute values $\acdot_1$, and $\acdot_2$ are equivalent iff they give rise to the same topology.
\end{x}

\begin{x}{\small\bf LEMMA} \ 
The topology induced by $\acdot$  is metrizable.\\

PROOF \  This is because $\acdot$ is equivalent to an absolute value satisfying the 
triangle inequality (cf. $\S1$, \ \#14), the underlying metric being 
\[
d(a,b) \ = \  \abs{a - b}.
\]
\end{x}

\begin{x}{\small\bf THEOREM} \ 
A field with a topology defined by an absolute value is a 
\underline{topological} \underline{field} 
\index{topological field}
i.e., the operations sum, product, and inversion are continuous.\\
\end{x}
\vspace{0.1cm}

Assume now that $\acdot$  is non-archimedean, hence that the ultrametric inequality 
\[
\abs{a - b} \ \le \  \sup (\abs{a},\abs{b})
\]
is in force.\\

\begin{x}{\small\bf LEMMA} \ 
$N_r(a)$ is closed (open is automatic).\\

\indent PROOF \   Let $p$ be a limit point of $N_r(a)$ $-$then $\forall$  $t > 0$,
\[
(N_t(p) - \{p\}) \cap N_r(a) \ne \emptyset
\]
Take $t = \ds\frac{r}{2}$ and choose $b \in N_r(a):$
\[
d(p,b) < \frac{r}{2}	\quad \text{$(p \ne b)$}.
\]
Then
\[
\begin{aligned}
d(a,p) 
&\le  \ \sup(d(a,b), d(b,p)) \\
&< \  r
\end{aligned}
\]
$\implies$
\[
p \in N_r(a).
\]
Therefore, $N_r(a)$ contains all its limit points, hence is closed.\\
\end{x}

\begin{x}{\small\bf LEMMA} \ 
If $a^{\prime} \in N_r(a)$, then $N_r(a^{\prime} ) = N_r(a)$.\\

PROOF \  E.g:
\[
b \in N_r(a) \implies \abs{b - a} < r
\]
$\implies$
\[
\begin{aligned}
\abs{b - a^{\prime}} 
&= \  \abs{(b - a) + (a - a^{\prime})}\\
&\le \  \sup(\abs{b - a},\abs{a - a^{\prime}})\\
&< r 
\end{aligned}
\]
$\indent\indent\implies$ 
\[
N_r(a) \subset N_r(a^{\prime}).
\]
\end{x}

\begin{x}{\small\bf REMARK} \ 
Put
\[
B_r(a) = \{b:\abs{b - a} \le r\}.
\]
Then a priori, $B_r(a)$ is closed.  But $B_r(a)$ is also open and if $a^{\prime} \in B_r(a)$, then  $B_r(a^{\prime}) = B_r(a)$.
\end{x}
\vspace{0.1cm}

\begin{x}{\small\bf LEMMA} \ 
If
\[
a_1 + a_2 + \dotsb + a_n = 0, 
\]
then $\exists$  $i \ne j$ such that
\[
\abs{a_i} = \abs{a_j} = \sup\abs{a_k}.
\]

PROOF \
Without loss of generality write $a_1 = \sup\limits_{1 \leq k \leq n}\abs{a_k}$.  Then
\allowdisplaybreaks
\begin{align*}
\abs{a_1} 	\ 
&=\  \abs{0 - a_1}\\	
&=\  \abs{a_1 + a_2 + \dotsb + a_n - a_1}\\
&=\  \abs{a_2 + \dotsb + a_n}\\	
&\leq\ \sup\limits_{2 \leq k \leq n}\abs{a_k}\\
&=\ \abs{a_j} \qquad (\exists \ j: 2 \leq j \leq n)\\
&\leq\ \sup\limits_{1\leq k \leq n}\abs{a_k}\\
&= \abs{a_1}.
\end{align*}

\end{x}

\chapter{
$\boldsymbol{\S}$\textbf{3}.\quad  COMPLETIONS}
\setlength\parindent{2em}
\setcounter{theoremn}{0}

\ \indent 

Let $\acdot$  be a an absolute value on a field $\F$ which satisfies the triangle inequality $-$then per $\acdot$, $\F$ might or might not be complete. 
(Recall, a metric space is \un{complete} iff every Cauchy sequence converges.)

\vspace{0.1cm}

\begin{x}{\small\bf EXAMPLE} \ 
Take $\F = \R$ or $\Q$ and let $\acdot = \acdot_\infty -$then $\R$ is complete but $\Q$ is not.
\end{x}

\vspace{0.1cm}

\begin{x}{\small\bf EXAMPLE} \ 
Take $\F =  \Q$ and let $\acdot = \acdot_p -$then $\Q$ is not complete.

\vspace{0.05cm}

[To illustrate this, choose $p = 5$ and starting with $x_1 = 2$, define inductively a sequence $\{x_n\}$ of integers subject to
\[\left\{
\begin{array}{l l}
{x_n}^2 + 1 \equiv 0	 \quad \text{ $\mod 5^n$ } \\
{x_{n+1}} \equiv x_n	 \quad \ \  \text{ $\mod 5^n$ }
\end{array}
.
\right.\]
Then
\[
\abs{x_m - x_n}_5  \ \le \ 5^{-n}		\quad (m > n),
\]
so $\{x_n\}$ is a Cauchy sequence and, to get a contradiction, assume that it has a limit x in $\Q$, thus
\[
\begin{aligned}
\abs{x_n^2 + 1}_5 \le \  5^{-n} 
&\implies \abs{x^2 + 1}_5 = 0 \\
&\implies x^2 + 1 = 0 \dots .] 
\end{aligned}
\]
\end{x}

\vspace{0.1cm}

\begin{x}{\small\bf DEFINITION} \ 
If an absolute value is not non-archimedean, then it is said to be \un{archimedean}.
\end{x}

\vspace{0.1cm}

\begin{x}{\small\bf FACT} \ 
Suppose that $\F$ is a field which is complete with respect to an archimedean absolute value $\acdot$ 
$-$then $\F$ is isomorphic to either $\R$ or $\C$ and $\acdot$ is equivalent to $\acdot_\infty.$
\end{x}

\vspace{0.1cm}

\begin{x}{\small\bf RAPPEL} \ 
Every metric space X has a completion 
$\compl{X}$
.  
Moreover, there is an isometry $\phi:X \ra \compl{X}$ such that $\phi$(X) is dense in $\compl{X}$ and $\compl{X}$ is unique up to isometric isomorphism.  
(Recall, an isometry is a distance preserving mapping.  An isometry is injective, indeed, is a homeomorphism onto its image.)
\end{x}

\vspace{0.1cm}

\begin{x}{\small\bf CONSTRUCTION} \ 
The standard model for $\compl{X}$ is the set of all Cauchy sequences in $X$ modulo the equivalence relation 
$\sim$, 
where
\[
\{x_n\} \sim \{y_n\} \Leftrightarrow d(x_n,y_n) \ra 0,
\]
the map $\phi:X \ra \compl{X}$ being the rule that sends $x \in$ X to the equivalence class of the constant sequence $x_n = x$.

[Note: The metric on $\compl{X}$ is specified by
\[
\ov{d}(\{x_n\},\{y_n\}) = \lim_{n \ra \infty} d(x_n,y_n).]
\]

Take $X = \F$ and 
\[
d(x,y) = \abs{x - y}.
\]
Then the claim is that $\compl{\F}$ is a field.  
E.g.: Let us deal with addition.  
Given $\ov{x}, \ov{y} \in \compl{\F}$, how does one define $\ov{x} + \ov{y}$ ?  
To this end, choose sequences
$
\begin{cases}
x_n  \\
y_n 
\end{cases}
$
in $\F$ such that
$
\begin{cases}
x_n \ra \ov{x} \\
y_n \ra \ov{y}
\end{cases}
$
$-$then
\[
\begin{aligned}
d(x_n + y_n,x_m + y_m) \ 
&= \  \abs{x_n + y_n - x_m - y_m}\\
&= \ \abs{(x_n - x_m ) + (y_n - y_m)}\\
&\leq \ \abs{x_n - x_m } + \abs{y_n - y_m}.
\end{aligned}
\]
Therefore  $\{x_n + y_n\}$ is a Cauchy sequence in $\F$, hence converges in $\compl{\F}$ to an element $\ov{z}$.  If
$
\begin{cases}
{x_n}^\prime  \\
{y_n}^\prime
\end{cases}
$
are sequences in $\F$ converging to 
$
\begin{cases}
\ov{x}  \\
\ov{y}
\end{cases}
$
as well, then $\{x_n^\prime + y_n^\prime\}$ converges in $\compl{\F}$ to an element $\ov{z}^\prime$.  
And
\[
\ov{z} = \ov{z}^\prime.
\]
Proof: \ Choose n $\in \N$ such that\\

$
\indent\indent\indent\indent
\begin{cases}
\abs{\ov{z} - (x_n + y_n)} \ \ < \ds\frac{\epsilon}{3}  \\
\abs{\ov{z}^\prime - (x_n^\prime  + y_n^\prime)} < \ds\frac{\epsilon}{3} 
\end{cases}
$
\\
and
\[
\abs{(x_n + y_n) - (x_n^\prime + y_n^\prime)}  \ 
\le \  \abs{x_n - x_n^\prime} + \abs{y_n - y_n^\prime} < \frac{\epsilon}{3}.
\]
Then
\[
\begin{aligned}
\abs{\ov{z} - \ov{z}^\prime} 
&\le \ \abs{\ov{z} - (x_n + y_n)} + \abs{\ov{z}^\prime - (x_n + y_n)}\\
&\le \ \abs{\ov{z} - (x_n + y_n)} + \abs{\ov{z}^\prime - 
(x_n^\prime + y_n^\prime)} + \abs{(x_n^\prime + y_n^\prime) - (x_n + y_n)} < \epsilon\\
&\implies \ov{z}= \ov{z}^\prime.
\end{aligned}
\]
Therefore addition in $\F$  extends to $\compl{\F}$.  
The same holds for multiplication and 
inversion.  
Bottom line: $\compl{\F}$ is a field.  
Furthermore, the prescription 
\[
\abs{\ov{x}} = \ov{d}(x,0)  	\quad (\ov{x} \in \compl{\F})
\]
is an absolute value on $\compl{\F}$ whose underlying topology is the metric topology.  
It thus follows that $\compl{\F}$ is a topological field (cf. $\S$2, $\#5$).
\end{x}

\vspace{0.1cm}

\begin{x}{\small\bf EXAMPLE} \ 
Take $\F = \Q$, $\acdot = \acdot_p -$then the completion $\compl{\F} = \compl{\Q}$ 
is denoted by $\Q_p$, the field of \un{ $p$-adic numbers}.
\end{x}

\vspace{0.1cm}

\begin{x}{\small\bf LEMMA} \ 
If $\acdot$ is non-archimedean per $\F$, then $\acdot$ is non-archimedean per $\compl{\F}$.

\vspace{0.1cm}

PROOF \   
Given
$
\begin{cases}
\ov{x}\\
\ov{y}
\end{cases}
\in \compl{\F},
$
choose
$
\begin{cases}
x_n\\
y_n
\end{cases}
\in \F 
$
such that
$
\begin{cases}
x_n \ra \ov{x}_n  \\
y_n \ra \ov{y}_n
\end{cases}
\text{in } \compl{\F}:
$
\allowdisplaybreaks
\[
\begin{aligned}
\abs{\ov{x} - \ov{y}} 
&\le \ \abs{\ov{x} - x_n + x_n - y_n + y_n - \ov{y}}\\
&\le \ \abs{\ov{x} - x_n} + \abs{x_n - y_n} + \abs{y_n - \ov{y}}.\\
& \qquad\quad \downarrow \qquad\qquad\qquad\qquad\quad  \downarrow\\
& \qquad\quad \  0 \qquad\qquad\qquad\qquad\quad  0
\end{aligned}
\]
And
\[
\begin{aligned}
\abs{x_n - y_n}  \ 
&\le \  \sup(\abs{x_n},\abs{y_n})\\
&= \  \frac{1}{2} (\abs{x_n} + \abs{y_n}) + \abs{x_n - y_n})\\
&\ra \ \frac{1}{2} (\abs{\ov{x}} + \abs{\ov{y}}) + \abs{\ov{x} - \ov{y}})\\
&= \  \sup(\abs{\ov{x}},\abs{\ov{y}}).
\end{aligned}
\]
\end{x}

\vspace{0.1cm}

\begin{x}{\small\bf LEMMA} \ 
If $\acdot$ is non-archimedean per $\acdot$, then
\[
\{\abs{\ov{x}}: \ov{x} \in \compl{\F}\} = \{\abs{x}: x \in \F\}.
\]

\vspace{0.1cm}

PROOF \    Take $\abs{\ov{x}} \in \compl{\F}:\ov{x} \ne 0$.  
Choose x $\in \F:\abs{\ov{x} - x} \ < \ \abs{\ov{x}}.$ 
Claim: $\abs{\ov{x}} = \abs{x}$.
\\
Thus, consider the other possibilities.
\[
\begin{aligned}
\text{\textbullet} \abs{x} < \abs{\ov{x}}:\\
&\abs{\ov{x} - x} = \abs{\ov{x} + (-x)} = \abs{\ov{x}}	\quad (\text{c.f.}  \ \S1, \ \# 18) < \abs{\ov{x}} \dots \ .\\
\text{\textbullet} \abs{\ov{x}} < \abs{x}:\\
&\abs{\ov{x} - x} = \abs{-x + \ov{x}} = \abs{-x}	 \quad (\text{c.f.} \   \S1, \ \# 18) = \abs{x} < \abs{\ov{x}} \dots \ .
\end{aligned}
\]
\end{x}
\vspace{0.1cm}

\begin{x}{\small\bf EXAMPLE} \ 
The image of $\Q_p$ under $\acdot_p$ is the same as the image of $\Q$ under $\acdot_p$, namely
\[
\{p^k:k \in\Z\} \cup \{0\}.
\]
\end{x}
\vspace{0.1cm}

Let $\K$ be a field, $\LL/\K$ a finite field extension.\\

\begin{x}{\small\bf EXTENSION PRINCIPLE} \ 
Let $\acdot_\K$ be a complete absolute value on $\K$ $-$then there is one and only one extension 
$\acdot_\LL$ of $\acdot_\K$  to $\LL$ and it is given by
\[
\abs{x}_\LL = \abs{N_{\LL/\K}(x)}_\K^{1/n},
\]
where n = $[\LL:\K]$.  In addition, $\LL$ is complete with respect to $\acdot_\LL$. 

\vspace{0.1cm}

[Note: $\acdot_\LL$ is non-archimedean if $\acdot_\K$ is non-archimedean.$]$
\end{x}
\vspace{0.1cm}

\begin{x}{\small\bf SCHOLIUM} \ 
There is a unique extension of $\acdot_\K$ to the algebraic closure 
$\K^{c\ell}$ 
of $\K$.\\
\indent [Note: It is not true in general that $\K^{c\ell}$  is complete.$]$\\

Suppose further that $\LL/\K$ is a Galois extension.  
Given $\sigma \in$ $\Gal(\LL/\K)$, define 
$\acdot_\sigma$ by 
$\abs{x}_\sigma = \abs{\sigma x}_\LL$
$-$then
\[
\restr{\acdot_\sigma}{\K} = \acdot_\K,
\]
so by uniqueness, $\acdot_\sigma$ = $\acdot_L$. 
But
\[
N_{\LL/\K}(x) = \prod_{\sigma\in \Gal(\LL/\K)} \sigma x
\]
\indent\indent$\implies$\\
\[
\begin{aligned}
\abs{N_{\LL/\K}(x)}_\K \ 
&= \ \abs{N_{\LL/\K}(x)}_\LL \\
&= \  \bigl |\prod_{\sigma\in \Gal(\LL/\K)} \sigma x \bigr |_\LL\\
&= \  \prod_{\sigma\in \Gal(\LL/\K)} \abs{\sigma x}_\LL\\
&= \  \prod_{\sigma\in \Gal(\LL/\K)} \abs{x}_\LL\\
&= \  \abs{x}_\LL^{\#(\Gal(\LL/\K))}\\
&= \  \abs{x}_\LL^{[\LL:\K]}\\
&= \  \abs{x}_\LL^n.
\end{aligned}
\]
\end{x}
\vspace{0.1cm}

\[
\textbf{APPENDIX}
\]
\setcounter{theoremn}{0}

\begin{x}{\small\bf APPROXIMATION PRINCIPLE} \ 
Let $\acdot_1, \dots , \acdot_N$ be pairwise inequivalent non-trivial absolute values on $\F$ .  
Fix elements $a_1, \dotsb , a_N$ in $\F$ $-$then $\forall$  $\epsilon > 0$, $\exists$ $a_\epsilon \in \F$:
\[
\abs{a_\epsilon - a_k}_k \  < \  \epsilon			\quad (k = 1, \dotsb, N). 
\]
\indent Let $\compl{\F}_1, \dotsb , \compl{\F}_N$ be the associated completions and let
\[
\Delta:\F \ra \prod_{k=1}^N \compl{\F}_k
\]
be the diagonal map $-$then the image $\Delta \F$ is dense (i.e., its closure is the whole of $\prod_{k=1}^N \compl{\F}_k$).

[Fix $ \epsilon > 0$ and elements $\ov{a}_1, \dots , \ov{a}_N$ in $\compl{\F}_1, \dots , \compl{\F}_N$ respectively 
$-$then there exist elements $a_k \in \F$:
\[
\abs{a_k - \ov{a}_k}_k < \epsilon			\quad (k = 1, \dots, N). 
\]
Choose $a_\epsilon \in \F$:
\[
\abs{a_\epsilon  - \ov{a}_k} < \epsilon	\quad (k = 1, \dots, N).
\]
Then
\[
\begin{aligned}
\abs{a_\epsilon  - \ov{a}_k}_k 
&= \ \abs{(a_\epsilon - a_k) + (a_k - \ov{a}_k)}_k\\
&\le \ \abs{a_\epsilon  - a_k}	 + \abs{a_k  - \ov{a}_k}_k\\
&< \ 2\epsilon.]
\end{aligned}
\]
\end{x}

\begin{x}{\small\bf \un{N.B.}} \ 
The product $\prod\limits_{k=1}^N \compl{\F}_k$ carries the product topology and the prescription
\[
\begin{aligned}
d((\ov{a}_1, \dots , \ov{a}_N), (\ov{b}_1, \dots , \ov{b}_N)) 
&= \ \sup_{1 \le k \le N} d_k(\ov{a}_k, \ov{b}_k) \\
&= \  \sup_{1 \le k \le N} \abs{\ov{a}_k - \ov{b}_k}_k
\end{aligned}
\]
metrizes the product topology.  Therefore 
\[
\begin{aligned}
d((a_\epsilon, \dots , a_\epsilon), (\ov{a}_1, \dots , \ov{a}_N)) \
&= \  \sup_{1 \le k \le N} d_k(a_\epsilon, \ov{a}_k)\\
&= \ \sup_{1 \le k \le N} \abs{a_\epsilon - \ov{a}_k}_k\\
&< \  2\epsilon.
\end{aligned}
\]
\end{x}


\chapter{
$\boldsymbol{\S}$\textbf{4}.\quad  p-ADIC STRUCTURE THEORY}
\setlength\parindent{2em}
\setcounter{theoremn}{0}

\ \indent 
Fix a prime $p$ and recall that 
$\Q_p$
\index{$\Q_p$}
is the completion of $\Q$ per the $p$-adic absolute value 
$\acdot_{p}$.\\

\begin{x}{\small\bf NOTATION} \  
Let
\[
\sA= \{0,1, ..., p-1\}.
\]
\index{$\sA$}
\end{x}
\vspace{0.1cm}
			
\begin{x}{\small\bf SCHOLIUM}  \ 
Structurally, $\Q_{p}$ is the set of all Laurent series in $p$ with coefficients in $\sA$ subject to the restriction that only finitely many of the negative powers of $p$ occur, thus generically a typical element x $\neq$ 0 of $\Q_{p} $ has the form 
\[
x =  \sum_{n=N}^\infty a_{n} p^n   \quad (a_{n}  \in \sA , \ N \in  \Z).
\]
\end{x}
\vspace{0.1cm}

\begin{x}{\small\bf \un{N.B.}}  \  
It follows from this that $\Q_{p}$ is uncountable, so $\Q$ is not complete per $\acdot_{p}$.\\
\end{x}
\vspace{0.1cm}

\indent The exact formulation of the algebraic rules (i.e., addition, multiplication, inversion) is elementary (but technically a bit of a mess) and will play no role in the sequel, hence can be omitted.\\

\begin{x}{\small\bf LEMMA} \ 
Every positive integer $N$ admits a base $p$ expansion:
\[
N =  a_{0} + a_{1} p + ... + a_{n} p^n, 	 
\]
where the $a_{n}  \in \sA$.
\end{x}
\vspace{0.1cm}

\begin{x}{\small\bf EXAMPLE} \ 
\[
1= 1 +  0p  + 0p^2 + \dots \ .
\]
\end{x}
\vspace{0.1cm}

\begin{x}{\small\bf EXAMPLE} \ 
Take $p = 3$ $-$then
\[\left\{
\begin{array}{l l}
24 = 0 +  2 \times 3 + 2 \times 3^2 = 2p + 2p^2\\
17 = 2 + 2 \times 3 + 1 \times 3^2 = 2 + 2p + p^2
\end{array}
\right.\]

\indent\indent $\implies$	
\[
\ds\frac{24}{17}  = \ds\frac{2p + 2p^2}{2 + 2p + p^2} = p + p^3 + 2p^5 + p^7 + p^8 + 2p^9 + \dots \ .
\]
\end{x}
\vspace{0.1cm}

\begin{x}{\small\bf LEMMA} \ 
\[
			-1= (p-1) + (p-1)p + (p-1)p^2 + \dots \ .
\]
\indent PROOF \ 
\[
\begin{aligned}
1 + &(p-1) + (p-1)p + (p-1)p^2 + (p-1)p^3 + ...\\
&= \ p + (p-1)p + (p-1)p^2 + (p-1)p^3 + ... \\
&= \ p^2 + (p-1)p^2 + (p-1)p^3 + ... \\
&= \ p^3 + (p-1)p^3 + ... \\
&= \ 0.
\end{aligned}
\]
\end{x}
\vspace{0.1cm}

\begin{x}{\small\bf APPLICATION} \ 
\[
\begin{aligned}
-N &= (-1) \cdot N \\
&= \ \bigl(\sum_{i=0}^\infty (p-1) p^i \bigr) (a_0 + a_1p + ... + a_np^n) \\
&= \ ...
\end{aligned}
\]
\end{x}
\vspace{0.1cm}

\begin{x}{\small\bf LEMMA} \  
A $p$-adic series
\[
\sum_{n=0}^\infty x_n  	 \quad (x_{n}  \in \Q_{p}) 
\]
is convergent iff 
$\abs{x_n}_{p} \ra 0  \quad (n  \ra  \infty)$.\\ 

\indent PROOF \ 
The usual argument establishes necessity.
So suppose that $\abs{x_n}_p  \ra  0$  $(n  \ra  \infty)$.  
Given $K  >  0$, $\exists \  N$:
\[
n > N \implies \abs{x_n}_{p} < p^{-K}.
\]
Let
\[
s_n = \sum_{k=1}^n x_k.
\]
Then
\[
\begin{aligned}
m > n > N 
&\implies \abs{s_m - s_n}_{p} = \abs{x_{n+1} + \dots + x_m}_{p}\\
&\le \  \sup(\abs{x_{n+1}}_{p}, \dots, \abs{x_m}_{p})\\
&<  \ p^{-K}.
\end{aligned}
\]
Therefore the sequence $\{s_n\}$ of partial sums is Cauchy, thus is convergent $(\Q_{p}$ being complete).\\
\end{x}

\vspace{0.1cm}

\begin{x}{\small\bf EXAMPLE} \ 
The p-adic series
\[
\sum_{i=0}^\infty p^i
\]
is convergent (to $\ds\frac{1}{1-p})$.
\end{x}

\vspace{0.1cm}

\begin{x}{\small\bf EXAMPLE} \ 
The $p$-adic series
\[
\sum_{n=0}^\infty n!
\]
is convergent.

\vspace{0.2cm}

[Note that 
\[
\abs{n!}_p = p^{-N},
\]
where
\[
N = [n/p] + [n/p^2] + \dots \ .]
\]
\end{x}

\vspace{0.1cm}

\begin{x}{\small\bf EXAMPLE} \ 
The $p$-adic series
\[
\sum_{n=0}^\infty n \cdot n!
\]
is convergent (to $-1$).
\end{x}

\vspace{0.1cm}

\begin{x}{\small\bf LEMMA} \ 
$\Q_p$ is a topological field (cf. $\S\ 2, \  \#5$).
\end{x}

\vspace{0.1cm}

\begin{x}{\small\bf LEMMA} \ 
$\Q_p$ is 0-dimensional, hence is totally disconnected.\\

\indent PROOF \  
A basic neighborhood $N_r(x)$ is open (by definition) and closed (cf. $\S2, \  \#6$).
\end{x}
\vspace{0.1cm}

\begin{x}{\small\bf NOTATION} \ 
\begin{align*} 
&\text{\textbullet} \quad \Z_p = \{x \in \Q_p:\abs{x}_p \le 1\} \qquad\qquad\qquad\qquad\qquad\qquad\qquad\qquad 
\index{$\Z_p$}\\
&\text{\textbullet} \quad p\Z_p = \{x \in \Q_p:\abs{x}_p < 1\} 
\index{$p\Z_p$}\\
&\text{\textbullet} \quad \Z_p^\times = \{x \in \Z_p:\abs{x}_p = 1\} 
\index{$\Z_p^\times$}
\end{align*} 
\end{x}
\vspace{0.1cm}

\begin{x}{\small\bf LEMMA} \ 
$\Z_p$ is a commutative ring with unit $($the ring of 
\underline{$p$-adic integers}, 
\index{$p$-adic integers})
in fact $\Z_p$ is an integral domain.\\
\end{x}

\vspace{0.1cm}

\begin{x}{\small\bf LEMMA} \ 
$p\Z_p$ is an ideal in $\Z_p$, in fact $p\Z_p$ is a maximal ideal in $\Z_p$, in fact $p\Z_p$ is the unique maximal ideal in $\Z_p$, hence $\Z_p$ is a local ring.\\
\end{x}

\vspace{0.1cm}

\begin{x}{\small\bf LEMMA} \ 
$\Z_p^\times$ is a group under multiplication, in fact $\Z_p^\times$ is the set of 
\underline{$p$-adic units}
\index{$p$-adic units} 
in  $\Z_p$, i.e., the set of elements in $\Z_p$ that have a multiplicative inverse in $\Z_p$.\\
\end{x}

\vspace{0.1cm}

Obviously, 
\[
\Z_p = \Z_p^\times \amalg (\Z_p - \Z_p^\times)
\]
or still, 
\[
\Z_p = \Z_p^\times \amalg p\Z_p.
\]
\begin{x}{\small\bf LEMMA} \ 
\[
\Z_p = \bigcup\limits_{0 \le k \le p-1} (k + p\Z_p).
\]

\indent PROOF \   
Let $x \in \Z_p$.  Matters being clear if $\abs{x}_p < 1$, (since in this case $x \in p\Z_p$), suppose that $\abs{x}_p = 1$.  Chose $q = \ds\frac{a}{b} \in \Q: \abs{q - x}_p < 1$, where $(a,b) = 1$ and 
$
\begin{cases}
(a,p) = 1\\
(b,p) = 1
\end{cases}
$
$-$then
\[
x + p\Z_p = q + p\Z_p.
\]
Choose $k$ with $ 0 < k \le p-1$ such that  $p$ divides $a - kb$, thus $\abs{a - kb}_p < 1$ and, moreover, 
$\abs{\ds\frac{a - kb}{b}}_p < 1$.  Therefore
\[
\begin{aligned}
\abs{k - \ds\frac{a}{b}}_p < 1 
&\implies k + p\Z_p = q + p\Z_p = x + p\Z_p \\
&\implies x \in k + p\Z_p.
\end{aligned}
\]

Consider a $p$-adic series
\[
\sum_{n=0}^\infty a_np^n		\qquad (a_n \in \sA).
\]
Then
\[
\begin{aligned}
\abs{\sum_{n=0}^\infty a_n p^n}_p 
&\le \ \sup_n \abs{a_np^n}_p \\
&\le \ \sup_n \abs{p^n}_p \\
&\le \ 1,
\end{aligned}
\]
so it converges to an element $x$ of $\Z_p $.  Conversely:
\end{x}
\vspace{0.1cm}

\begin{x}{\small\bf THEOREM} \ 
Every $x \in \Z_p$ admits a unique representation
\[
x = \sum_{n=0}^\infty a_np^n		\qquad \text{$a_n \in \sA$}.
\]

\indent PROOF \   
Let $x \in \Z_p$ be given.  
Choose uniquely $a_0 \in \sA$ such that $\abs{x - a_0}_p < 1$, hence $x = a_0 + px_1$ for some $x_1 \in \Z_p$.  
Choose uniquely $a_1 \in \sA $ such that $\abs{x_1 - a_1}_p < 1$, hence $x_1 = a_1 + p x_2$ for some $x_2 \in \Z_p$.  Continuing: $\forall$ $N$,
\[
x = a_0 + a_1p + \dots + a_Np^N + x_{N+1}p^{N+1},
\]
where $a_n \in \sA$ and $x_{N+1} \in \Z_p$.  But
\[
x_{N+1}p^{N+1} \ra 0.
\]
\end{x}
\vspace{0.1cm}

\begin{x}{\small\bf APPLICATION} \ 
$\Z$ is dense in $\Z_p$.
\end{x}
\vspace{0.1cm}

\begin{x}{\small\bf EXAMPLE} \ 
Let $x \in \Z_p -$then $\forall$ n $\in \N,$
\[
\binom{x}{n} = \frac{x(x-1) \dots (x-n+1)}{n!} \in \Z_p.
\]
\end{x}
\vspace{0.1cm}

\begin{x}{\small\bf LEMMA} \ 
\[
\Z_p^\times = \bigcup\limits_{1 \le k \le p-1} 	(k + p\Z_p).
\]
\end{x}
\vspace{0.1cm}

Consequently, if
\[
x = \sum_{n=0}^\infty a_np^n		\qquad \text{($a_n \in\sA$)}
\]
and if $x \in \Z_p^\times$, then $a_0 \ne 0$.\\
\vspace{0.1cm}
\indent [In fact, there is a unique $k$ $(1 \le k \le p-1)$ such that $x \in k + p\Z_p$ and this "k" is $a_0.]$\\

\begin{x}{\small\bf THEOREM} \ 
An element
\[
x = \sum_{n=0}^\infty a_np^n		\qquad (a_n \in\sA)
\]
in $\Z_p$ is a unit iff $a_0 \ne 0$.\\

\indent PROOF \  
To establish the characterization, construct a multiplicative inverse $y$ for $x$ as follows.  
First choose uniquely $b_0$  $(1 \le b_0 \le p-1)$ such that $a_0 b_0 \equiv 1 \mod p$.  
Proceed from here by recursion and assume that $b_1, \dots, b_M$ between 0 and $p-1$ have already been found subject to
\[
x\bigl(\sum_{0 \le m \le M} b_mp^m\bigr) \equiv 1 \mod p^{M+1}.
\]
Then there is exactly one $0 \le b_{M+1} \le p-1$ such that
\[
x \bigl(\sum_{0 \le m \le M+1} b_mp^m\bigr) \equiv 1 \mod p^{M+2}.
\]
Now put $y = \sum\limits_{m=0}^\infty b_mp^m$, thus $x y = 1$.\\
\end{x}
\vspace{0.1cm}

\begin{x}{\small\bf EXAMPLE} \ 
$1 - p$ is invertible in $\Z_p$ but $p$ is not invertible in $\Z_p$.\\
\end{x}
\vspace{0.1cm}

\begin{x}{\small\bf REMARK} \ 
The arrow
\[
\epsilon: \Z_p \ra \Z / p\Z
\]
that sends
\[
x = \sum_{n=0}^\infty a_np^n		\qquad \text{($a_n \in\sA$)}
\]
to $a_0 \ \modx p$ is a homomorphism of rings called 
\underline{reduction mod $p$}\index{reduction mod $p$}.  
It is surjective with kernel $p\Z_p$, hence $[\Z_p:p\Z_p] = p$.\\

Consider now the topological aspects of $\Z_p$:\\
\indent\textbullet \quad $\Z_p$ is totally disconnected.\\
\indent\textbullet \quad $\Z_p$ is closed, hence complete.\\
\indent\textbullet \quad $\Z_p$ is open.\\

\indent [As regards the last point, observe that
\[
\Z_p = \{x \in \Q_p: \abs{x}_p < r\} \equiv N_r(0) 	\qquad (1 < r < p).]
\]
\end{x}
\vspace{0.1cm}

\begin{x}{\small\bf THEOREM} \ 
$\Z_p$ is compact.\\

\indent PROOF \  
Since $\Z_p$ is a metric space, it suffices to show that $\Z_p$ is sequentially compact.  So let $x_1, x_2, \dots $ be an infinite sequence in $\Z_p$.  Choose $a_0 \in \sA$ such that $a_0 + p\Z_p$ contains infinitely many of the $x_n$.  Write
\[
\begin{aligned}
a_0 + p\Z_p  \ 
&= \ a_0 + p( \bigcup\limits_{a \in \sA} (a + p\Z_p))\\
&= \ a_0 +  \bigcup\limits_{a \in \sA} (ap + p^2\Z_p)\\
&= \ \bigcup\limits_{a \in \sA} (a_0 + ap + p^2\Z_p).
\end{aligned}
\]
Choose $a_1 \in \sA$ such that $a_0 + a_1p + p^2\Z_p$ contains infinitely many of the $x_n$.  Etc.  
The construction thus produces a descending sequence of cosets of the form
\[
A_j + p^j\Z_p,
\] 
each of which contains infinitely many of the $x_n$.  But
\[
\begin{aligned}
A_j + p^j\Z_p \  
&= \  \{x \in \Z_p: \abs{x - A_j}_p \le p^{-j}\}\\
&\equiv \  B_{p^{-j}}(A_j),
\end{aligned}
\]
a closed ball in the p-adic metric of radius $p^{-j} \ra 0$ $(j \ra \infty)$, hence by the completeness of $\Z_p$,
\[
\bigcap\limits_{j=1}^\infty B_{p^{-j}}(A_j) = \{A\}.
\]
Finally choose
\[
x_{n_1} \in B_{p^{-1}}(A_1),  \ x_{n_2} \in B_{p^{-2}}(A_2), \dots \ .
\]
Then
\[
\lim_{j \ra \infty} x_{n_j} = A.
\]
\end{x}
\vspace{0.1cm}

\begin{x}{\small\bf APPLICATION} \ 
$\Q_p$ is locally compact.\\

\indent [Since $\Q_p$ is Hausdorff, it is enough to prove that each $x \in \Q_p$ has a compact neighborhood.  But $\Z_p$ is a compact neighborhood of 0, so $x + \Z_p$ is a compact neighborhood of x.$]$
\end{x}
\vspace{0.1cm}

\indent The set $p^{-n}\Z_p$ $(n \ge 0)$ is the set of all $x \in \Q_p$ such that $\abs{x}_p \le p^n$.  
Therefore
\[
\Q_p = \bigcup\limits_{n=0}^\infty  p^{-n}\Z_p.
\]
Accordingly, $\Q_p$ is $\sigma$-compact (the $p^{-n}\Z_p$ being compact).\\
\vspace{0.1cm}

\begin{x}{\small\bf SCHOLIUM} \ 
A subset of $\Q_p$  is compact off it is closed and bounded.\\
\end{x}

\vspace{0.1cm}

\begin{x}{\small\bf LEMMA} \ 
Given $n,m \in \Z$,
\[
p^n\Z_p \subset p^m\Z_p \Leftrightarrow m \le n.
\]
\end{x}

\vspace{0.1cm}

\begin{x}{\small\bf REMARK} \ 
Take $ n\ge 1$ $-$then the $p^n\Z_p$ are principal ideals in $\Z_p$ and, apart from $\{0\}$, these are the only ideals in $\Z_p$, thus $\Z_p$ is a principal ideal domain.
\end{x}

\vspace{0.1cm}

\begin{x}{\small\bf LEMMA} \ 
For every $x_0 \in \Q_p$ and $r > 0$, there is an integer $n$ such that 
\[
\begin{aligned}
N_r(x_0) \ 
&= \{x \in \Q_p: \abs{x - x_0}_p < r\}\\
&= \ N_{p-n}(x_0) \\
&= \ \{x \in \Q_p: \abs{x - x_0}_p < p^{-n}\}\\
&= \ x_0 + p^{n+1}\Z_p
\end{aligned}
\]
\end{x}

\vspace{0.1cm}

\begin{x}{\small\bf SCHOLIUM} \ 
The basic open sets in  $\Q_p$  are the cosets of some power of $p\Z_p$.\\

\indent [Note: It is a corollary that every nonempty open subset of $\Q_p$ can be written as a disjoint union of cosets of the $p^n\Z_p$ $(n \in \Z).]$ \\
\end{x}

\vspace{0.1cm}

\begin{x}{\small\bf LEMMA} \ 
\[
p^n\Z_p^\times = p^n\Z_p - p^{n+1}\Z_p.
\]
\end{x}

\vspace{0.1cm}

\begin{x}{\small\bf DEFINITION} \ 
The $p^n\Z_p^\times$ are called 
\un{shells}\index{shells}.
\end{x}

\vspace{0.1cm}

\begin{x}{\small\bf \un{N.B.}} \ 
There is a disjoint decomposition
\[
\Q_p^\times = \bigcup\limits_{n \in \Z} p^n \Z_p^\times,
\]
where
\[
p^n\Z_p^\times = \bigcup\limits_{1 \le k \le p-1}(p^nk + p^{n+1}\Z_p).
\]

[Note: \   
For the record, $\Q_p^\times$ is totally disconnected and, being open in $\Q_p$, is Hausdorff and locally compact.  
Moreover, $\Z_p^\times$ is open-closed (indeed, open-compact).]\\

Let $x \in \Q_p^\times -$then there is a unique $v(x) \in \Z$  and a unique $u(x) \in \Z_p^\times$ 
such that $x = p^{v(x)}u(x)$.  
Consequently, 
\[
\Q_p^\times \approx \langle p \rangle \times \Z_p^\times
\] 
or still, 
\[
\Q_p^\times \approx \Z \times \Z_p^\times.
\] 
\end{x}

\vspace{0.1cm}

\begin{x}{\small\bf NOTATION} \ 
For $ n = 1, 2, \dots$, put
\[
U_{p,n} = 1 + p^n\Z_p.
\]

[Note: \ 
\[
1 + p^n\Z_p = \{x \in \Z_p^\times: \abs{1 - x}_p \le p^{-n}\}.]
\]

The $U_{p,n}$ are open-compact subgroups of $\Z_p^\times$ and 
\[
\Z_p^\times \supset U_{p,1} \supset U_{p,2} \supset \dots \ .
\]
\end{x}

\vspace{0.1cm}

\begin{x}{\small\bf LEMMA} \ 
The collection $\{U_{p,n}:n \in \N\}$ is a neighborhood basis at 1.\\
\end{x}

\vspace{0.1cm}

\begin{x}{\small\bf DEFINITION} \ 
$U_{p,1} = 1 + p\Z_p$ is called the group of 
\un{principal units}
\index{principal units} 
of $\Z_p$.\\
\end{x}

\vspace{0.1cm}

\begin{x}{\small\bf LEMMA} \ 
The quotient $\Z_p^\times/U_{p,1}$ is isomorphic to $\F_p^\times$ and the index of $U_{p,1}$ in $\Z_p^\times$ is $p-1$.
\end{x}
\vspace{0.1cm}

A generator of $\F_p^\times$ can be "lifted" to $\Z_p^\times$.\\

\begin{x}{\small\bf THEOREM} \ 
There exists a $\zeta \in \Z_p^\times$ such that $\zeta^{p-1} = 1$ and $\zeta^k \ne 1$ $(0 < k < p-1)$.

[This is a straightforward application of Hensel's lemma.]
\end{x}

\vspace{0.1cm}

\begin{x}{\small\bf \un{N.B.}} \ 
$\zeta \notin U_{p,1}$ ($p$ odd).

[If $x\in \Z_p$ and if for some $n \ge 1$,
\[
(1 + px)^n = 1,
\]
then using the binomial theorem one finds that $x = 0$.  This said, suppose that $\zeta \in U_{p,1}$:
\[
\zeta = 1 + pu \ (u \in \Z_p) \implies (1 + pu)^{p-1} = 1 \implies u = 0,
\]
a contradiction.]
\end{x}

\vspace{0.1cm}

\begin{x}{\small\bf SCHOLIUM} \ 
$\Z_p^\times$ can be written as a disjoint union
\[
\Z_p^\times = U_{p,1} \cup \zeta U_{p,1} \cup \zeta^2 U_{p,1} \cup \dots \cup \zeta^{p-2}U_{p,1}.
\]
Therefore
\[
\Q_p^\times \approx \Z \times \Z_p^\times \approx \Z \times \Z/(p-1)\Z \times U_{p,1}.
\]
\end{x}

\vspace{0.1cm}

\begin{x}{\small\bf LEMMA} \ 
Any root of unity in $\Q_p$ lies in $\Z_p^\times$.

PROOF \  
If $x = p^{v(x)}u(x)$ and if $x^n = 1$, then $n v(x) = 0$, so $v(x) = 0$, thus $x \in \Z_p^\times$.
\end{x}

\indent The roots of unity in $\Z_p^\times$ are a subgroup (as in any abelian group), call it $T_p$.  
If, on the other hand, $G_{p-1}$ is the cyclic subgroup of 
$\Z_p^\times$ generated by $\zeta$, 
then $G_{p-1}$ consists of  $(p-1)^{st}$ roots of unity, hence $G_{p-1} \subset T_p$.

\vspace{0.1cm}

\begin{x}{\small\bf LEMMA} \ 
If $p \ne 2$, then $G_{p-1} = T_p$ but if $p = 2$, then $T_p = \{\pm 1\}$.
\end{x}

\vspace{0.1cm}

\begin{x}{\small\bf APPLICATION} \ 
If $p_1$, $p_2$ are distinct primes, then $\Q_{p_1}$ is not field isomorphic to $\Q_{p_2}$.
\end{x}

\vspace{0.1cm}

\begin{x}{\small\bf REMARK} \ 
$\Q_p$ is not a field isomorphic to $\R$.

\vspace{0.1cm}

[$\Q_p$ has algebraic extensions of arbitrarily large linear degree which is not the case of $\R$ (cf. $\S5$, $\#26$).]
\end{x}

\vspace{0.1cm}

\begin{x}{\small\bf LEMMA} \ 
Let $x \in \Q_p^\times$ $-$then $x \in \Z_p^\times$ iff $x^{p-1}$ possesses $n^{th}$ roots for infinitely many $n$.

\vspace{0.1cm}

PROOF \   
If $x \in \Z_p^\times$ and if $n$ is not a multiple of $p$, then one can use Hensel's lemma to infer the existence of a $y_n \in \Z_p$ such that $y_n^n = x^{p-1}$.  Conversely, if $y_n^n = x^{p-1}$, then
\[
nv(y_n) = (p-1)v(x),
\]
thus $n$ divides $(p-1)v(x)$.  
But this can happen for infinitely many $n$ only if $v(x) = 0$, implying thereby that $x$ is a unit.
\end{x}

\vspace{0.1cm}

\begin{x}{\small\bf APPLICATION} \ 
Let $\phi: \Q_p \ra \Q_p$ be a field automorphism $-$then $\phi$ preserves units.

\vspace{0.1cm}

[In fact, if $x \in  \Z_p^\times$, then
\[
y_n^n = x^{p-1} \implies \phi(y_n)^n = (\phi(x))^{p-1}.]
\]
\end{x}
\vspace{0.1cm}

\begin{x}{\small\bf THEOREM} \ 
The only field automorphism $\phi$ of $\Q_p$ is the identity.

\indent PROOF \  
Given $x \in \Q_p^\times$, write $x = p^{v(x)}u(x)$, hence
\[
\begin{aligned}
\phi(x)  \ 
&= \  \phi(p^{v(x)}u(x))\\
&= \ \phi(p^{v(x)})\phi(u(x)) \\
&=   p^{v(x)} \phi(u(x)),
\end{aligned}
\]
hence
\[
v(\phi(x)) = v(x)			\qquad (\phi(u(x)) \in \Z_p^\times).
\]
Therefore $\phi$ is continuous.  
Since $\Q$ is dense in $\Q_p$, it follows that $\phi = id_{\Q_p}$.
\end{x}
\vspace{0.01cm}

[Note: 
\[
\begin{aligned}
x_k \ra 0 
&\implies \abs{x_k}_p \ra 0 \\
&\implies p^{-v(x_k)} \ra 0\\
&\implies p^{-v(\phi(x_k))} \ra 0 \\
&\implies \abs{\phi(x_k)}_p \ra 0 \\
&\implies \phi(x_k) \ra 0.]
\end{aligned}
\]

\indent The final structural item to be considered is that of quadratic extensions and to this end it is necessary to explicate $(\Q_p^\times)^2$, bearing in mind that
\[
\Q_p^\times \approx \Z \times \Z_p^\times \approx \Z \times \Z / (p-1)\Z \times U_{p,1}.
\]
\begin{x}{\small\bf LEMMA} \ 
If $p \ne 2$, then $U_{p,1}^2 = U_{p,1}$ but if $p = 2$, then $U_{2,1}^2 = U_{2,3}$.
\end{x}

\vspace{0.1cm}

\begin{x}{\small\bf APPLICATION} \ 
If $p \ne 2$, then
\[
(\Q_p^\times)^2 \approx 2\Z \times 2(\Z/(p-1)\Z) \times U_{p,1}
\]
but if $p = 2$, then
\[
(\Q_p^\times)^2 \approx 2\Z \times U_{2,3}.
\]
\end{x}

\vspace{0.1cm}

\begin{x}{\small\bf THEOREM} \ 
If $p \ne 2$, then
\[
[\Q_p^\times: (\Q_p^\times)^2] = 4
\]
but if $p = 2$, then
\[
[\Q_2^\times: (\Q_2^\times)^2] = 8.
\]
\end{x}

\vspace{0.1cm}

\begin{x}{\small\bf REMARK} \ 
If $p \ne 2$, then
\[
\Q_p^\times / (\Q_p^\times)^2 \approx \Z/2\Z \times \Z/2\Z
\]
but if $p = 2$, then
\[
\Q_2^\times / (\Q_2^\times)^2 \approx \Z/2\Z \times \Z/2\Z \times \Z/2\Z.
\]
\end{x}

\vspace{0.1cm}

\begin{x}{\small\bf CRITERION} \ 
Suppose that $p \ne 2$.

\vspace{0.1cm}

\indent\indent \textbullet \quad $p$ is not a square.

[If $p = x^2$, write $x = p^{v(x)}u(x)$ to get 
\[
1 = v(p) = v(x^2) = 2v(x),
\] 
an untenable relation.]\\

\indent\indent \textbullet \quad $\zeta$ is not a square.

[Assume that $\zeta = x^2 -$then
\[
\zeta^{p-1} = 1 \implies x^{2(p-1)} = 1,
\]
thus $x$ is a root of unity, thus $x\in T_p$, 
thus $x \in G_{p-1}$ (cf. $\# 45$), 
thus $x = \zeta^k$ $(0 < k < p-1)$, 
thus $\zeta = (\zeta^k)^2 = \zeta^{2k}$, thus $1 = \zeta^{2k-1}$.  
But
\[
2k < 2p - 2 \implies 2k-1 < 2p -1.
\]
And
\[
\qquad \qquad\qquad
\begin{cases}
2k - 1 = p -1 \implies 2k = p \implies \quad \text{$p$ even $\dots$} \\
2k - 1 = 2p -2 \implies 2k - 1 = 2(p-1) \implies   \quad \text{$2k - 1$ even $\dots \  .]$}
\end{cases}
.\]

\vspace{0.1cm}

\indent\indent \textbullet \quad $p\zeta$ is not a square.

\vspace{0.1cm}

[For if $p\zeta = p^{2n}u^2$ $(n \in \Z)$, then
\[
\begin{aligned}
\zeta = p^{2n-1}u^2 
&\implies 1 = \abs\zeta_p = \abs{p^{2n-1}}_p = p^{1 - 2n} \\
&\implies 1 - 2n = 0,
\end{aligned}
\]
an untenable relation.]
\end{x}
\vspace{0.1cm}

\begin{x}{\small\bf THEOREM} \ 
If $p \ne 2$, then up to isomorphism, $\Q_p$ has three quadratic extensions, viz.
\[
\Q_p(\sqrt{p}), \  \Q_p(\sqrt{\zeta}), \  \Q_p(\sqrt{p\zeta})
\]

[Note: if $\tau_1 = p$, $\tau_2 = \zeta$, $\tau_3 = p\zeta$, then these extensions of $\Q_p$ are inequivalent 
since $\tau_i\tau_j^{-1} (i \ne j)$ is not a square in $\Q_p.]$
\end{x}
\vspace{0.1cm}

\begin{x}{\small\bf REMARK} \ 
Another choice for the three quadratic extensions of $\Q_p$ when $p \ne 2$ is
\[
\Q_p( \sqrt{p} ), \  \Q_p( \sqrt{a} ), \ \Q_p( \sqrt{pa} ),
\]
where $1 < a < p$ is an integer that is not a square mod $p$.\\
\end{x}
\vspace{0.1cm}

\begin{x}{\small\bf REMARK} \ 
It can be shown that up to isomorphism, $\Q_2$ has seven quadratic extensions, viz.
\[
\Q_2(\sqrt{-1}),  \ \Q_2(\sqrt{\pm2}), \  \Q_2(\sqrt{\pm5}), \  \Q_2(\sqrt{\pm10}).
\]
\end{x}
\vspace{0.01cm}

\begin{x}{\small\bf EXAMPLE} \ 
Take $p = 5$ $-$then $2 \notin (\Q_5^\times)^2, 3 \notin (\Q_5^\times)^2$, but 6 $\in (\Q_5^\times)^2$.  
And
\[
\Q_5(\sqrt{2}) = \Q_5(\sqrt{3}).
\]
\indent [Working within $\Z_5^\times$, consider the equation $x^2 = 2$ and expand $x$ as usual:
\[
x = \sum_{n=0}^\infty a_n5^n			\qquad (a_n \in \sA).
\]
Then
\[
a_0^2 \equiv 2 \mod 5.
\]
But the possible values of $a_0$ are 0, 1, 2, 3, 4, thus the congruence is impossible, 
so $2 \notin (\Q_5^\times)^2$. Analogously, $3 \notin (\Q_5^\times)^2$.  On the other hand, $6 \in (\Q_5^\times)^2$ $($by direct verification or Hensel's lemma$)$, hence $6 = \gamma^2$  $(\gamma \in \Q_5)$.  
Finally, to see that
\[
\Q_5(\sqrt{2}) = \Q_5(\sqrt{3}),
\]
it need only be shown that $\sqrt{2} = a + b\sqrt{3}$ for certain $a, b \in \Q_5$.  
To this end, note that $\sqrt{2} \ \sqrt{3} = \pm \gamma$, from which
\[
\sqrt{2} \ =\  \pm \frac{\gamma}{\sqrt{3}} \ = \  \pm \frac{\gamma}{3}\sqrt{3}.]
\]
\end{x}
\vspace{0.01cm}

\begin{x}{\small\bf EXAMPLE} \ 
If $p$ is odd, then $p - 1$ is even and $-1 \in G_{p-1}$.  In addition, $-1 \in (\Q_2^\times)^2$ iff (p-1)/2 is even, i.e. iff $p \equiv 1 \mod 4$.  
Accordingly, to start $\sqrt{-1}$ exists in $\Q_5, \  \Q_{13}, \dots$ .
\end{x}
\vspace{0.1cm}

\indent [Note: $\sqrt{-1}$ does not exist in $\Q_2.]$\\

\[
\textbf{APPENDIX}
\]
\setcounter{theoremn}{0}

\indent Let $\Q_p^{c\ell}$ be the algebraic closure of $\Q_p$ 
$-$then $\acdot_p$ extends uniquely to $\Q_p^{c\ell}$ (cf. $\S 3$,  $\# 12$) (and satisfies the ultrametric inequality).  Furthermore, the range of $\acdot_p$ per $\Q_p^{c\ell}$ is the set of all rational powers of $p$ (plus 0).\\

\begin{x}{\small\bf THEOREM} \ 
$\Q_p^{c\ell}$ is not second category.\\
\end{x}
\vspace{0.1cm}

\begin{x}{\small\bf APPLICATION} \ 
The metric space $\Q_p^{c\ell}$ is not complete.\\
\end{x}
\vspace{0.1cm}

\begin{x}{\small\bf APPLICATION} \ 
The Hausdorff space $\Q_p^{c\ell}$ is not locally compact (cf. $\S5, \  \#5$).\\
\end{x}
\vspace{0.1cm}

\begin{x}{\small\bf NOTATION} \ 
Put
\[
\complement_p = \overline{\bigl(\Q_p^{c\ell}\bigr)}, \index{$\complement_p$}
\]
the completion of $\Q_p^{c\ell}$ per $\acdot_p$.\\
\end{x}
\vspace{0.1cm}

\begin{x}{\small\bf THEOREM} \ 
$\complement_p$ is algebraically closed.\\
\end{x}
\vspace{0.1cm}

\begin{x}{\small\bf \un{N.B.}} \ 
The metric space $\complement_p$ is separable but the Hausdorff space 
$\complement_p$ is not locally compact (cf. $\S5, \  \#5$).\\
\end{x}


\chapter{
$\boldsymbol{\S}$\textbf{5}.\quad  LOCAL FIELDS}
\setlength\parindent{2em}
\setcounter{theoremn}{0}

\ \indent Let $\K$ be a field of characteristic 0 equipped with a non-archimedean absolute value $\acdot$.

\vspace{0.1cm}

\begin{x}{\small\bf NOTATION} \ 
Let
\index{$R$}
\index{$R^\times$}
\[
\begin{cases}
R = \{a \in \K: \abs{a} \le 1\}\\
R^\times = \{a \in \K: \abs{a} = 1\}
\end{cases}
.
\]
\end{x}

\begin{x}{\small\bf LEMMA} \ 
$R$ is a commutative ring with unit and $R^\times$ is its multiplicative group of invertible elements.
\end{x}

\vspace{0.1cm}

\begin{x}{\small\bf NOTATION} \ 
Let
\index{$P$}
\[
P = \{a \in \K: \abs{a} < 1\}.
\]
\end{x}

\vspace{0.1cm}

\begin{x}{\small\bf LEMMA} \ 
$P$ is a maximal ideal.
\end{x}

\vspace{0.1cm}

Therefore the quotient $R/P$ is a field, the 
\un{residue field}
\index{residue field} of $\K$.

\vspace{0.2cm}

\begin{x}{\small\bf THEOREM} \ 
$\K$ is locally compact iff the following conditions are satisfied.

\vspace{0.1cm}

\indent 1. \ $\K$ is a complete metric space.\\
\indent 2. \  $R/P$ is a finite field.\\
\indent 3. \  $\abs{\K^\times}$ is a nontrivial discrete subgroup of $\R_{>0}$.
\end{x}

\vspace{0.1cm}

\begin{x}{\small\bf DEFINITION} \ 
A 
\un{local field}
\index{local field} is a locally compact field of characteristic 0.
\end{x}

\vspace{0.1cm}

\begin{x}{\small\bf EXAMPLE} \ 
$\R$ and $\C$ are local fields.
\end{x}

\vspace{0.1cm}

\begin{x}{\small\bf EXAMPLE} \ 
$\Q_{p}$ is a local field.
\end{x}

\vspace{0.1cm}

Assume that $\K$ is a non-archimedean local field.

\vspace{0.1cm}

\begin{x}{\small\bf LEMMA} \ 
$R$ is compact.
\end{x}

\vspace{0.1cm}

\begin{x}{\small\bf LEMMA} \ 
$P$ is principal, say $P = \pi$R, and 
\[
\abs{\K^\times}  = \abs{\pi}^{\Z}, \quad \text{where } 0 < \abs{\pi} < 1.
\]
\end{x}

\vspace{0.1cm}

[Note: Such a $\pi$ is said to be a \un{prime element}\index{prime element}.]

\vspace{0.1cm}

\begin{x}{\small\bf REMARK} \ 
A nontrivial discrete subgroup $\Gamma$ of $\R_{>0}$ is free on one generator $ 0 < \gamma < 1:$
\[
\Gamma = \{\gamma^n:n \in \Z\}.
\]
This said, choose $\pi$ with the largest absolute value $<$ 1, thus $\pi \in P \subset  R \Rightarrow \pi R \subset P.$
In the other direction, 
\[
a \in P \Rightarrow \abs{a} \le \abs{\pi} \Rightarrow \frac{a}{\pi} \in R.
\]
And
\[
a = \pi \cdot \frac{a}{\pi}  \Rightarrow a \in \pi R.
\]
\end{x}

\vspace{0.1cm}

\begin{x}{\small\bf FACT} \ 
A locally compact topological vector space over a local field is necessarily finite dimensional.
\end{x}

\vspace{0.1cm}

\begin{x}{\small\bf THEOREM} \ 
$\K$ is a finite extension of $\Q_{p}$ for some $p$.

\vspace{0.1cm}

PROOF \ 
First, $\K\supset \Q$ (since char $\K = 0$).  
Second, the restriction of $\acdot$ to $\Q$ is equivalent to $\acdot_{p}$ $(\exists \ p$)
(cf. \S1, \#20), 
hence the closure of $\Q$ in $\K$ "is" $\Q_{p}$ $($since $\K$ is complete$)$. Third, $\K$ is finite dimensional over $\Q_{p}$ $($since $\K$ is locally compact).
\end{x}
\indent There is also a converse.
\begin{x}{\small\bf THEOREM} \ 
Let $\K$ be a finite extension of $\Q_{p}$ $-$then $\K$ is a local field.

\vspace{0.1cm}

\indent PROOF \ 
In view of $\#5$, it suffices to equip $\K$ with a non-archimedean absolute value subject to the conditions 1, 2, 3. 
But, by the extension principle (cf. $\S$3, $\#$11), $\acdot_{p}$ extends uniquely to $\K$. 
This extension is non-archimedean and points 1, 3 are manifest. 
As for point 2, it suffices to observe that the canonical arrow
\[
\Z_{p}/p\Z_{p} \ra R/P
\]
is injective and
\[
[R/P:\F_p] \le [\K:\Q_{p}] < \infty.
\] 
\indent [Details: To begin with,
\[
\Q_{p} \cap P = p\Z_{p},
\]
thus the inclusion $\Z_{p} \ra$ R induces an injection
\[
\Z_{p}/p\Z_{p} \ra R/P.
\] 
Put now $n = [\K:\Q_{p}]$ and let $A_{1}, ..., A_{n+1} \in R$ 
$-$then the claim is that the residue classes $\ov{A}_{1}, ..., \ov{A}_{n+1}  \in$ 
$R/P$ are linearly dependent over $\Z_{p}/p\Z_{p}$.  
In any event, 
there are elements $x_{1}, ..., x_{n+1} \in \Q_{p}$ such that 
\[
\sum_{i=1}^{n+1} x_iA_i = 0, 
\]
matters being arranged in such a way that 
\[
\max\abs{x_i}_{p} = 1.
\]
Therefore the $ x_{i} \in \Z_{p}$ and not every residue class $\ov{x}_{i} \in \Z_{p}/p\Z_{p}$ is zero. But then
\[
\sum_{i=1}^{n+1} \ov{x}_i \ov{A}_i = 0
\]
is a nontrivial dependence relation.]\\
\end{x}

\vspace{0.1cm}

\begin{x}{\small\bf SCHOLIUM} \ 
A non-archimedean field of characteristic zero is a local field iff it is a finite extension of $\Q_{p}$ $(\exists$ $p)$.
\end{x}

\vspace{0.1cm}

Let $\K/\Q_{p}$ be a finite extension of degree $n$ $-$then the 
\un{canonical absolute value}
\index{canonical absolute value}  
on $\K$ is given by
\[
\abs{a}_{p} = \abs{N_{\K/\Q_{p}}(a)}_{p}^{1/n}.
\]

\vspace{0.1cm}

[Note: The 
\un{normalized absolute value}
\index{normalized absolute value}  
on $\K$ is given by
\[
\abs{a}_\K = \abs{a}_{p}^{n}.
\]
Its intrinsic significance will emerge in due course but for now observe that $\acdot_{\K}$ is equivalent to $\acdot_{p}$ and is non-archimedean (cf. $\S 1, \  \#23$).]

\vspace{0.1cm}

\begin{x}{\small\bf LEMMA} \ 
The range of $\restr{\acdot_{p}}{\K^\times}$  is $\abs{\pi}_{p}^\Z.$
\end{x}

\vspace{0.1cm}

\begin{x}{\small\bf DEFINITION} \ 
The 
\un{ramification index}
\index{ramification index} 
 of $\K$ over $\Q_{p}$ is the positive integer
\[
e = [\abs{\K^\times}_{p}:\abs{Q^\times_p}_{p}].
\]
I.e.,
\[
e = [\abs{\pi}_{p}^\Z:\abs{p}_{p}^\Z].
\]
Therefore
\[
\abs{\pi}_{p}^{e} = \abs{p}_{p} \qquad (= \frac{1}{p}).
\]

[Consider $\Z$ and e$\Z$  $-$then the generator 1 of $\Z$  is related to the generator e of e$\Z$ by the triviality 
$1 + \cdots + 1 = e\cdot 1 = e$.]
\end{x}

\vspace{0.1cm}

\begin{x}{\small\bf \un{N.B.} } \ 
If $\pi^\prime$ has the property that $\abs{\pi^\prime}_p^{e}$ = ${\abs{p}}_{p}$ then $\pi^\prime$ is a prime element.
\\
\indent [Using obvious notation, write $\pi^\prime$ = ${\pi^{v(\pi)}}u$, thus
\[
\begin{aligned}
\abs{p}_{p} \ 
&= \  \abs{\pi^\prime}_p^{e}\\ 
&= \  {(\abs{\pi}^{v(\pi)}_p)}^e \\
&= \ {(\abs{\pi}^{e}_p)}^{v(\pi)} \\
&= \ {\abs{p}_{p}}^{v(\pi)},
\end{aligned}
\]
thus v($\pi$) = 1.]\\
\end{x}

\begin{x}{\small\bf NOTATION} \ 
\[
q \equiv \card R/P = (\card \F_p)^f = p^f,
\]
so
\[
f = [R/P:\F_p],
\]
the 
\un{residual index}
\index{residual index}  
of $\K$ over $\Q_{p}$.\\
\end{x}

\vspace{0.1cm}

\begin{x}{\small\bf THEOREM} \ 
Let $\K/\Q_{p}$ be a finite extension of degree $n$ $-$then
\[
n = [\K:\Q_{p}] = ef.
\]
\end{x}

\begin{x}{\small\bf APPLICATION} \ 
\[
\begin{aligned}
{\abs{\pi}}_\K \ 
= \ \abs{\pi}_p^n \\
= \ \abs{p}_{p}^{n/e} \\
= \ \left(\frac{1}{p}\right)^{n/e} \\
= \ \left({\frac{1}{p}}\right)^f \\
= \ \frac{1}{p^f} \\
= \ \frac{1}{q}.
\end{aligned}
\]
\end{x}

\vspace{0.1cm}


View $p$ as an element of $\K$:

\vspace{0.1cm}

\qquad \textbullet \quad $\abs{p}_{p} =  \abs{N_{\K/\Q_{p}}(p)}_{p}^{1/n}= \abs{p^n}_{p}^{1/n} = \abs{p}_{p}$.

\qquad \textbullet \quad $\abs{p}_{\K} 
= \abs{N_{\K/\Q_{p}}(p)}_{p} 
= \abs{p^n}_{p} 
= \ds\frac{1}{p^n} 
= \ds\frac{1}{p^{ef}} 
= {\left(\ds\frac{1}{p^f}\right)}^e 
= q^{-e}.$

\vspace{0.3cm}

\begin{x}{\small\bf DEFINITION} \  
A finite extension $\K/\Q_{p}$ is

\vspace{0.1cm}

\qquad \textbullet\quad \un{unramified}\index{unramified} if $e = 1$

\qquad \textbullet\quad \un{ramified}\index{ramified} if $f = 1$.

\end{x}

\vspace{0.1cm}

Take the case $\K = \Q_{p}$ $-$then $e = 1$, hence $\K$ is unramified, and $f = 1$, hence $\K$ is ramified.

\vspace{0.1cm}

\begin{x}{\small\bf LEMMA} \ 
If $\K/\Q_{p}$ is is unramified, then $p$ is a prime element.
\end{x}

\vspace{0.1cm}

\begin{x}{\small\bf THEOREM} \ 
$\forall$ $n = 1, 2, \ldots$,  there is up to isomorphism one unramified extension $\K/\Q_{p}$ of degree $n$.
\end{x}

\vspace{0.1cm}

Let $\K/\Q_{p}$ be a finite extension.

\vspace{0.1cm}

\begin{x}{\small\bf LEMMA} \ 
The group $M^\times$ of roots of unity of order prime to $p$ in $\K$ is cyclic of order 
\[
p^f - 1 \quad \text{$(= q-1)$}.
\]
\end{x}

\vspace{0.1cm}

\begin{x}{\small\bf LEMMA} \ 
The set M = $M^\times$ $\cup$ \{0\} is  a set of coset representatives for $R/P$. Therefore (cf. $\S4$, $\#$43)
\[
\K^\times \approx \Z \times R^\times \approx \Z \times \Z/(q-1)\Z \times1 + P.
\]
\end{x}

\vspace{0.1cm}

\begin{x}{\small\bf NOTATION} \ 
Let
\index{$\K_{ur}$}
\[
\K_{ur} = \Q_p(M^\times).
\]
\end{x}

\begin{x}{\small\bf LEMMA} \ 
$\K_{ur}$ is the maximal unramified extension of $\Q_p$ in $\K$ and
\[
[\K_{ur}:\Q_p] = f.
\]
\end{x}

\begin{x}{\small\bf REMARK} \ 
The maximal unramified extension 
$(\Q_p^{c\ell})_{ur} \subset \Q_p^{c\ell}$ 
is the field extension generated by all roots of unity of order prime to $p$.
\end{x}

\vspace{0.1cm}

\begin{x}{\small\bf QUADRATIC EXTENSIONS} \ \ 
\index{Quadratic extensions}
(cf. $\S 4$, $\#56$)  \ Suppose that $p \ne$ 2, 
let $\tau \in \Q_p^\times - (\Q_p^\times)^2$, and form the quadratic extension
\[
\Q_p(\tau) = \{x + y\sqrt{\tau} : x, y \in \Q_p\}.
\]
Then the canonical absolute value on $\Q_p(\sqrt{\tau})$ is given by
\[
\begin{aligned}
\abs{x + y\sqrt{\tau}}_p \ 
&= \ \abs{N_{\Q_p(\sqrt{\tau})/\Q_p}   (x + y\sqrt{\tau}}_p^{1/2} \\
&= \ \abs{x^2 - \tau y^2}_p^{1/2}.
\end{aligned}
\]
\end{x}

\vspace{0.1cm}

\begin{x}{\small\bf CLASSIFICATION} \ 
Consider the three possibilities
\[
\Q_p(\sqrt{p}), \ \Q_p(\sqrt{\tau}), \ \Q_p(\sqrt{p \tau}),
\]
thus here $ef = 2$.
\begin{itemize}
\item $\Q_p(\sqrt{p})$ is ramified or still, $e = 2$.

\vspace{0.1cm}

[Note that
\[
\abs{\sqrt{p}}_p^2 = \abs{0^2 - (p)1^2}_p = \abs{p}_p = \frac{1}{p}.]
\]
\item $\Q_p(\sqrt{p \zeta})$ is ramified or still, $e = 2$.\\

\vspace{0.1cm}

[Note that
\[
\abs{\sqrt{p \zeta}}^2 
= \abs{0^2 - (p \zeta)1^2}_p 
= \abs{p \zeta}_p 
= \abs{p}_p \cdot \abs{\zeta}_p 
= \abs{p}_p = \frac{1}{p}.]
\]

If $e = 1$, then in either case, the value group would be $p^{\Z}$, an impossibility since 
$\ds\frac{1}{\sqrt{p}} \notin p^{\Z}$, so $e = 2$.
\item $\Q_p(\sqrt{\zeta})$ is unramified or still, $e = 1$.
\end{itemize}

\vspace{0.1cm}

[There is up to isomorphism one unramified extension $\K$ of $\Q_p$ of degree 2 (cf. $\# 24$)].

\vspace{0.1cm}

[Instead of quoting theory, one can also proceed directly, it being simplest to work instead with $\Q_p(\sqrt{a})$, 
where $1 < a < p$ is an integer that is not a square mod $p$ (cf. $\S 4$, $\#57$) 
$-$then the residue field of $\Q_p(\sqrt{a})$ is $\F_p(\sqrt{a})$, hence $f = 2$, hence $e = 1$ (since $n = 2$).]
\end{x}

\vspace{0.1cm}

The preceding developments are absolute, i.e., based at $\Q_p$.  
It is also possible to relativize the theory.  
Thus let $\LL/\K$, $\K/\Q_p$ be finite extensions.  
Append subscripts to the various quantities involved:
\[
\begin{cases}
R_\K \supset P_{\K}, \ R_\K/P_{\K}, \ e_{\K}, \ f_{\K}, \ M_\K^\times \\
R_\LL \supset P_{\LL}, \ R_\LL/P_{\LL}, \ e_{\LL}, \ f_{\LL}, \ M_\LL^\times
\end{cases}
.
\]
Introduce
\[
\begin{cases}
e(\LL/\K) = [\abs{\LL^\times}:\abs{\K^\times}] \\
f(\LL/\K) = [R_\LL/P_\LL:R_\K/P_\K]
\end{cases}
.
\]
\vspace{0.2cm}

\begin{x}{\small\bf LEMMA} \ 
\[
[\LL:\K] = e(\LL/\K)f(\LL/\K).
\]

\indent PROOF \   We have
\[
\begin{cases}
[\LL:\Q_p ] = e_\LL f_\LL \\
[\K:\Q_p] = e_\K f_\K
\end{cases}
\qquad (\text{ cf. $\# 20$}).
\]
Therefore
\[
[\LL:\K] = \frac{[\LL:\Q_p] }{[\K:\Q_p]} = \frac{e_\LL f_\LL}{e_\K f_\K} = e(\LL/\K) f(\LL/\K).
\]
\end{x}

\vspace{0.1cm}

\begin{x}{\small\bf THEOREM} \ 
Let $\LL/\K$, $\K/\Q_p$ be finite extensions $-$then there exists a unique maximal intermediate extension 
$\K \subset \K_{ur} \subset \LL$ that is unramified over $\K$.

\vspace{0.1cm}

[In fact, 
\[
\K_{ur} = \K(M_\LL^\times) \subset \LL.]
\]

[Note: The extension $\LL/\K_{ur}$ is ramified.]
\end{x}


\chapter{
$\boldsymbol{\S}$\textbf{6}.\quad  HAAR MEASURE}
\setlength\parindent{2em}
\setcounter{theoremn}{0}

\ \indent 
Let $X$ be a locally compact Hausdorff space.

\vspace{0.25cm}

\begin{x}{\small\bf DEFINITION} \ 
A 
\underline{Radon measure}
\index{Radon measure} 
is a measure $\mu$ defined on the Borel $\sigma$-algebra of $X$ subject to the following conditions.

\vspace{0.1cm}

1.  $\mu$ is finite on compacta, i.e., for every compact set $K$ $\subset$ X, $\mu(K) < \infty$.

\vspace{0.1cm}

2.  $\mu$ is outer regular, i.e., for every Borel set $A \subset$ $X$, 
\[
\mu(A) = \inf_{U\supset A} \mu(U),	 \quad \text{where $U \subset$ X is open.}
\]

3.  $\mu$ is inner regular, i.e., for every open set $A \subset$ X, 
\[
\mu(A) = \sup_{K\subset A} \mu(K),	 \quad \text{where $K \subset X$ is compact.}
\]
\end{x}

\vspace{0.1cm}

Let $G$ be a locally compact abelian group.

\vspace{0.2cm}

\begin{x}{\small\bf DEFINITION} \ 
A 
\underline{Haar measure}
\index{Haar measure} 
on $G$ is a Radon measure $\mu_{G}$ which is translation invariant: 
$\forall$ Borel set $A$, $\forall \ x \in G$,
\[
\mu_G(x+A) = \mu_G(A) = \mu_G(A+x)
\]
or still, $\forall \ f \in C_c(G), \forall \ y \in G$,
\[
\int_Gf(x+y)d\mu_G(x) \ = \  \int_Gf(x)d\mu_g(x).
\]
\end{x}

\vspace{0.1cm}

\begin{x}{\small\bf THEOREM} \ 
$G$ admits a Haar measure and for any two Haar measures $\mu_G$,  $\nu_G$ differ by a positive constant: 
$\mu_G  = c\nu_G  \ (c > 0)$.
\end{x}

\vspace{0.1cm}

\begin{x}{\small\bf LEMMA} \ 
Every nonempty open subset of $G$ has positive Haar measure.
\end{x}

\vspace{0.1cm}

\begin{x}{\small\bf LEMMA} \ 
$G$ is compact iff $G$ has finite Haar measure.
\end{x}

\vspace{0.1cm}

\begin{x}{\small\bf LEMMA} \ 
$G$ is discrete iff every point of $G$ has positive Haar measure.
\end{x}

\vspace{0.1cm}

\begin{x}{\small\bf EXAMPLE}   \ 
Take $G = \R$ $-$then $\mu_{ \R} = dx$ ($dx$ = Lebesgue measure) is a Haar measure 
($\mu_{\R} ([0,1]) = \ds\int_0^1dx = 1$).
\end{x}

\vspace{0.1cm}

\begin{x}{\small\bf EXAMPLE} \ 
Take $G = \R^\times -$then $\mu_{R^\times} = \ds\frac{dx}{\abs{x}}$ ($dx$ = Lebesgue measure) 
is a Haar measure ($\mu_{\R^\times} ([1,e]) = \ds\int_1^e\frac{dx}{\abs{x}} = 1$).
\end{x}

\begin{x}{\small\bf EXAMPLE} \ 
Take $G = \Z$ $-$then $\mu_{ \Z} =$ counting measure is a Haar measure.
\end{x}

\vspace{0.1cm}

\begin{x}{\small\bf LEMMA} \ 
Let $G^\prime$ be a closed subgroup of $G$ and put $G^{\prime \prime} = G/G^\prime$.  
Fix Haar measures $\mu_{G}$, $\mu_{G^{\prime}}$ on $G$, $G^\prime$ respectively 
$-$then there is a unique determination of the Haar measure 
$\mu_{G^{\prime \prime}}$ on $G^{\prime \prime}$ such that $\forall \ f \in C_c(G)$, 
\[
\int_Gf(x)d\mu_G(x) 
= \int_{G^{\prime \prime}} \left( \int_{G^\prime} f(x + x^\prime)d\mu_{G^\prime}(x^\prime)\right) 
d\mu_{G^{\prime \prime}} (x^{\prime \prime}).
\]

\vspace{0.1cm}

[Note: The function
\[
x \ra \int_{G^\prime} f(x + x^\prime) d\mu_{G^\prime} (x^\prime).
\]
is $G^{\prime}$-invariant, hence is a function on $G^{\prime \prime}.]$
\end{x}

\vspace{0.1cm}

\begin{x}{\small\bf EXAMPLE}
Take $G = \R$,  $G^\prime = \Z$ with the usual choice of Haar measures.  
Determine $\mu_{\R/\Z}$ per $\# 10$ $-$then $\mu_{\R/\Z}(\R/\Z) = 1$.

\vspace{0.1cm}

[Let $\chi$ be the characteristic function of $[0,1[$ $-$then
\[
\sum_{n \in \Z} \chi(x+n)
\]
is $\equiv$ 1, hence when integrated over $\R/\Z$ gives the volume of $\R/\Z$. 
On the other hand, 
\[
\int_{\R} \chi = 1.]
\]
\end{x}

\vspace{0.1cm}

Let $\K$ be a local field (cf. $\S$5, $\#6$).  
Given $a \in \K^\times$, let $M_a:\K \ra \K$ be the automorphism that sends $x$ to $ax = xa$ 
$-$then for any Haar measure $\mu_\K$ on $\K$, the composite $\mu_\K \circ M_a$ is again a Haar measure on  $\K$, 
hence there exists a positive constant $\modxs_\K(a)$ such that for every Borel set $A$,
\[
\mu_\K(M_a(A)) = \modxs_\K(a)\mu_\K(A)
\]
or still, $\forall \ f \in C_c(\K)$,
\[
\int_\K f(a^{-1}x)d\mu_\K(x) = \modxs_\K(a)\int_\K f(x)d\mu_\K(x).
\]

\vspace{0.1cm}

[Note: $\modxs_\K(a)$ is independent of the choice of $\mu_\K$.]

\vspace{0.2cm}

Extend $\modxs_\K$ to all of $\K$ by setting $\modxs_\K(0)$ equal to 0.

\vspace{0.1cm}

\begin{x}{\small\bf LEMMA} \ 
Let $\K$, $\LL$ be local fields, where $\LL/\K$ is a finite field extension $-$then $\forall$ $x \in \LL$,
\[
\begin{aligned}
\modxs_\LL(x) 
&= \modxs_\K(N_{\LL/\K}(x)) \\
&\equiv \ \modxs_\K(\det(M_x))
\end{aligned}
\]

[Let $n = [\LL:\K]$, view $\LL$ as a vector space of dimension $n$, and identify $\LL$ with $\K^n$ by choosing a basis.  
Proceed from here by breaking $M_x$ into a product of $n$ 
"elementary" transformations.]
\end{x}

\vspace{0.1cm}

\begin{x}{\small\bf EXAMPLE} \ 
Take $\K = \R$, $\LL = \R$ $-$then $\forall$ a $\in \R$,
\[
\modxs_{\R}(a) = \abs{a}.
\]
$[\forall$ f $\in C_c(\R)$,
\[
\int_{\R} f(a^{-1}x)dx = \abs{a} \int_{\R} f(x)dx.]
\]
\end{x}

\vspace{0.2cm}

\begin{x}{\small\bf EXAMPLE} 
Take $\K = \C$, $\LL = \C$ $-$then $\forall$ $a\in \C$,
\[
\allowdisplaybreaks
\begin{aligned}
\modxs_{\C}(z) \ 
&= \  \modxs_{\R}(N_{\C/\R}(z)) \\
&= \ \abs{z\ov{z}} \\
&=\  \abs{z}^2.
\end{aligned}
\]
\end{x}

\vspace{0.1cm}

\begin{x}{\small\bf LEMMA} \ 
\[
\modxs_{\Q_p} = \acdot_{p}
\]
To prove this we need a preliminary.
\end{x}

\begin{x}{\small\bf LEMMA} \ 
The arrow
\[
\epsilon_k : \Z_p  \ra \Z/p^k\Z
\]
that sends
\[
x = \sum_{n=0}^\infty a_np^n		\qquad(a_{n}  \in \sA)
\]
to
\[
\sum_{n=0}^{k-1} a_np^n \modxs p^k
\]
is a homomorphism of rings.  
It is surjective with kernel $p^k\Z_p$, so $[\Z_p:p^k\Z_p] =p^k$ 
(cf. $\S 4$,   $\#26$), thus there is a disjoint decomposition of $\Z_p:$
\[
\Z_p = \bigcup\limits_{j=1}^{p^k} (x_j + p^k\Z_p).
\]

Normalize the Haar measure on $\Q_p$ by stipulating that
\[
\mu_{\Q_p} (\Z_p) = 1.
\]

\vspace{0.1cm}

[Note: In this connection, recall that $\Z_p$ is an open-compact set.$]$

\vspace{0.2cm}

The claim now is that for every Borel set $A$,
\[
\mu_{\Q_p} (M_x(A)) = \abs{x}_p\mu_{\Q_p} (A).
\]
Since the Borel $\sigma$-algebra is generated by the open sets, it is enough to take $A$ open.  
But any open set can be written as the disjoint union of cosets of the subgroups 
$p^k\Z_p$ $($cf. $\S 4$, $\# 33$), hence thanks to translation invariance, it suffices to deal with these alone:
\[
\begin{aligned}
\mu_{\Q_p} (p^k\Z_p) \ 
&= \ \modxs_{\Q_p}(p^k)\mu_{\Q_p}(\Z_p) \\
&= \ \modxs_{\Q_p}(p^k) \\
&= \ |p^k|_p.
\end{aligned}
\]

1. $k \ge 0$:
\[
\begin{aligned}
1  \
&= \ \mu_{\Q_p}(\Z_p)  \\
&= \ \mu_{\Q_p}(\bigcup\limits_{j=1}^{p^k} (x_j + p^k\Z_p)) \\
&= \ p^k\mu_{\Q_p}(p^k\Z_p)
\end{aligned}
\]
\qquad\qquad\qquad\qquad$\implies$
\[
\begin{aligned}
\mu_{\Q_p}(p^k\Z_p) \ 
&= \ p^{-k} \\
&= \ |p^k|_p.
\end{aligned}
\]

2. $k < 0$:
\[
\begin{aligned}
1 \ 
&=\  \mu_{\Q_p}(\Z_p) \\
&= \ \mu_{\Q_p}(p^{-k}p^k\Z_p)\\
&= \ \modxs_{\Q_p} (p^{-k}) \mu_{\Q_p} (p^k\Z_p)\\
&= \ |p^{-k}|_p \mu_{\Q_p}(p^k\Z_p)
\end{aligned}
\]
\qquad\qquad\qquad\qquad$\implies$
\[
\begin{aligned}
\mu_{\Q_p}(p^k\Z_p) \ 
&= \ |p^{-k}|_p^{-1} \\
&= \ |p^k|_p.
\end{aligned}
\]
\end{x}

\vspace{0.1cm}

\begin{x}{\small\bf SCHOLIUM}  \ 
If $\K$ is a finite field extension of $\Q_p$, then $\forall$ a $\in \K$, 
\[
\modxs_\K(a) = \abs{N_{\K/ \Q_p}(a)}_{p},
\]
the normalized absolute value on $\K$ mentioned in $\S$ 5:
\[
\modxs_\K(a) = \abs{a}_\K \quad (= \abs{a}_p^n, \ n = [\K:\Q_p]).
\]
\end{x}

\vspace{0.1cm}

\begin{x}{\small\bf CONVENTION}  \ 
Integration w.r.t. $\mu_{\Q_p}$ will be denoted by dx:
\[
\int_{\Q_p}f(x)d\mu_{\Q_p}(x) = \int_{\Q_p}f(x)dx.
\]

[Note: Points are of Haar measure zero:
\[
\{0\} = \bigcap\limits_{k=1}^{\infty} p^k\Z_p
\]
\qquad\qquad\qquad\qquad$\implies$
\[
\begin{aligned}
\mu_{\Q_p} (\{0\}) \ 
&= \  \lim_{k \ra \infty} \mu_{\Q_p}(p^k\Z_p)\\
&= \  \lim_{k\ra \infty} p^{-k} = 0.]
\end{aligned}
\]
\end{x}

\vspace{0.1cm}

\begin{x}{\small\bf EXAMPLE} \ 
\[
\Z_p^\times = \bigcup\limits_{1 \le k \le p-1} (k+ p\Z_p)	\qquad \text{ $($cf. $\S 4$, $\#23$)}.
\]
Therefore
\[
\begin{aligned}
\vol_{dx}(\Z_p^\times)  \
&= \ (p-1)\vol_{dx}(p\Z_p ) \\
&= \ \frac{p-1}{p}.
\end{aligned}
\]
\end{x}

\vspace{0.1cm}

\begin{x}{\small\bf EXAMPLE}  \ 
\[
\begin{aligned}
\vol_{dx}(p^n\Z_p^\times)  \
&= \ \vol_{dx}(p^n\Z_p - p^{n+1}\Z_p)		\qquad \text{(cf.} \  \S 4, \ \# 34)\\
&= \ \vol_{dx}(p^n\Z_p) - \vol_{dx}(p^{n+1}\Z_p)	\\	
&= \ \abs{p^n}_p \vol_{dx}(\Z_p) - \abs{p^{n+1}}_p \vol_{dx}(\Z_p)\\
&= \ p^{-n} - p^{-n-1}.
\end{aligned}
\]
\end{x}

\vspace{0.1cm}

\begin{x}{\small\bf EXAMPLE}  \ 
Write
\[
\Z_p - \{0\} = \bigcup\limits_{n \ge 0}p^n \Z_p^\times.
\]
Then
\begin{align*}
\allowdisplaybreaks
\int_{\Z_p - \{0\}} \log \abs{x}_p dx \ 
&= \ \sum_{n=0}^\infty \  \int_{p^n \Z_p^\times} \log \abs{x}_p dx\\
&= \ \sum_{n=0}^\infty \log p^{-n} \vol_{dx}(p^n\Z_p^\times)\\
&= \ -\log p \ \sum_{n=0}^\infty n(p^{-n} - p^{-n-1})\\
&= \ -\log p \  \bigl(\sum_{n=0}^\infty \frac{n}{p^n} - \frac{1}{p} \sum_{n=0}^\infty \frac{n}{p^n}\bigr)\\
&= \ - (1 - \frac{1}{p}) \log p \ \sum_{n=0}^\infty \frac{n}{p^n} \\
&= \ - (1 - \frac{1}{p}) \log p \ \frac{p}{(p-1)^2}\\
&= \ - \frac{\log p}{p-1}.
\end{align*}
\end{x}


\begin{x}{\small\bf EXAMPLE}   \ 
\[
\int_{\Z_p^\times} \log \abs{1 - x}_p dx = -\frac{\log p}{p-1}.
\]

[Break $\Z_p^\times$ up via the scheme
\[
(\Z_p^\times:a_0 \ne 1) \cup (\Z_p^\times:a_0 = 1, a_1 \ne 0) \cup (\Z_p^\times:a_0 = 1, a_1 = 0, a_2 \ne 0) \cup \dotsb .]
\]
\end{x}


\begin{x}{\small\bf LEMMA} \ 
The measure $\ds\frac{dx}{\abs{x}_p}$ is a Haar measure on the multiplicative group $\Q_p^\times$.


PROOF \  $\forall \  y \in \Q_p^\times$,

\allowdisplaybreaks
\begin{align*}
\allowdisplaybreaks
\int_{\Q_p^\times} f(y^{-1}x) \frac{dx}{\abs{x}_p} \ 
&= \ \abs{y}_p^{-1} \int_{\Q_p^\times} f(y^{-1}x) \frac{1}{\abs{y^{-1}x}_p} dx\\
&= \ \abs{y}_p^{-1} \modxs_{\Q_p}(y) \int_{\Q_p^\times} f(x) \frac{dx}{\abs{x}_p}\\
&= \ \abs{y}_p^{-1} \abs{y}_p \int_{\Q_p^\times} f(x) \frac{dx}{\abs{x}_p}\\
&= \ \int_{\Q_p^\times} f(x) \frac{dx}{\abs{x}_p}.
\end{align*}
\end{x}
\vspace{0.1cm}

\begin{x}{\small\bf EXAMPLE} \ 
\[
\allowdisplaybreaks
\begin{aligned}
\ds\vol_{\frac{dx}{\abs{x}_p}}(p^n\Z_p^\times) \
&= \ \vol_{\frac{dx}{\abs{x}_p}}(\Z_p^\times)\\
&= \ \int_{\Z_p^\times} \frac{dx}{\abs{x}_p} \\
&= \ \int_{\Z_p^\times} dx \\
&= \ \vol_{dx}(\Z_p^\times) \\
&= \ \frac{p-1}{p}.
\end{aligned}
\]
\end{x}

\vspace{0.1cm}

\begin{x}{\small\bf DEFINITION} \ 
The 
\underline{normalized Haar measure}
\index{normalized Haar measure} 
on the multiplicative group $\Q_p^\times$ is given by
\[
d^\times x = \frac{p}{p-1}\frac{dx}{\abs{x}_p}.
\]
Accordingly,
\[
\vol_{d^\times x} (\Z_p^\times) = 1,
\]
this condition characterizing $d^\times x$.
\end{x}

\vspace{0.1cm}

\begingroup
\allowdisplaybreaks
\begin{x}{\small\bf EXAMPLE} 
Let $s$ be a complex variable with $\Re(s) > 1$.  Write 
\[
\Z_p - \{0\} = \bigcup\limits_{n \geq 0} p^n \Z_p^\times.
\]
Then
\[
\allowdisplaybreaks
\begin{aligned}
\allowdisplaybreaks
\int_{\Z_p-\{0\}} \abs{x}_p^s d^\times x  \
&= \ \sum_{n = 0}^\infty p^{-ns} \int_{\Z_p^\times} d^\times x \\
&= \ \sum_{n = 0}^\infty p^{-ns} \\
&= \ \frac{1}{1 - p^{-s}},
\end{aligned}
\]
the $p^{th}$ factor in the Euler product for the Riemann zeta function.

\vspace{0.1cm}

Let $\K/\Q_p$ be a finite extension. 
Given a Haar measure $da$ on $\K$, put
\[
d^\times a = \frac{q}{q-1} \frac{da}{\abs{a}_\K}.
\]
Then $\ds\frac{da}{\abs{a}_\K}$ is a Haar measure on $\K^\times$ and we have
\allowdisplaybreaks
\begin{align*}
\allowdisplaybreaks
\vol_{d^\times a} (R^\times) \ 
&= \ \int_{R^\times} \frac{q}{q-1} \frac{da}{\abs{a}_\K}\\
&= \ \frac{q}{q-1} \int_{R^\times}  da\\
&= \ \sum_{n = 0}^\infty q^{-n} \int_{R^\times}  da\\
&= \ \sum_{n = 0}^\infty \int_{R^\times}  q^{-n} da\\
&= \ \sum_{n = 0}^\infty \int_{\pi^nR^\times}  da\\
&= \ \int_{\bigcup\limits_{n \ge 0}\pi^nR^\times}  da\\
&= \ \int_R da \\
&= \  \vol_{da}(R).
\end{align*}
\end{x}
\endgroup


\chapter{
$\boldsymbol{\S}$\textbf{7}.\quad  HARMONIC ANALYSIS}
\setlength\parindent{2em}
\setcounter{theoremn}{0}

\ \indent 
Let $G$ be a locally compact abelian group.

\vspace{0.25cm}

\begin{x}{\small\bf DEFINITION} \ 
A 
\un{character}
\index{character} 
of $G$ is a continuous homomorphism $\chi:G \ra \C^\times$.
\end{x}

\vspace{0.1cm}

\begin{x}{\small\bf NOTATION} \ 
Write $\widetilde{G}$ for the group whose elements are the characters of $G$.
\end{x}

\vspace{0.1cm}

\begin{x}{\small\bf DEFINITION} \ 
A 
\un{unitary character}
\index{unitary character} 
of $G$ is a continuous homomorphism $\chi: G \ra \T$.
\end{x}

\vspace{0.1cm}

\begin{x}{\small\bf NOTATION} \ 
Write $\widehat{G}$ for the group whose elements are the unitary characters of $G$.
\end{x}

\vspace{0.1cm}

\begin{x}{\small\bf LEMMA} \ 
There is a decomposition
\[
\widetilde{G} \approx \widetilde{G}_+ \times \widehat{G},
\]
where $\widetilde{G}_+$ is the group of positive characters of $G$.

\vspace{0.1cm}

PROOF \  The only positive unitary character is trivial, so $\widetilde{G}_+ \cap \widehat{G} = \{1\}$.  On the other hand, if $\chi$ is a character, then $\abs{\chi}$ is a positive character, $\chi / \abs{\chi}$ is a unitary character, and $\chi = \abs{\chi} (\ds\frac{\chi}{\abs{\chi}}$).
\end{x}

\vspace{0.1cm}

\begin{x}{\small\bf LEMMA} \ 
Every bounded character of $G$ is a unitary character.

\vspace{0.1cm}

PROOF \  The only compact subgroup of $\R_{>0}$ is the trivial subgroup $\{1\}$.
\end{x}

\vspace{0.1cm}

\begin{x}{\small\bf APPLICATION} \ 
If $G$ is compact, then every character of $G$ is unitary.
\end{x}

\begin{x}{\small\bf EXAMPLE} \ 
Take $G = \Z$ $-$then $\widetilde{G} \approx \C^\times$, the isomorphism being given by the map $\chi \ra \chi (1)$.
\end{x}

\vspace{0.1cm}

\begin{x}{\small\bf EXAMPLE} \ 
Take $G = \R$ $-$then $\widetilde{G} \approx \R \times \R $ and every character has the form 
$\chi (x) = e^{zx}$  $(z \in \C)$.
\end{x}

\vspace{0.1cm}

\begin{x}{\small\bf EXAMPLE} \ 
Take $G = \C$ $-$then $\widetilde{G} \approx \C \times \C $ and every character has the form 
$\chi (x) = \exp(z_1 \Re(x) + z_2 \Im(x))$  \ $(z_1, z_2 \in \C)$.
\end{x}
\vspace{0.1cm}

\begin{x}{\small\bf EXAMPLE} \ 
Take $G = \R^\times$ $-$then $\widetilde{G} \approx \Z/ 2 \Z \times \C$, and every character has the form 
$\chi (x) = (\sgn x)^\sigma \abs{x}^s$ $(\sigma \in \{0,1\}, \ s \in \C)$.
\end{x}
\vspace{0.1cm}

\begin{x}{\small\bf EXAMPLE} \ 
Take $G = \C^\times$ $-$then $\widetilde{G} \approx \Z \times \C$, 
and every character has the form $ \chi (x) = \exp(\sqrt{-1} $ $n \arg x) \abs{x}^s$ $(n \in \Z, s \in \C)$.
\end{x}

\vspace{0.1cm}

\begin{x}{\small\bf DEFINITION} \ 
The 
\un{dual group}
\index{dual group} 
of $G$ is $\widehat{G}$.
\end{x}

\vspace{0.1cm}

\begin{x}{\small\bf RAPPEL} \ 
Let $X$, $Y$ be topological spaces and let $F$ be a subspace of $C(X,Y)$.  
Given a compact set $K \subset X$ and an open subset $V \subset Y$, 
let $W(K,V)$ be the set of all $f \in F$ such that $f(K) \subset V$ 
$-$then the collection $\{W(K,V)\}$ is a subbasis for the 
\un{compact open topology}
\index{compact open topology} 
on $F$.

\vspace{0.1cm}
$[$Note: The family of finite intersections of sets of the form $W(K,V)$ is then a basis for the compact open topology: 
Each member has the form $\bigcap\limits_{i= 1}^n W(K_i,V_i)$, where the $K_i \subset X$ 
are compact and the $V_i \subset Y$ are open.]
\end{x}

\vspace{0.1cm}

Equip $\widehat{G}$ with the compact open topology.
\vspace{0.1cm}
\begin{x}{\small\bf FACT} \ 
The compact open topology on $\widehat{G}$ coincides with the topology of uniform convergence on compact subsets of $G$.
\end{x}

\begin{x}{\small\bf LEMMA} \ 
$\widehat{G}$ is a locally compact abelian group.
\end{x}

\begin{x}{\small\bf REMARK} \ 
$\widetilde{G}$  is also a locally compact abelian group and the decomposition
\[
\widetilde{G} \approx \widetilde{G}_+ \times \widehat{G}
\]
is topological.
\end{x}
\vspace{0.1cm}

\begin{x}{\small\bf EXAMPLE} \ 
Take $G = \R$ and given a real number $t$, let $\chi_t(x) = e^{\sqrt{-1}\ tx}$ $-$then $\chi_t$ 
is a unitary character of $G$ and for any $\chi \in \widehat{G}$, there is a unique $t \in \R$ such that $\chi = \chi_t$, 
hence $G$ can be identified with $\widehat{G}$.
\end{x}

\vspace{0.1cm}

\begin{x}{\small\bf EXAMPLE} \ 
Take $G = \R^2$ and given a point $(t_1,t_2)$, let $\chi_{(t_1,t_2)}(x_1,x_2) = e^{\sqrt{-1}(t_1x_1 + t_2x_2)}$ 
$-$then $\chi_{(t_1,t_2)}$ is a unitary character of $G$ and for any $\chi \in \widehat{G}$, 
there is a unique $(t_1,t_2) \in \R^2$ such that $\chi = \chi_{(t_1,t_2)}$, 
hence $G$ can be identified with $\widehat{G}$.
\end{x}

\vspace{0.1cm}

\begin{x}{\small\bf EXAMPLE} \ 
Take $G = \Z/ n\Z$ and given an integer $m = 0, 1, \dotsb, n-1$, let 
$
\chi_m(k) = \exp\bigl(2\pi \sqrt{-1} \ \ds\frac{km}{n}\bigr)
$
$-$then $\chi_0, \chi_1, \dotsb, \chi_{n-1}$ are characters of $G$, hence $G$ can be identified with $\widehat{G}$.
\end{x}

\vspace{0.1cm}

\begin{x}{\small\bf LEMMA} \ 
If $G$ is compact, then $\widehat{G}$ is discrete.
\end{x}

\vspace{0.1cm}

\begin{x}{\small\bf EXAMPLE} \ 
Take $G = \T$ and given $n \in$ $\Z$, let 
$\chi_n\bigl(e^{\sqrt{-1}\ \theta}\bigr ) = e^{\sqrt{-1} \ n\theta}$ 
$-$then $\chi_n$ is a unitary character of $G$ and all such have this form, so $\T$ $\approx \Z$.
\end{x}

\vspace{0.1cm}

\begin{x}{\small\bf LEMMA} \ 
If $G$ is discrete, then $\widehat{G}$ is compact.
\end{x}

\vspace{0.1cm}

\begin{x}{\small\bf EXAMPLE} \ 
Take $G = \Z$ and given $e^{\sqrt{-1}\ \theta}  \in \T$, let $\chi_\theta(n) = e^{\sqrt{-1} \ \theta n}$ 
$-$then $\chi_\theta$ is unitary character of $G$ and all such have this form, so $\widehat{\Z} \approx \T$.
\end{x}

\vspace{0.1cm}

\begin{x}{\small\bf LEMMA} \ 
If $G_1$, $G_2$ are locally compact abelian groups, then 
$\reallywidehat{G_1 \times G_2}$
is topologically isomorphic to $\widehat{G_1} \times \widehat{G_2}$.
\end{x}

\vspace{0.1cm}

\begin{x}{\small\bf EXAMPLE} \ 
Take $G = \R^\times$ $-$then 
\[
G \approx \Z / 2 \Z \times \R_{>0}^\times \ \approx \ \Z / 2 \Z \times \R, 
\]
thus $\widehat{G}$ is topologically isomorphic to $\Z / 2 \Z \times \R:$
\[
(u,t) \ra \chi_{(u,t)}			\quad (u \in \Z / 2 \Z, t \in  \R),
\]
where
\[
\chi_{(u,t)}(x) = \left (\ds\frac{x}{\abs{x}}\right)^u \abs{x}^{\sqrt{-1}\  t}.
\]
\end{x}

\vspace{0.2cm}

\begin{x}{\small\bf EXAMPLE} \ 
Take $G = \C^\times$ $-$then 
\[
G \approx \T \times \R_{>0}^\times \ \approx \T \ \times \R, 
\]
thus $\widehat{G}$ is topologically isomorphic to $\Z \times \R:$
\[
(n,t) \ra \chi_{n,t}			\quad (n \in \Z , t \in  \R),
\]
where
\[
\chi_{(n,t)}(z) = \left(\ds\frac{z}{\abs{z}}\right)^n \abs{z}^{\sqrt{-1} \  t}.
\]
\end{x}

\vspace{0.5cm}

Denote by $\ev_G$ the canonical arrow 
$G \ra \reallywidehat{\reallywidehat{G}}:$
\[
\ev_G(x) (\chi) = \chi (x).
\]

\vspace{0.1cm}

\begin{x}{\small\bf REMARK} \ 
If $G$, $H$ are locally compact abelian groups and if $\phi : G \ra H$ is a continuous homomorphism, then there is a commutative diagram
\[
\begin{tikzcd}[sep=large]
G \ar{d}[swap]{\phi} \ar{rr}{\text{ev}_G}	
&&\reallywidehat{\reallywidehat{G}} \ar{d}{\reallywidehat{\reallywidehat{\phi}}}\\
H \ar{rr}{\text{ev}_H} 		
&&\reallywidehat{\reallywidehat{H}}
\end{tikzcd}
\quad .
\]
\end{x}

\vspace{0.1cm}

\index{Pontryagin duality}
\begin{x}{\small\bf PONTRYAGIN DUALITY} \ 
$\ev_G$ is an isomorphism of groups and a homeomorphism of topological spaces.
\end{x}

\vspace{0.1cm}

\begin{x}{\small\bf SCHOLIUM} \ 
Every compact abelian group is the dual of a discrete abelian group and every discrete abelian group is the dual of a compact abelian group.
\end{x}

\vspace{0.1cm}

\begin{x}{\small\bf REMARK} \ 
Every finite abelian group $G$ is isomorphic to its dual $\widehat{G}:$ $G \approx \widehat{G}$  (but the isomorphism is not "functorial").
\end{x}

\vspace{0.1cm}

Let $H$ be a closed subgroup of $G$.

\vspace{0.2cm}

\index{$H^\perp$}
\begin{x}{\small\bf NOTATION} \ 
Put
\[
H^\perp = \{\chi \in \widehat{G}: \restr{\chi}{H} = 1\}.
\]
\end{x}
\vspace{0.1cm}

\begin{x}{\small\bf LEMMA} \ 
$H^\perp$ is a closed subgroup of $\widehat{G}$ and H = $H^{\perp \perp}$.
\end{x}

\vspace{0.1cm}

Let $\pi_H:G \ra G/H$ be the projection and define
\[\left\{
\begin{array}{l l}
\Phi : \widehat{G/H} \ra H^\perp\\
\Psi : \widehat{G}/H^\perp \ra \widehat{H}
\end{array}
\right.\]
by
\[\left\{
\begin{array}{l l}
\Phi(\chi) = \chi \circ \pi_H \\
\Psi(\chi H^\perp) = \restr{\chi}{H}.
\end{array}
\right.\]

\vspace{0.1cm}

\begin{x}{\small\bf LEMMA} \ 
$\Phi$ and $\Psi$ are isomorphisms of topological groups.
\end{x}

\vspace{0.1cm}

\begin{x}{\small\bf APPLICATION} \ 
Every unitary character of $H$ extends to a unitary character of $G$.
\end{x}

\vspace{0.1cm}

\begin{x}{\small\bf EXAMPLE} \ 
Let $G$ be a finite abelian group and let $H$ be subgroup of $G$ $-$then $G$ contains a subgroup isomorphic to $G/H$.

[In fact, 
\[
G/H  \approx \  \widehat{G/H}\  \approx \  H^\perp \ \subset \  \widehat{G} \ \approx \ G.]
\]
\end{x}

\vspace{0.1cm}

\begin{x}{\small\bf REMARK} \ 
Denote by $\textbf{LCA}$ the category whose objects are the locally compact abelian groups and whose morphisms are the continuous homomorphisms $-$then
\[
\widehat{ } : \mathbf{LCA} \ra \mathbf{LCA}
\]
is a contravariant functor.  This said, consider the short exact sequence
\[
\begin{tikzcd}
1 	\ar{r}	
&H	\ar{r}	
&G 	\ar{rr}{\pi_H}	
&&{G/H}	\ar{r}  
&{1}
\end{tikzcd}
\]
and apply  $ $   $\widehat{ }:$
\[
\begin{tikzcd}
1 	\ar{r}	
&{\widehat{G/H} \  \approx \ H^\perp} 	\ar{r}	
&\widehat{G} 	\ar{r}	
&{\widehat{H} \ \approx \  \widehat{G}/H^\perp}	\ar{r}  
&{1}
\end{tikzcd}
,
\]
which is also a short exact sequence.
\end{x}

\vspace{0.1cm}

Given f $\in L^1(G)$, its 
\un{Fourier transform}
\index{Fourier transform} 
is the function 
\[
\widehat{f} : \widehat{G} \ra \C
\]
defined by the rule
\[
\widehat{f}(\chi) = \int_G f(x) \chi(x) d \mu_G(x).
\]

\vspace{0.1cm}

\begin{x}{\small\bf EXAMPLE} \ 
Take $G = \R$ $-$then $\widehat{\R} \approx \R$ and
\[
\widehat{f}(\chi_t) \equiv \widehat{f}(t) = \int_{-\infty}^\infty  f(x) e ^{\sqrt{-1}\  tx} dx.
\]
\end{x}

\vspace{0.1cm}

\begin{x}{\small\bf EXAMPLE} \ 
Take $G = \R^2$ $-$then $\widehat{\R}^2 \approx \R^2$ and
\[
\widehat{f}(\chi_{(t_1,t_2)}) \equiv \widehat{f}(t_1 , t_2) \ 
= \ \int_{-\infty}^\infty  \int_{-\infty}^\infty  f(x_1, x_2) e ^{\sqrt{-1} \ (t_1 x_1 + t_2 x_2)} dx_1 dx_2.
\]
\end{x}

\vspace{0.1cm}

\begin{x}{\small\bf EXAMPLE} \ 
Take $G = \T$ $-$then $\widehat{\T} \approx \Z$ and
\[
\widehat{f}(\chi_n) \equiv \widehat{f}(n) 
=  \int_{0}^{2\pi}  f(\theta) e^{\sqrt{-1} \ n \theta} d\theta
\]
\end{x}

\vspace{0.1cm}

\begin{x}{\small\bf EXAMPLE} \ 
Take $G = \Z$ $-$then $\widehat{Z} \approx \T$ and
\[
\widehat{f}(\chi_\theta) \equiv \widehat{f}(\theta) =  \sum_{n = -\infty}^\infty f(n) e^{\sqrt{-1} \ n \theta}.
\]
\end{x}

\vspace{0.1cm}

\begin{x}{\small\bf EXAMPLE} \ 
Take $G = \Z/n\Z$ $-$then $\widehat{\Z/n\Z} \approx \Z/n\Z$ and
\[
\widehat{f}(\chi_m) \equiv \widehat{f}(m) =  \sum_{k = 0}^{n - 1}f(k) \exp(2 \pi \sqrt{-1} \ \frac{km}{n}). 
\]
\end{x}

\vspace{0.1cm}

\begin{x}{\small\bf LEMMA} \ 
$\widehat{f} : \widehat{G} \ra \C$ is a continuous function on $\widehat{G}$ that vanishes at infinity and
\[
\| \widehat{f} \|_\infty \le \| f  \|_1.
\]
\end{x}

\vspace{0.1cm}

\index{$\mathbf{INV}$}
\begin{x}{\small\bf NOTATION} \ 
$\mathbf{INV}(G)$ is the set of continuous functions $f \in L^1(G)$ with the property that $\widehat{f} \in L^1(\widehat{G})$.
\end{x}

\vspace{0.1cm}

\index{Fourier inversion}
\begin{x}{\small\bf FOURIER INVERSION} \ 
Given a Haar measure $\mu_G$ on $G$, 
there exists a unique Haar measure $\mu_{\widehat{G}}$ on $\widehat{G}$ such that $\forall f \  \in \mathbf{INV}(G)$,
\[
f(x) = \int_{\widehat{G}} \widehat{f}(\chi) \ov{\chi (x)} d\mu_{\widehat{G}} (\chi).
\]
\end{x}

\vspace{0.1cm}

If $G$ is compact, then it is customary to normalize $\mu_G$ by the requirement $\int_G 1 d\mu_G = 1$.

\vspace{0.1cm}

\begin{x}{\small\bf LEMMA} \ 
\[ 
\int_G \chi(x) d\mu_G(x) \ = \  
\begin{cases}
\ 1  \quad \text{if} \  \chi = 0\\
\ 0 \quad \text{if} \ \chi \ne 0
\end{cases}
.
\]

\vspace{0.1cm}

PROOF \  The case $\chi = 0$ is clear.  On the other hand, if $\chi \ne 0$, then there exists $x_0 : \chi (x_0) \ne 1$, hence
\begin{align*}
\int_G \chi(x) d\mu_G(x) \ 
&= \ \int_G \chi (x - x_0 + x_0) d\mu_G (x)\\
&= \ \chi (x_0) \int_G \chi (x - x_0) d\mu_G (x)\\
&= \ \chi (x_0) \int_G \chi (x) d\mu_G (x)
\end{align*}
$\implies$
\[
\int_G \chi (x) d\mu_G (x) = 0 .
\]
\end{x}

\vspace{0.1cm}

Assuming still that $G$ is compact ($\implies \widehat{G}$ is discrete), take $f \equiv 1:$
\[
\widehat{f}(0) = 1, \   \widehat{f}(\chi) = 0 	 \quad (\chi \ne 0). 
\]
I.e.$:$	$\widehat{f}$ is the characteristic function of \{0\}, hence is integrable, thus $f \in \mathbf{INV}(G)$.  
Accordingly, if $\mu_{\widehat{G}}$ is the Haar measure on $\widehat{G}$ per Fourier inversion, then
\begin{align*}
1 \ 
&= \ f(0) \\
&= \ \int_{\widehat{G}} \widehat{f}(\chi) d\mu_{\widehat{G}} (\chi) \\
&= \ \mu_{\widehat{G}} (\{0\}), 
\end{align*}
so $\forall$ $\chi \in \widehat{G}$,
\[
\mu_{\widehat{G}} (\{\chi\}) = 1.
\]

\vspace{0.1cm}

\begin{x}{\small\bf EXAMPLE} \ 
Let $G = \T$ $-$then $d\mu_G = \ds\frac{d\theta}{2\pi}$, so for $f \in \mathbf{INV}(G)$,
\[
f(\theta) = \sum_{n = -\infty}^\infty \widehat{f}(n) e^{-\sqrt{-1} \ n \theta},
\]
where
\[
\widehat{f}(n) = \int_0^{2\pi} f(\theta) e^{\sqrt{-1} \  n \theta} \frac{d\theta}{2\pi}.
\]

\vspace{0.1cm}

If $G$ is discrete, then it is customary to normalize $\mu_G$ by stipulating that singletons are assigned measure 1.
\end{x}

\vspace{0.1cm}

\begin{x}{\small\bf REMARK} \ 
There is a conflict if $G$ is both compact and discrete, i.e., if $G$ if finite.
\end{x}

\vspace{0.1cm}

\allowdisplaybreaks
Assuming still that $G$ is discrete ( $\implies \widehat{G}$ is compact), take $f(0) = 1, f(x) = 0$ $(x \ne 0):$
\begin{align*}
\widehat{f}(\chi)	 \ 
&= \  \int_G f(x) \chi(x) d\mu_G(x)\\	
&= \ f(0) \chi(0) \mu_G(\{0\})\\		
&= \ 1.
\end{align*}
I.e.$:$ $\widehat{f}$ is the constant function 1, hence is integrable, thus $f \in \mathbf{INV}(G)$.  Accordingly, if $\mu_{\widehat{G}}$ is the Haar measure on $\widehat{G}$ per Fourier inversion, then
\begin{align*}
\mu_{\widehat{G}}(\widehat{G}) 	 \ 
&= \ \int_{\widehat{G}} 1 d \mu_{\widehat{G}}(\chi)  \\		 						
&= \  \int_{\widehat{G}} \widehat{f}(\chi) d \mu_{\widehat{G}}(\chi)  \\
&= \ \int_{\widehat{G}} \widehat{f}(\chi) \chi(0)d \mu_{\widehat{G}}(\chi)  \\	
&= \ f(0)\\
&=	\ 1.
\end{align*}

\begin{x}{\small\bf EXAMPLE} \ 
Take $G = \Z / n\Z$ and let $\mu_G$ be the counting measure $($thus here $\mu_G(G) = n)$ $-$then $\mu_{\widehat{G}}$ is the counting measure divided by $n$ and for $f \in \mathbf{INV}(G)$, 
\[
f(k) = \frac{1}{n} \sum_{m = 0} ^{n - 1} \widehat{f}(m) \exp(-2\pi\sqrt{-1} \ \frac{km}{n}),
\]
where
\[
\widehat{f}(m) = \sum_{k = 0} ^{n - 1} f(k) \exp(2\pi\sqrt{-1} \ \frac{km}{n}).
\]
\end{x}

\begin{x}{\small\bf EXAMPLE} \ 
Take $G = \R$ and let $\mu_G = \alpha dx$ $(\alpha > 0)$, hence $\mu_{\widehat{G}} = \beta dt $ $(\beta > 0)$ 
and we claim that
\[
1 = 2\alpha \beta \pi.
\]
To establish this, recall first that the formalism is
\[
\begin{cases}
\widehat{f}(t) 	&= \ds\int_{-\infty}^\infty f(x) e^{\sqrt{-1}\  tx} \alpha dx \\	
\\
f(x) 	&= \ds\int_{-\infty}^\infty \widehat{f}(t) e^{-\sqrt{-1} \ tx} \beta dx 
\end{cases}
.
\]
Let $f(x) = e^{-\abs{x}}$ $-$then
\[
\frac{2\alpha}{1 + t^2} = \int_{-\infty}^\infty e^{-\abs{x}} e^{\sqrt{-1}\  tx} \alpha dx,
\]
so $f \in \mathbf{INV}(G)$, thus
\[
\begin{aligned}
e^{-\abs{x}} \ 
&=\  \int_{-\infty}^\infty \frac{2\alpha}{1 + t^2} e^{-\sqrt{-1} \ tx} \beta dt \\
&= \  2\alpha \beta \int_{-\infty}^\infty \frac{e^{-\sqrt{-1}\  tx }}{1 + t^2}dt.
\end{aligned}
\]
Now put $x = 0:$
\[
1 = 2\alpha \beta \int_{-\infty}^\infty \frac{dt}{1 + t^2}dt = 2\alpha \beta  \pi,
\]
as claimed.  One choice is to take
\[
\alpha = \beta = \frac{1}{\sqrt{2\pi}},
\]
the upshot being that the Haar measure of $[0,1]$ is not 1 but rather $\ds\frac{1}{\sqrt{2\pi}}$.
\end{x}

\vspace{0.1cm}

\begin{x}{\small\bf NOTATION} \ 
Given $f \in L^1(\R)$, let
\[
\sF_\R f(t) = \int_{-\infty}^\infty f(x) e^{2\pi \sqrt{-1}\   tx} dx.
\]
Therefore
\begin{align*}
\sF_\R f(t) \ 
&=\  \sqrt{2\pi}\   \frac{1}{\sqrt{2\pi}} \int_{-\infty}^\infty f(x) e^{2\pi \sqrt{-1}\   tx} dx \\
&=\   \sqrt{2\pi}\   \widehat{f}(2 \pi t).
\end{align*}
\end{x}

\vspace{0.1cm}

\begin{x}{\small\bf STANDARDIZATION} \ 
$(G = \R)$ Let $f \in \mathbf{INV}(\R)$, $-$then
\[
\sF_\R \sF_\R f(x) = f(-x).
\]

[In fact,
\begin{align*}
\sF_\R \sF_\R f(x)	 \ 
&=\  \int_{-\infty}^\infty \sF_\R f(t) e^{2 \pi \sqrt{-1}\  tx} dx \\	
&=\  \ \int_{-\infty}^\infty \sqrt{2\pi} \widehat{f}(2 \pi t) e^{2 \pi \sqrt{-1}\  tx} dx \\	
&=\  \ \sqrt{2\pi}  \int_{-\infty}^\infty  \widehat{f}(u) e^{\sqrt{-1}\  ux} \frac{du}{2\pi} \\	
&=\  \ \frac{1}{\sqrt{2\pi}}  \int_{-\infty}^\infty  \widehat{f}(t) e^{\sqrt{-1}\  tx} dt \\	
&=\  \ f(-x).]
\end{align*}

\vspace{0.1cm}

Fourier inversion in the plane takes the form
\[
\begin{cases}
\widehat{f}(t_1, t_2)  \
= \ \ds\frac{1}{2\pi} \int_{-\infty}^\infty \int_{-\infty}^\infty f(x_1, x_2) e^{\sqrt{-1}\  (t_1x_1 + t_2x_2)} dx_1 dx_2\\
\\
f(x_1, x_2)  \ 
= \ \ds\frac{1}{2\pi} \int_{-\infty}^\infty \int_{-\infty}^\infty \widehat{f}(t_1, t_2) e^{-\sqrt{-1}\  (t_1x_1 + t_2x_2)} dt_1 dt_2
\end{cases}
.
\]
One may then introduce
\begin{align*}
\sF_{\R^2}f(t_1,t_2) \ 
&= \ \int_{-\infty}^\infty \int_{-\infty}^\infty f(x_1, x_2) e^{2\pi\sqrt{-1}\  (t_1x_1 + t_2x_2)} dx_1 dx_2 \\
&= \ 2\pi \widehat{f} (2\pi t_1, 2\pi t_2)
\end{align*}
and proceeding as above we find that
\[
\sF_{\R^2} \sF_{\R^2} f(x_1, x_2) = f(-x_1, -x_2).
\]

Now identify $\R^2$ with $\C$ and recall that $\trs_{\C /\R} (z) = z + \bar{z}$.  
Write
\[
\begin{cases}
w = a + \sqrt{-1} \ b\\
z = x + \sqrt{-1}  \ y
\end{cases}
.
\]

Then
\[
wz + \ov{wz} = 2\Re(wz) = 2(ax -by). 
\]
Therefore
\[
\frac{1}{2\pi} \int_{-\infty}^\infty \int_{-\infty}^\infty f(x,y) e^{2 \sqrt{-1}\ (ax-by)}dx dy = \widehat{f}(2a, -2b).
\]

\vspace{0.1cm}

[Note: \  Let $\chi_w(z) = \exp(\sqrt{-1}(wz + \ov{wz}))$ 
$-$then $\chi_w$ is a unitary character of $\C$ and for any $\chi \in \widehat{\C}$, 
there is a unique $w \in \C$ such that $\chi = \chi_w$, hence $\widehat{\C} = \C$.]\\
\end{x}

\begin{x}{\small\bf NOTATION} \ 
Given $f \in L^1(\R^2)$, let
\begin{align*}
\sF_\C f(w) \ 
&= \  \sF_\C f(a,b)\\
&= \  2 \sF_{\R^2} f(2a,-2b)\\
&= \  4 \pi \widehat{f} (4 \pi a, -4\pi b)\\
&= \  \int_{-\infty}^\infty  \int_{-\infty}^\infty \ f(x,y)e^{4 \pi \sqrt{-1}\  (ax - by)} 2dxdy\\
\end{align*}
\end{x}

\begin{x}{\small\bf STANDARDIZATION} \ 
$(G = \C)$ Let $f \in \mathbf{INV}(\C)$, $-$then
\[
\sF_\C \sF_\C f(x,y) = f(-x, -y).
\]

[In fact,
\begin{align*}
\sF_\C \sF_\C f(x,y) \ 
&= \  \int_{-\infty}^\infty \int_{-\infty}^\infty \sF_\C f(a,b) e^{4 \pi \sqrt{-1} \ (ax - by)} \ 2dadb \\	
&= \  \int_{-\infty}^\infty \int_{-\infty}^\infty 4\pi \widehat{f}(4\pi a, -4\pi b) e^{4\pi\sqrt{-1}\  (ax - by)}\  2dadb \\	
&= \  \frac{4\pi}{(4\pi)^2} \int_{-\infty}^\infty \int_{-\infty}^\infty  \widehat{f}(u,-v) e^{\sqrt{-1}\  (ux - vy)}\  2dudv\\	
&= \  \frac{1}{2\pi}  \int_{-\infty}^\infty \int_{-\infty}^\infty  \widehat{f}(u,-v) e^{\sqrt{-1}\  (ux - vy)}\  dudv\\	
&= \  \frac{1}{2\pi}  \int_{-\infty}^\infty \int_{-\infty}^\infty  \widehat{f}(u,-v) e^{-\sqrt{-1}\  (-ux + vy)}\  dudv\\
&= \  \frac{1}{2\pi}  \int_{-\infty}^\infty \int_{-\infty}^\infty  \widehat{f}(u,v) e^{-\sqrt{-1}\  (-ux - vy)}\  dudv\\
&= \  f(-x, -y).]
\end{align*}
\end{x}

\vspace{0.1cm}

\index{Plancherel theorem}
\begin{x}{\small\bf PLANCHEREL THEOREM} \ 
The Fourier transform restricted to $L^1(G) \cap L^2(G)$ is an isometry 
(with respect to $L^2$ norms$)$ onto a dense linear subspace of $L^2(\widehat{G})$, 
hence can be extended uniquely to an isometric isomorphism $L^2(G) \ra L^2(\widehat{G})$.\\
\end{x}

\vspace{0.1cm}

\index{Parseval theorem}
\begin{x}{\small\bf PARSEVAL FORMULA} \ 
$\forall$ $f, g \in  L^2(G)$,
\[
\int_G f(x) \ov{g(x)} d_G(x) = \int_{\widehat{G}} \widehat{f}(\chi) \ov{\widehat{g}(\chi)} d_{\widehat{G}} (\chi).
\]
\end{x}

\vspace{0.1cm}

\begin{x}{\small\bf \un{N.B.}} \ 
In both of these results, the Haar measure on $\widehat{G}$ is per Fourier inversion.
\end{x}

\chapter{
$\boldsymbol{\S}$\textbf{8}.\quad  ADDITIVE p-ADIC CHARACTER THEORY}
\setlength\parindent{2em}
\setcounter{theoremn}{0}


\begin{x}{\small\bf FACT} \ 
Every proper closed subgroup of $\T$ is finite.
\end{x}

\vspace{0.1cm}

Suppose that $G$ is compact abelian and totally disconnected.

\vspace{0.1cm}

\begin{x}{\small\bf LEMMA} \ 
If $\chi \in \widehat{G}$, then the image $\chi(G)$ is a finite subgroup of $\T$.

\vspace{0.1cm}

PROOF \  $\ker \chi$ is closed and 
\[
\chi(G) \thickapprox G / \ker \chi.
\]
But the quotient $G / \ker \chi$ is 0-dimensional, hence totally disconnected.  
Therefore $\chi(G)$ is totally disconnected.  
Since $\T$ is connected, it follows that $\T \ne \chi(G)$, thus $\chi(G)$ is finite.
\end{x}

\vspace{0.1cm}

\begin{x}{\small\bf \un{N.B.}} \ 
The torsion of $\R/\Z$ is $\Q/\Z$, so $\chi$ factors through the inclusion
\[
\Q/\Z \hookrightarrow \R/\Z, \quad \text{i.e., } \chi(G) \subset \Q/\Z.
\]
\end{x}

\vspace{0.1cm}

The foregoing applies in particular to $G = \Z_p$.

\vspace{0.1cm}

\begin{x}{\small\bf LEMMA} \ 
Every character of $\Q_p$ is unitary.

\vspace{0.1cm}

PROOF \   This is because
\[
\Q_p \ = \  \bigcup_{n \in \Z} p^n \Z_p,
\]
where the $p^n \Z_p$ are compact, thus $\S7$, $\# 7$ is applicable.
\end{x}

\vspace{0.1cm}

\begin{x}{\small\bf LEMMA} \ 
If $\chi \in \widehat{\Q}_p$ is nontrivial, then there exists an $n \in \Z$ such that $\chi \equiv 1$ on $p^n\Z_p$ but $\chi \not\equiv 1$ on $p^{n-1}\Z_p$.

\vspace{0.1cm}

PROOF \  Consider a ball $B$ of radius $\frac{1}{2}$ about 1 in 
$\C^\times$ $-$then the only subgroup of  $\C^\times$ contained in $B$ is the trivial subgroup and, by continuity, $\chi(p^n \Z_p)$ must be inside $B$ for all sufficiently large $n$, thus must be 
identically 1 there.
\end{x}
\vspace{0.1cm}

\begin{x}{\small\bf DEFINITION} \ 
The 
\un{conductor}
\index{conductor} 
$\con \chi$ of a nontrivial $\chi \in \widehat{\Q}_p$ is the largest subgroup $p^n \Z_p$ 
on which $\chi$ is trivial (and $n$ is the minimal integer with this property).

\vspace{0.2cm}

A typical $x \ne 0$ of $\Q_p$ has the form
\begin{align*}
x \ 
&=\  \sum_{n = v(x)}^\infty a_n p^n \qquad ( a_n \in \sA, v(x) \in \Z) \\
&=\  f(x) + [x].
\end{align*}
Here the 
\un{fractional part}
\index{fractional part} 
f(x) of $x$ is defined by the prescription
\[
f(x) \ = \ 
\begin{cases}
\sum\limits_{n = v(x)}^{-1} a_np^n \quad \ \   \text{if} \  v(x) < 0\\
0  \qquad\qquad\qquad \text{if} \ v(x) \geq 0
\end{cases}
\]
and the 
\un{integral part} 
\index{integral part} 
$[x]$ of $x$ is defined by the prescription
\[
[x] = \sum_{n = 0 }^\infty a_np^n,
\]
with $f(0) = 0$, $[0] = 0$ by convention.
\end{x}
\vspace{0.1cm}

\begin{x}{\small\bf \un{N.B.}} \ 
\[
f(x) \in \Z\bigl[\frac{1}{p}\bigr] \subset \Q,
\]
where
\[
\Z\bigl[\frac{1}{p}\bigr] = \{\frac{n}{p^k}: n \in \Z, k \in \Z\},
\]
while $[x] \in \Z_p$.
\end{x}
\vspace{0.1cm}

\begin{x}{\small\bf OBSERVATION} \ 
\begin{align*}
0 \ 
&\leq \ f(x) \\
&= \ \sum\limits_{1 \le j \le -v(x)} \frac{a_{-j}}{p^j} \\
&<\  (p-1) \sum_{j = 1}^\infty \frac{1}{p^j}\\ 
&=\  1 
\end{align*}
\qquad\qquad\qquad\qquad$\implies$
\[
f(x) \in [0,1[ \ \cap \  \Z\bigl[\frac{1}{p}\bigr].
\]

Let $\mu_{p^\infty}$ stand for the group of roots of unity in $\C^\times$ having order a power of $p$, 
thus $\mu_{p^\infty}$ is a $p$-group and there is an increasing sequence of cyclic groups
\[
\begin{cases}
\mu_p \subset \mu_{p^2} \subset \cdots \subset \mu_{p^k} \subset  \cdots \\
\mu_{p^\infty} = \bigcup\limits_{k \ge 0}  \mu_{p^k} 
\end{cases}
,
\]
where
\[
 \mu_{p^k}  = \{z \in \C^\times : z^{p^k} = 1\}.
\]
\end{x}

\vspace{0.1cm}

\begin{x}{\small\bf REMARK} \ 
Denote by $\mu$ the group of all roots of unity in $\C^\times$, hence
\[
\mu = \bigcup\limits_{m \ge 1}  \mu_m, \quad \mu_m = \{z \in \C^\times : z^m = 1\}.
\]
Then $\mu$ is an abelian torsion group and $\mu_{p^\infty}$ is the $p$-Sylow subgroup of $\mu$, 
i.e., the maximal $p$-subgroup of $\mu$.
\vspace{0.2cm}

Put
\[
\chi_p(x) = \exp(2\pi \sqrt{-1} \ f(x)) \qquad ( x \in \Q_p).
\]
Then
\[
\chi_p: \Q_p \ra \T
\]
and $\Z_p \subset \ker \chi_p$.
\end{x}

\vspace{0.1cm}

\begin{x}{\small\bf EXAMPLE} \ 
Suppose that $v(x) = -1$, so $x = \ds\frac{k}{p} + y$ with $0 < k \le p-1$ and $y \in \Z_p:$
\[
\chi_p(x) = \exp(2\pi \sqrt{-1} \ \frac{k}{p}) = \zeta^k,
\]
where $\zeta = \exp(2\pi \sqrt{-1}/p)$ is a primitive $p^{th}$ root of unity.
\end{x}

\vspace{0.1cm}

\begin{x}{\small\bf LEMMA} \ 
$\chi_p$ is a unitary character

\vspace{0.1cm}

PROOF \  Given $x, y \in \Q_p$, write
\begin{align*}
f(x+y) - f(x) - f(y) \ 
&=\  x + y - [x+y] - (x - [x]) - (y - [y]) \\
&=\  [x] + [y] - [x+y] \in \Z_p.
\end{align*}
But at the same time
\[
f(x+y) - f(x) - f(y) \in \Z\bigl[\frac{1}{p}\bigr].
\]
Thus
\[
f(x+y) - f(x) - f(y) \in \Z\bigl[\frac{1}{p}\bigr] \cap \Z_p = \Z
\]
and so
\[
\exp(2\pi \sqrt{-1} \ (f(x + y) - f(x) - f(y)) = 1
\]
or still, 
\[
\chi_p(x + y) = \chi_p(x) \chi_p(y).
\]
Therefore $\chi_p : \Q_p \ra \T$ is a homomorphism.  
As for continuity, it suffice to check this at 0, matters then being clear $($since $\chi_p$ is trivial in a neighborhood of 0$)$ 
$(\Z_p$ is open and $0 \in  \Z_p)$.
\end{x}

\vspace{0.1cm}

\begin{x}{\small\bf LEMMA} \ 
The kernel of $\chi_p$ is $\Z_p$.

\vspace{0.1cm}

[A priori, the kernel of $\chi_p$ consists of those $x \in \Q_p$ such that $f(x) \in \Z$.  
Therefore
\[
\con  \chi_p = \Z_p.]
\]
\end{x}

\vspace{0.1cm}

\begin{x}{\small\bf LEMMA} \ 
The image of $\chi_p$ is $\mu_{p^\infty}$.

\vspace{0.1cm}

[A priori, the image of $\chi_p$ consists of the complex numbers of the form
\[
\exp( 2\pi \sqrt{-1} \ \frac{k}{p^m}) \ = \ \exp(2\pi \sqrt{-1}/p^m)^k.
\]
Since $\exp(2\pi \sqrt{-1}/p^m)$ is a root of unity of order $p^m$, these roots generate $\mu_{p^\infty}$ 
as $m$ ranges over the positive integers.]
\end{x}

\vspace{0.1cm}

\begin{x}{\small\bf SCHOLIUM} \ 
$\chi_p$ implements an isomorphism
\[
\Q_p/\Z_p \thickapprox \mu_{p^\infty}.
\]
\end{x}

\vspace{0.1cm}

\begin{x}{\small\bf REMARK} \ 
\begin{align*}
x \in p^{-k}\Z_p	&\Leftrightarrow p^{k}x \in \Z_p\\	
								&\Leftrightarrow \chi_p(p^{k}x) = 1\\
								&\Leftrightarrow \chi_p(x)^{p^k} = 1\\
								&\Leftrightarrow \chi_p(x) \in \mu_{p^k}.
\end{align*}
\end{x}

\vspace{0.1cm}

\begin{x}{\small\bf RAPPEL} \ 
Let $p$ be a prime $-$then a group is 
\un{$p$-primary}
\index{$p$-primary group} 
if every element has order a power of $p$.
\end{x}
\vspace{0.1cm}

\begin{x}{\small\bf RAPPEL} \ 
Every abelian torsion group $G$ is a direct sum of its $p$-primary subgroups $G_p$.
\end{x}

\vspace{0.1cm}

[Note: \ The $p$-primary component of $G_p$ is the $p$-Sylow subgroup of $G$.]

\vspace{0.1cm}

\begin{x}{\small\bf NOTATION} \ 
$\Z(p^\infty)$ is the $p$-primary component of  $\Q /  \Z$.  

Therefore
\[
\Q / \Z \ \thickapprox \ \bigoplus_p \ \Z(p^\infty).
\]
\end{x}
\vspace{0.1cm}

\begin{x}{\small\bf LEMMA} \ 
$\Z(p^\infty)$ is isomorphic to $\mu_{p^\infty}$.

\vspace{0.1cm}

[$\Z(p^\infty)$ is generated by the $1/p^n$ in $\Q /  \Z$.]

\vspace{0.2cm}

Therefore
\[
\Q /  \Z \ \thickapprox \ \bigoplus_p \  \mu_{p^\infty} \ \thickapprox \ \bigoplus_p \  \Q_p /  \Z_p.
\]

[Note: Consequently, 
\begin{align*}
\End(\Q / \Z)	\ 
&\thickapprox \ \End \bigl(\bigoplus\limits_p \ \Q_p /  \Z_p\bigr)\\	
&\thickapprox \ \prod\limits_p \End(\Q_p /  \Z_p)\\	
&\thickapprox \ \prod \limits_p \Z_p.]
\end{align*}
\end{x}

\vspace{0.1cm}

\begin{x}{\small\bf REMARK} \ 
$\widehat{\Z}_p$ is isomorphic to $\mu_{p^\infty}$ (c.f. $\# 26$ infra).
\end{x}

\vspace{0.1cm}

Given $t \in \Q_p$, let $L_t$ be left multiplication by $t$ and put $\chi_{p,t} = \chi_p \circ L_t$ 
$-$then $\chi_{p,t}$ is continuous and $\forall \ x \in \Q_p$,
\[
\chi_{p,t}(x) = \chi_p(tx).
\]

[Note: Trivially, $\chi_{p,0} \equiv 1$.  And $\forall$ $t \ne 0$, 
\[
\con \chi_{p,t} = p^{-v(t)}\Z_p.
\]

Proof: \ 
\begin{align*}
x \in \text{ con }\chi_{p,t}\ 	
&\Leftrightarrow\  tx \in  \Z_p\\			
&\Leftrightarrow\  \abs{tx}_p \le 1\\
&\Leftrightarrow\  \abs{x}_p \le \frac{1}{\abs{t}_p} = p^{v(t)}\\
&\Leftrightarrow\  x \in p^{-v(t)} \Z_p.]	
\end{align*}
Next
\begin{align*}
\chi_{p,t}(x+y)	\  
&=\  \chi_p(t(x+y))\\	
&=\  \chi_p(tx+ty)\\	
&=\  \chi_p(tx)\chi_p(ty)\\		
&=\   \chi_{p,t}(x)\chi_{p,t}(y).
\end{align*}
Therefore $\chi_{p,t} \in \widehat{\Q}_p$.

Next
\begin{align*}
\chi_{p,t+s}(x)	\ 
&= \ \chi_p((t+s)x)\\	
&= \ \chi_p(tx+sx)\\	
&= \ \chi_p(tx)\chi_p(sx)\\		
&= \ \chi_{p,t}(x)\chi_{p,s}(x).
\end{align*}
Therefore the arrow
\begin{align*}
\Xi_p:\Q_p	&\ra \widehat{\Q}_p\\	
t &\mapsto \chi_{p,t}
\end{align*}
is a homomorphism.

\vspace{0.2cm}
											
\begin{x}{\small\bf LEMMA} \ 
If $t\ne s$, then $\chi_{p,t} \ne \chi_{p,s}$.

\vspace{0.1cm}

PROOF \  If to the contrary, $\chi_{p,t} = \chi_{p,s}$, then $\forall$ $x \in \Q_p$, $\chi_p(tx) = \chi_p(sx)$ or still, 
$\forall$ $x \in \Q_p$, $\chi_p((t-s)x) = 1$. 
But $L_{t-s}: \Q_p \ra \Q_p$ is an automorphism, hence $\chi_p$ is trivial, which it isn't.
\end{x}

\vspace{0.1cm}

\begin{x}{\small\bf LEMMA} \ 
The set
\[
\Xi_p(\Q_p) \ = \  \{\chi_{p,t}: t \in \Q_p\}
\]
is dense in $\widehat{\Q}_p$.

\vspace{0.1cm}

PROOF \  Let $H$ be the closure in $\widehat{\Q}_p$ of the $\chi_{p,t}$.  
Consider the quotient $\widehat{\Q}_p/H$.  
To get a contradiction, assume that $H \ne\widehat{\Q}_p$, 
thus that there is a nontrivial $\xi \in \widehat{\widehat{\Q}}_p$ which is trivial on $H$.   
By definition, $H^\perp$ is computed in $\widehat{\widehat{\Q}}_p$, which by Pontryagin duality, is 
identified with $\Q_p$, so spelled out
\[
H^\perp = \{x \in \Q_p: \ev_{\Q_p}\restr{(x)}{H} = 1\}.
\]
Accordingly, for some $x$, $\xi=  \ev_{\Q_p}(x)$, hence $\forall \ t$, 
\begin{align*}
\xi(\chi_{p,t}) \ 
&=\  \ev_{\Q_p}(x)(\chi_{p,t})  \\
&=\  \chi_{p,t}(x) \\
&=\  \chi_p(tx) \\
&= 1,
\end{align*}
which is possible only if $x = 0$ and this implies that $\xi$ is trivial.
\end{x}

\vspace{0.1cm}

\begin{x}{\small\bf LEMMA} \ 
The arrows
\[
\begin{cases}
\Q_p	\ra \Xi_p(\Q_p)\\	
\Xi_p(\Q_p)\ra  \Q_p\
\end{cases}
\]
are continuous.
\end{x}

Therefore $\Xi(\Q_p)$ is a locally compact subgroup of $\widehat{\Q}_p$.  
But a locally compact subgroup of a locally compact group is closed.  
Therefore $\Xi_p(\Q_p) = \widehat{\Q}_p$. 

In summary$:$
\begin{x}{\small\bf THEOREM} \ 
$\widehat{\Q}_p$ is topologically isomorphic to $\Q_p$ via the arrow
\[
\Xi_p:\Q_p		\ra  \widehat{\Q}_p.
\]
\end{x}

\vspace{0.1cm}

\begin{x}{\small\bf LEMMA} \ 
Fix $t$ $-$then $\restr{\chi_{p,t}}{\Z_p} = 1$ iff $t \in \Z_p$.

\vspace{0.1cm}

PROOF \  Recall that the kernel of $\chi_p$ is $\Z_p$.
\[
\begin{aligned}
&\text{\textbullet} \quad t \in \Z_p, \ x \in \Z_p \implies tx \in \Z_p \implies \chi_p(tx) = 1 \implies \restr{\chi_{p,t}}{\Z_p} = 1.\\
&\text{\textbullet} \quad \restr{\chi_{p,t}}{\Z_p} = 1 \implies \chi_{p,t}(1) = 1 \implies \chi_p(t) = 1 \implies t \in \Z_p.
\end{aligned}
\]
\end{x}

\vspace{0.1cm}

\begin{x}{\small\bf APPLICATION} \ 
$\widehat{\Z}_p$ is isomorphic to $\mu_{p^\infty}$.

\vspace{0.1cm}

$[\widehat{\Z}_p$ can be computed as $\widehat{\Q}_p/\Z_p^\perp$.  
But $\Z_p^\perp$, when viewed as a subset of $\Q_p$, consists of those $t$ such that $\restr{\chi_{p, t}}{\Z_p} = 1.$ 
Therefore 
\[
\widehat{\Z}_p \ 
\thickapprox \ \widehat{\Q}_p/\Z_p \ 
\thickapprox \ \Q_p /\Z_p \ 
\thickapprox \ \mu_{p^\infty}.]
\]
\end{x}

\vspace{0.1cm}

\begin{x}{\small\bf NOTATION} \ 
Let
\[
x_\infty(x) = \exp(-2\pi \sqrt{-1} \  x) \qquad (x \in \R).
\]
\end{x}

\vspace{0.1cm}

\index{Product principle}
\begin{x}{\small\bf PRODUCT PRINCIPLE} \ 
$\forall$ $x \in \Q$,
\[
\prod_{p \le \infty} \chi_p(x) \ = \ 1.
\]

\vspace{0.1cm}

PROOF \  Take $x$ positive $-$then there exist primes $p_1, \cdots, p_n$ such that $x$ admits a representation
\[
x \ = \ 
\frac{N_1}{p_1^{\alpha_1}} + 
\frac{N_2}{p_2^{\alpha_2}} + \cdots + 
\frac{N_n}{p_n^{\alpha_n}} + M,
\]
where the $\alpha_k$ are positive integers, the $N_k$ are positive integers 
$(1 \le N_k < p_k^{\alpha_k} - 1)$, and $M \in \Z$.  
Appending a subscript to $f$, we have
\[
f_{p_k}(x) \ = \ 
\frac{N_k}{p_k^{\alpha_k}}, \quad f_p(x) = 0 \quad (p \ne p_k, \ k = 1, 2, \ldots, n).
\]
Therefore
\begin{align*}
\prod_{p < \infty} \chi_p(x) \ 	
&= \  \prod_{1 \le k \le n} \chi_{p_k}(x)\\	
&= \  \prod_{1 \le k \le n} \exp(2\pi\sqrt{-1}\ 	 f_{p_k}(x))\\
&= \  \exp(2\pi\sqrt{-1}\ 	 \sum_{k = 1}^n f_{p_k}(x))\\
&= \  \exp(2\pi\sqrt{-1}\ 	 (x - M))\\
&= \  \exp(2\pi\sqrt{-1}\ 	 x)
\end{align*}
\qquad\qquad$\implies$
\begin{align*}
\prod_{p \le \infty} \chi_p(x) \ 	
&=\ \prod_{p < \infty} \chi_p(x) \chi_\infty (x) \\	
&=\ \exp(2\pi\sqrt{-1}\ 	 x) \exp(-2\pi\sqrt{-1}\ 	 x)\\							
&=\ 1.
\end{align*}
\end{x}

\vspace{0.1cm}

\[
\textbf{APPENDIX}
\]
\setcounter{theoremn}{0}

Let $\K$ be a finite extension of $\Q_p$.
\vspace{0.25cm}

\begin{x}{\small\bf THEOREM} \ 
The topological groups $\K$ and $\widehat{\K}$ are topologically isomorphic.

\vspace{0.1cm}

[Put
\begin{align*}
\chi_{\K,p}(a) \ 
&=\  \exp(2\pi \sqrt{-1} \ f(\tr_{\K/\Q_p}(a)))\\
&=\ \chi_p(\tr_{\K/\Q_p}(a))
\end{align*}
and given $b \in \K$, put
\[
\chi_{\K,p,b}(a) = \chi_{\K,p}(ab).
\]
Proceed from here as above.]
\end{x}

\vspace{0.1cm}

\begin{x}{\small\bf REMARK} \ 
Every character of $\K$ is unitary.
\end{x}

\vspace{0.1cm}

\begin{x}{\small\bf LEMMA} \ 
\[
\begin{cases}
\ a \in R 	&\implies \tr_{\K/ \Q_p}(a) \in \Z_p\\	
\ a \in P &\implies \tr_{\K/ \Q_p}(a) \in p\Z_p	
\end{cases}
.
\]
\end{x}

\vspace{0.1cm}

\begin{x}{\small\bf DEFINITION} \ 
The 
\un{differential of $\K$}
\index{differential of $\K$} 
is the set
\[
\Delta_\K = \{b \in \K:\tr_{\K/ \Q_p}(R b) \subset \Z_p\}.
\]
\end{x}

\vspace{0.1cm}

\begin{x}{\small\bf LEMMA} \ 
$\Delta_\K$ is a proper $R$-submodule of $\K$ containing $R$.
\end{x}

\vspace{0.1cm}

\begin{x}{\small\bf LEMMA} \ 
There exists a unique nonnegative integer $d$ 
$-$\un{the differential exponent} \un{of $\K$}
\index{the differential exponent of $\K$} 
$-$characterized by the condition that
\[
\pi^{-d}R = \Delta_\K.
\]

[This follows from the theory of "fractional ideals" $($details omitted$)$.]
\end{x}

\vspace{0.1cm}

[Note: $\chi_{\K,p}$ is trivial on $\pi^{-d}R$ but is nontrivial on $\pi^{-d-1}R.]$

\vspace{0.2cm}

\begin{x}{\small\bf LEMMA} \ 
Let $e$ be the ramification index of $\K$ over $\Q_p$ (cf. $\S5$, $\#17)$ $-$ then
\[
a \in P^{-e+1} \implies \tr_{\K/ \Q_p}(a) \in \Z_p.
\]

\vspace{0.1cm}

PROOF \  Let
\[
a \in P^{-e+1} = \pi^{-e+1}R = \pi^{-e}(\pi R) =  \pi^{-e}P,
\]
so $a = \pi^{-e}b \ (b \in P)$.  
Write $p = \pi^eu$ and consider $pa$:
\[
pa = \pi^eu\pi^{-e}b = ub.
\]
But
\begin{align*}
\abs{u} = 1, \  \abs{b} < 1 	\ 
&\implies\  \abs{ub} < 1\\	
&\implies\  ub \in P\\
&\implies\  \tr_{\K/ \Q_p}(ub) \in p\Z_p\\
&\implies\  \tr_{\K/ \Q_p}(pa) \in p\Z_p\\
&\implies\  p\tr_{\K/ \Q_p}(a) \in p\Z_p\\
&\implies\  \tr_{\K/ \Q_p} \in \Z_p.
\end{align*}
\end{x}

\vspace{0.1cm}

\begin{x}{\small\bf APPLICATION} \ 
\[
d \ge e-1.
\]

[It suffices to show that
\[
P^{-e+1} \subset \Delta_\K \quad  (\equiv \pi^{-d}R).
\]
Thus let $a \in  P^{-e+1}$, say $a = \pi^eb$  $(b \in P)$, and let $r \in R$ $-$then the claim is that
\[
\tr_{\K/ \Q_p}(ar) \in \Z_p.
\]
But
\[
ar = \pi^{-e}br \in \pi^e P \quad (\abs{br} < 1)
\]
or still,
\[
ar \in P^{-e+1} \implies \tr_{\K/ \Q_p}(ar) \in \Z_p.]
\]
\end{x}

\vspace{0.1cm}

\begin{x}{\small\bf REMARK} \ 
Therefore $d = 0 \implies e = 1$, hence in this situation, $\K$ is unramified.

[Note: There is also a converse, viz. if $\K$ is unramified, then $d = 0.$]
\end{x}

\vspace{0.1cm}

\begin{x}{\small\bf \un{N.B.}} \ 
It can be shown that
\[
\tr_{\K/ \Q_p}(R) = \Z_p \ \text{ iff } \ d = e-1.
\]
\end{x}

\vspace{0.1cm}

\begin{x}{\small\bf CRITERION} \ 
Fix $b \in \K$ $-$then
\[
b \in \Delta_\K \Leftrightarrow \forall \  a \in R, \ \chi_{\K,p}(ab) = 1.
\]

\vspace{0.1cm}

PROOF \ 

\indent\indent \textbullet \quad 
$a \in R, b \in \Delta_\K  \ \implies\  ab \in \Delta_\K$

$\indent\indent\indent\indent\indent\implies\  \tr_{\K/ \Q_p}(ab) \in \Z_p$

$\indent\implies$
\[
\chi_{\K,p}(ab) = \chi_p(\tr_{\K/ \Q_p}(ab)) = 1.
\]
\indent\indent\textbullet \quad 
$\forall \ a \in R, \ \chi_{\K,p}(ab) = 1$ 

$\indent\implies\ \forall \ a \in R, \  \tr_{\K/ \Q_p}(ab) \in \Z_p$

$\indent\implies\  b \in \Delta_\K$.

Normalize Haar measure on $\K$ by the condition
\[
\mu_\K(R) = \int_{R} da = q^{-d/2}.
\]
Let $\chi_R$ be the characteristic function of $R$ $-$then
\[
\int_\K \chi_R(a) \chi_{\K,p}(ab) da = \int_R \chi_{\K,p}(ab)da.
\]

\indent\indent\textbullet \quad 
$b \in \Delta_\K \implies \chi_{\K,p}(ab) = 1 \quad (\forall$ $a \in R)$

$\hspace{2.65cm} \implies \int_R\chi_{\K,p}(ab)da = \mu_\K(R) = q^{-d/2}$.

\indent\indent\textbullet \quad 
$b \notin \Delta_\K \implies \chi_{\K,p}(ab) \ne 1$ $(\exists$ $a \in R)$  

$\hspace{2.65cm}\implies \int_R\chi_{\K,p}(ab)da = 0$.

\text{}\\
Consequently, as a function of $b$,
\[
\int_R\chi_{\K,p}(ab)da =q^{-d/2} \chi_{\Delta_\K} (b),
\]
$\chi_{\Delta_\K}$ the characteristic function of $\Delta_\K$.
\end{x}

\vspace{0.1cm}

\begin{x}{\small\bf LEMMA} \ 
\[
[\pi^{-d}R:R] = q^d.
\]
Therefore
\begin{align*}
\mu_\K(\Delta_\K) \ 
&=\  \mu_\K(\pi^{-d}R)\\	
&=\  q^d\mu_\K(R)\\
&= q^d q^{-d/2}\\								
&=\  q^{d/2}.
\end{align*}
\end{x}

\vspace{0.1cm}

\begin{x}{\small\bf LEMMA} \ 
$\forall$ $a \in \K$,
\[
\int_\K q^{-d/2} \chi_{\Delta_\K} (b) \chi_{\K,p} (ab) db = \chi_R(a).
\]

\vspace{0.1cm}

PROOF \  The left hand side reduces to
\[
q^{-d/2} \int_{\Delta_\K} \chi_{\K,p} (ab) db
\]
and there are two possibilities

\indent\indent\textbullet \quad 
$a \in R  \implies ab \in \Delta_\K \quad (\forall \ b \in \Delta_\K)$

$\hspace{2.4cm}\implies \tr_{\K/ \Q_p}(ab) \in \Z_p$

$\hspace{2.4cm}\implies \chi_{\K,p}(ab) = 1$

$\indent\implies$
\begin{align*}
q^{-d/2} \int_{\Delta_\K} \chi_{\K,p}(ab)db \ 
&=\  q^{-d/2} \mu_\K(\Delta_\K) \\
&=\  q^{-d/2} q^{d/2}  \\
&=\  1.
\end{align*}

\indent\indent\textbullet \quad 
$a \notin R: \chi_{K,p}(ab) \ne 1 \quad (\exists \  b \in \Delta_\K)$

$\qquad\qquad\implies$
\[
q^{-d/2} \int_{\Delta_\K} \chi_{\K,p}(ab)db = 0.
\]


To detail the second point of this proof, work with the normalized absolute value 
(cf. $\S 6, \ \# 18$) and recall that 
$\abs{\pi}_K = \ds\frac{1}{q}$ 
(cf. $\S 5, \ \# 21$).  
Accordingly, 
\[
x \in \pi^n R \Leftrightarrow \abs{x}_\K \le q^{-n}.
\]
Fix $a \notin R$ $-$then the claim is that $b \ra \chi_{\K,p}(ab)$ $(b \in \Delta_\K)$ is nontrivial.  
For
\begin{align*}
\chi_{\K,p}(ab) = 1 	\ 
&\Leftrightarrow\  ab \in \pi^{-d}R\\
&\Leftrightarrow\  \abs{ab}_\K \le q^d\\
&\Leftrightarrow\  \abs{a}_\K \abs{b}_\K \le q^d\\
&\Leftrightarrow\   \abs{b}_\K \le \frac{q^d}{\abs{a}_\K} = q^{d+v(a)}.
\end{align*}
But
\begin{align*}
a \notin R 	
&\implies\  v(a) < 0\\
&\implies\  -v(a) > 0\\
&\implies\  -d-v(a) > -d\\
&\implies\ \pi^{-d-v(a)}R \subsetneqq \pi^{-d}R,
\end{align*}
a proper containment.
\end{x}


\chapter{
$\boldsymbol{\S}$\textbf{9}.\quad  MULTIPLICATIVE p-ADIC CHARACTER THEORY}
\setlength\parindent{2em}
\setcounter{theoremn}{0}

\ \indent 
Recall that
\[
\Q_p^\times \ \thickapprox \  \Z \times \Z_p^\times,
\]
the abstract reflection of the fact that for ever $x \in \Q_p^\times$, there is a unique $v(x) \in \Z$ and a unique $u(x) \in \Z_p^\times$ such that $x = p^{v(x)}u(x)$.  
Therefore
\[
\widehat{(\Q_p^\times)} \ 
\thickapprox \ \widehat{\Z} \times \widehat{(\Z_p^\times)}  \ 
\thickapprox \ \T \times \widehat{(\Z_p^\times)}.
\]

\vspace{0.1cm}

\begin{x}{\small\bf \un{N.B.}} \ 
A character of $\Q_p$ is necessarily unitary (cf. $ \S8, \  \# 4$) but this is definitely not the case for $\Q_p^\times$ $($cf. infra$)$.
\end{x}

\vspace{0.1cm}

\begin{x}{\small\bf DEFINITION} \ 
A character $\chi: \Q_p^\times \ra \C^\times$ is \un{unramified} if it is trivial on $\Z_p^\times$.
\end{x}

\vspace{0.1cm}

\begin{x}{\small\bf EXAMPLE} \ 
Given any complex number $s$, the arrow $x \ra \abs{x}_p^s$ is an unramified character of $\Q_p^\times$.
\end{x}

\vspace{0.1cm}

\begin{x}{\small\bf LEMMA} \ 
If $\chi: \Q_p^\times \ra \C^\times$ is an unramified character, then there exists a complex number $s$ such that 
$\chi = \acdot_p^s$.

\vspace{0.1cm}

PROOF \ Such a $\chi$ factors through the projection $\Q_p^\times \ra p^\Z$ defined by $x \ra \abs{x}_p$, 
hence gives rise to a character $\widetilde{\chi}: p^\Z \ra \C^\times$ 
which is completely determined by its value on $p$, say $\widetilde{\chi}(p) = p^s$ for the complex number
\[
s = \frac{\log \widetilde{\chi}(p)}{\log p},
\]
itself determined up to an integral multiple of
\[
\frac{2\pi \sqrt{-1}}{\log p}.
\]
Therefore
\begin{align*}
\chi(x) \ 		
&=\ \widetilde{\chi}(\abs{x}_p)\\	
&=\ \widetilde{\chi}(p^{-v(x)})\\		
&=\  (\widetilde{\chi}(p))^{-v(x)}\\	
&=\  (p^s)^{-v(x)}\\
&=\  (p^{-v(x)})^s\\
&=\  \abs{x}_p^s.	
\end{align*}

[Note: \  For the record, 
\begin{align*}
\abs{x}_p^{2 \pi \sqrt{-1}/ \log p}	\ 	
&=\  (p^{-v(x)})^{2 \pi \sqrt{-1}/ \log p}\\	
&=\  (e^{-v(x) \log p})^{2 \pi \sqrt{-1} / \log p}\\		
&=\  e^{-v(x) 2 \pi \sqrt{-1}}\\	
&=\  1.]\\	
\end{align*}
\end{x}

\vspace{0.1cm}

Suppose that $\chi: \Q_p^\times \ra \C^\times$ is a character $-$then $\chi$ can be written as
\[
\chi(x) = \abs{x}_p^s \un{\chi}(u(x)),
\]
where $s \in \C$ and $\un{\chi} \equiv \restr{\chi}{\Z_p^\times} \in \widehat{(\Z_p^\times)}$, 
thus $\chi$ is unitary iff $s$ is pure imaginary.
\vspace{0.2cm}

\begin{x}{\small\bf LEMMA} \ 
If $\un{\chi} \in \widehat{(\Z_p^\times)}$ is nontrivial, then there is an $n \in \N$ such that $\un{\chi} \equiv 1$ on $U_{p,n}$ but $\chi \not\equiv 1$ on $U_{p,n-1}$ (cf. $\S 8, \  \# 5$).
\end{x}

\vspace{0.1cm}

Assume again that $\chi: \Q_p^\times \ra \C^\times$ is a character.

\vspace{0.2cm}

\begin{x}{\small\bf DEFINITION} \ 
$\chi$ is 
\un{ramified of degree $n \ge 1$}
\index{ramified of degree $n \ge 1$} 
if 
$\restr{\un{\chi}}{U_{p,n}} \equiv 1$ 
and 
$\restr{\un{\chi}}{U_{p,n-1}} \not\equiv 1$.
\end{x}

\vspace{0.1cm}

\begin{x}{\small\bf DEFINITION} \ 
The 
\un{conductor}
\index{conductor} 
$\con\chi$ of $\chi$ is $\Z_p^\times$ if $\chi$ is unramified and $U_{p,n}$ if 
$\chi$ is ramified of degree $n$.
\end{x}

\vspace{0.1cm}

\begin{x}{\small\bf RAPPEL} \ 
If $G$ is a finite abelian group, then the number of unitary characters of $G$ is card $G$.
\end{x}

\vspace{0.1cm}

\begin{x}{\small\bf LEMMA} \ 
\[
[\Z_p^\times:U_{p,1}] = p-1  \qquad (\text{cf. } \S 4, \  \#40)
\]
and 
\[
[U_{p,1}:U_{p,n}] = p^{n-1} .
\]
\end{x}

\vspace{0.1cm}

If $\chi$ is ramified of degree $n$, then $\un{\chi}$ can be viewed as a unitary character of $\Z_p^\times / U_{p,n}$.   
But the quotient $\Z_p^\times / U_{p,n}$ is a finite abelian group, thus has 
\[
\card \ \Z_p^\times / U_{p,n} = [\Z_p^\times:U_{p,n}] 
\]
unitary characters.  
And
\begin{align*}
[\Z_p^\times:U_{p,n}] \ 
&=\  [\Z_p^\times:U_{p,1}] \cdot [U_{p,1}:U_{p,n}]  \\
&=\  (p-1) p^{n-1},
\end{align*}
this being the number of unitary characters of $\Z_p^\times$ of degree $\le n$.  
Therefore the 
group $\Z_p^\times$ has p-2 unitary characters of degree 1 and for $n \ge 2$, the group $\Z_p^\times$ has
\[
(p-1) p^{n-1} - (p-1) p^{n-2} = p^{n-2}(p-1)^2
\]
unitary characters of degree $n$.

\vspace{0.2cm}

\begin{x}{\small\bf LEMMA} \ 
Let $\chi \in \widehat{\Q_p^\times}$ $-$then
\[
\chi(x) = \abs{x}_p^{\sqrt{-1} \  t} \un{\chi}(u(x)),
\]
where $t$ is real and 
\[
-(\pi / \log p) < t \le \pi / \log p.
\]
\end{x}

\newpage


\[
\textbf{APPENDIX}
\]
\setcounter{theoremn}{0}

Suppose that $p \ne 2$, let $\tau \in \Q_p^\times - (\Q_p^\times)^2$, and form the quadratic extension
\[
\Q_p(\tau) = \{x+y \sqrt{\tau} : x,y \in \Q_p\}.
\]

\begin{x}{\small\bf NOTATION} \ 
Let $\Q_{p,\tau}$ be the set of points of the form $x^2 - \tau y^2$ $(x \ne 0$, $y \ne 0)$.
\end{x}

\vspace{0.1cm}

\begin{x}{\small\bf LEMMA} \ 
$\Q_{p,\tau}$ is a subgroup of $\Q_p^\times$ containing $(\Q_p^\times)^2$.
\end{x}

\vspace{0.1cm}

\begin{x}{\small\bf LEMMA} \ 
\[
[\Q_p^\times:\Q_{p,\tau}] = 2 \text{ and } [\Q_{p,\tau}:(\Q_p^\times)^2] = 2.
\]

[Note:
\[
[\Q_p^\times:(\Q_p^\times)^2] = 4		\qquad (\text{cf. } \S 4, \ \# 53).]
\]
\end{x}

\vspace{0.1cm}

\begin{x}{\small\bf DEFINITION} \ 
Given $x \in \Q_p^\times$, let
\[
\sgn_\tau(x) \ = \ 
\begin{cases}
\  \ 1	\quad \text{if }x \in \Q_{p,\tau}\\
-1 	\quad \text{if } x \notin \Q_{p,\tau}
\end{cases}
.
\]
\end{x}

\vspace{0.1cm}

\begin{x}{\small\bf LEMMA} \ 
$\sgn_\tau$ is a unitary character of $\widehat{\Q}_p$.
\end{x}

\chapter{
$\boldsymbol{\S}$\textbf{10}.\quad  TEST FUNCTIONS}
\setlength\parindent{2em}
\setcounter{theoremn}{0}
 
\ 
\indent The 
\un{Schwartz space}
\index{Schwartz space} 
$\sS(\R^n)$ consists of those complex valued $\sC^\infty$ functions which, together with all their derivatives, vanish at infinity faster than any power of $\norm{\cdot}$.

\vspace{0.2cm}

\begin{x}{\small\bf DEFINITION} \ 
The elements $f$ of $\sS(\R^n)$ are the 
\un{test functions}
\index{test functions} 
on $\R^n$.
\end{x}

\vspace{0.1cm}

\begin{x}{\small\bf EXAMPLE} \ 
Take $n = 1$ $-$then
\[
f(x) = C x^A \exp(-\pi x^2),			
\]
where $A = 0$ or 1, is a test function, said to be 
\un{standard}. 
\index{standard (test function)} 
Here
\[
\int_\R x^A \exp(-\pi x^2) e^{2 \pi \sqrt{-1} \ tx} dx = (\sqrt{-1})^A t^A \exp(- \pi t^2),
\]
thus $\sF_\R$ of a standard function is again standard (c.f. $\S 7$, $51$).
\end{x}

\vspace{0.1cm}

[Note: \  Henceforth, by definition, the Fourier transform of an $f \in  L^1(\R)$ will be the function

\[
\widehat{f}:\R \lra \C
\]
defined by the rule
\begin{align*}
\widehat{f}(t) \ \
&= \ \sF_\R f(t) \\
&= \  \int_{\R} f(x) e^{2 \pi \sqrt{-1}\  tx} dx.]
\end{align*}

\vspace{0.1cm}

\begin{x}{\small\bf EXAMPLE} \ 
Take $n = 2$ and identify $\R^2$ with $\C$ $-$then
\[
f(z) = Cz^A \ov{z}^{B} \exp( -2 \pi \abs{z}^2),
\]
where $A,B \in \Z_{\ge 0}$ $\&$ $AB = 0$, is a test function, said to be 
\un{standard}.   
\index{standard (test function)}
Here
\[
\int_{\C} z^A \ov{z}^B \exp( -2 \pi \abs{z}^2) e^{2 \pi \sqrt{-1} \ (wz + \ov{w} \ov{z}}) \abs{dz \wedge d \ov{z}} 
\ = \ \sqrt{-1}^{A+B} \  w^B \ov{w}^A \exp( -2 \pi \abs{w}^2),
\]
thus $\sF_\C$ of a standard function is again standard ( c.f. $\S 7$, $\# 53$ ).
\end{x}

\vspace{0.1cm}

[Note: Henceforth, by definition, the Fourier transform of an $f \in L^1(\C)$ will be the function 
\[
\widehat{f} : \C \lra \C
\]
defined by the rule
\begin{align*}
\widehat{f}(w) \ 
&= \ \sF_\C f(w) \\
&= \ \int_\C f(z) e^{2\pi \sqrt{-1} \  (wz + \ov{w}\ov{z})} \abs{dz \wedge d\ov{z}} .]
\end{align*}

\vspace{0.2cm}

\begin{x}{\small\bf DEFINITION} \ 
Let $G$ be a totally disconnected locally compact group $-$then a function $f:G \ra \C$ is said to be 
\un{locally constant} if for any $x \in G$, 
there is an open subset $U_x$ of $G$ containing $x$ such that $f$ is constant on $U_x$.
\end{x}

\vspace{0.1cm}

\begin{x}{\small\bf LEMMA} \ 
A locally constant function $f$ is continuous.

\vspace{0.1cm}

PROOF \ 
Fix $x \in G$ and suppose that $\{x_i\}$ is a net converging to $x$ $-$then $x_i$ is eventually in $U_x$, hence there $f(x_i) = f(x)$.
\end{x}

\vspace{0.1cm}

\begin{x}{\small\bf DEFINITION} \ 
The 
\un{Bruhat space}
\index{Bruhat space} 
$\sB(G)$
\index{$\sB(G)$} 
consists of those complex valued locally constant functions whose support is compact.
\end{x}

\vspace{0.1cm}

[Note: \  $\sB(G)$ carries a "canonical topology" but I shall pass in silence as regards to its precise formulation$]$.\\

\begin{x}{\small\bf DEFINITION} \ 
The elements $f$ of $\sB(G)$ are the 
\un{test functions}
\index{test functions (on \mG)} 
on $G$.
\end{x}

\vspace{0.1cm}

\begin{x}{\small\bf LEMMA} \ 
Given a test function $f$, there exists an open-compact subgroup $K$ of $G$, and integer 
$n \ge 0$, elements $x_1, \ldots, x_n$ in $G$ and elements $c_1, \ldots, c_n$ in $\C$ such that the union
$\bigcup\limits_{k=1}^n K x_k K	$ is disjoint and 
\[
f = \sum_{k=1}^n c_k \chi_{K x_k K},	
\]
$\chi_{K x_k K}$ the characteristic function of $K x_k K$.
\vspace{0.1cm}

\vspace{0.1cm}

PROOF \ 
Since $f$ is locally constant, for every $z \in \C$ the pre image $f^{-1}(z)$ is an open subset of $G$.  
Therefore $X = \{x:f(x) \ne 0\}$ is the support of $f$.  
This said, given $x \in X$, define a map
\[
\begin{aligned}
\phi_x: G \times\ G 		&\ra \C\\	
(x_1, x_2) 				&\mapsto f(x_1 x x_2)\\	
\end{aligned}
,
\]
thus $\phi_x(e,e) = f(x) $ and $\phi_x$ is continuous if $\C$ has the discrete topology. 
Consequently, one can find an open-compact subgroup $K_x$ of $G$ such that $\phi_x$ is constant on $K_x \times K_x$.  
Put $U_x = K_x \times K_x$ $-$then $U_x$ is open-compact and $f$ is constant on $U_x$.  
But $X$ is covered by the $U_x$, hence, being compact, is covered by finitely many of them.  
Bearing in mind that distinct double cosets are disjoint, consider now the intersection $K$ of the finitely many $K_x$ that occur.

\end{x}

\vspace{0.1cm}

Specialize and let $G = \Q_p$.
\\
\begin{x}{\small\bf EXAMPLE} \ 
If $K \subset \Q_p$ is open-compact, then its characteristic function $\chi_K$ is a test function on $\Q_p$ .
\end{x}

\vspace{0.1cm}

\begin{x}{\small\bf LEMMA} \ 
Every $f \in \sB(\Q_p)$ is a finite linear combination of functions of the form
\[
\chi_{x+p^n \Z_p}		\qquad (x \in \Q_p, \ n \in \Z).
\]

[This is an instance of $\# 8$ or argue directly (c.f. $\S 4$, $\#33$).]
\end{x}

\vspace{0.1cm}

\begin{x}{\small\bf DEFINITION} \ 
Given $f \in L^1(\Q_p)$, its 
\un{Fourier transform}
\index{Fourier transform} 
is the function 
\[
\widehat{f}:\Q_p \lra \C
\]
defined by the rule
\begin{align*}
\widehat{f}(t) 
&= \   \int_{\Q_p} f(x) \chi_{p,t}(x) dx\\
&= \ \int_{\Q_p} f(x) \chi_p(tx) dx.	
\end{align*}
\end{x}

\vspace{0.1cm}

\begin{x}{\small\bf LEMMA} \ 
$\forall$ $f \in L^1(\Q_p)$,
\[
\widehat{\ov{f\hspace{0.07cm}}}(t) = \overline{\widehat{f}(-t)}.
\]

\vspace{0.1cm}

PROOF \ 
\begin{align*}
\widehat{\ov{f\hspace{0.07cm}}}(t)   \ 	
&=\  \int_{\Q_p} \overline{f(x)} \chi_p(tx) dx \\
&=\  \int_{\Q_p} \overline{f(x) \chi_p(-tx)} dx \\
&=\  \int_{\Q_p} \overline{f(x) \chi_p((-t)x)} dx \\
&=\  \overline{\int_{\Q_p)} f(x) \chi_p((-t)x) dx} \\
&=\  \overline{\widehat{f}(-t)}. 
\end{align*}
\end{x}

\vspace{0.1cm}

\begin{x}{\small\bf SUBLEMMA} \ 

\[
\int_{p^n\Z_p} \chi_p(x) dx \ = \ 
\begin{cases}
p^{-n} \qquad (n \geq 0)\\
0 \qquad\quad \ (n < 0)
\end{cases}
.
\]

\vspace{0.1cm}

[Recall that 
\[
\mu_{\Q_p}(p^n \Z_p) \ = \  p^{-n}
\]
and apply $\S 7$, $\# 46$ and $\S 8$, $\# 12$.]
\end{x}

\vspace{0.1cm}

\begin{x}{\small\bf LEMMA} \ 
Take $f = \chi_{p^n\Z_p}$  $-$then
\[
\widehat{\chi}_{p^n\Z_p} = p^{-n} \chi_{p^{-n}\Z_p}.
\]

\vspace{0.1cm}

PROOF \ 
\begin{align*}
\widehat{\chi}_{p^n\Z_p}(t) \ 	
&=\  \int_{\Q_p} \chi_{p^n\Z_p}(x) \chi_{p,t}(x) dx\\
&=\  \int_{\Q_p} \chi_{p^n\Z_p}(x) \chi_p(tx) dx\\
&=\  \abs{t}_p^{-1}\int_{\Q_p} \chi_{p^n\Z_p}(t^{-1}x) \chi_p(x) dx\\
&=\  \abs{t}_p^{-1}\int_{p^{n+v(t)}\Z_p} \chi_p(x) dx.
\end{align*}
The last integral equals 
\[
p^{-n - v(t)}
\]
if $n + v(t) \geq 0$ and equals 0 if $n + v(t) < 0$ (cf. $\# 13$).  
But
\[
t \in p^{-n}\Z_p \Leftrightarrow v(t) \ge -n \Leftrightarrow n+v(t) \ge 0.
\]
Since
\[
\abs{t}^{-1}_pp^{v(t)} \ =\  1,
\]
it therefore follows that
\[
\widehat{\chi}_{p^n\Z_p} \ =\  p^{-n} \chi_{p^{-n}\Z_p}.
\]
In particular, 
\[
\widehat{\chi}_{\Z_p} \ =\  \chi_{\Z_p}.
\]
\end{x}

\vspace{0.1cm}

\begin{x}{\small\bf THEOREM} \ 
Take $f = \chi_{x+p^n\Z_p}$ $-$then
\[
\widehat{\chi}_{x+p^n\Z_p}(t)  = \ 
\begin{cases}
\chi_p(tx)p^{-n} \qquad (\abs{t}_p \le p^n)\\
0 \qquad \qquad\qquad (\abs{t}_p > p^n)
\end{cases}
.
\]

\vspace{0.1cm}

PROOF \ 
\begin{align*}
\widehat{\chi}_{x+p^n\Z_p}(t) \ 
&=\  \int_{\Q_p} \chi_{x+p^n\Z_p}(y) \chi_{p,t}(y) dy\\
&=\  \int_{\Q_p} \chi_{x+p^n\Z_p}(y) \chi_p(ty) dy\\
&=\  \int_{x+p^n\Z_p} \chi_p(ty) dy\\
&=\  \int_{p^n\Z_p} \chi_p(t(x+y)) dy\\
&=\  \int_{p^n\Z_p} \chi_p(tx+ty) dy\\
&=\  \int_{p^n\Z_p} \chi_p(tx) \chi_p(ty) dy\\
&=\  \chi_p(tx) \int_{p^n\Z_p} \chi_p(ty) dy\\
&=\  \chi_p(tx) \int_{\Q_p} \chi_{p^n\Z_p}(y) \chi_p(ty) dy\\										
&=\  \chi_p(tx) \int_{\Q_p} \chi_{p^n\Z_p}(y) \chi_{p,t}(y) dy\\
&=\  \chi_p(tx) \widehat{\chi}_{{p^n}\Z_p} (t)\\
&=\  \chi_p(tx) p^{-n}\chi_{p^{-n}\Z_p}(t).
\end{align*}
\end{x}

\vspace{0.1cm}

\begin{x}{\small\bf APPLICATION} \ 
Taking into account $\# 10$, 
\[
f \in \sB(\Q_p) \Rightarrow \widehat{f} \in \sB(\Q_p). 
\]
\end{x}

\vspace{0.1cm}

\begin{x}{\small\bf THEOREM} \ 
$\forall$ $f \in \mathbf{INV}(\Q_p)$, 
\[
\widehat{\widehat{f}\hspace{0.1cm}} = f(-x)	\qquad (x \in \Q_p).
\]

\vspace{0.1cm} 

PROOF \ 
It suffices to check this for a single function, so take $f = \chi_{Z_p}$ $-$then as noted above,
\[
\widehat{\chi}_{\Z_p} = \chi_{\Z_p},
\]
thus $\forall$ x,
\[
\widehat{\widehat{\chi}\hspace{0.05cm}}_{\Z_p}(x) = \chi_{\Z_p} (x) =  \chi_{\Z_p} (-x).
\]
\end{x}
\vspace{0.1cm}

\begin{x}{\small\bf \un{N.B.}} \ 
It is clear that
\[
\sB(\Q_p) \subset \mathbf{INV}(\Q_p).
\]
\end{x}

\vspace{0.1cm}

\begin{x}{\small\bf SCHOLIUM} \ 
The arrow $f \ra \widehat{f}$ is a linear bijection of $\sB(\Q_p)$ onto itself.

[Injectivity is manifest.  
As for surjectivity, the arrow $f \ra \check{f}$, where
\[
\widecheck{f\hspace{0.05cm}} = f(-x),
\]
maps $\sB(\Q_p)$ into itself.  And
\[
f \ 
=\ \widecheck{\widecheck{f\hspace{0.05cm}}} 
=\  (\widecheck{f\hspace{0.05cm}}) \raisebox{0.22cm}{$\text{ }\widecheck{}$}
=\  (\widecheck{f\hspace{0.05cm}})    \raisebox{0.11cm}{$\text{ }\widehat{} \hspace{0.15cm}\text{ }\widehat{}$}
=\  ((\widecheck{f\hspace{0.05cm}}) 
\raisebox{0.08cm}{$\text{ }\widehat{}$} \hspace{0.1cm} )
\raisebox{0.08cm}{$\text{ }\widehat{}$}.]
\]
\end{x}

\begin{x}{\small\bf REMARK} \ 
As is well-known, the same conclusion obtains if $\Q_p$ is replaced by $\R$ or $\C$.
\end{x}

\vspace{0.1cm}

Pass now from  $\Q_p$ to  $\Q_p^\times$.

\vspace{0.2cm}

\begin{x}{\small\bf LEMMA} \ 
Let $ f \in \sB(\Q_p^\times)$ $-$then $\exists$ $n \in \N:$
\[
\begin{cases}
\ \abs{x}_p < p^{-n} 	\implies f(x) = 0\\
\ \abs{x}_p > p^{n} \ \ 	\implies f(x) = 0
\end{cases}
.
\]

Therefore an element $f$ of $\sB(\Q_p^\times)$ can be viewed as an element of $\sB(\Q_p)$ with the property that $f(0) = 0$.
\end{x}

\vspace{0.1cm}

\begin{x}{\small\bf DEFINITION} \ 
Given $f \in L^1(\Q_p^\times, d^\times x)$, its 
\un{Mellin transform}
\index{Mellin transform} 
$\widetilde{f}$ 
is the Fourier transform of $f$ per $\Q_p^\times:$
\[
\widetilde{f}(\chi) = \int_{\Q_p^\times} f(x) \chi(x) d^\times x.
\]

[Note: \  
By definition, 
\[
d^\times x = \frac{p}{p-1} \frac{dx}{\abs{x}_p}		\qquad \text{(c.f. } \S 6, \ \# 26),
\]
so
\[
\vol_{d^\times x}(\Z_p^\times) \ =\  \vol_{dx}(\Z_p) \ =\  1.]
\]
\end{x}

\vspace{0.1cm}

\begin{x}{\small\bf EXAMPLE} \ 
Take $f = \chi_{\Z_p^\times}$ $-$then
\begin{align*}
\widetilde{\chi}_{\Z_p^\times} (\chi) \ 
&= \  \int_{\Q_p^\times} \chi_{\Z_p^\times} (x) \chi(x) d^\times x \\
&= \ \int_{\Z_p^\times} \chi(x) d^\times x.
\end{align*}
Decompose $\chi$ as in $\S 9$, $\# 10$, hence
\begin{align*}
\int_{\Z_p^\times} \chi(x) d^\times x \ 
&=\  \int_{\Z_p^\times} \abs{x}_p^{\sqrt{-1}\ t} \un{\chi}(p^{-v(x)}x) d^\times x\\
&=\  \int_{\Z_p^\times}  \un{\chi}(x) d^\times x\\
&=\ 
\begin{cases}
 \ 0  \qquad (\un{\chi} \not\equiv 1)\\
\ 1  \qquad (\un{\chi} \equiv 1)
\end{cases}
.
\end{align*}

According to $\S 9$, $\# 2$, a unitary character $\chi \in \widehat{(\Q_p^\times)}$ is unramified if its restriction 
$\un{\chi}$ to $\Z_p^\times$ is trivial.  
Therefore the upshot is that the Mellin transform of $\chi_{\Z_p^\times}$ 
is the characteristic function of the set of unramified elements of $\widehat{(\Q_p^\times)}$.
\end{x}

\vspace{0.1cm}


\[
\textbf{APPENDIX}
\]
\setcounter{theoremn}{0}

Let $\K$ be a finite extension of $\Q_p$ $-$then
\[
\K^\times \thickapprox \Z \times R^\times
\]
and the generalities developed in \S9 go through with but minor changes when $\Q_p$ is replace by $\K$.

In particular$:$ $\forall$ $\chi \in \widehat{K}^\times$, there is a splitting
\[
\chi(a) = \abs{a}_\K^{\sqrt{-1}\  t} \un{\chi}(\pi^{-v(a)}a),
\]
where $t$ is real and 
\[
-(\pi/\log q) < t \le \pi / \log q.
\]

[Note: \  $\chi$ is \un{unramified} if it is trivial on $R^\times$.$]$

\begin{x}{\small\bf \un{N.B.}} \ 
The $"\pi"$ in the first instance is a prime element (c.f. $\S 5$, $\# 10$) and $\abs{\pi}_\K = \ds\frac{1}{q}$.  
On the other hand, the $"\pi"$ in the second instance is $3.14\ldots$ .

\vspace{0.4cm}

The extension of the theory  from $\sB(\Q_p)$ to $\sB(\K)$ is straightforward, the point of departure being the observation that
\[
\int_{\pi^nR} \chi_{\K,p}(a) da = \mu_\K(R) \ 
\begin{cases}
\ q^{-n} \qquad (n = -d, -d+1, \ldots)\\
\  0 \qquad \quad (n = -d-1, -d-2, \ldots)
\end{cases}
.
\]
\end{x}

\vspace{0.1cm}

\begin{x}{\small\bf CONVENTION} \ 
Normalize the Haar measure on $\K$ by stipulating that $\ds\int_R da = q^{-d/2}$.
\end{x}

\vspace{0.1cm}

\begin{x}{\small\bf DEFINITION} \ 
Given $f \in L^1(\K)$, its 
\un{Fourier transform}
\index{Fourier transform} 
is the function 
\[
\widehat{f}:\K \lra \C
\]
defined by the rule
\[
\begin{aligned}
\widehat{f}(b)	
=&  \int_{\K} f(a) \chi_{\K,p,b} (a) da\\	
=&  \int_{\K} f(a) \chi_{\K,p} (ab) da.
\end{aligned}
\]
\end{x}

\vspace{0.2cm}

\begin{x}{\small\bf THEOREM} \ 
$\forall$ $f \in \mathbf{INV}(\K)$,
\[
\widehat{\widehat{f}\hspace{0.1cm}}(a) = f(-a)	\qquad (a \in \K).
\]

\vspace{0.1cm}

PROOF \  
It suffices to check this for a single function, so take $f = \chi_R$, 
in which case the work has already been done in the Appendix to $\S 8$.  
To review:

\vspace{0.4cm}

\indent\indent$\text{\textbullet} \quad \widehat{\chi}_R(b)$
\indent\indent\indent\indent\quad \  \ $=\  \ds\int_{\K} \chi_R(a) \chi_{\K,p} (ab)da$\\
\indent\indent\indent\indent\indent\indent\indent\indent\indent $=\  \ds\int_{R} \chi_{\K,p} (ab)da$\\
\indent\indent\indent\indent\indent\indent\indent\indent\indent $=\  q^{-d/2} \chi_{\Delta_\K}(b)$.\\

\vspace{0.1cm}

\indent\indent$\text{\textbullet} \quad \int_{\K} q^{-d/2} \chi_{\Delta_\K}(b) \chi_{\K,p} (ab)db \quad $
$ =\  q^{-d/2} \ds\int_{\Delta_\K} \chi_{\K,p} (ab)db$\\
\indent\indent\indent\indent\indent\indent\indent\indent\indent $=\  \chi_R(a) \qquad \  \text{(loc. cit., $\# 13$)}$\\
\indent\indent\indent\indent\indent\indent\indent\indent\indent $=\  \chi_R(-a)$.

\end{x}

\vspace{0.1cm}

\begin{x}{\small\bf \un{N.B.}} \ 
It is clear that 
\[
\sB(k) \subset \mathbf{INV}(\K).
\]
\end{x}

\vspace{0.1cm}

\begin{x}{\small\bf SCHOLIUM} \ 
The arrow $f \ra \widehat{f}$ is a linear bijection of $\sB(k)$ onto itself.
\end{x}

\vspace{0.1cm}

\begin{x}{\small\bf CONVENTION} \ 
Put
\[
d^\times a = \frac{q}{q-1} \frac{da}{\abs{a}_\K}.
\]
Then $d^\times a$ is a Haar measure on $\K^\times$ and 
\[
\vol_{d^\times a} (R^\times) = \vol_{da}(R) = q^{-d/2}.
\]
\end{x}

\vspace{0.1cm}

\begin{x}{\small\bf DEFINITION} \ 
Given $f \in L^1(\K^\times, d^\times a)$, its 
\un{Mellin transform}
\index{Mellin transform} 
$\widetilde{f}$ is the Fourier transform of $f$ per $\K^\times:$
\[
\widetilde{f}(\chi) = \int_{\K^\times} f(a) \chi(a) d^\times a.
\]
\end{x}

\vspace{0.1cm}

\begin{x}{\small\bf EXAMPLE} \ 
Take $f = \chi_{R^\times}$ $-$then

\[
\widetilde{\chi}_{R^\times}(\chi) = \ 
\begin{cases}
0 \quad\qquad \ \  (\un{\chi} \ne 1)\\
q^{-d/2} \qquad (\chi \equiv 1)
\end{cases}
.\]

\end{x}

\chapter{
$\boldsymbol{\S}$\textbf{11}.\quad  LOCAL ZETA FUNCTIONS: $\R^\times$ or $\C^\times$}
\setlength\parindent{2em}
\setcounter{theoremn}{0}
\ \indent 
We shall first consider $\R^\times$, hence $\widetilde{\R}^\times \thickapprox \Z / 2 \Z \times \C$ and every character has the form
\[
\chi(x) \equiv \chi_{\sigma, s}(x) = (\sgn x)^\sigma \abs{x}^s \quad 
(\sigma \in \{0, 1\}, \ s \in \C) \quad 
(\text{cf. } \S7, \  \#11).
\]

\vspace{0.1cm}

\begin{x}{\small\bf DEFINITION} \ 
Given $f \in \sS(\R^n)$ and a character $\chi: \R^\times \ra \C^\times$, 
the 
\un{local zeta} \un{function}
\index{local zeta function} 
attached to the pair $(f, \chi)$ is
\[
Z (f,\chi) = \int_{\R^\times} f(x) \chi(x) d^\times x, \qquad \text{where $d^\times x = \frac{dx}{\abs{x}}$.}
\]

[Note: The parameters $\sigma$ and $s$ are implicit:
\[
Z(f,x) \equiv Z(f,\chi_{\sigma,s}).]
\]
\end{x}

\vspace{0.1cm}

\begin{x}{\small\bf LEMMA} \ 
The integral defining $Z(f,\chi)$ is absolutely convergent for $\Re(s) > 0$.

\vspace{0.1cm}

PROOF \   
Since $f$ is Schwartz, there are no issues at infinity.  
As for what happens at the origin, let $I = ]-1,1[ \ -\  \{0\}$ and fix $C > 0$ such that $\abs{f(x)} \le C$ $(x \in I)$. $-$then
\begin{align*}
\abs{Z(f,\chi)} \ 	
&\le\  \int_{\R- \{0\}} \abs{f(x)} \abs{x}^{\Re(s) - 1} dx\\	
&\le\   \bigl(\int_{\R - I} + \int_I\bigr) \abs{f(x)} \abs{x}^{\Re(s) - 1} dx\\	
&\le\   M + C \int_I \abs{x}^{\Re(s) - 1} dx,
\end{align*}
a finite quantity.
\end{x}

\vspace{0.1cm}

\begin{x}{\small\bf LEMMA} \ 
$Z(f,\chi)$ is a holomorphic function of $s$ in the strip $\Re(s) > 0$.

[Formally, 
\[
\frac{d}{ds}Z(f,\chi) \ =\  \int_{\R^\times} f(x) (sgn x)^\sigma (\log \abs{x}) \abs{x}^s d^\times x,
\]

and while correct, "differentiation under the integral sign" does require a formal proof $\ldots$ .]
\end{x}

\vspace{0.1cm}

\begin{x}{\small\bf NOTATION} \ 
Put
\[
\widecheck{\chi} = \chi^{-1} \| \cdot \|.
\]

The integral defining $Z(f,\widecheck{\chi} )$ is absolutely convergent if $\Re(1-s) > 0$, i.e., if $1-\Re(s) > 0$ 
or still, if $\Re(s) < 1$.
\end{x}

\vspace{0.1cm}

\begin{x}{\small\bf LEMMA} \ 
Let $f, g \in \sS(\R)$ and suppose that $0 < \Re(s) < 1$ $-$then
\[
Z(f,\chi) Z(\widehat{g},\widecheck{\chi}) =  Z(\widehat{f},\widecheck{\chi}) Z(g,\chi) 
\]

\vspace{0.1cm}

PROOF \ 
Write
\[
Z(f,\chi) Z(\widehat{g},\widecheck{\chi})  
= \int \int_{\R^\times \times \R^\times} f(x) \widehat{g}(y) \chi(xy^{-1}) \abs{y} d^\times x \hspace{0.05cm} d^\times y
\]
and make the substitution $t = yx^{-1}$ to get
\[
Z(f,\chi) Z(\widehat{g},\widecheck{\chi})  
= \int_{\R^\times} \bigl(\int_{\R^\times} f(x) \widehat{g}(tx) \abs{x} d^\times x\bigr) 
\chi(t^{-1}) \abs{t} d^\times t.
\]
The claim now is that the inner integral is symmetric in $f$ and $g$ (which then implies that
\[
Z(f,\chi) Z(\widehat{g},\widecheck{\chi})  = Z(g,\chi)   Z(\widehat{f},\widecheck{\chi}), 
\]
the desired equality$)$.   To see this is so, observe first that
\[
\abs{x} du \cdot d^\times x = \abs{u} dx \cdot d^\times u.
\]
Since $\R^\times$ and $\R$ differ by a single element, it therefore follows that
\begin{align*}
\int_{\R^\times} f(x) \widehat{g}(tx) \abs{x} d^\times x \ 
&=\  \int_{\R^\times} f(x) \abs{x} \bigl(\int_\R g(u) e^{2\pi\sqrt{-1}\  txu} du\bigr) d^\times x\\	
&=\   \int \int_{\R \times \R^\times} f(x) g(u) \abs{x} e^{2 \pi \sqrt{-1}\  txu} du d^\times x\\
&=\  \int_{\R^\times} g(u) \abs{u} \bigl(\int_\R f(x) e^{2 \pi \sqrt{-1}\  txu} dx \bigr) d^\times u\\
&=\  \int_{\R^\times} g(u) \widehat{f}(tu) \abs{u} d^\times u.
\end{align*}

Fix $\phi \in \sS(\R)$ and put
\[
\rho(\chi) = \frac{Z(\phi,\chi)}{Z(\widehat{\phi},\widecheck{\chi})}
\]
Then $\rho(\chi)$ is independent of the choice of $\phi$ and $\forall$ $f \in \sS(\R)$, the 
\un{functional equation}
\index{functional equation}
\[
Z(f,\chi) = \rho(\chi) Z(\widehat{f}, \widecheck{\chi})
\]
obtains.
\end{x}

\vspace{0.1cm}

\begin{x}{\small\bf LEMMA} \ 
$\rho(\chi)$ is a meromorphic function of $s$ (cf. infra).
\end{x}

\vspace{0.1cm}

\begin{x}{\small\bf APPLICATION} \ 
$\forall$ $f \in \sS(\R), Z(f,\chi)$ admits a meromorphic continuation to the whole $s$-plane.
\end{x}

\vspace{0.1cm}

\begin{x}{\small\bf NOTATION} \ 
Set
\[
\Gamma_\R(s) = \pi^{-s/2} \Gamma(s/2).
\]
\end{x}

\vspace{0.1cm}

\begin{x}{\small\bf DEFINITION} \ 
Write
\[L(\chi) =\ 
\begin{cases}
\Gamma_\R(s) &\quad (\sigma = 0)\\
\Gamma_\R(s+1) &\quad (\sigma = 1)
\end{cases}
.\]
\end{x}

\vspace{0.1cm}

Proceeding to the computation of $\rho(\chi)$, distinguish two cases.
\vspace{0.2cm}

\qquad \textbullet \quad \un{$\sigma = 0$} \quad Take $\phi_0(x)$ to be $e^{-\pi x^2}$ $-$then
\begin{align*}
Z(\phi_0, \chi) \ 
&=\  \int_{\R^\times} e^{-\pi x^2} \abs{x}^s d^\times x\\	
&=\  2\int_0^\infty e^{-\pi x^2} x^{s-1} dx\\	
&=\  \pi^{-s/2} \Gamma(s/2) \\
&=\  \Gamma_\R(s)\\
&=\  L(\chi).
\end{align*}
Next $\widehat{\phi}_0 = \phi_0$ ( cf. \S10, \#2) so by the above argument,
\[
Z(\widehat{\phi}_0,\widecheck{\chi}) = L(\widecheck{\chi}),
\]
from which
\begin{align*}
\rho(\chi) \ 	
&=\  \frac{L(\chi)}{L(\widecheck{\chi})}\\	
&=\  \frac{\pi^{-s/2} \Gamma\bigl(\ds\frac{s}{2}\bigr)}{\pi^{-(1-s)/2} \Gamma\bigl(\ds\frac{1-s}{2}\bigr)}\\	
&=\  2^{1-s} \pi^{-s} \cos\bigl( \ds\frac{\pi s}{2}\bigr) \Gamma(s).
\end{align*}

\qquad \textbullet \quad \un{$\sigma = 1$} \quad Take $\phi_1(x)$ to be $xe^{-\pi x^2}$ $-$then
\allowdisplaybreaks
\begin{align*}
Z(\phi_1, \chi) \ 
&=\  \int_{\R^\times} xe^{-\pi x^2} \frac{x}{\abs{x}} \abs{x}^s d^\times x\\	
&=\  \int_{\R^\times} e^{-\pi x^2}  \abs{x}^{s+1} d^\times x\\
&=\  2\int_0^\infty e^{-\pi x^2} x^s dx\\	
&=\  \pi^{-(s+1)/2} \Gamma\bigl(\frac{s+1}{2}\bigr)\\
&=\  \Gamma_\R(s+1)\\
&=\  L(\chi).
\end{align*}
Next
\[
\widehat{\phi}_1(t) = \sqrt{-1}\ t \exp(- \pi t^2) \quad \text{(cf. \S10, \#2)}.
\]
Therefore
\begin{align*}
Z(\widehat{\phi_1}, \widecheck{\chi}) \ 	
&=\  \sqrt{-1}\   \int_{\R^\times} xe^{-\pi x^2} \frac{x}{\abs{x}} \abs{x}^{1-s} d^\times x\\	
&=\  \sqrt{-1}\   \int_{\R^\times} e^{-\pi x^2}  \abs{x}^{2-s} d^\times x\\
&=\  \sqrt{-1}\   2\int_0^\infty e^{-\pi x^2} x^{1-s} dx\\	
&=\  \sqrt{-1}\  \pi^{-(2-s)/2} \Gamma\bigl(\frac{2-s}{2}\bigr)\\
&=\  \sqrt{-1}\  \Gamma_\R(2-s)\\
&=\  \sqrt{-1}\  L(\widecheck{\chi}).
\end{align*}
Accordingly
\begin{align*}
\rho(\chi) \ 
&=\  - \sqrt{-1}\  \ds\frac{L(\chi)}{L(\widecheck{\chi})}\\	
&=\  - \sqrt{-1}\   \frac{\pi^{-(s+1)/2} \Gamma(\ds\frac{s+1}{2})}{\pi^{(s-2)/2} \Gamma\bigl(\ds\frac{2-s}{2}\bigr)}\\	
&=\  - \sqrt{-1}\    2^{1-s} \pi^{-s} \sin(\ds\frac{\pi s}{2}) \Gamma(s).
\end{align*}

\vspace{0.1cm}

\begin{x}{\small\bf FACT} \ 
\[
\begin{cases}
\ds\frac{\zeta(1-s)}{\zeta(s)} 	&= 2^{1-s} \pi^{-s} \cos\bigl(\ds\frac{\pi s}{2}\bigr) \Gamma(s)\\
\\
\ds\frac{\zeta(s)}{\zeta(1-s)} 	&= 2^{s} \pi^{s-1} \sin\bigl(\ds\frac{\pi s}{2}\bigr) \Gamma(1-s)
\end{cases}
.
\]

\vspace{0.1cm}

To recapitulate: $\rho(\chi)$ is a meromorphic function of $s$ and 
\[
\rho(\chi) 	= \epsilon(\chi) \frac{L(\chi)}{L(\widecheck{\chi})},
\]
where
\[\epsilon(\chi) =\ 
\begin{cases}
\  1 				&\quad \text{$(\sigma = 0)$}\\
\  - \sqrt{-1} 	&\quad \text{$(\sigma = 1)$}
\end{cases}
.\]
\end{x}

\vspace{0.1cm}

Having dealt with $\R^\times$, let us now turn to $\C^\times$, hence 
$\widetilde{\C}^\times$ $\thickapprox \Z \times \C$ and every character has the form
\[
\chi(x) \equiv \chi_{n,s}(x) = \exp(\sqrt{-1} \  n \  \arg x) \abs{x}^s 
\quad (n \in \Z, \ s \in \C) \quad (\text{cf. } \S 7, \ \#12).
\]
Here, however, it will be best to make a couple of adjustments.

1. Replace $x$ by $z$.

2. Replace $\acdot$ by $\acdot_\C$, the normalized absolute value, so
\[
\abs{z}_\C = \abs{z  \bar{z}} = \abs{z}^2	\qquad ( \text{cf. } \S6, \  \#15).
\]

\begin{x}{\small\bf DEFINITION} \ 
Given $f \in \sS(\C)$ $(= \sS(\R^2))$ and a character $\chi : \C^\times \ra \C^\times$, the \un{local zeta function} attached to the pair $(f, \chi)$ is
\[
Z(f,\chi) = \int_{\C^\times} f(z) \chi(z) d^\times z,
\]
where $d^\times z = \ds\frac{\abs{dz \wedge d\ov{z}}}{\abs{z}_\C }$.
\vspace{0.2cm}

[Note: \  The parameters $n$ and $s$ are implicit:
\[
Z(f,\chi) \equiv Z(f, \chi_{n,s}).]
\]
\end{x}

\vspace{0.1cm}

\begin{x}{\small\bf NOTATION} \ 
Put
\[
\widecheck{\chi} = \chi^{-1} \acdot_\C.
\]

The analogs of \#2 and \#3 are immediate, as is the analog of $\#5$ 
(just replace $\R^\times$ by $\C^\times$ and $\acdot$ by $\acdot_\C$), the crux then being the analog of $\#6$.
\end{x}

\vspace{0.1cm}

\begin{x}{\small\bf NOTATION} \ 
Set
\[
\Gamma_\C(s) = (2\pi)^{1-s} \Gamma(s).
\]
\end{x}

\vspace{0.1cm}

\begin{x}{\small\bf DEFINITION} \ 
Write
\[
L(\chi) = \Gamma_\C(s + \frac{\abs{n}}{2}).
\]

To determine $\rho(\chi)$ via a judicious choice of $\phi$ per the relation
\[
\rho(\chi) = \frac{Z(\phi,\chi)}{Z(\widehat{\phi},\widecheck{\chi})},
\]
let
\[
\phi_n(z) =\ 
\begin{cases}
\ \ov{z}^{n} e^{-2\pi\abs{z}^2}		&\quad (n \ge 0)\\
\ z^{-n} e^{-2\pi\abs{z}^2}	 		&\quad (n < 0)\\
\end{cases}
.\]
Then
\[
\widehat{\phi}_n = \sqrt{-1}^{\abs{n}} \phi_{-n} 	\qquad (\text{cf. } \S 10, \  \#3).
\]
\end{x}

\vspace{0.1cm}

\begin{x}{\small\bf \un{N.B.}} \ 
In terms of polar coordinates $z = r e^{\sqrt{-1} \text{ } \theta}$,\\

\qquad\textbullet \quad$\phi_n(z) = r^{\abs{n}} \exp (-2\pi r^2 - \sqrt{-1} \ n \theta)$\\

\qquad\textbullet \quad$d^\times z =\ds\frac{2 r dr d \theta}{r^2} = \ds\frac{2}{r} dr d \theta$\\

\qquad\textbullet \quad$\chi(z) = e^{\sqrt{-1} \  n \theta} \abs{z}_\C^s = e^{\sqrt{-1} \ n \theta} r^{2s}$.\\

Therefore

\begin{align*}
Z(\phi_n,\chi) \ 	
&=\  \int_0^{2\pi} \int_0^\infty r^{\abs{n}} \exp(-2\pi r^2 - \sqrt{-1} \  n \theta) e^{\sqrt{-1} \ n \theta} r^{2s} \frac{2}{r} dr d\theta\\	
&=\  \int_0^{2\pi} \int_0^\infty r^{2(s-1) + \abs{n}} \exp(-2\pi r^2) 2r dr d\theta\\		
&=\  2\pi \int_0^\infty t^{(s-1) + \abs{n}/2} \exp(-2\pi t) dt\\
&=\  (2\pi)^{1 - s - \abs{n}/2} \Gamma\bigl(s + \frac{\abs{n}}{2}\bigr)\\ 
&=\  \Gamma_\C\bigl(s + \frac{\abs{n}}{2}\bigr)\\
&=\  L(\chi)\\  
\end{align*}
and
\begin{align*}
Z(\widehat{\phi}_n,\widecheck{\chi}) \ 	
&=\  Z((\sqrt{-1})^{\abs{n}} \phi_{-n}, \widecheck{\chi})\\	
&=\  (\sqrt{-1})^{\abs{n}} (2\pi)^{1 - (1 - s)- \abs{n}/2} \Gamma\bigl(1 - s + \frac{\abs{n}}{2}\bigr)\\		
&=\  (\sqrt{-1})^{\abs{n}} (2\pi)^{s - \abs{n}/2} \Gamma\bigl(1 - s + \frac{\abs{n}}{2}\bigr)\\	
&=\  (\sqrt{-1})^{\abs{n}}  \Gamma_\C\bigl(1 - s + \frac{\abs{n}}{2}\bigr)\\
&=\  (\sqrt{-1})^{\abs{n}}  L(\widecheck{\chi}).
\end{align*}
Consequently, 
\begin{align*}
\rho(\chi)	\ 	
&=\  \frac{Z(\phi_n,\chi)}{Z(\widehat{\phi}_n,\widecheck{\chi}) }\\	
&=\  (\sqrt{-1})^{-\abs{n}} \frac{L(\chi)}{L(\widecheck{\chi})}\\		
&=\  \epsilon(\chi)\frac{L(\chi)}{L(\widecheck{\chi})},		
\end{align*}
where
\[
\epsilon(\chi) = (\sqrt{-1})^{-\abs{n}}.
\]
And
\[
\frac{L(\chi)}{L(\widecheck{\chi})} = (2\pi)^{1 - 2s} 
\frac{\Gamma\bigl(s + \ds\frac{\abs{n}}{2}\bigr)}{\Gamma\bigl(1-s + \ds\frac{\abs{n}}{2}\bigr)}.
\]
\end{x}


\chapter{
$\boldsymbol{\S}$\textbf{12}.\quad  LOCAL ZETA FUNCTIONS: $\Q^\times_p$}
\setlength\parindent{2em}
\setcounter{theoremn}{0}

\ \indent 
The theory set forth below is in the same spirit as that of \S11 but matters are technically more complicated due to the presence of ramification.
\vspace{0.1cm}

\begin{x}{\small\bf DEFINITION} \ 
Given $f \in \sB(\Q_p)$ and a character $\chi: \Q_p^\times \ra \C^\times$, the 
\un{local zeta} \un{function}
\index{local zeta function} 
attached to the pair $(f,\chi)$ is
\[
Z(f, \chi) = \int_{\Q_p^\times} f(x) \chi(x) d^\times x,
\]
where $d^\times x = \ds\frac{p}{p-1} \ds\frac{dx}{\abs{x}_p}$ \ (cf. $\S 6$, $\#26$).
\end{x}

\vspace{0.1cm}

[Note: \  There are two parameters associated with $\chi$, viz. $s$ and $\un{\chi}$ (cf. $\S9$).]
\vspace{0.2cm}

\begin{x}{\small\bf LEMMA} \ 
The integral defining $Z(f, \chi)$ is absolutely convergent for $\Re(s) > 0$.

\vspace{0.01cm}

PROOF \   
It suffices to check the absolute convergence for $f = \chi_{p^n \Z_p}$ (cf. $\S 10$, $\#10$) 
and then we might just as well take $n = 0:$
\begin{align*}
\abs{Z(f,\chi)} \ 	
&\le \int_{\Q_p^\times} \abs{f(x)} \abs{x}_p^{\Re(s)} d^\times x\\	
&=\  \int_{\Q_p^\times} \chi_{\Z_p}(x) \abs{x}_p^{\Re(s)} d^\times x\\		
&=\  \int_{\Z_p - \{0\}} \abs{x}_p^{\Re(s)} d^\times x\\	
&=\  \frac{1}{1 - p^{-\Re(s)}} \qquad \text{(cf. $\S 6$, $\#27$)}.
\end{align*}
\end{x}

\vspace{0.1cm}

\begin{x}{\small\bf LEMMA} \ 
$Z(f, \chi)$ is a holomorphic function of $s$ in the strip $\Re(s) > 0$.
\end{x}

\vspace{0.1cm}

\begin{x}{\small\bf NOTATION} \ 
Put
\[
\widecheck{x} = x^{-1}  \acdot_p.
\]

The integral defining $Z(f, \widecheck{\chi})$ is absolutely convergent if  $\Re(1-s) > 0$, i.e., if $1 - \Re(s) > 0$ or still, if $\Re(s) < 1$.
\end{x}

\vspace{0.1cm}

\begin{x}{\small\bf LEMMA} \ 
Let $f, g \in \sB(\Q_p)$ and suppose that $0 < \Re(s) < 1$ $-$then
\[
Z(f, \chi)  Z(\widehat{g}, \widecheck{\chi}) = Z(\widehat{f}, \widecheck{\chi}) Z(g, \chi).
\]

[Simply follow verbatim the argument employed in \S11, \#5.]
\end{x}

\vspace{0.1cm}

Fix $\phi \in \sB(\Q_p)$ and put
\[
\rho(\chi) \ = \ \frac{Z(\phi,\chi)}{Z(\widehat{\phi},\widecheck{\chi}}).
\]
Then $\rho(\chi)$ is independent of the choice of $\phi$ and $\forall$ $f \in \sB(\Q_p)$, the 
\index{functional equation}
\un{functional equation} 
\[
Z(f, \chi) = \rho(\chi) Z(\widehat{f}, \widecheck{\chi})
\]
obtains.

\vspace{0.1cm}

\begin{x}{\small\bf LEMMA} \ 
$\rho(\chi)$ is a meromorphic function of $s$ $($cf. infra$)$.
\end{x}

\vspace{0.1cm}

\begin{x}{\small\bf APPLICATION} \ 
$\forall$ $f \in \sB(\Q_p)$, $Z(f, \chi)$ admits a meromorphic continuation to the whole $s$-plane
\end{x}

\vspace{0.1cm}

\begin{x}{\small\bf DEFINITION} \ 
Write
\[L(\chi) =\ 
\begin{cases}
\ (1 - \chi(p))^{-1} &\quad (\chi  \ \text{unramified})\\
\ 1 &\quad (\chi \ \text{ramified}
\end{cases}
.\]
\end{x}

\vspace{0.1cm}

There remains the computation of $\rho(\chi)$, the simplest situation being when $\chi$ is unramified, say 
$\chi = \acdot_p^s$, in which case we take $\phi_0(x) = \chi_p(x) \chi_{\Z_p}(x):$
\begin{align*}
Z(\phi_0,\chi)	\ 
&=\  \int_{\Q_p^\times} \phi_0(x) \chi(x) d^\times x\\	
&=\  \int_{\Q_p^\times}  \chi_p(x) \chi_{\Z_p}(x) \abs{x}_p^s d^\times x\\		
&=\  \int_{\Z_p - \{0\}} \chi_p(x) \abs{x}_p^s d^\times x\\	
&=\  \int_{\Z_p - \{0\}} \abs{x}_p^s d^\times x\\
&=\  \frac{1}{1 - p^{-s}} \qquad \text{(cf. \S6, \#27)}\\
&=\  \frac{1}{1 - \abs{p}_p^s} \\
&=\  \frac{1}{1 - \chi(p)}\\
&=\  L(\chi).
\end{align*}

To finish the determination, it is necessary to explicate the Fourier transform 
$\widehat{\phi}_0$ of $\phi_0$ $($cf. \S10, \#11$):$
\begin{align*}
\widehat{\phi}_0(t)\ 
&=\  \int_{\Q_p} \phi_0(x) \chi_p(tx) dx\\	
&=\  \int_{\Q_p} \chi_p(x) \chi_{\Z_p}(x) \chi_p(tx) dx\\	
&=\  \int_{\Z_p} \chi_p(x)  \chi_p(tx) dx\\	
&=\  \int_{\Z_p} \chi_p((1+t)x) dx\\
&=\  \chi_{\Z_p}(t).
\end{align*}
Therefore
\begin{align*}
Z(\widehat{\phi}_0, \widecheck{\chi})	
&=\  \int_{\Q_p^\times} \widehat{\phi}_0(x)  \widecheck{\chi}(x) d^\times x\\	
&=\  \int_{\Q_p^\times} \chi_{\Z_p}(x) \abs{x}_p^{1-s} d^\times x\\	
&=\  \int_{\Z_p-\{0\}} \abs{x}_p^{1-s} d^\times x\\	
&=\  \frac{1}{1 - p^{-(1-s)}}	\qquad \text{(cf. \S6, \#27)}\\
&=\  \frac{1}{1 - \abs{p}_p^{1-s}}\\
&=\  \frac{1}{1 - \widecheck{\chi}(p)}\\
&=\  L(\widecheck{\chi}).
\end{align*}
And finally
\[
\rho(\chi) = \frac{Z(\phi_0, \chi)}{Z(\widehat{\phi}, \widecheck{\chi})	} = \frac{L(\chi)}{L(\widecheck{\chi})}
\]
or still, 
\[
\rho(\chi) = \frac{1 - p^{-(1-s)}}{1 - p^{-s}}.
\]

\vspace{0.1cm}

\begin{x}{\small\bf REMARK} \ 
The function
\[
\frac{1 - p^{-(1-s)}}{1 - p^{-s}}
\]
has a simple pole at $s = 0$ with residue
\[
\frac{p-1}{p}\log p
\]
and there are no other singularities. 
\end{x}

\vspace{0.1cm}

Suppose now that $\chi$ is ramified of degree $n \ge 1 : \chi = \acdot_p^s$ $\un{\chi}$ \ 
$($cf. \S9, \#6 $)$ and take $\phi_n(x) = \chi_p(x) \chi_{p^{-n} \Z_p}(x):$
\begin{align*}
Z(\phi_n, \chi)	
&=\  \int_{\Q_p^\times} \phi_n(x)  \chi(x) d^\times x\\	
&=\  \int_{\Q_p^\times} \chi_p(x) \chi_{p^{-n} \Z_p}(x) \abs{x}_p^s\un{\chi}(x) d^\times x\\
&=\  \int_{p^{-n}\Z_p-\{0\}} \chi_p(x)  \abs{x}_p^s\un{\chi}(x) d^\times x\\	
&=\  \sum_{k = -n}^\infty \int_{\Z_p^\times} \chi_p(p^ku) \abs{p^ku}_p^s \un{\chi}(u) d^\times u\\
&=\  \sum_{k = -n}^\infty p^{-ks}  \int_{\Z_p^\times} \chi_p(p^ku) \un{\chi}(u) d^\times u.					
\end{align*}

\vspace{0.1cm}

\begin{x}{\small\bf LEMMA} \ 
If $\abs{v}_p \ne p^n$, then
\[
\int_{\Z_p^\times} \chi_p(vu) \un{\chi}(u) d^\times u = 0.
\]

Since $\abs{p^k}_p = p^{-k}$, $Z(\phi_n, \chi)$ reduces to
\[
p^{ns} \int_{\Z_p^\times} \chi_p(p^{-n}u) \un{\chi}(u) d^\times u.
\]
Let $E = \{e_i: i \in I\}$ be a system of coset representatives for 
$\Z_p^\times  / U_{p,n}$ $-$then by assumption, $\un{\chi}$ is constant on the cosets mod $U_{p,n}$, hence
\[
\int_{\Z_p^\times} \chi_p(p^{-n}u) \un{\chi}(u) d^\times u 
= \sum_{i = 1}^r \un{\chi}(e_i) \int_{e_i U_{p,n}} \chi_p(p^{-n}u) d^\times u.
\]
But
\[
u \in e_i U_{p,n} \implies p^{-n} u \in p^{-n} e_i + \Z_p
\]
$\implies$
\begin{align*}
\chi_p(p^{-n} u) 	
&=\   \chi_p(p^{-n} e_i + x) \qquad ( x \in \Z_p)\\	
&=\ \chi_p(p^{-n} e_i).						
\end{align*}
Therefore
\begin{align*}
\int_{\Z_p^\times} \chi_p(p^{-n} u)  \un{\chi}(u) d^\times u 	
&=\  \sum_{i = 1}^r \un{\chi}(e_i)  \chi_p(p^{-n}e_i) \int_{e_i U_{p,n}}  d^\times u\\	
&=\  \tau(\chi) \int_{U_{p,n}}  d^\times u							
\end{align*}
if
\[
 \tau(\chi) = \sum_{i = 1}^r \un{\chi}(e_i)  \chi_p(p^{-n}e_i).
\]
And
\begin{align*}
\int_{U_{p,n}} d^\times u  \ 
&=\  \int_{1 + p^n \Z_p} d^\times u \\	
&=\  \frac{p}{p-1}\int_{1 + p^n \Z_p} \frac{du}{\abs{u}_p} \\
&=\  \frac{p}{p-1}\int_{1 + p^n \Z_p} du \\
&=\  \frac{p}{p-1}\int_{p^n \Z_p} du \\		
&=\  \frac{p}{p-1} p^{-n } \\
&=\  \frac{p^{1-n}}{p-1}.				
\end{align*}
So in the end
\[
Z(\phi_n, \chi) = \tau(\chi) \frac{p^{1 + n(s - 1)}}{p-1}.
\]
Next
\begin{align*}
\widehat{\phi}_n(t) \ 
&=\  \int_{\Q_p} \phi_n(x) \chi_p(tx) dx \\	
&=\  \int_{\Q_p} \chi_p(x) \chi_{p^{-n}\Z_p (x)}  \chi_p(tx) dx \\
&=\  \int_{p^{-n}\Z_p}  \chi_p(x)  \chi_p(tx) dx \\
&=\  \int_{p^{-n}\Z_p} \chi_p((1+t)x) dx \\	
&=\  \vol_{dx} (p^{-n} \Z_p) \chi_{p^n \Z_p-1} (t)\\	
&=\  p^n\chi_{p^n \Z_p-1} (t).			
\end{align*}

Therefore 
\begin{align*}
Z(\widehat{\phi}_n, \widecheck{\chi}) \ 
&=\  \int_{\Q_p^\times} \widehat{\phi}_n(x) \widecheck{\chi}(x) d^\times x \\		
&=\  \int_{\Q_p^\times} p^n \chi_{p^n \Z_p-1} (x)  \chi^{-1}(x) \abs{x}_p d^\times x \\	
&=\  p^n\int_{p^n\Z_p-1} \ov{\un{\chi}(x)} \abs{x}_p^{1-s}d^\times x \\	
&=\  p^n\int_{p^n\Z_p-1} \ov{\un{\chi}(x)} d^\times x \\
&=\  p^n\int_{1 + p^n\Z_{p}} \ov{\un{\chi}(-x)} d^\times x \\	
&=\  p^n \ov{\chi(-1)} \int_{1 + p^n\Z_{p}} \ov{\un{\chi}(x)} d^\times x \\
&=\  p^n \chi(-1) \int_{U_{p,n}} d^\times x \\	
&=\  p^n \chi(-1) \frac{p^{1-n}}{p-1} \\
&=\  \frac{p}{p-1} \chi(-1).
\end{align*}

\vspace{0.1cm}

[Note: \ $\chi(-1) = \pm 1:$
\[
1 = (-1)(-1) \implies 1 = \chi(-1)\chi(-1) = \chi(-1)^2.]
\]
Assembling the data then gives
\begin{align*}
\rho(\chi) \ 
&=\  \frac{Z(\phi_n, \chi)}{Z(\widehat{\phi}_n, \widecheck{\chi})} \\		
&=\  \frac{\tau(\chi) \ds\frac{p^{1 + n(s-1)}}{p-1} }    {\ds\frac{p}{p-1}\chi(-1)}\\
&=\  \tau(\chi) \ds\frac{p^{1 + n(s-1)}}{p-1}  \frac{p-1}{p\chi(-1)}\\
&=\  \tau(\chi) \chi(-1) p^{n(s-1)}\\
&=\  \tau(\chi) \chi(-1) p^{n(s-1)} \ds\frac{1}{1}\\
&=\  \tau(\chi) \chi(-1) p^{n(s-1)} \ds\frac{L(\chi)}{L(\widecheck{\chi})}.
\end{align*}
\end{x}

\vspace{0.1cm}

\begin{x}{\small\bf THEOREM} \ 
\[\rho(\chi) = \epsilon(\chi) \frac{L(\chi)}{L(\widecheck{\chi})}, \quad \text{ where } \epsilon(\chi)  \ =\left\{
\begin{array}{l l}
1 \quad \text{ }  \text{ if $\chi$ is unramified}\\
\rho(\chi)  \text{ if $\chi$ is ramified of degree $n \ge 1$.}\\
\end{array}
\right.\]
\end{x}

\vspace{0.2cm}

\begin{x}{\small\bf LEMMA} \ 
Suppose that $\chi$ is ramified of degree $n \ge 1$ $-$then
\[
 \epsilon(\chi)  \epsilon(\widecheck{\chi}) = \chi(-1).
\]

\vspace{0.1cm}

PROOF \  
$\forall$ $f \in \sB(\Q_p)$, 
\begin{align*}
Z(f, \chi) 	\ 
&=\   \epsilon(\chi) Z(\widehat{f}, \widecheck{\chi}) \\		
&=\   \epsilon(\chi)  \epsilon(\widecheck{\chi}) Z(\widehat{\widehat{f}\hspace{.125cm}}, \widecheck{\widecheck{\chi}}).
\end{align*}
But $\widecheck{\widecheck{\chi}} = \chi$, hence
\begin{align*}
Z(\widehat{\widehat{f}\hspace{.125cm}}, \widecheck{\widecheck{\chi}}) \ 	
&=\   \int_{\Q_p^\times} \widehat{\widehat{f}\hspace{.125cm}}(x) \chi(x) d^\times x \\		
&=\   \int_{\Q_p^\times} f(-x)\chi(x) d^\times x \\	
&=\   \int_{\Q_p^\times} f(x)\chi(-x) d^\times x \\
&=\   \chi(-1) \int_{\Q_p^\times} f(x)\chi(x) d^\times x \\
&=\   \chi(-1) Z(f, \chi).
\end{align*}
\end{x}

\vspace{0.1cm}

\begin{x}{\small\bf APPLICATION} \ 
\[
\tau(\chi) \tau(\widecheck{\chi}) = p^n \chi(-1).
\]

[In fact,
\begin{align*}
\epsilon(\chi)  \epsilon(\widecheck{\chi}) \ 
&= \  \tau(\chi)p^{n(s-1)} \chi(-1) \tau(\widecheck{\chi}) p^{n(1 - s - 1)} \widecheck{\chi}(-1) \\
&= \  \tau(\chi) \tau(\widecheck{\chi}) p^{-n} \\
&= \  \chi(-1)
\end{align*}
\qquad\qquad $\implies$
\[
\tau(\chi) \tau(\widecheck{\chi}) = p^n \chi(-1).]
\]
\end{x}

\vspace{0.1cm}

\begin{x}{\small\bf LEMMA} \ 
Suppose that $\chi$ is ramified of degree $n \ge 1$ $-$then
\[
\epsilon(\ov{\chi}) = \chi(-1) \ov{\epsilon(\chi)}.
\]

\vspace{0.1cm}

PROOF \ 
$\forall$ $f \in \sB(\Q_p)$, 
\begin{align*}
Z(\widehat{\ov{f}}, \chi) \ 
&=\   \int_{\Q_p^\times} \widehat{\ov{f}}(x) \chi(x) d^\times x \\		
&=\   \int_{\Q_p^\times} \ov{\widehat{f}(-x)}\chi(x) d^\times x  \qquad (\text{cf.} \ \S 10, \ \#12)\\	
&=\   \int_{\Q_p^\times} \ov{\widehat{f}(x)}\chi(-x) d^\times x \\
&=\   \chi(-1) \int_{\Q_p^\times} \ov{\widehat{f}(x)}\chi(x) d^\times x \\
&=\   \chi(-1) Z(\ov{\widehat{f}}, \chi).
\end{align*}

But $\widecheck{\ov{\chi}} = \ov{\widecheck{\chi}}$, hence
\begin{align*}
\ov{Z(f, \chi)}	\ 
&=\   Z(\ov{f}, \ov{\chi}) \\		
&=\   \epsilon(\ov{\chi}) Z(\widehat{\ov{f}}, \widecheck{\ov{\chi}}) \\	
&=\   \epsilon(\ov{\chi}) Z(\widehat{\ov{f}}, \ov{\widecheck{\chi}}) \\
&=\   \epsilon(\ov{\chi}) \chi(-1) Z(\ov{\widehat{f}}, \ov{\widecheck{\chi}}) \\	
&=\   \epsilon(\ov{\chi}) \chi(-1) \ov{Z(\widehat{f}, \widecheck{\chi})}.	
\end{align*}
On the other hand,
\begin{align*}
\ov{Z({f}, {\chi})} \ 
&=\   \ov{\epsilon(\chi) Z(\widehat{f}, \widecheck{\chi})}\\		
&=\   \ov{\epsilon(\chi)} \ov{Z(\widehat{f}, \widecheck{\chi})}.
\end{align*}
Therefore 
\[
\epsilon(\ov{\chi}) \chi(-1) = \ov{\epsilon(\chi)}
\]
\qquad\qquad\qquad$\implies$
\[
\epsilon(\ov{\chi}) = \chi(-1)\ov{\epsilon(\chi)}.
\]
\end{x}

\vspace{0.1cm}

\begin{x}{\small\bf APPLICATION} \ 
\[
\tau(\ov{\chi}) = \chi(-1)\ov{\tau(\chi)}.
\]

[In fact,
\begin{align*}
\epsilon(\ov{\chi})  	\ 
&=\   \tau(\ov{\chi}) p^{n(\ov{s} - 1)} \ov{\chi}(-1)\\	
&=\   \chi(-1)\ov{\epsilon(\chi)}\\
&=\   \chi(-1)\ov{\tau(\chi)}p^{n(\ov{s} - 1)} \ov{\chi(-1)}\\
&=\   \chi(-1)\ov{\tau(\chi)}p^{n(\ov{s} - 1)} \ov{\chi}(-1)
\end{align*}
\qquad\qquad\qquad$\implies$
\[
\tau(\ov{\chi}) = \chi(-1)\ov{\tau(\chi)}.]
\]
\end{x}

\vspace{0.1cm}

\begin{x}{\small\bf DEFINITION} \ 
Let $\un{\chi} \in \widehat{\Z_p^\times}$ be a nontrivial unitary character $-$then its \un{root number} 
$W(\un{\chi}) $ is prescribed by the relation
\[
W(\un{\chi}) = \epsilon(\acdot_p^{1/2} \un{\chi}).
\]

[Note: \  If $\un{\chi}$ is trivial, then $W(\un{\chi}) = 1.]$
\end{x}

\vspace{0.1cm}

\begin{x}{\small\bf LEMMA} \ 
\[
\abs{W(\un{\chi})} = 1.
\]

\vspace{0.1cm}

PROOF \ 
Put $\chi = \acdot_p^{1/2} \un{\chi}$ $-$then
\[
\epsilon(\chi)\epsilon(\widecheck{\chi}) = \chi(-1) 	\qquad \text{(cf. } \# 12)
\]
$\implies$
\begin{align*}
\epsilon(\chi)^{-1}  \ 	
&=\   \epsilon(\widecheck{\chi}) \chi(-1)^{-1}\\	
&=\   \epsilon(\widecheck{\chi}) \chi(-1)\\
&=\   \epsilon(\ov{\chi}) \chi(-1) \qquad (\widecheck{\chi} = \ov{\chi})\\
&=\   \chi(-1) \ov{\epsilon(\chi)} \chi(-1) \qquad (\text{cf.} \  \#14)\\
&=\   \chi(-1)^2 \ \ov{\epsilon(\chi)}\\
&=\   \ov{\epsilon(\chi)}.
\end{align*}
$\implies$
\[
\abs{\epsilon(\chi)} = 1 \implies \abs{W(\un{\chi})} = 1.
\]
\end{x}

\vspace{0.1cm}

\begin{x}{\small\bf APPLICATION} \ 
\[
\abs{\tau(\acdot_p^{1/2} \chi)} = p^{n/2}.
\]

[In fact,
\[
1 = \abs{W(\un{\chi}) } = \abs{\tau(\acdot_p^{1/2} \un{\chi}) p^{n(\frac{1}{2} - 1)} }.
\]
\end{x}

\vspace{0.2cm}

\begin{x}{\small\bf EXERCIZE AD LIBITUM} \ 
Show that the theory expounded above for $\Q_p$ can be carried over to any finite extension $\K$ of $\Q_p$.
\end{x}

\chapter{
$\boldsymbol{\S}$\textbf{13}.\quad  RESTRICTED PRODUCTS}
\setlength\parindent{2em}
\setcounter{theoremn}{0}

\ \indent 
Recall:

\begin{x}{\small\bf FACT} \ 
Suppose that $X_i$ $(i \in I)$ is a nonempty Hausdorff space $-$then the product 
$\prod\limits_{i \in I} X_i$ is locally compact iff each 
$X_i$ is locally compact and all but a finite number of the $X_i$ are compact.
\end{x}

\vspace{0.1cm}

Let $X_i \  (i \in I)$ be a family of nonempty locally compact Hausdorff spaces and for each 
$i \in I$, let $K_i \subset X_i$ be an open-compact subspace.
\begin{x}{\small\bf DEFINITION} \ 
The 
\underline{restricted product}
\index{restricted product}
\[
\prod\limits_{i \in I} (X_i : K_i)
\]
consists of those $x = \{x_i\}$ in $\prod\limits_{i \in I} X_i$ such that $x_i \in K_i$ for all but a finite number of $i \in I$. 
\end{x}

\vspace{0.1cm}

\begin{x}{\small\bf \un{N.B.}} \ 
\[
\prod_{i \in I} (X_i : K_i) = \bigcup_{S \subset I} \ \prod_{i \in S} X_i \times \prod_{i \notin S} K_i,
\]
where $S \subset I$ is finite.
\end{x}

\vspace{0.1cm}

\begin{x}{\small\bf DEFINITION} \ 
A 
\underline{restricted open rectangle}
\index{restricted open rectangle} 
is a subset of $\prod\limits_{i \in I} (X_i : K_i)$ of the form
\[
\prod_{i \in S} U_i \times \prod_{i \notin S} K_i,
\]
where $S \subset I$ is finite and $U_i \subset X_i$ is open.
\end{x}

\vspace{0.1cm}

\begin{x}{\small\bf LEMMA} \ 
The intersection of two restricted open rectangles is a restricted open rectangle.
\end{x}

\vspace{0.1cm}

Therefore the collection of restricted open rectangles is a basis for a topology on $\prod\limits_{i \in I} (X_i : K_i)$, the \underline{restricted product topology}.  
\index{restricted product topology}

\vspace{0.2cm}

\begin{x}{\small\bf LEMMA} \ 
If $I$ is finite, then
\[
\prod_{i \in I} X_i = \prod_{i \in I} (X_i : K_i)
\]
and the restricted product topology coincides with the product topology.
\end{x}

\vspace{0.1cm}

\begin{x}{\small\bf LEMMA} \ 
If $I = I_1 \cup I_2$, with $I_1 \cap I_2 = \emptyset$, then
\[
 \prod_{i \in I} (X_i : K_i) \ \approx \  \bigl( \prod_{i \in I_1} (X_i : K_i) \bigr) \times \bigl(\prod_{i \in I_2} (X_i : K_i)\bigr),
\]
the restricted product topology on the left being the product topology on the right.
\end{x}

\vspace{0.1cm}

\begin{x}{\small\bf LEMMA} \ 
The inclusion $\prod\limits_{i \in I} (X_i : K_i)  \hookrightarrow \prod\limits_{i \in I} X_i$ is continuous but the restricted product topology coincides with the relative topology only if $X_i = K_i$ for all but a finite number of $i \in I$.
\end{x}
\vspace{0.1cm}

\begin{x}{\small\bf LEMMA} \ 
$\ds\prod\limits_{i \in I} (X_i : K_i)$ is a Hausdorff space.

\vspace{0.1cm}

PROOF \   
Taking into account \#8, this is because

1. A subspace of a Hausdorff space is Hausdorff;

2. Any finer topology on a Hausdorff space is Hausdorff.
\end{x}

\vspace{0.1cm}

\begin{x}{\small\bf LEMMA} \ 
$\prod_{i \in I} (X_i : K_i)$ is a locally compact Hausdorff space.

\vspace{0.1cm}

PROOF \  
Let $x \in \prod\limits_{i \in I} (X_i:K_i)$ $-$then there exists a finite set $S \subset I$ such that $x_i \in K_i$ if $i \notin S$.  
Next, for each $i \in S$, choose a compact neighborhood $U_i$ of $x_i$.  
This done, consider
\[
\prod_{i \in S} U_i \times \prod_{i \notin S} K_i,
\]
a compact neighborhood of $x$.
\end{x}

\vspace{0.1cm}

From this point forward, it will be assumed that $X_i \equiv G_i$ is a locally compact abelian group and $K_i \subset G_i$ is an open-compact subgroup.
\vspace{0.2cm}

\begin{x}{\small\bf NOTATION} \ 
\[
G = \prod_{i \in I} \ (G_i:K_i).
\]
\end{x}

\vspace{0.1cm}

\begin{x}{\small\bf LEMMA} \ 
$G$ is a locally compact abelian group.  
\end{x}

\vspace{0.1cm}

Given $i \in I$, there is a canonical arrow
\begin{align*}
\ins_i : G_i	&\ra G\\	
x			&\mapsto (\cdots, 1, 1, x, 1, 1, \cdots).\\	
\end{align*}

\begin{x}{\small\bf LEMMA} \ 
$\ins_i$ is a closed embedding.

\vspace{0.1cm}

PROOF \ 
Take $S = \{i\}$ and pass to
\[
G_i \times \prod_{j \ne i} K_j,
\]
an open, hence closed subgroup of G.  The image $\ins_i(G_i)$ is a closed subgroup of
\[
G_i \times \prod_{j \ne i} K_j
\]
in the product topology, hence in the restricted product topology.\\

Therefore $G_i$ can be regarded as a closed subgroup of $G$.
\end{x}

\vspace{0.1cm}

\begin{x}{\small\bf LEMMA} \ 

1. Let $\chi \in \widetilde{G}$ 
$-$then $\chi_i = \chi \circ \ins_i = \restr{\chi}{G_i} \in \widetilde{G}_i$ and 
$\restr{\chi}{K_i} \equiv 1$ for all but a finite number of $i \in I$, so for each $x \in G$, 
\[
\chi(x) \ =\  \chi(\{x_i\}) \ =\  \prod_{i \in I} \chi_i(x_i).
\]

2.  Given $i \in I$, let $\chi_i \in \widetilde{G}_i$ and assume that $\restr{\chi}{K_i} \equiv 1$ 
for all but a finite number of $i \in I$ $-$then the prescription 
\[
\chi(x) \ =\  \chi(\{x_i\}) \ =\  \prod_{i \in I} \chi_i(x_i)
\]
defines a $\chi \in \widetilde{G}$.
\end{x}

\vspace{0.1cm}

These observations also apply if $\widetilde{G}$ is replaced by $\widehat{G}$, in which case more can be said.
\vspace{0.2cm}

\begin{x}{\small\bf THEOREM} \ 
As topological groups,
\[
\widehat{G} \thickapprox \prod_{i \in I} \  (\widehat{G}_i:K_i^\perp).
\]

[Note: \   Recall that
\[
K_i^\perp = \{\chi_i \in \widehat{G}_i:\chi{|K_i} \equiv 1\} \qquad (\text{cf.} \  \S7, \  \#32)
\]
and a tacit claim is that $K_i^\perp$ is an open-compact subgroup of $\widehat{G}$.  
To see this, 
quote \S7, \#34 to get
\[
\widehat{K}_i \thickapprox \widehat{G} / K_i^\perp, \quad K_i^\perp \thickapprox \widehat{G / K_i}.
\]
Then

\qquad \textbullet \quad
$K_i$ compact $\implies \widehat{K}_i$ discrete $\implies \widehat{G}/K_i^\perp$ discrete $\implies K_i^\perp$ open.

\qquad \textbullet \quad
$K_i$ open $\implies G/K_i$ discrete $\implies \widehat{G/K_i}$ compact $\implies K_i^\perp$ compact.$]$
\end{x}

\vspace{0.1cm}

Let $\mu_i$ be the Haar measure on $G_i$ normalized by the condition
\[
\mu_i(K_i) = 1.
\]
\begin{x}{\small\bf LEMMA} \ 
There is a unique Haar measure $\mu_G$ on $G$ such that for every finite subset $S \subset I$, the restriction of $\mu_G$ to
\[
G_S \ \equiv \ \prod_{i \in S} \ G_i \ \times  \ \prod_{i \notin S}\  K_i
\]
is the product measure.
\end{x}

\vspace{0.1cm}

Suppose that $f_i$ is a continuous, integrable function on $G_i$ such that $\restr{f_i}{K_i} = 1$ 
for all $i$ outside some finite set and let $f$ be the function on $G$ defined by
\[
f(x) = f(\{x_i\}) = \prod_i f_i(x_i).
\]
Then $f$ is continuous.
Proof: The $G_S$ are open and cover $G$ and on each of them $f$ is continuous.

\vspace{0.2cm}

\begin{x}{\small\bf LEMMA} \ 
Let $S \subset I$ be a finite subset of $I$ $-$then
\[
\int_{G_S} f(x) d\mu_{G_S} (x) \ =\  \prod_{i \in S} \  \int_{G_i} f_i(x_i)d\mu_{G_i}(x_i).
\]
\end{x}

\vspace{0.1cm}

\begin{x}{\small\bf APPLICATION} \ 
If
\[
\sup_S \ \prod_{i \in S} \  \int_{G_i} \  \abs{f_i(x_i)} d\mu_{G_i}(x_i) \ < \  \infty,
\]
then $f$ is integrable on $G$ and
\[
\int_G f(x) d\mu_{G} (x) \ = \ \prod_{i \in I} \ \int_{G_i} f_i(x_i)d\mu_{G_i}(x_i).
\]
\end{x}

\begin{x}{\small\bf EXAMPLE} \ 
Take  $f_i = \chi_{K_i}$ (which is continuous, $K_i$ being open-compact) $-$then $\widehat{f}_i = \chi_{K_i^\perp}.$  
Setting
\[
f \ = \ \prod_{i \in I} f_i,
\]
it thus follow that $\forall$ $\chi \in \widehat{G}$,
\[
\widehat{f}(\chi) \ = \  \prod_{i \in I} \widehat{f}_i(\chi_i).
\]
\end{x}

\vspace{0.1cm}

Working within the framework of \S7, \#45, let $\mu_{\widehat{G}_i}$ be the Haar measure on $\widehat{G}_i$ per Fourier inversion.

\vspace{0.1cm}

\begin{x}{\small\bf LEMMA} \ 
\[
\mu_{\widehat{G}_i} (K_i^\perp) = 1.
\]

\vspace{0.1cm}

PROOF \ 
Since $\chi_{K_i} \in \mathbf{INV}(G_i)$, $\forall$ $x_i \in G_i$,
\begin{align*}
\chi_{K_i} (x_i) \ 
&= \  \int_{\widehat{G}_i } \widehat{\chi}_{K_i} (x_i) \overline{\chi_i (x_i)} d \mu_{\widehat{G}_i} (\chi_i) \\
&= \  \int_{K_i^\perp} \overline{\chi_i (x_i)} d \mu_{\widehat{G}_i} (\chi_i).
\end{align*}
Now set $x_i = 1$ to get
\begin{align*}
1 
&=\  \int_{K_i^\perp} d\mu_{\widehat{G}_i} (\chi_i) \\
&=\  \mu_{\widehat{G}_i}(K_i^\perp).
\end{align*}
\end{x}

\vspace{0.1cm}

Let $\mu_{\widehat{G}}$ be the Haar measure on $\widehat{G}$ constructed as in \#16  (i.e., replace $G$ by 
$\widehat{G}$, bearing in mind \#20).

\vspace{0.2cm}

\begin{x}{\small\bf LEMMA} \ 
$\mu_{\widehat{G}}$ is the Haar measure on $\widehat{G}$ figuring in the Fourier inversion per $\mu_G$.

\vspace{0.1cm}

PROOF \  Take
\[
f \ =\  \prod_{i \in I} f_i,
\]
where $f_i = \chi_{K_i}$ (cf. \#19 ) $-$then
\begin{align*}
\int_{\widehat{G}} \widehat{f}(\chi) \overline{\chi(x)} d\mu_{\widehat{G}} (\chi) \ 
&=\  \prod_{i \in I} \ \int_{\widehat{G}_i} \widehat{f}_i(\chi_i) \overline{\chi_i(x_i)} d\mu_{\widehat{G}_i} (\chi_i) \\	
&=\  \prod_{i \in I}\  f_i(x_i) \\	
&=\  f(\{x_i\}) \\
&=\  f(x).
\end{align*}
\end{x}

\chapter{
$\boldsymbol{\S}$\textbf{14}.\quad  ADELES AND IDELES}
\setlength\parindent{2em}
\setcounter{theoremn}{0}

\begin{x}{\small\bf DEFINITION} \ 
The set of 
\underline{finite adeles}
\index{finite adeles} 
is the restricted product
\[
\A_{\fin} \ = \ \prod\limits_p (\Q_p: \Z_p).
\]
\end{x}

\vspace{0.1cm}

\begin{x}{\small\bf DEFINITION} \ 
The set of 
\underline{adeles}
\index{adeles} 
is the product
\[
\A \ = \ \A_{\fin} \times \R.
\]
\end{x}

\vspace{0.1cm}

\begin{x}{\small\bf LEMMA} \ 
$\A$ is a locally compact abelian group (under addition).
\end{x}

\vspace{0.1cm}

\begin{x}{\small\bf \un{N.B.}} \ 
$\A$ is a subring of \quad $\prod\limits_p \Q_p \times \R$.
\end{x}

\vspace{0.1cm}

The image of the diagonal map
\[
\Q \ra \prod\limits_p \Q_p \times \R
\]
lies in $\A$, so $\Q$ can be regarded as a subring of $\A$.

\vspace{0.1cm}

\begin{x}{\small\bf LEMMA} \ 
$\Q$ is a discrete subspace of $\A$.

\vspace{0.1cm}

PROOF \ 
To establish the discreteness of $\Q \subset \A$, one need only exhibit a neighborhood \mU of 0 in $\A$ such that $\Q \cap U = \{0\}$.  
To this end, consider
\[
U \ = \  \prod\limits_p \\ Z_p \ \times \  ] -\frac{1}{2},\frac{1}{2}[.
\]
If $x \in \Q \cap U$, then $\abs{x}_p \le 1$ $\forall \ p$.  
But \ $\bigcap\limits_p (\Q \cap \Z_p) = \Z$, so $x \in \Z$.  
And further, $\abs{x}_{\infty} < \frac{1}{2}$, hence finally $x = 0$.
\end{x}

\vspace{0.1cm}


\begin{x}{\small\bf FACT} \ 
Let G be a locally compact group and let $\Gamma \subset G$ be a discrete 
subgroup $-$then $\Gamma$ is closed in G and $G/\Gamma$ is a locally compact Hausdorff space.
\end{x}



\begin{x}{\small\bf THEOREM} \ 
The quotient $\A/\Q$ is a compact Hausdorff space.

\vspace{0.01cm}

PROOF \  
Since $\Q \subset \A$ is a discrete subgroup, $\Q$ must be closed in $\A$ and the quotient  $\A/\Q$ must be Hausdorff.  
As for compactness, it suffices to show that the compact set $\prod\limits_p \Z_p \times [0,1]$ contains a set of representatives of  
$\A/\Q$ because this implies that the projection
\[
\prod\limits_p \Z_p \times [0,1] \ra \A/\Q
\]
is surjective, hence that  $\A/\Q$ is the continuous image of a compact set.  
So let $x \in \A$ $-$then there is a finite set $S$ of primes such that $p \notin S \implies x_p \in \Z_p$.  
For $p \in S$, write
\[
x_p = f(x_p) + [x_p], 
\]
thus $[x_p] \in \Z_p$ and if $q \ne p$ is another prime,
\begin{align*}
\abs{f(x_p)}_q \ 
&=\  \abs{\sum_{n = v(x_p)}^{-1} a_n p^n}_q \\
&\le\  \sup\{\abs{a_np^n}_q\} \\
&\le\  1.
\end{align*}
Agreeing to denote $f(x_p)$ by  $r_p$, write
\[
x = (x - r_p) + r_p. 
\]
Then $r_p$ is a rational number and per $x - r_p$, $S$ reduces to $S - \{p\}$.  
Proceed from here by iteration to get 
\[
\allowdisplaybreaks
x = y + r, 
\]
where $\forall$ $p$, $y_p \in \Z_p$, and $r \in \Q$.
At infinity,
\[
\allowdisplaybreaks
x_\infty = y_\infty + r \quad (r_\infty = r)
\]
and there is a unique $k \in \Z$ such that
\[
y_\infty = (y_\infty - k) + k 
\]
with $ 0 \le y_\infty - k < 1$.  
Accordingly,
\[
y = y + r = (y - k) + k + r.
\]
And
\[
\forall \ p, \quad (y - k)_p = y_p - k_p = y_p - k \in \Z_p, 
\]
while
\[
x_\infty = (y_\infty - k) + k + r.
\]
It therefore follows that $x$ can be written as the sum of an element in $\prod\limits_p \Z_p \times [0,1]$ and a rational number, the contention.
\end{x}

\vspace{0.1cm}

\begin{x}{\small\bf DEFINITION} \ 
The topological group $\A/\Q$ is called the 
\underline{adele class group}.
\index{adele class group}
\end{x}

\vspace{0.1cm}

\begin{x}{\small\bf DEFINITION} \ 
Let $G$ be a locally compact group and let $\Gamma \subset G$ be a discrete subgroup $-$then a 
\underline{fundamental domain}
\index{fundamental domain} 
for $G/\Gamma$ is a Borel measurable subset $D \subset G$ which is a system of representatives for $G/\Gamma$.
\end{x}

\vspace{0.1cm}

\begin{x}{\small\bf LEMMA} \ 
The set
\[
D = \prod\limits_p \Z_p \times [0,1[
\]
is a fundamental domain for $\A/\Q$.

\vspace{0.1cm}

PROOF \  
The claim is that every $x \in \A$ can be written uniquely as $d + r$, where $d \in D, r \in \Q$.  
The proof of \#7 settles existence, thus the remaining issue is uniqueness:
\[
d_1 + r_1 = d_2 + r_2 \implies d_1 = d_2,\  r_1 = r_2
\]
To see this, consider
\[
\rho = d_1 - d_2 = r_2 - r_1 \in (D-D) \cap \Q.
\]

\qquad\qquad \textbullet \quad $\forall$ p, $\rho = \rho_p \in D_p - D_p = D_p = \Z_p$ 

\qquad\qquad\qquad $\implies \rho \in \bigcap\limits_p \ (\Q \cap \Z_p) = \Z$.

\qquad\qquad \textbullet \quad $\rho = \rho_\infty \in D_\infty - D_\infty = \  ]-1,1[.$

Therefore
\[
\rho \in  \Z \  \cap \   ]-1,1[ \   \implies \   \rho = 0.
\]
\end{x}

\vspace{0.1cm}

\begin{x}{\small\bf REMARK} \ 
$\Q$ is dense in $\A_\fin$.

\vspace{0.1cm}

[The point is that $\Z$ is dense in $\ds\prod\limits_p \Z_p$.$]$
\end{x}

\vspace{0.1cm}

\begin{x}{\small\bf DEFINITION} \ 
The set of 
\underline{finite ideles}
\index{finite ideles} 
is the restricted product
\[
\I_{\fin}\ = \ \prod\limits_p (\Q_p^\times:\Z_p^\times).
\]
\end{x}

\vspace{0.1cm}

\begin{x}{\small\bf DEFINITION} \ 
The set of 
\underline{ideles}
\index{ideles} 
is the product
\[
\I = \I_{\fin} \times \R^\times.
\]
\end{x}

\vspace{0.1cm}

\begin{x}{\small\bf LEMMA} \ 
$\I$ is a locally compact abelian group $($under multiplication$)$.
\end{x}

\vspace{0.1cm}

Algebraically, $\I$  can be identified with $\A^\times$ but there is a topological issue since when endowed with the relative topology, $\A^\times$ is not a topological group: Multiplication is continuous but inversion is not continuous.

\begin{x}{\small\bf LEMMA} \ 
Equip $\A \times \A$ with the product topology and define
\begin{align*}
\phi : \I 	&\ra \A \times \A\\	
x 		&\mapsto \bigl(x, \frac{1}{x}\bigr).
\end{align*}
Endow the image $\phi(\I)$ with the relative topology $-$then $\phi$ is a topological isomorphism of $\I$ onto $\phi(\I)$.
\end{x}

\vspace{0.1cm}

The image of the diagonal map
\[
\Q^\times \lra \prod\limits_p \Q_p \times \R^\times
\]
lies in $\I$, so $\Q^\times$ can be regarded as a subgroup of $\I$.

\vspace{0.2cm}

\begin{x}{\small\bf LEMMA} \ 
$\Q^\times$ is a discrete subspace of $\I$.

\vspace{0.1cm}

PROOF \  
$\Q$ is a discrete subspace of $\A$ (cf. \#5), hence  $\Q \times \Q$ is a discrete subspace of $\A \times \A$, 
hence $\phi(\Q^\times)$ is a discrete subspace of $\phi(\I)$.
\end{x}

\vspace{0.1cm}

Consequently, $\Q^\times$ is a closed subgroup of $\I$ and the quotient $\I/\Q^\times$ is a locally compact Hausdorff space but, as opposed to the adelic situation, it is not compact (see below).

\vspace{0.2cm}

\begin{x}{\small\bf DEFINITION} \ 
The topological group $\I/\Q^\times$ is called the 
\underline{idele class group}.
\index{idele class group}
\end{x}

\vspace{0.1cm}

\begin{x}{\small\bf NOTATION} \ 
Given $x \in \I$, put
\[
\abs{x}_\A = \prod_{p \le \infty} \abs{x_p}_p.
\]
Extend the definition of $\acdot_\A$ to all of $\A$ by setting $\abs{x}_\A = 0$ if $x \in \A  - \A^\times$.
\end{x}

\vspace{0.1cm}

\begin{x}{\small\bf LEMMA} \ 
$\forall$ $x \in \Q^\times$, $\abs{x}_\A = 1$ (cf. \S1, \  \#21 ).
\end{x}

\vspace{0.1cm}

\begin{x}{\small\bf LEMMA} \ 
The homomorphism
\[
\acdot_\A : \I \ra \R_{> 0}^\times
\]
is continuous and surjective.

\vspace{0.1cm}

PROOF \ 
Omitting the verification of continuity, fix $t \in \R_{> 0}^\times$ and let $x$ be the idele specified by
\[x_p = \ 
\begin{cases}
\ 1 \quad \text{$(p < \infty)$}\\
\ t \quad \text{$(p = \infty)$}
\end{cases}
.\]
Then $\abs{x}_\A = t.$
\end{x}

\vspace{0.1cm}

\begin{x}{\small\bf SCHOLIUM} \ 
The idele class group $\I/\Q^\times$ is not compact.
\end{x}

\vspace{0.1cm}

\begin{x}{\small\bf NOTATION} \ 
Let
\[
\I^1  \ = \ \ker \acdot_\A.
\]
\end{x}

\vspace{0.1cm}

\begin{x}{\small\bf \un{N.B.}} \ 
$x \in \I^1 \implies x_\infty \in \Q^\times$.
\end{x}
\vspace{0.1cm}

\begin{x}{\small\bf THEOREM} \ 
The quotient $\I^1/\Q^\times$ is a compact Hausdorff space, in fact
\[
\I^1/\Q^\times \ \thickapprox \  \prod\limits_p \Z_p^\times,
\]
hence
\[
\prod\limits_p \Z_p^\times \times \{1\}
\]
is a fundamental domain for $\I^1/\Q^\times$.

\vspace{0.1cm}

PROOF \  
The arrow 
\[
\prod\limits_p \Z_p^\times \ra \I^1/\Q^\times
\]
that sends $x$ to $(x,1)\Q^\times$ is an isomorphism of topological groups.

[In obvious notation, the inverse is the map
\[
x \ = \  (x_{\fin}, x_\infty) \ra \frac{1}{x_\infty} x_{\fin}.]
\]
\end{x}
\vspace{0.1cm}

\begin{x}{\small\bf REMARK} \ 
$\forall$ p, $\Z_p^\times$ is totally disconnected.  
But a product of totally disconnected spaces is totally disconnected, thus $\prod\limits_p \Z_p^\times$ is totally disconnected, thus $\I^1/\Q^\times$ is totally disconnected.
\end{x}

\vspace{0.1cm}

\begin{x}{\small\bf \un{N.B.}} \ 
 $\prod\limits_p \Z_p^\times \times \R_{>0}^\times$ is a fundamental domain for $\I^1/\Q^\times$.

\vspace{0.1cm}

[Note: \ If $r \in \Q$ and if $\abs{r}_p = 1$ $\forall$ p, then $r = \pm1.]$
\end{x}

\vspace{0.1cm}

\begin{x}{\small\bf LEMMA} \ 
\[
\I \ \thickapprox \ \I^1 \times  \R_{>0}^\times.
\]

\vspace{0.1cm}

PROOF \  
The arrow
\[
\I \ra \I^1 \times  \R_{>0}^\times
\]
that sends x to $(\widetilde{x}, \abs{x}_\A)$, where 

\[(\widetilde{x})_p = \ 
\begin{cases}
x_p \quad \quad (p < \infty)\\
\ds\frac{x_\infty}{\abs{x}_\A} \  \quad (p = \infty)
\end{cases}
,\]
is an isomorphism of topological groups.
\end{x}

\vspace{0.1cm}

\begin{x}{\small\bf LEMMA} \ 
There is a disjoint decomposition
\[
\I_{\fin} = \coprod_{q \in \Q_{>0}^\times} q\bigl(\prod\limits_p \Z_p^\times\bigr).
\]

\vspace{0.1cm}

PROOF \  The right hand side is obviously contained in the left hand side.  
To go the other way, fix an $x \in \I_{\fin}$ $-$then $\abs{x}_\A \in  \Q_{>0}^\times$.  
Moreover, 
$\abs{x}_\A x \in \I_{\fin}$ and $\forall \ p$, $\abs{ \abs{x}_\A x_p}_p = 1$ 
(for $x_p = p^ku \ (u \in \Z_p^\times) \implies \abs{x}_\A = p^{-k}r  \ (r \in \Q_p^\times$, $r$ coprime to $p$)), hence
\[
\abs{x}_\A x \in \prod\limits_p \Z_p^\times.
\]
Now write
\[
x = \abs{x}_\A^{-1} (\abs{x}_\A x)
\]
to conclude that
\[
x \in q \prod\limits_p \Z_p^\times \qquad (q = \abs{x}_\A^{-1} ).
\]
\end{x}

\vspace{0.1cm}

\begin{x}{\small\bf LEMMA} \ 
There is a disjoint decomposition
\[
\I_{\fin} \ \cap \ \prod\limits_p \Z_p \ = \ \coprod_{n \in \N} \ n\bigl(\prod\limits_p \Z_p^\times\bigr).
\]
\end{x}

\vspace{0.1cm}

Normalize the Haar measure $d^\times x$ on $\I_{\fin}$ by assigning the open-compact subgroup $\prod\limits_p \Z_p^\times$ total volume 1.
\vspace{0.2cm}

\begin{x}{\small\bf EXAMPLE} \ 
Suppose that $\Re (s) > 1$ $-$then
\allowdisplaybreaks
\begin{align*}
\int_{\I_{\fin} \hspace{0.03cm}\cap \hspace{0.03cm} \prod\limits_p \Z_p} \abs{x}_\A^s d^\times x 	\ 
&= \sum_{n \in \N} \ \int_{n(\prod\limits_p \Z_p^\times)} \abs{x}_\A^s d^\times x \\	
&= \sum_{n \in \N} \ \int_{\prod\limits_p \Z_p^\times} \abs{nx}_\A^s d^\times x \\
&= \sum_{n \in \N} n^{-s} \vol_{d^\times x} \bigl(\prod\limits_p \Z_p^\times\bigr) \\
&= \sum_{n \in \N} n^{-s} \\
&= \zeta(s). 
\end{align*}

\vspace{0.1cm}

[Note: \  Let $x \in \prod\limits_p \Z_p^\times:$\\

\qquad\qquad $\implies  \  \abs{x_p}_p \ = \ 1 	\quad  \forall p$, 
\begin{align*}
\implies \abs{nx}_\A \quad \ 	
&= \prod\limits_p \abs{nx_p}_p \\	
&= \prod\limits_p \abs{n}_p \abs{x_p}_p \\	
&= \prod\limits_p \abs{n}_p \\
&= \prod\limits_p \abs{n}_p \cdot n \cdot \frac{1}{n} \\
&= 1\cdot \frac{1}{n}\\
&= n^{-1}.]
\end{align*}
\end{x}

\vspace{0.1cm}

The idelic absolute value $\acdot_\A$ can be interpreted measure theoretically.
\vspace{0.1cm}

\begin{x}{\small\bf NOTATION} \ 
Write
\[
dx_\A = \prod_{p \le \infty} dx_p
\]
for the Haar measure $\mu_\A$ on $\A$ (cf. \S13, \#16).
\end{x}

\vspace{0.1cm}

Consider a function of the form $f = \prod\limits_{p \le \infty} f_p$, where $\forall$ $p, f_p$ is a continuous, 
integrable function on $\Q_p$ and for all but a finite number of $p$, 
$f_p = \chi_{\Z_p}$ 
$-$then
\[
\int_\A f(x)dx_\A = \prod_{p \le \infty}  \int_{\Q_p}	f_p(x_p) dx_p	\qquad (\text{cf.} \ \S13, \  \#18),
\]
it being understood that $\Q_\infty = \R$.

\vspace{0.1cm}

\begin{x}{\small\bf LEMMA} \ 
Let $M \subset \A$ be a Borel set with $0 < \mu_\A (M) < \infty$ $-$then $\forall$ $x \in \I$,
\[
\frac{\mu_\A (xM)}{\mu_\A (M)} = \abs{x}_\A.
\]

\vspace{0.1cm}

PROOF \   Take $M = D = \prod\limits_p \Z_p \times [0,1[$ \  (cf. \#10):

\allowdisplaybreaks
\begin{align*}
\mu_\A (xM) 	
&= \prod\limits_p \mu_{\Q_p} (x_p\Z_p) \times \mu_\R(x_\infty [0,1[)\\	
&= \prod\limits_p \abs{x_p}_p \mu_{\Q_p} (\Z_p) \times \abs{x_\infty} \mu_\R([0,1[)\\	
&= \prod\limits_p \abs{x_p}_p \times \abs{x_\infty}_\infty\\
&= \prod_{p \le \infty} \abs{x_p}_p \\
&= \abs{x}_\A.
\end{align*}

[Note: \  Needless to say, multiplication by an idele $x$ is an automorphism of $\A$, thus transforms $\mu_\A$ into a positive constant multiple of itself, the multiplier being $\abs{x}_\A.]$
\end{x}

\chapter{
$\boldsymbol{\S}$\textbf{15}.\quad  GLOBAL ANALYSIS}
\setlength\parindent{2em}
\setcounter{theoremn}{0}

\ \indent 
By definition,
\[
\A \ =\  \A_{\fin} \times \R.
\]
Therefore
\[
\widehat{\A} \ \thickapprox \  \widehat{\A}_{\fin} \times \widehat{\R}.
\]
And
\[
\A_{\fin} \ =\  \prod_p  \  (\Q_p: \Z_p)
\]
\qquad\qquad $\implies$
\[
\widehat{\A}_{\fin} \ \thickapprox \   \prod_p \  (\widehat{\Q}_p: \Z_p^\perp)	\qquad (\text{cf.} \  \S13, \  \#15).
\]
Put
\[
\chi_\Q \ =\  \prod_{p \le \infty} \chi_{p},
\]
where
\[
\chi_\infty \ =\  \exp( -2\pi \sqrt{-1}  \ x)  \qquad (x \in \R)	\qquad (\text{cf.}  \ \S8, \  \#27).
\]
Then
\[
\chi_\Q \in \widehat{\A}.
\]
Given $t \in \A$, define $\chi_{\Q, t} \in \widehat{\A}$ by the rule
\[
\chi_{\Q, t} (x) = \chi_\Q (tx).
\]
Then the arrow
\[
\Xi_\Q : \A \ra \widehat{\A}
\]
that sends $t$ to $\chi_{\Q, t}$ is an isomorphism of topological groups (cf.  \S8, \  \#24).

Recall now that $\forall$ $q \in \Q$,
\[
\chi_\Q(q) = 1	\qquad (\text{cf.} \  \S8, \  \#28).
\]
Accordingly, $\chi_\Q$ passes to the quotient and defines a unitary character of the adele class group $\A/\Q$.  
So, $\forall$ $q \in \Q$, $\chi_{\Q, q}$ is constant on the cosets of $\A/\Q$, 
thus it too determines an element of $\widehat{\A/\Q}$.

Equip $\Q$ with the discrete topology.
\vspace{0.1cm}

\begin{x}{\small\bf THEOREM} \ 
The induced map
\begin{align*}
\restr{\Xi_\Q}{\Q}: \Q  	
&\ra \widehat{\A/\Q}\\
q 		
&\mapsto \chi_{\Q, q}
\end{align*}
is an isomorphism of topological groups.

\vspace{0.1cm}

PROOF \  
Form $\Q^\perp \subset \widehat{\A}$, the closed subgroup of $\widehat{\A}$ consisting of those $\chi$ 
that are trivial on $\Q$ $-$then $\Q \subset \Q^\perp$ and $\widehat{\A/\Q} \thickapprox \Q^\perp.$  
But $\A/\Q$ is compact, thus its unitary dual $\widehat{\A/\Q}$ is discrete, thus $\Q^\perp$ is discrete.  
The quotient $\Q^\perp/\Q \subset \A/\Q$ $(\A \thickapprox \widehat{\A})$ is therefore discrete and closed, 
hence discrete and compact, hence finite.  
But $\Q^\perp/\Q$ is a $\Q$-vector space, so $\Q^\perp/\Q = \{0\}$ or still, $\Q^\perp = \Q$, 
which implies that $\Q \thickapprox \widehat{\A/\Q}$.
\end{x}

\vspace{0.1cm}

\begin{x}{\small\bf \un{N.B.}} \ 
There are two points of detail that have been tacitly invoked in the foregoing derivation.

\qquad \textbullet \quad $\Q^\perp/\Q$ in the quotient topology is discrete.  
Reason$:$  Let $S$ be an arbitrary nonempty subset of $\Q^\perp/\Q$, say 
$S = \{x\Q: x \in U\}$, \mU a subset of $\Q^\perp$ $-$then $U$ is automatically open $(\Q^\perp$ being discrete$)$, 
thus by the very definition of the quotient 
topology, $S$ is an open subset of $\Q^\perp/\Q$.

\qquad \textbullet \quad The quotient $\Q^\perp/\Q$ is closed in $\A/\Q$.  
Reason: $\Q^\perp$ is a closed subgroup of $\A$ containing $\Q$, so the following generality is applicable:  
If $G$ is a topological group, if $H$ is a subgroup of $G$, if $F$ is a closed subgroup of $G$ containing $H$, 
then $\pi(F)$ is closed in $G/H$ $(\pi : G \ra G/H$ the projection$)$.
\end{x}

\vspace{0.1cm}

\begin{x}{\small\bf SCHOLIUM} \ 
\[
\Q \thickapprox \widehat{\A/\Q} \implies \widehat{\Q} \thickapprox \widehat{\widehat{\A/\Q}} \thickapprox \A/\Q.
\]

\vspace{0.1cm}

[Note: \  Bear in mind that $\Q$ carries the discrete topology.]
\end{x}

\vspace{0.1cm}

\begin{x}{\small\bf DISCUSSION} \ 
Explicated, if $\chi \in \widehat{\Q}$, then there exists a $t \in \A$ such that $\chi = \chi_{\Q, t}$ and $\chi_{\Q, t_1} = \chi_{\Q, t_2}$ iff $ t_1 - t_2 \in \Q$.
\end{x}

\vspace{0.1cm}

\begin{x}{\small\bf DEFINITION} \ 
The 
\underline{Bruhat space}
\index{Bruhat space} 
$\sB(\A_{\fin})$
\index{$\sB(\A_{\fin})$} 
consists of all finite linear combinations of functions of the form
\[
f = \prod_p f_{p},
\]
where $\forall$ $p$, $f_p \in \sB(\Q_p)$ and $f_p = \chi_{\Z_p}$ for all but a finite number of $p$.
\end{x}

\vspace{0.1cm}

\begin{x}{\small\bf DEFINITION} \ 
The 
\underline{Bruhat-Schwartz space}
\index{Bruhat-Schwartz space} 
$\sB_\infty(\A)$
\index{$\sB_\infty(\A)$}  
consists of all finite linear combinations of functions of the form
\[
f = \prod_p f_{p} \times f_\infty,
\]
where 
\[
\prod_p f_p = \sB(\A_{\fin}) \text{ and } f_\infty \in \sS(\R).
\]
\end{x}

\vspace{0.1cm}

Given an $f \in \sB_\infty(\A)$, its Fourier transform is the function:
\begin{align*}
\widehat{f}:\A 	
&\ra \C\\	
t 		
&\mapsto \int_\A f(x) \chi_{\Q,t}(x) d\mu_\A(x) = \int_\A f(x) \chi_\Q(tx) d\mu_\A(x).	
\end{align*}

\vspace{0.1cm}

\begin{x}{\small\bf LEMMA} \ 
If
\[
f = \prod_p f_p \times f_\infty
\]
is a Bruhat-Schwartz function, then
\[
\widehat{f} = \prod_p \widehat{f}_p \times \widehat{f}_\infty.
\]
\end{x}
\vspace{0.1cm}

\begin{x}{\small\bf REMARK} \ 
$\widehat{f}_p$ is computed per \S10, \  \#11 but $\widehat{f}_\infty$ is computed per
\[
\chi_\infty(x) = \exp(-2\pi\sqrt{-1} \ x),
\]
meaning that the sign convention here is the opposite of that laid down in \S10 (a harmless deviation).
\end{x}

\vspace{0.1cm}

\begin{x}{\small\bf APPLICATION} \ 
\[
f \in \sB_\infty(\A) \implies \widehat{f} \in \sB_\infty(\A) \qquad (\text{cf.} \  \S10, \  \#16).
\]
\end{x}

\vspace{0.1cm}

\begin{x}{\small\bf \un{N.B.}} \ 
It is clear that
\[
\sB_\infty(\A) \subset \bINV(\A)
\]
and $\forall$ $f \in \sB_\infty(\A)$,
\[
\widehat{\widehat{f}\hspace{0.1cm}} = f(-x) \quad \text{$( x \in \A)$}.
\]
\end{x}

\vspace{0.11cm}

\begin{x}{\small\bf LEMMA} \ 
Given $f \in \sB_\infty(\A)$, the series
\[
\sum_{r \in \Q} f(x+r),  \qquad \sum_{q \in \Q} \widehat{f}(x+q)
\]
are absolutely and uniformly convergent on compact subsets of $\A$.
\end{x}

\vspace{0.1cm}

\begin{x}{\small\bf POISSON SUMMATION FORMULA} \ 
Given $f \in \sB_\infty(\A)$,
\[
\sum_{r \in \Q} f(r) \ = \  \sum_{q \in \Q} \widehat{f}(q).
\]

The proof is not difficult but there are some measure theoretic issue to be dealt with first.

On general grounds, 
\[
\int_\A \ =\  \int_{\A/\Q} \ \sum_\Q	\qquad (\text{cf.} \ \S6,  \ \#11).
\]
Here the integral $\ds\int_\A$ is with respect to the Haar measure $\mu_\A$ on $\A$ (cf. \S14, \  \#31).  
Taking $\mu_\Q$ to be counting measure, this choice of data fixes the Haar measure $\mu_{\A/\Q}$ on $\A/\Q$.

\vspace{0.1cm}

[Note: \  The restriction of $\mu_\A$ to the fundamental domain
\[
D = \prod_p \Z_p \times [0,1[
\]
for $\A/\Q$ ( cf. \S14, \  \#10 ) determines $\mu_{\A/\Q}$ and
\[
1 \ =\  \mu_\A(D) \ =\  \mu_{\A/\Q} (\A/\Q).]
\]
If $\phi:\Q \ra \C$, then $\widehat{\phi}:\widehat{\Q} \ra \C$, i.e. $\widehat{\phi}:
\A/\Q \ra \C$ or still,
\[
\widehat{\phi}(\chi) \ = \ \sum_{r \in \Q} \phi(r) \chi(r).
\]

Specialize and suppose that $\phi$ is the characteristic function of $\{0\}$, so $\forall$ $\chi$,
\[
\widehat{\phi}(\chi) \ =\  \chi(0) \ =\  1.
\]
Therefore $\widehat{\phi}$ is the constant function 1 on $\A/\Q$.  
Pass now to $\widehat{\widehat{\phi}\hspace{0.059cm}}$, thus 
$\widehat{\widehat{\phi}\hspace{0.05cm}}: \widehat{\A/\Q} \ra \C$ or still,
\begin{align*}
\widehat{\widehat{\phi}\hspace{0.05cm}}:(\chi_{\Q,q}) \ 
&=\  \int_{\A/\Q} \widehat{\phi}(x) \chi_{\Q,q} (x) d\mu_{\A/\Q}(x) \\
&=\  \int_{\A/\Q} \chi_{\Q,q} (x) d\mu_{\A/\Q}(x)
\end{align*}
which is 1 if $q = 0$ and is 0 otherwise (cf. \S7, \  \#46 $(\A/\Q$ is compact)), 
hence $\widehat{\widehat{\phi}\hspace{0.05cm}} = \phi$.  
But $\phi(r) = \phi(-r)$, thereby leading to the conclusion that the Haar measure 
$\mu_{\A/\Q}$ on $\A/\Q$ is the one singled out by Fourier inversion ( cf. \S7, \  \#45).

\vspace{0.2cm}

Summary: Per Fourier inversion,

\qquad \textbullet \quad $\mu_\Q$ is paired with $\mu_{\A/\Q}$. 

\qquad \textbullet \quad $\mu_{\A/\Q}$ is paired with $\mu_\Q$.
\vspace{0.2cm}

Given $f \in \sB_\infty(\A)$, put
\[
F(x) = \sum_{r \in \Q} f(x+r).
\]
Then $F$ lives on $\A/\Q$, so $\widehat{F}$ lives on $\widehat{\A/\Q} \thickapprox \Q:$
\begin{align*}
\widehat{F}(q) \ 
&=\  \int_{\A/\Q} F(x) \chi_{\Q, q}(x) d\mu_{\A/\Q}(x) \\
&=\  \int_{\A/\Q} F(x) \chi_\Q(qx) d\mu_{\A/\Q}(x).
\end{align*}
On the other hand,
\begin{align*}
\widehat{f}(q)  	
&= \int_\A f(x) \chi_{\Q, q}(x) d\mu_\A(x) \\
&= \int_\A f(x) \chi_\Q (qx) d\mu_\A(x) \\
&= \int_{\A/\Q} \bigl(\sum_{r \in \Q} f(x+r)  \chi_\Q (q(x+r))\bigr) d\mu_{\A/\Q} (x) \\
&= \int_{\A/\Q} \bigl(\sum_{r \in \Q} f(x+r)  \chi_\Q (qx+qr)\bigr) d\mu_{\A/\Q} (x) \\
&= \int_{\A/\Q} \bigl(\sum_{r \in \Q} f(x+r)  \chi_\Q (qx) \chi_\Q (qr)  \bigr) d\mu_{\A/\Q} (x) \\
&= \int_{\A/\Q} \bigl(\sum_{r \in \Q} f(x+r)\bigr)  \chi_\Q (qx)   d\mu_{\A/\Q} (x) \\
&= \int_{\A/\Q} F(x)  \chi_\Q (qx)   d\mu_{\A/\Q} (x) \\
&= \widehat{F}(q).
\end{align*}

To finish the proof, per Fourier inversion, write
\[
F(x) \ =\  \sum_{q \in \Q} \widehat{F}(q) \overline{\chi_\Q(qx)}
\]
and then put $x = 0$:
\[
F(0) \ =\  \sum_{r \in \Q} f(r) \ =\  \sum_{q \in \Q} \widehat{F}(q) \ =\  \sum_{q \in \Q} \widehat{f}(q).
\]
\end{x}

\vspace{0.1cm}

\begin{x}{\small\bf THEOREM} \ 
Let $x \in \I$ $-$then $\forall$ $f \in \sB_\infty(\A)$,
\[
\sum_{r \in \Q} f(rx) \ =\  \frac{1}{\abs{x}_\A} \  \sum_{q \in \Q} \  \widehat{f}(qx^{-1}).
\]

\vspace{0.1cm}

PROOF \  
Work with $f_x \in \sB_\infty(\A)$ $(f_x(y) = f(xy)):$
\[
\sum_{r \in \Q} f_x(r) \ =\  \sum_{q \in \Q} \widehat{f}_x(q).
\]
But
\begin{align*}
\widehat{f}_x(q)  \ 	
&=\  \int_\A f_x(y) \chi_{\Q, q}(y) d\mu_\A(y) \\
&=\  \int_\A f_x(y) \chi_{\Q}(qy) d\mu_\A(y) \\
&=\  \int_\A f(xy) \chi_{\Q}(qxx^{-1}y) d\mu_\A(y) \\
&=\  \frac{1}{\abs{x}_\A} \int_\A f(y) \chi_{\Q}(qx^{-1}y) d\mu_\A(y) \\
&=\  \frac{1}{\abs{x}_\A} \widehat{f}(qx^{-1}).
\end{align*}\\
\end{x}

\chapter{
$\boldsymbol{\S}$\textbf{16}.\quad  FUNCTIONAL EQUATIONS}
\setlength\parindent{2em}
\setcounter{theoremn}{0}

\ \indent 
Let
\[
\zeta(s) = \sum_{n = 1}^\infty \frac{1}{n^s}	\qquad \text{$(\Re(s) > 1)$}
\]
be the Riemann zeta function $-$then $\zeta(s)$ can be meromorphically continued into the whole 
$s$-plane with a simple pole at $s = 1$ and satisfies there the functional equation
\[
\pi^{-s/2}\Gamma(s/2)\zeta(s) \ =\  \pi^{-(1-s)/2}\Gamma((1-s)/2)\zeta(1-s).
\]

\vspace{0.2cm}

\begin{x}{\small\bf REMARK} \ 
The product $\pi^{-s/2}\Gamma(s/2)$ was denoted by $\Gamma_\R(s)$ in \S11, \#8.
\end{x}

\vspace{0.1cm}

There are many proofs of the functional equation satisfied by $\zeta(s)$.  Of these, we shall single out two, one "classical", the other "modern".

\vspace{0.25cm}

To proceed in the classical vein, start with 
\[
\Gamma(s) = \int_0^\infty e^{-x} x^s \frac{dx}{x}			\qquad (\Re(s) > 1).
\]
Then by change of variable, 
\[
\pi^{-s/2}\Gamma(s/2)n^{-s} = \int_0^\infty e^{-n^2\pi x}x^{s/2} \frac{dx}{x}.
\]
So, upon summing from $n = 1$ to $\infty$:
\[
\pi^{-s/2}\Gamma(s/2)\zeta(s) = \int_0^\infty \psi(x) x^{s/2} \frac{dx}{x},
\]
where
\[
\psi(x) = \sum_{n=1}^\infty e^{-n^2\pi x}.
\]
Put now
\[
\theta(x) \ =\  1 + 2\psi(x) \ =\  \sum_{n \in \Z} e^{-n^2\pi x}.
\]

\vspace{0.1cm}

\begin{x}{\small\bf LEMMA} \ 
\[
\theta\bigl(\frac{1}{x}\bigr) \ =\  \sqrt{x}\  \theta(x).
\]
Therefore
\begin{align*}
\psi\bigl(\frac{1}{x}\bigr) 	\ 
&=\  -\frac{1}{2} + \frac{1}{2} \ \theta\bigl(\frac{1}{x}\bigr)\\	
&=\  -\frac{1}{2} + \frac{\sqrt{x}}{2} \ \theta\bigl(x\bigr)\\
&=\  -\frac{1}{2} + \frac{\sqrt{x}}{2} + \sqrt{x}\ \psi(x).		
\end{align*}
One may then write
\begin{align*}
\pi^{-s/2}\Gamma(s/2)\zeta(s) \ 
&=\vsx\  \int_0^\infty \psi(x) x^{s/2} \frac{dx}{x}\\	
&=\vsx\  \int_0^1 \psi(x) x^{s/2} \frac{dx}{x} +  \int_1^\infty \psi(x) x^{s/2} \frac{dx}{x}\\
&=\vsx\  \int_1^\infty  \psi\bigl(\frac{1}{x}\bigr) x^{-s/2} \frac{dx}{x} +  \int_1^\infty \psi(x) x^{s/2} \frac{dx}{x}\\
&=\vsx\  \int_1^\infty  \bigl(-\frac{1}{2} + \frac{\sqrt{x}}{2} + \sqrt{x} \ \psi(x)\bigr) x^{-s/2} \frac{dx}{x} + 
\int_1^\infty \psi(x) x^{s/2} \frac{dx}{x}\\	
&= \frac{1}{s-1} - \frac{1}{s} + \int_1^\infty \psi(x) \bigl(x^{s/2} + x^{(1-s)/2}\bigr) \frac{dx}{x}.			
\end{align*}

The last integral is convergent for all values of $s$ and thus defines a holomorphic function.  
Moreover, the last expression is unchanged if $s$ is replaced by $1-s$.  I.e.:
\[
\pi^{-s/2}\Gamma(s/2)\zeta(s) = \pi^{-(1-s)/2}\Gamma((1-s)/2)\zeta(1-s).
\]

\vspace{0.1cm}

The modern proof of this relation uses the adele-idele machinery.  

Thus let
\[
\Phi(x) = e^{-\pi x_\infty^2} \prod_p \chi_{\Z_p} (x_p) 	\qquad (x \in \A).
\]
Then if $\Re(s) > 1$,
\begin{align*}
\int_\I \Phi(x) \abs{x}_\A^s d^\times x	\ 
&=\vsx\  \int_{\R^\times} e^{-\pi t^2} \abs{t}^s \frac{dt}{\abs{t}} \cdot \prod_p \  \int_{\Q_p^\times} \chi_{\Z_p} (x_p) \abs{x_p}_p^s d^\times x_p\\													
&=\vsx\  \pi^{-s/2}\Gamma(s/2) \cdot \prod_p \ \int_{\Z_p-\{0\}}  \abs{x_p}_p^s d^\times x_p\\	
&=\vsx\  \pi^{-s/2}\Gamma(s/2) \cdot \prod_p \frac{1}{1 - p^{-s}} \qquad (\text{cf.} \  \S6, \ \#26) \\	
&=\vsx\  \pi^{-s/2}\Gamma(s/2) \zeta(s).
\end{align*}

To derive the functional equation, we shall calculate the integral
\[
\int_\I \Phi(x) \abs{x}_\A^s d^\times x
\]
in another way.  
To this end, put
\[
D^\times \ =\  \prod_p \Z_p^\times \times \R_{>0}^\times,
\]
a fundamental domain for $\I/\Q^\times$ (cf. \S14, \  \# 26), so
\[
\I \ =\  \coprod_{r \in \Q^\times} r D^\times \qquad \text{(disjoint union)}.
\]
Therefore
\allowdisplaybreaks
\begin{align*}
\int_\I \Phi(x) \abs{x}_\A^s d^\times x	\ 
&=\vsx\  \sum_{r \in \Q^\times} \int_{rD^\times} \Phi(x) \abs{x}_\A^s d^\times x \\	
&=\vsx\   \int_{D^\times} \sum_{r \in \Q^\times} \Phi(rx) \abs{rx}_\A^s d^\times x \\	
&=\vsx\   \int_{D^\times : \abs{x}_\A \le 1} \sum_{r \in \Q^\times} \Phi(rx) \abs{x}_\A^s d^\times x  +  
\int_{D^\times : \abs{x}_\A \ge 1} \sum_{r \in \Q^\times} \Phi(rx) \abs{x}_\A^s d^\times x.
\end{align*}
To proceed further, recall that $\widehat{\Phi} = \Phi$ $(\implies \widehat{\Phi}(0) = \Phi(0) = 1)$, hence (cf. \S15, \#13) 
\[
1 + \sum_{r \in \Q^\times} \Phi(rx) = \frac{1}{\abs{x}_\A} + \frac{1}{\abs{x}_\A} \sum_{q \in \Q^\times} \Phi(qx^{-1}). 
\]
Accordingly, $\vsx$

\vspace{0.1cm}

$
\ds\int_{D^\times : \abs{x}_\A \le 1} \sum_{r \in \Q^\times} \Phi(rx) \abs{x}_\A^s d^\times x	\\
\vsx\text{\qquad\qquad\qquad} = \int_{D^\times : \abs{x}_\A \le 1} (-1 + \frac{1}{\abs{x}_\A} + 
\frac{1}{\abs{x}_\A} \sum_{q \in \Q^\times} \Phi(qx^{-1}))  \abs{x}_\A^s d^\times x\\	
\vsx\text{\qquad\qquad\qquad} = \int_{D^\times : \abs{x}_\A \le 1}  (\abs{x}_\A^{s-1}  - \abs{x}_\A^s) d^\times x   +  
\int_{D^\times : \abs{x}_\A \ge 1} \sum_{q \in \Q^\times} \Phi(qx) \abs{x}_\A^{1-s} d^\times x.\\	
\\	
$
But 
\begin{align*}
\int_{D^\times : \abs{x}_\A \le 1}  (\abs{x}_\A^{s-1}  - \abs{x}_\A^s) d^\times x \ 
&=\vsx\  \int_0^1 (t^{s - 1} - t) \frac{dt}{t} \\
&=\vsx\  \frac{1}{s-1} - \frac{1}{s}.
\end{align*}
So, upon assembling the data, we conclude that 
\[
\int_\I \Phi(x) \abs{x}_\A^s d^\times x = \frac{1}{s-1} - \frac{1}{s} + \int_{D^\times : \abs{x}_\A \ge 1} \sum_{q \in \Q^\times} \Phi(qx)) (\abs{x}_\A^s  + \abs{x}_\A^{1-s}) d^\times x.
\]
Since the second expression is invariant under the transformation $s \ra 1-s$, the functional equation for $\zeta(s)$ follows once again.
\end{x}

\begin{x}{\small\bf REMARK} \ 
Consider
\[
 \int_{D^\times : \abs{x}_\A \ge 1} \sum_{q \in \Q^\times} \Phi(qx)) \ldots \ .
\]
Then from the definitions, 
\begin{align*}
x \in D^\times
&\implies x_p \in \Z_p^\times \ \& \  qx_p \in \Z_p \\
&\implies q \in \Z.
\end{align*}
Matters thus reduce to $\vsx$
\[
2 \int_1^\infty \ \sum_{n=1}^\infty e^{-n^2 \pi t^2} (t^s + t^{1-s}) \frac{dt}{t}
\]
or still,$\vsx$
\[
\int_1^\infty \psi(t) (t^{s/2} + t^{(1 - s)/2}) \frac{dt}{t},
\]
\end{x}
the classical expression.

\chapter{
$\boldsymbol{\S}$\textbf{17}.\quad  GLOBAL ZETA FUNCTIONS}
\setlength\parindent{2em}
\setcounter{theoremn}{0}

\ \indent 
Structurally, there is a short exact sequence
\[
1 \ra \I^1/ \Q^\times \ra \I/ \Q^\times \ra \R_{>0}^\times \ra 1 	\qquad 
(\text{cf.} \  \S14, \ \#27)
\]
and $\I^1/ \Q^\times$ is compact (cf. \S14, \ \#24).

\vspace{0.3cm}

\begin{x}{\small\bf DEFINITION} \ 
Given $f \in \mathcal{B}_\infty(\A)$ and a unitary character $\omega:\I/ \Q^\times \ra \T$, the 
\underline{global zeta function}
\index{global zeta function} 
attached to the pair $(f,\omega)$ is 
\[
Z (f,\omega,s) \ =\  \int_\I f(x)\omega(x)\abs{x}_\A^s d^\times x		\qquad (\Re(s) > 1).
\]
\end{x}

\vspace{0.1cm}

\begin{x}{\small\bf EXAMPLE} \ 
In the notation of \S16, take
\[
f(x) \ = \ \Phi(x) \ = \ e^{-\pi x_\infty ^2} \ \prod_p \chi_{\Z_p} (x_p)		\qquad (x \in \A)
\]
and let $\omega = 1$ $-$then as shown there
\[
Z(f,1,s) = \pi ^{-s/2} \Gamma(s/2) \zeta(s).
\]
\end{x}

\vspace{0.1cm}

\begin{x}{\small\bf LEMMA} \ 
$Z(f,\omega,s)$ is a holomorphic function of $s$ in the strip $\Re(s) > 1.$
\end{x}

\vspace{0.1cm}

\begin{x}{\small\bf THEOREM} \ 
$Z(f,\omega,s)$ can be meromorphically continued into the whole $s$-plane and satisfies the functional equation
\[
Z(f,\omega,s) = Z(\widehat{f},\ov{\omega},1-s).
\]

[Note:
\[
f \in \mathcal{B}_\infty(\A) \implies \widehat{f} \in \mathcal{B}_\infty(\A)	\qquad (\text{cf.} \  \S15, \  \#9).]
\]
The proof is a computation, albeit a lengthy one.

To begin with,
\[
\I \thickapprox \R_{>0}^\times \times \I^1	\qquad (\text{cf.} \   \S14, \  \#27).
\]
Therefore
\begin{align*}
Z (f,\omega,s) \ 	
&=\vsx \int_\I f(x)\omega(x)\abs{x}_\A^s d^\times x\\	
&=\vsx \int_{\R_{>0}^\times \times \I^1} f(tx)\omega(tx)\abs{tx}_\A^s \frac{dt}{t} d^\times x\\	
&=\vsx \int_0^\infty \bigl(\int_{\I^1} f(tx)\omega(tx)\abs{tx}_\A^s  d^\times x\bigr) \frac{dt}{t}.	
\end{align*}
\end{x}

\vspace{0.1cm}

\begin{x}{\small\bf NOTATION} \ 
Put
\[
Z_t(f,\omega,s) = \int_{\I^1} f(tx)\omega(tx)\abs{tx}_\A^s  d^\times x.
\]
\end{x}

\vspace{0.1cm}

\begin{x}{\small\bf LEMMA} \ 
\ \qquad\\
\vspace{0.1cm}

$
\text{\qquad\qquad} Z_t(f,\omega,s) + f(0) \ds\int_{\I^1/\Q^\times} \omega(tx)\abs{tx}_\A^s  d^\times x\\
\text{\qquad\qquad\qquad\qquad} = Z_{t^{-1}}(\widehat{f},\ov{\omega},1-s) + 
\widehat{f}(0) \int_{\I^1/\Q^\times} \ov{\omega}(t^{-1}x)\abs{t^{-1}x}_\A^{1-s}  d^\times x.\\
$

\vspace{0.1cm}

PROOF \ 
Write
\begin{align*}
\int_{\I^1} f(tx)\omega(tx)\abs{tx}_\A^s  d^\times x \ 
&=\vsx \int_{\I^1/\Q^\times} \bigl(\sum_{r \in \Q^\times} f(rtx) \omega(rtx) \abs{rtx}_\A^s \bigr) d^\times x\\	
&=\vsx \int_{\I^1/\Q^\times} \bigl(\sum_{r \in \Q^\times} f(rtx) \omega(tx) \abs{tx}_\A^s \bigr) d^\times x.	
\end{align*}
Then

\allowdisplaybreaks
$
Z_t(f,\omega,s) + f(0) \ds\int_{\I^1/\Q^\times} \omega(tx)\abs{tx}_\A^s  d^\times x
$
\begin{align*}
\quad \quad     
&=\vsx  \int_{\I^1/\Q^\times} (\sum_{q \in \Q} f(rtx) \omega(tx) \abs{tx}_\A^s d^\times x\\	
&=\vsx \int_{\I^1/\Q^\times} \bigl(\frac{1}{\abs{tx}_\A} \ \sum_{q \in \Q} \widehat{f}(qt^{-1}x^{-1})\bigr) \omega(tx) \abs{tx}_\A^s d^\times x \qquad (\text{cf.} \  \S15, \ \#13)\\
&=\vsx \int_{\I^1/\Q^\times} \bigl(\sum_{q \in \Q} \widehat{f}(qt^{-1}x)\bigr) \abs{t^{-1}x}_\A\omega(tx^{-1}) \abs{tx^{-1}}_\A^s d^\times x \qquad \text{$(x \ra x^{-1})$}\\
&=\vsx \int_{\I^1/\Q^\times} \bigl(\sum_{q \in \Q} \widehat{f}(qt^{-1}x)\bigr) \omega^{-1}(t^{-1}x) \abs{t^{-1}x}_\A^{1-s} d^\times x \\
&=\vsx \int_{\I^1/\Q^\times} \bigl(\sum_{q \in \Q} \widehat{f}(qt^{-1}x)\bigr) \ov{\omega}(t^{-1}x) \abs{t^{-1}x}_\A^{1-s} d^\times x \\
&=\vsx \int_{\I^1/\Q^\times} \bigl(\sum_{q \in \Q^\times} \widehat{f}(qt^{-1}x) \ov{\omega}(qt^{-1}x) \abs{qt^{-1}x}_\A^{1-s} \bigr)d^\times x \\
&\hspace{3.5cm} + \widehat{f}(0) \int_{\I^1/\Q^\times} \ov{\omega}(t^{-1}x)\abs{t^{-1}x}_\A^{1-s}  d^\times x \\
&=\vsx \int_{\I^1} \widehat{f}(t^{-1}x)) \ov{\omega}(t^{-1}x) \abs{t^{-1}x)}_\A^{1-s} d^\times x \\
&\hspace{3.5cm} + \widehat{f}(0) \int_{\I^1/\Q^\times} \ov{\omega}(t^{-1}x)\abs{t^{-1}x}_\A^{1-s}  d^\times x\\
&=\vsx Z_{t^{-1}}(\widehat{f},\ov{\omega},1-s) + \widehat{f}(0) \int_{\I^1/\Q^\times} \ov{\omega}(t^{-1}x)\abs{t^{-1}x}_\A^{1-s}  d^\times x.
\end{align*}

Return to $Z(f,\omega,s)$ and break it up as follows:
\[
Z(f,\omega,s) = \int_0^1 Z_t(f,\omega,s) \frac{dt}{t} + \int_1^\infty Z_t(f,\omega,s) \frac{dt}{t}.
\]

\vspace{0.1cm}

\begin{x}{\small\bf LEMMA} \ 
The integral
\[
 \int_1^\infty Z_t(f,\omega,s) \frac{dt}{t}
\]
is a holomorphic function of $s$.

[It can be expressed as
\[
\int_{\I :\abs{x}_\A \ge 1} f(x)\omega (x) \abs{x}_\A^s d^\times x.]
\]
\end{x}

\vspace{0.1cm}

This leaves
\[
\int_0^1 Z_t(f,\omega,s) \frac{dt}{t},
\]
which can thus be represented as
\[
\int_0^1 (Z_{t^{-1}}(\widehat{f},\ov{\omega},1-s)  \  - \   f(0) \int_{\I^1/\Q^\times} \omega(tx)  \abs{tx}_\A^s d^\times x \  + \  \widehat{f}(0) \int_{\I^1/\Q^\times} \ov{\omega}(t^{-1}x)  \abs{t^{-1}x}_\A^{1-s} d^\times x )\frac{dt}{t}.
\]

To carry out the analysis, subject
\[
\int_0^1 Z_{t^{-1}}(\widehat{f},\ov{\omega},1-s) \frac{dt}{t}
\]
to the change of variable $t \ra t^{-1}$, thereby leading to
\[
\int_1^\infty Z_{t}(\widehat{f},\ov{\omega},1-s) \frac{dt}{t},
\]
a holomorphic function of $s$ (cf. \#7 supra).

\vspace{0.2cm}

It remains to discuss
\begin{align*}
R(f,\omega,s)     	
&=\vsx \int_0^1 (-f(0) \int_{\I^1/\Q^\times} \omega(tx)  \abs{tx}_\A^s d^\times x \  + \  \widehat{f}(0) \int_{\I^1/\Q^\times} \ov{\omega}(t^{-1}x)  \abs{t^{-1}x}_\A^{1-s} d^\times x )\frac{dt}{t}\\
&=\vsx \int_0^1 \bigl(-f(0)\omega(t) \abs{t}^s \int_{\I^1/\Q^\times} \omega(x) d^\times x \  + \  \widehat{f}(0) \ov{\omega}(t^{-1})\abs{t^{-1}}^{1-s} \int_{\I^1/\Q^\times} \ov{\omega}(x) d^\times x \bigr)\frac{dt}{t},
\end{align*}
there being two cases.

\vspace{0.1cm}

1. \quad $\omega$ is nontrivial on $\I^1$.  
Since $\I^1/\Q^\times$ is compact (cf. \S14, \#24), the integrals
\[
\int_{\I^1/\Q^\times} \omega(x) d^\times x, \qquad   \int_{\I^1/\Q^\times} \ov{\omega}(x) d^\times x 
\]
must vanish (cf. \S7, \#46).  
Therefore $R(f,\omega,s) = 0$, hence
\[
Z(f,\omega,s) = \int_1^\infty Z_t(f,\omega,s) \frac{dt}{t} + \int_1^\infty Z_t(\widehat{f},\ov{\omega},1-s) \frac{dt}{t},
\]
is a holomorphic function of $s$.

\vspace{0.1cm}

2. \quad $\omega$ is trivial on $\I^1$.  Let $\phi: \R_{>0}^\times \ra \I/\I^1$ be the isomorphism per \S14, \#27
$-$then $\omega \circ \phi: \R_{>0}^\times \ra \T$ is a unitary character of $\R_{>0}^\times$, 
thus for some $w \in \R$, $\omega \circ \phi = \acdot^{-\sqrt{-1}\ w}$, so
\[
\omega = \acdot^{-\sqrt{-1}\ w} \circ \phi^{-1} \implies \omega(x) = \abs{x}_\A^{-\sqrt{-1}\ w}.
\]
Therefore

\begin{align*}
R(f,\omega,s)  \ 	
&=\  -f(0) \vol(\I^1/\Q^\times) \int_0^1 t^{-\sqrt{-1} \ w + s - 1} dt + \widehat{f}(0) \vol(\I^1/\Q^\times) \int_0^1 t^{-\sqrt{-1}\ w + s - 2} dt \\
&=\  -f(0) \frac{\vol(I^1/\Q^\times)}{-\sqrt{-1}\ w + s} + \widehat{f}(0) \frac{\vol(\I^1/\Q^\times)}{-\sqrt{-1}\ w + s - 1},
\end{align*}
a meromorphic function that has a simple pole at
\\
\[
\begin{cases}
\ s = \sqrt{-1}\ w \quad \quad \quad \text{ with residue } \quad -f(0) \ \vol(\I^1/\Q^\times) \quad \text{if} \ f(0) \ne 0\\
\ s = \sqrt{-1}\ w+1 \quad \ \text{ with residue } \qquad  \widehat{f}(0) \ \vol(\I^1/\Q^\times) \quad \text{if} \ \widehat{f}(0) \ne 0\\
\end{cases}
.\]

\vspace{0.1cm}

\begin{x}{\small\bf \un{N.B.}} \ 
To explicate $\vol(\I^1/\Q^\times)$ use the machinery of \S16:  In the notation of \#2 above, 
\[
Z(f,1,s) = -\frac{1}{s} + \frac{1}{s-1} + \dotsb 
\]
\[
\implies \vol(\I^1/\Q^\times) = 1.
\]

[Note: \ Here, $w = 0$ and $f(0) = 1$, $\widehat{f}(0) = 1.]$
\end{x}

That $Z(f,\omega,s)  $ can be meromorphically continued into the whole $s$-plane is now manifest.  
As for the functional equation, we have
\\
\begin{align*}
Z(f,\omega,s)      
&= \vsx \int_1^\infty Z_t(f,\omega,s) \frac{dt}{t} 
+ \int_1^\infty Z_t(\widehat{f},\ov{\omega},1-s)  \frac{dt}{t} 
+ R(f,\omega,s)  \\
&= \vsx\int_1^\infty \bigl(\int_{\I^1} f(tx) \omega(tx) \abs{tx}_\A^s d^\times x\bigr)  \frac{dt}{t} 
+ \int_1^\infty \bigl(\int_{\I^1} \widehat{f}(tx) \ov{\omega}(tx) \abs{tx}_\A^{1-s} d^\times x\bigr)  \frac{dt}{t} 
+ R(f,\omega,s).
\end{align*}
And we also have
\\
\begin{align*}
Z(\widehat{f},\ov{\omega},1-s)      	
&=\vsx \int_1^\infty Z_t(\widehat{f},\ov{\omega},1-s) \frac{dt}{t} + \int_1^\infty Z_t(\widehat{\widehat{f}\hspace{0.125cm}},\ov{\ov{\omega}},1-(1-s))  \frac{dt}{t} + R(\widehat{f},\ov{\omega},1-s)  \\
&=\vsx \int_1^\infty Z_t(\widehat{f},\ov{\omega},1-s) \frac{dt}{t} 
+ \int_1^\infty Z_t(\widehat{\widehat{f}\hspace{0.125cm}},\omega,s)\bigr)  \frac{dt}{t} 
+ R(\widehat{f},\ov{\omega},1-s)  \\						
&=\vsx \int_1^\infty \bigl(\int_{\I^1} \widehat{f}(tx) \ov{\omega}(tx) \abs{tx}_\A^{1-s} d^\times x\bigr)  \frac{dt}{t} 
+ \int_1^\infty \bigl(\int_{\I^1} \widehat{\widehat{f}\hspace{0.125cm}}(tx) \omega(tx) \abs{tx}_\A^s d^\times x\bigr)  \frac{dt}{t} \\
&\hspace{7.5cm} +  R(\widehat{f},\ov{\omega},1-s).
\end{align*}
The first of these terms can be left as is (since it already figures in the formula for $Z(f,\omega,s)).$  
Recalling that
\[
\widehat{\widehat{f}\hspace{0.125cm}}(x) = f(-x) \quad (x \in \A) \qquad (\text{cf.} \  \S15, \ \#10),
\]
The second term becomes
\[
\int_1^\infty \bigl(\int_{\I^1} f(-tx) \omega(tx) \abs{tx}_\A^s d^\times x\bigr) \frac{dt}{t}
\]
or still,
\[
\int_1^\infty \bigl(\int_{\I^1} f(tx) \omega(-tx) \abs{-tx}_\A^s d^\times x\bigr) \frac{dt}{t} 
= \int_1^\infty \bigl(\int_{\I^1} f(tx) \omega(-tx) \abs{tx}_\A^s d^\times x\bigr) \frac{dt}{t}.
\]
But by hypothesis, $\omega$ is trivial on $\Q^\times$, hence
\[
\omega(-tx) = \omega((-1)tx) = \omega(-1) \omega(tx) = \omega(tx),
\]
and we end up with 
\[
 \int_1^\infty \bigl(\int_{\I^1} f(tx) \omega(tx) \abs{tx}_\A^s d^\times x\bigr) \frac{dt}{t}
\]
which likewise figures in the formula for $Z(f,\omega,s)$.  
Finally, if $\omega$ is trivial on $\I^1$, then
\\
\begin{align*}
R(\widehat{f},\ov{\omega},1-s) \   	
&= -\frac{\widehat{f}(0)}{\sqrt{-1} \ w + 1 - s} +  \frac{\widehat{\widehat{f}\hspace{0.125cm}}(0)}{\sqrt{-1} \ w + (1 - s) - 1}\\
&= \frac{f(0)}{\sqrt{-1}\  w - s} -  \frac{\widehat{f}(0)}{\sqrt{-1} \ w + 1 - s}\\
&= -\frac{f(0)}{-\sqrt{-1}\  w + s} +  \frac{\widehat{f}(0)}{-\sqrt{-1}\  w + s - 1}\\
&= R(f,\omega,s).
\end{align*}
On the other hand, if $\omega$ is nontrivial on $\I^1$, then $\ov{\omega}$ is nontrivial on $\I^1$ and 
\[
R(f,\omega,s) = 0, \quad R(\widehat{f},\ov{\omega},1-s) = 0.
\]
\end{x}

\chapter{
$\boldsymbol{\S}$\textbf{18}.\quad  LOCAL ZETA FUNCTIONS (BIS)}
\setlength\parindent{2em}
\setcounter{theoremn}{0}

\ \indent 
To be in conformity with the global framework laid down in \S17, we shall reformulate the local theory of \S11 and \S12.

\vspace{0.25cm}

\begin{x}{\small\bf DEFINITION} \ 
Given $f \in \mathcal{S}(\R)$ and a unitary character $\omega : \R^\times \ra \T$, the 
\un{local zeta function}
\index{local zeta function} 
attached to the pair $(f,\omega)$ is
\[
Z(f,\omega,s) = \int_{\R^\times} f(x) \omega(x) \abs{x}^s d^\times x 		\qquad (\Re (s) > 0).
\]
\end{x}

\vspace{0.1cm}

\begin{x}{\small\bf THEOREM} \ 
There exists a meromorphic function $\rho(\omega,s)$ such that $\forall \ f$,
\[
\rho(\omega,s) = \frac{Z(f,\omega,s) }{Z(\widehat{f},\ov{\omega},1-s) }.
\]

Decompose $\omega$ as a product:
\[
\omega(x) = (\sgn\  x)^\sigma \abs{x}^{-\sqrt{-1} \ w} \qquad (\sigma \in \{0,1\}, w \in \R).
\]
\end{x}
\vspace{0.1cm}

\begin{x}{\small\bf DEFINITION} \ 
Write (cf. \S11, \#9)
\[
L(\omega,s) =\ 
\begin{cases}
\Gamma_\R (s - \sqrt{-1} \ w) \quad \quad \quad (\sigma = 0)
\vspace{0.25cm}\\
\Gamma_\R (s - \sqrt{-1} \ w + 1 )  \ \quad (\sigma = 1)
\end{cases}
.
\]
\end{x}
\vspace{0.1cm}

\begin{x}{\small\bf FACT} \ 
\[
\rho(\omega,s) = \ 
\begin{cases}
\ds\frac{L(\omega,s) }{L(\omega,1-s)} \quad \quad \quad \quad (\sigma = 0)
\vspace{0.25cm}\\
-\sqrt{-1} \ \ds\frac{L(\omega,s) }{L(\ov{\omega},1-s)}  \  \quad (\sigma = 1)
\end{cases}
.\]
\end{x}
\vspace{0.1cm}

\begin{x}{\small\bf REMARK} \ 
The complex case can be discussed analogously but it will not be needed in the sequel.
\\
\end{x}

\vspace{0.1cm}

\begin{x}{\small\bf DEFINITION} \ 
Given $f \in \sB(\Q_p)$ and a unitary character $\omega:\Q_p^\times \ra \T$, 
the 
\un{local zeta function}
\index{local zeta function} 
attached to the pair $(f,\omega)$ is
\[
Z(f,\omega,s) = \int_{\Q_p^\times} f(x) \omega(x) \abs{x}_p^s d^\times x \qquad (\Re (s) > 0).
\]
\end{x}

\vspace{0.1cm}

\begin{x}{\small\bf THEOREM} \ 
There exists a meromorphic function $\rho(\omega,s)$ such that $\forall$ $f$,
\[
\rho(\omega,s) = \frac{Z(f,\omega,s) }{Z(\widehat{f},\ov{\omega},1-s) }.
\]

Decompose $\omega$ as a product:
\[
\omega(x) = \un{\omega}(x) \abs{x}_p^{-\sqrt{-1} \ w} \qquad (\un{\omega} \in \widehat{\Z_p^\times}, \ w \in \R).
\]
\end{x}

\vspace{0.1cm}

\begin{x}{\small\bf DEFINITION} \ 
Write (cf. \S12, \ \#8)
\[
L(\omega,s) =\ 
\begin{cases}
(1 - \omega(p)p^{-s})^{-1}   \quad (\un{\omega} = 1)
\vspace{0.25cm}\\
\ 1  \  \quad \quad \quad \quad \quad \quad \quad (\un{\omega} \ne 1)
\end{cases}
.\]

[Note: \  if $\un{\omega}$ = 1, then
\[
\omega(p) = \abs{p}_p^{-\sqrt{-1}\ w} = p^{\sqrt{-1}\ w}.]
\]
\end{x}

\vspace{0.1cm}

\vspace{0.15cm}
\begin{x}{\small\bf FACT} \ 
$(\un{\omega} = 1)$
\[
\rho(\omega,s) \ = \  \frac{L(\omega,s) }{L(\ov{\omega},1-s)} \ = \  \frac{1 - \ov{\omega}(p)p^{-(1-s)}}{1 - \omega(p)p^{-s}}.
\]
\end{x}

\vspace{0.1cm}

\begin{x}{\small\bf FACT} \ 
$(\un{\omega} \ne 1$)
\[
\rho(\omega,s) = \tau(\omega) \hspace{0.05cm} \un{\omega}(-1) \hspace{0.05cm} p^{n(s \ + \ \sqrt{-1} \  -1)},
\]
where
\[
\tau(\omega) = \sum_{i = 1}^r \un{\omega}(e_i) \chi_p(p^{-n}e_i)
\]
and $\deg \omega = n \ge 1.$
\end{x}


\[
\textbf{APPENDIX}
\]
\setcounter{theoremn}{0}

It can happen that
\\
\[
Z(f,\omega,s)  \ \equiv \  0.
\]
To illustrate, suppose that $\omega(-1) = -1$ and $f(x) = f(-x)$.  Working with $\Q_p^\times$ $($the story for $\R^\times$ being the same$)$, we have
\\
\begin{align*}
Z(f,\omega,s)	\ 	
&=\vsx\  \int_{\Q_p^\times} f(x) \omega(x) \abs{x}_p^s d^\times x\\
&=\vsx\  \int_{\Q_p^\times} f(-x) \omega(-x) \abs{-x}_p^s d^\times x\\
&=\vsx\ \omega(-1)  \int_{\Q_p^\times} f(x) \omega(x) \abs{x}_p^s d^\times x\\	
&=\vsx\ \omega(-1) Z(f,\omega,s)\\
&=\vsx\ -Z(f,\omega,s).
\end{align*}


\chapter{
$\boldsymbol{\S}$\textbf{19}.\quad  L-FUNCTIONS}
\setlength\parindent{2em}
\setcounter{theoremn}{0}

\ \indent 
Let $\omega:\I/\Q^\times \ra \T$ be a unitary character.

\vspace{0.25cm}

\begin{x}{\small\bf LEMMA} \ 
There is a unique unitary character $\un{\omega}$ of $\I/\Q^\times $ of finite order and a unique real number $w$ such that
\[
\omega \hsx = \hsx \un{\omega}\acdot_\A^{-\sqrt{-1}\ w}.
\]

[Note: \  To say that $\un{\omega}$ is of finite order means that there exists a positive integer $n$ such that 
$\un{\omega}(x)^n = 1$ $\forall$ $x \in \I$.]
\end{x}

\vspace{0.1cm}

\begin{x}{\small\bf \un{N.B.}} \ 
\[
\omega \hsx =\hsx  \prod_p \omega_p \times \omega_\infty,
\]
where
\[
\omega_p \hsx =\hsx  \un{\omega}_p \acdot_p^{-\sqrt{-1}\ w}
\]
and
\[
\omega_\infty \hsx =\hsx  (\sgn)^\sigma \acdot_\infty^{-\sqrt{-1}\ w}.
\]
\end{x}

\vspace{0.1cm}

\begin{x}{\small\bf DEFINITION} \ 
\[
L(\omega,s) = \prod_p L(\omega_p,s) \times L(\omega_\infty,s).
\]
\end{x}

\begin{x}{\small\bf RAPPEL} \ 
\[L(\omega_p,s) =\ 
\begin{cases}
\ (1 - \omega_p(p)p^{-s})^{-1}   \quad (\un{\omega}_p = 1)\\
\ 1  \qquad \qquad\qquad \qquad  (\un{\omega}_p \ne 1)
\end{cases}
\quad (\text{cf.} \ \S18, \ \#8).
\]

\vspace{0.1cm}

[Note: \  The set $S_\omega$ of primes for which $\un{\omega}_p \ne 1$ is finite.$]$
\end{x}

\vspace{0.1cm}

\begin{x}{\small\bf SUBLEMMA} \ 
\[
\abs{x} < 1 \implies \log (1 - x) = -\sum_{k = 1}^\infty \frac{x^{k}}{k}.
\]
Therefore
\begin{align*}
\abs{x} > 1 \implies  \log \frac{1}{1 - x^{-1}}
&=  \log 1 - \log (1 - x^{-1})\\
&=  -\bigl( -\sum_{k = 1}^\infty \frac{x^{-k}}{k}\bigr)\\
&=  \sum_{k = 1}^\infty \frac{x^{-k}}{k}.
\end{align*}
\end{x}

\vspace{0.1cm}

\begin{x}{\small\bf \un{N.B.}} \ 
\[
\log f(z) = \log \abs{f(z)} + \sqrt{-1} \  \arg f(z) 
\]
\qquad \qquad\qquad \qquad$\implies$
\[
\Re \log f(z) = \log \abs{f(z)}.
\]
\end{x}

\vspace{0.1cm}

\begin{x}{\small\bf LEMMA} \ 
The product
\[
\prod_p L(\omega_p,s)
\]
is absolutely convergent provided $\Re (s) > 1$.

\vspace{0.1cm}

PROOF \ 
Ignoring $S_\omega$ $($a finite set$)$, it is a question of estimating
\[
\prod \frac{1}{\abs{ 1 - \omega_p(p) p^{-s}}}.
\]
So take its logarithm and consider
\allowdisplaybreaks
\begin{align*}
\sum \log \bigl(\frac{1}{\abs{ 1 - \omega_p(p) p^{-s}}}\bigr) \ 	
&=\  \sum \Re \log \bigl(\frac{1}{ 1 - \omega_p(p) p^{-s}}\bigr)\\
&\\
&=\  \Re \sum \log \bigl(\frac{1}{ 1 - \omega_p(p) p^{-s}}\bigr)\\
&\\	
&=\  \Re \sum \sum_{k = 1}^\infty \frac{\omega_p(p)^k p^{-ks}}{k}.	
\end{align*}
The claim then is that the series
\[
\sum \sum_{k = 1}^\infty \frac{\omega_p(p)^k p^{-ks}}{k}
\]
is absolutely convergent.  
But
\[
\sum \sum_{k = 1}^\infty \abs{\frac{\omega_p(p)^k p^{-ks}}{k}} \ =\   \sum \sum_{k = 1}^\infty \frac{p^{-k\Re(s)}}{k}
\]
which is bounded by

\allowdisplaybreaks
\begin{align*}
\sum\limits_p \sum_{k = 1}^\infty \frac{p^{-k\Re(s)}}{k} \ 
&=\vsx\  \sum\limits_p \sum_{k = 1}^\infty \frac{p^{-k(1 + \delta)}}{k}	\qquad (\Re(s) = 1 + \delta)\\	
&\le\vsx \sum\limits_p \sum_{k = 1}^\infty p^{-k(1 + \delta)}\\	
&=\vsx\  \sum\limits_p \frac{p^{-(1 + \delta)}}{1 - p^{-(1 + \delta)}}\\	
&=\vsx\  \sum\limits_p \frac{1}{p^{1 + \delta}(1 - p^{-(1 + \delta)})}\\	
&=\vsx\  \sum\limits_p \frac{1}{p^{(1 + \delta)} - 1}\\	
&\le\vsx\  2 \ \sum\limits_p \frac{1}{p^{1 + \delta}}\\	
&<\vsx \  \infty.
\end{align*}
\end{x}

\begin{x}{\small\bf EXAMPLE} \ 
Take $\omega = 1$ $-$then
\begin{align*}
L(\omega,s) \ 
&=\  \prod_p \frac{1}{1 - p^{-s}} \times \Gamma_\R(s)\\ 
&=\  \pi^{-s/2} \Gamma(s/2) \zeta(s).
\end{align*}
\end{x}

\vspace{0.1cm}

\begin{x}{\small\bf LEMMA} \ 
$L(\omega,s)$ is a holomorphic function of $s$ in the strip $\Re(s) > 1.$
\end{x}

\vspace{0.1cm}

\begin{x}{\small\bf LEMMA} \ 
$L(\omega,s)$ admits a meromorphic continuation to the whole $s$-plane (see below).
\end{x}

\vspace{0.1cm}

\allowdisplaybreaks

Owing to \S17, \#4, $\forall$ $f \in \sB_\infty(\A)$, 
\allowdisplaybreaks
\[
Z(f,\omega,s) \hsx = \hsx Z(\hat{f},\ov{\omega},1-s).
\]
To exploit this, assume that 
\[
f \ = \ \prod_p \ f_p \times f_\infty,
\]
where $\forall \ p$, $f_p \in \sB(\Q_p)$ and $f_p = \chi_{\Z_p}$ for all but a finite number of $p$, while 
$f_\infty \in \sS(\R)$ $-$then 
\begin{align*}
Z(f,\omega,s) \ 
&= \ \int_\I f(x) \omega(x) \abs{x}_\A^s d^\times x\\
&= \ \prod_p \ \int\limits_{\Q_p^\times} f_p(x_p) \omega_p(x_p) \abs{x_p}_p^s d^\times x_p 
\times 
\int\limits_{\R^\times} f_\infty(x_\infty)\omega_\infty(x_\infty) \abs{x_\infty}_\infty^s d^\times x_\infty\\
&= \ \prod_p Z(f_p,\omega_p,s) \times Z(f_\infty,\omega_\infty,s)
\end{align*}
and analogously for $Z(\widehat{f},\ov{\omega}, 1 - s)$.

\vspace{0.1cm}

Therefore
\allowdisplaybreaks
\begin{align*}
\allowdisplaybreaks
1 \ 
&=\vsx\ \frac{Z(f,\omega,s)}{Z(\widehat{f},\ov{\omega}, 1 - s)} \\
&=\vsx \ \prod_p 
\frac{Z(f_p,\omega_p,s)}{Z(\widehat{f}_p,\ov{\omega}_p, 1 - s)} \times 
\frac{Z(f_\infty,\omega_\infty,s)}{Z(\widehat{f}_\infty,\ov{\omega}_\infty, 1 - s)}\\
&=\vsx \ \prod_p \rho(\omega_p,s) \times \rho(\omega_\infty,s)\\
&=\vsx \ \prod_{p \notin S_\omega} \rho(\omega_p,s) \times 
\prod_{p \in S_\omega} \rho(\omega_p,s) \times \rho(\omega_\infty,s)\\
&=\vsx \ \prod_{p \notin S_\omega} \frac{L(\omega_p,s)}{L(\ov{\omega}_p,1-s)} \times 
\prod_{p \in S_\omega} \rho(\omega_p,s) \times \frac{L(\omega_\infty,s)}{L(\ov{\omega}_\infty,1-s)}\\
&=\vsx \ 
\prod_{p \in S_\omega} \rho(\omega_p,s) \times 
\prod_{p \notin S_\omega} \frac{L(\omega_p,s)}{L(\ov{\omega}_p,1-s)} \times
\prod_{p \in S_\omega} \frac{L(\omega_p,s)}{L(\ov{\omega}_p,1-s)} \times \frac{L(\omega_\infty,s)}{L(\ov{\omega}_\infty,1-s)}\\
&=\vsx \ \prod_{p \in S_\omega} \rho(\omega_p, s) \times 
\prod_p 
\frac{L(\omega_p,s)}{L(\ov{\omega}_p,1-s)} \times 
\frac{L(\omega_\infty,s)}{L(\ov{\omega}_\infty,1-s)}\\
&=\vsx \ \prod_{p \in S_\omega} \rho(\omega_p, s) \times 
\ds\frac
{\ds\prod\limits_p L(\omega_p,s) \times L(\omega_\infty,s)}
{\ds\prod\limits_p L(\ov{\omega}_p,1-s) \times L(\ov{\omega}_\infty,1-s)}\\
&=\vsx \ \prod_{p \in S_\omega} \rho(\omega_p, s) \times \frac{L(\omega,s)}{L(\ov{\omega},1-s)}\\
&=\vsx \ \prod_{p \in S_\omega} \varepsilon(\omega_p, s) \times \frac{L(\omega,s)}{L(\ov{\omega},1-s)}
\qquad\qquad (\text{cf.} \ \S12, \ \#11)\\
&=\vsx \ \varepsilon(\omega, s) \times \frac{L(\omega,s)}{L(\ov{\omega},1-s)},
\end{align*}
where 
\[
\varepsilon(\omega, s) \ = \ \prod\limits_{p \in S_\omega} \ \varepsilon(\omega_p, s).
\]
\vspace{0.1cm}

\begin{x}{\small\bf THEOREM} \ 
\[
L(\ov{\omega}, 1 - s) \ = \ \varepsilon(\omega, s) \ L(\omega,s).
\]
\end{x}
\vspace{0.1cm}

\begin{x}{\small\bf EXAMPLE} \ 
Take $\omega = 1$ (cf. \# 8) $-$then $\varepsilon(\omega,s) = 1$ and
\[
L(\ov{\omega}, 1 - s) \ = \ L(\omega,s)
\]
\end{x}
\vspace{0.1cm}
translates into 
\[
\pi^{-(1-s)/2} \Gamma((1-s)/2) \zeta(1-s) \ = \ \pi^{-s/2}\Gamma(s/2) \zeta(s) \qquad (\text{cf.} \ \# 16).
\]
\vspace{0.1cm}

Make the following explicit choice for 
\[
f \ = \ \prod_p \ f_p \times f_\infty.
\]

\qquad \textbullet \quad If $\un{\omega}_p = 1$, let
\[
f_p(x_p) \ = \ \chi_p(x_p) \chi_{\Z_p}(x_p).
\]
Then
\[
Z(f_p,\omega_p,s) \ = \ L(\omega_p,s).
\]

\qquad \textbullet \quad If $\un{\omega}_p \neq 1$ and deg $\omega_p = n \geq 1$, let
\[
f_p(x_p) \ = \ \chi_p(x_p) \chi_{p^{-n}\Z_p} (x_p).
\]

Then 
\[
Z(f_p, \omega_p, s) \ = \ \tau(\omega_p) \ \frac{p^{1 + n(s + \sqrt{-1} \ w - 1)}}{p-1} \ L(\omega_p,s).
\]

At infinity, take 
\[
f_\infty(x_\infty)  \ = \ e^{-\pi x_\infty^2} \ (\sigma = 0) \ \ \text{or} \ \  
f_\infty(x_\infty)  \ = \ x_\infty e^{-\pi x_\infty^2} \ (\sigma = 1).
\]
Then
\[
Z(f_\infty,x_\infty,s) \ = \ L(\omega_\infty,s).
\]
\vspace{0.1cm}

\begin{x}{\small\bf NOTATION} \ 
Put
\[
H(\omega,s) \ = \ \prod\limits_{p \in S_\omega} \ \tau(\omega_p) \frac{p^{1 + n(s + \sqrt{-1} \ w - 1)}}{p-1}.
\]
\end{x}

\vspace{0.1cm}

\begin{x}{\small\bf \un{N.B.}} \ 
$H(\omega,s)$ is a never zero entire function of $s$.
\end{x}

\vspace{0.1cm}

\begin{x}{\small\bf LEMMA} \ 
\[
Z(f,\omega,s) \ = \  H(\omega,s) L(\omega,s).
\]
\end{x}

\vspace{0.1cm}

Since $Z(f,\omega,s)$ is a meromorphic function of $s$ (cf. \S17, \#4), it therefore follows that $L(\omega,s)$ is a meromorphic function of $s$.

\vspace{0.1cm}

Working now within the setting of \S17, we distinguish two cases per $\omega$.

\vspace{0.1cm}

1. \quad $\omega$ is nontrivial on $\I^1$, hence $\un{\omega} \neq 1$ and in this situation,  
$Z(f,\omega,s)$ is a holomorphic function of $s$, hence the same is true of $L(\omega,s)$.

\vspace{0.1cm}

2. \quad $\omega$ is trivial on $\I^1$ $-$then $\omega = \acdot_\A^{-\sqrt{-1} \ w}$ 
and there are simple poles at 
\[
\begin{cases}
\ s = \sqrt{-1} \  w \quad \text{with residue} \  -f(0) \ \text{if} \ f(0) \neq 0\\
\ s = \sqrt{-1} \  w + 1 \quad \text{with residue} \  \widehat{f}(0) \ \text{if} \ \widehat{f}(0) \neq 0
\end{cases}
.
\]
But $\forall \ p$, $\omega_p = \acdot_p^{-\sqrt{-1} \ w}$ ($\implies \un{\omega}_p = 1$), so 
$f_p(0) = 1$.  
And likewise $f_\infty(0) = 1$ $(\sigma = 0)$.    
Conclusion: $f(0) = 1$.  
As for the Fourier transforms, $\widehat{f}_p = \chi_{\Z_p}$ $\implies$ $\widehat{f}_p(0) = 1$.  
Also $\widehat{f}_\infty = f_\infty$ $(\sigma = 0)$ $\implies$ $\widehat{f}_\infty(0) = 1$.  
Conclusion: $\widehat{f}(0) = 1$.  
The respective residues are therefore $-1$ and 1.
\vspace{0.2cm}

\begin{x}{\small\bf THEOREM} \ 
Suppose that 
$\omega_{1,p} = \omega_{2,p}$ for all but finitely many $p$ and 
$\omega_{1,\infty} = \omega_{2,\infty}$ $-$then 
$\omega_1 = \omega_2$.  

\vspace{0.1cm}

PROOF \ 
Put $\omega = \omega_1 \omega_2^{-1}$, thus $\omega_p = 1$ for all $p$ outside a finite set $S$ of primes, so 
\begin{align*}
\allowdisplaybreaks
L(\omega,s) \ 
&=\vsx \ \prod_p L(\omega_p,s) \times L(\omega_\infty,s) \\
&=\vsx \ \prod_{p \in S} L(\omega_p,s) \prod_{p \notin S} L(1_p,s) \times L(1_\infty,s) \\
&=\vsx\ L(1,s) \ \prod_{p \in S} \frac{L(\omega_p,s)}{L(1_p,s)}\\
&=\vsx\ L(1,s) \ \prod_{p \in S} \frac{1 - p^{-s}}{1 - \alpha_p p^{-s}},
\end{align*}
where 
$\alpha_p = \omega_p(p)$ if $\un{\omega}_p = 1$ and 
$\alpha_p = 0$ if $\un{\omega}_p \neq 1$, and each factor 
\[
\frac{1 - p^{-s}}{1 - \alpha_pp^{-s}}
\]
is nonzero at $s = 0$ and $s = 1$.  
Therefore $L(\omega,s)$ has a simple pole at $s = 0$ and $s = 1$.  
Consider the decomposition
\[
\omega = \un{\omega} \acdot_\A^{-\sqrt{-1} \ w} \qquad (\text{cf.} \  \S19, \ \#1).
\]
Then $\un{\omega} = 1$ since otherwise $L(\omega,s)$ would be holomorphic, which it isn't.  
But then from the theory, $L(\omega,s)$ has simple poles at
\[
\begin{cases}
\ s = \sqrt{-1} \  w \quad \text{with residue} \  -1\\
\ s = \sqrt{-1} \  w + 1 \quad \text{with residue} \  1
\end{cases}
,
\]
thereby forcing $w = 0$, which implies that $\omega = 1$, i.e., $\omega_1 = \omega_2$.   

\vspace{0.1cm}

[Note: \ In the end, $\omega_p = 1$ $\forall \ p$, hence 
\[
\prod\limits_{p \in S} \ \frac{1 - p^{-s}}{1 - \alpha_p p^{-s}} \ = \ 
\prod\limits_{p \in S} \ \frac{1 - p^{-s}}{1 - p^{-s}}  \ = \ 1,
\]
as it has to be.]
\end{x}

\chapter{
$\boldsymbol{\S}$\textbf{20}.\quad  FINITE CLASS FIELD THEORY}
\setlength\parindent{2em}
\setcounter{theoremn}{0}

\ \indent 
Given a finite field $\F_q$ of characteristic $p$ (thus $q$ is an integral power of $p$), then in $\F_p^{\cl}$,
\[
\F_q \ = \ \{x:x^q = x\}.
\]
\vspace{0.1cm}

\begin{x}{\small\bf LEMMA} \ 
The multiplicative group 
\[
\F_q^\times \ = \ \{x:x^{q-1} = 1\}
\]
is cyclic of order $q-1$.
\end{x}
\vspace{0.1cm}

\begin{x}{\small\bf NOTATION} \ 
\[
\F_{q^n} \ = \ \{x:x^{q^n} = x\} \qquad (n \geq 1).
\]
\end{x}
\vspace{0.1cm}

\begin{x}{\small\bf LEMMA} \ 
$\F_{q^n}$ is a Galois extension of $\F_q$ of degree $n$.
\end{x}
\vspace{0.1cm}

\begin{x}{\small\bf LEMMA} \ 
$\Gal(\F_{q^n}/F_q)$ is a cyclic group of order $n$ generated by the element $\sigma_{q,n}$, where
\[
\sigma_{q,n}(x) \ = \ x^q \qquad (x \in \F_{q^n}).
\]
\end{x}
\vspace{0.1cm}

\begin{x}{\small\bf LEMMA} \ 
The $\F_{q^n}$ are finite abelian extensions of $\F_q$ and they comprise all the finite extensions of $\F_q$, hence the algebraic closure of $\bigcup\limits_n \F_{q^n}$ is $\F_q^\ab$.
\end{x}
\vspace{0.1cm}

\begin{x}{\small\bf THEOREM} \ 
There is a 1-to-1 correspondence between the finite abelian 
extensions of $\F_q$ and the subgroups of $\Z$ of finite index which is given by 
\[
\F_{q^n} \longleftrightarrow n\Z \qquad (n \geq 1).
\]
\end{x}
\vspace{0.1cm}

Schematically:\\

\[
\begin{tikzcd}[sep=tiny]
{\F_q} &{\subset} &{\F_{q^2}} &{\subset} &{\F_{q^4}}\\
{\cap} &&{\cap}\\
{\F_{q^3}} &{\subset} &{\F_{q^6}}\\
{\cap}\\
{\F_{q^9}}
\end{tikzcd}
\begin{tikzcd}[sep=tiny]
\\
\\
{ } \ar{rrrrr} &&&&&{} \ar{lllll}
\\
\\
\end{tikzcd}
\begin{tikzcd}[sep=tiny]
{\Z} &{\supset} &{2\Z} &{\supset} &{4\Z}\\
{\cup} &&{\cup}\\
{3\Z} &{\supset} &{6\Z}\\
{\cup}\\
{9\Z}
\end{tikzcd}
.
\]

The ``class field'' aspect of all this is the existence of a canonical homomorphism
\[
\rec_q:\Z \lra \Gal(\F_q^\ab / \F_q).
\]

\vspace{0.1cm}

\begin{x}{\small\bf NOTATION} \ 
Define
\[
\sigma_q \in \Gal(\F_q^\ab / \F_q)
\]
by 
\[
\sigma_q(x) \ = \ x^q.
\]
\end{x}

\vspace{0.2cm}

\begin{x}{\small\bf \un{N.B.}} \ 
Under the arrow of restriction
\[
\Gal(\F_q^\ab / \F_q) \lra \Gal(\F_{q^n} / \F_q),
\]
$\sigma_q$ is sent to $\sigma_{q,n}$.
\end{x}

\vspace{0.1cm}

\begin{x}{\small\bf DEFINITION} \ 
\[
\rec_q(k) \ = \ \sigma_q^k \qquad (k \in \Z).
\]
\index{$\rec_q$}
\end{x}

\vspace{0.1cm}

\begin{x}{\small\bf LEMMA} \ 
The identification
\[
\Z/n\Z \ \approx \ \Gal(\F_{q^n} / \F_q). 
\]
is the arrow $k \ra \sigma_{q,n}^k$.
\end{x}

\vspace{0.1cm}

On general grounds, 
\[
\Gal(\F_q^\ab / \F_q) \ = \ \lim\limits_{\lla} \Gal(\F_{q^n} / \F_q).
\]

[Note: \  The open subgroups of $\Gal(\F_q^\ab / \F_q)$ are the $\Gal(\F_{q}^\ab / \F_{q^n})$ and 
\[
\Gal(\F_q^\ab / \F_q) /  \Gal(\F_q^\ab / \F_{q^n}) \ \approx \ \Gal(\F_{q^n} / \F_q).]
\]

Therefore 
\[
\Gal(\F_q^\ab / \F_q)  \ \approx  \ \lim\limits_{\lla} \Z / n\Z,
\]
another realization of the RHS being $\prod\limits_p \Z_p$ which if invoked leads to 
\[
\sigma_q \longleftrightarrow (1, 1, 1, \ldots).
\]

\begin{x}{\small\bf \un{N.B.}} \ 
The composition
\[
\begin{tikzcd}[sep=small]
{\Z} \ar{rr}{\rec_q} &&{\Gal(\F_q^\ab / \F_q) \ \approx \ \lim\limits_{\lla} \Z / n\Z} 
\end{tikzcd}
\]
coincides with the canonical map 
\[
k \ra (k  \ \modx n)_n.
\]
\end{x}

\vspace{0.2cm}

\begin{x}{\small\bf REMARK} \ 
Give $\Z$ the discrete topology $-$then 
\[
\rec_q: \Z \lra \Gal(\F_q^\ab / \F_q) 
\]
is continuous and injective but it is not a homeomorphism ($\Gal(\F_q^\ab / \F_q)$ is compact).

\vspace{0.1cm}

[Note: \ 
The image $\rec_q(\Z)$ is the cyclic subgroup $\langle \sigma_q \rangle$ generated by $\sigma_q$.  And:

\vspace{0.2cm}

\qquad \textbullet \quad $\langle \sigma_q \rangle \ \neq \  \Gal(\F_q^\ab / \F_q)$\\

\qquad \textbullet \quad $\ov{\langle \sigma_q \rangle} \ = \  \Gal(\F_q^\ab / \F_q)$.]
\end{x}
\vspace{0.1cm}

\begin{x}{\small\bf SCHOLIUM} \ 
The finite abelian extensions of $\F_q$ correspond 1-to-1 with the open subgroups of $\Gal(\F_q^\ab / \F_q)$.

\vspace{0.1cm}

[Quote the appropriate facts from infinite Galois theory.]
\end{x}

\vspace{0.1cm}

\begin{x}{\small\bf SCHOLIUM} \ 
The open subgroups of $\Gal(\F_q^\ab / \F_q)$ correspond 1-to-1 with the open subgroups of $\Z$ of finite index.

\vspace{0.1cm}

[Given an open subgroup $U \subset \Gal(\F_q^\ab / \F_q)$, send it to $\rec_q^{-1}(U) \subset \Z$ (discrete topology).  
Explicated: 
\[
\rec_q^{-1}(\Gal(\F_q^\ab / \F_{q^n})) \ = \ n\Z.]
\]
\end{x}
\vspace{0.1cm}


\[
\textbf{APPENDIX}
\]
\setcounter{theoremn}{0}

The norm map 
\[
\tN_{\F_{q^n} / \F_q}: \F_{q^n}^\times \lra \F_q^\times
\]
is surjective.

\vspace{0.1cm}

[Let $x \in \F_{q^n}^\times$:
\begin{align*}
\tN_{\F_{q^n} / \F_q}(x) \ 
&= \ \prod\limits_{i = 0}^{n-1} \ (\sigma_{q,n})^{i_x}\\
&= \ \prod\limits_{i = 0}^{n-1} \ x^{q^i}\\
&= \ x^{\sum\limits_{i = 0}^{n-1} q^i}\\
&= \ x^{(q^n - 1) / (q-1)}.
\end{align*}
Specialize now and take for $x$ a generator of $\F_{q^n}^\times$, hence $x$ is of order $q^n - 1$, hence 
$\tN_{\F_{q^n} / \F_q}(x)$ is of order $q - 1$, hence is a generator of $\F_q$.]

\chapter{
$\boldsymbol{\S}$\textbf{21}.\quad  LOCAL CLASS FIELD THEORY}
\setlength\parindent{2em}
\setcounter{theoremn}{0}

\ \indent 
Let $\K$ be a local field $-$then there exists a unique continuous homomorphism 
\[
\rec_\K : \K^\times \lra \Gal(\K^\ab /\K),
\]
the so-called \un{reciprocity map}, that has the properties delineated in the results that follow.

\vspace{0.3cm}

\begin{x}{\small\bf CHART} \ 

\[
\begin{tikzcd}[sep=small]
{}     &{} \ar[dash]{ddd} &{}    &{} \ar[dash]{ddd} &{}\\
{\text{finite field} \quad \K}  &{}                      &{\Z} &{}                      &{\Gal(\K^\ab /\K)}\\
{\text{local field} \quad \K}  &{}                      &{\K^\times} &{}                      &{\Gal(\K^\ab /\K)}\\
{}     &{}  &{}    &{}  &{}
\end{tikzcd}
.
\]

\end{x}

\vspace{0.1cm}

\begin{x}{\small\bf CONVENTION} \ 
An \un{abelian extension} is a Galois extension whose Galois group is abelian.
\end{x}

\vspace{0.1cm}

\begin{x}{\small\bf SCHOLIUM} \ 
The finite abelian extensions $\LL$ of $\K$ correspond 1-to-1 with the open subgroups of $\Gal(\K^\ab/\K)$: 
\[
\LL \ \longleftrightarrow \ \Gal(\K^\ab/\LL).
\]

[Note: \ $\Gal(\LL/\K)$ is a homomorphic image of $\Gal(\K^\ab/\K)$:
\[
\Gal(\LL/\K) \ \approx \ \Gal(\K^\ab/\K) / \Gal(\K^\ab/\LL).]
\]
\end{x}

\vspace{0.1cm}

\begin{x}{\small\bf LEMMA} \ 
Suppose that $\LL$ is a finite extension of $\K$ $-$then
\[
\tN_{\LL/\K}: \LL^\times  \ra \K^\times
\]
is continuous, sends open sets to open sets, and closed sets to closed sets.
\end{x}

\vspace{0.1cm}


\begin{x}{\small\bf LEMMA} \ 
Suppose that $\LL$ is a finite extension of $\K$ $-$then 
\[
[\K^\times : \tN_{\LL/\K} (\LL^\times)] \ \leq \ [\LL:\K].
\]
\end{x}

\vspace{0.1cm}

\begin{x}{\small\bf LEMMA} \ 
Suppose that $\LL$ is a finite extension of $\K$ $-$then 
\[
[\K^\times : \tN_{\LL/\K} (\LL^\times)] \ = \ [\LL:\K].
\]
iff $\LL / \K$ is abelian.
\end{x}

\vspace{0.1cm}

\begin{x}{\small\bf NOTATION} \ 
Given a finite abelian extension $\LL$/$\K$, denote the composition 
\[
\begin{tikzcd}[sep=small]
{\K^\times} \ar{rr}{\rec_\K} &&{\Gal(\K^\ab/\K)} \ar{rr}{\pi_{\LL/\K}} &&{\Gal(\K/\LL)}
\end{tikzcd}
\]
by $(.,\LL/\K)$, the 
\un{norm residue symbol}.
\index{norm residue symbol}
\end{x}

\vspace{0.1cm}

\begin{x}{\small\bf THEOREM} \ 
Suppose that $\LL$ is a finite extension of $\K$ $-$then the kernel of $(.,\LL/\K)$ is 
$\tN_{\LL/\K}(\LL^\times)$, hence 
\[
\K^\times / \tN_{\LL/\K}(\LL^\times) \ \approx \ \Gal(\LL/\K).
\]
\end{x}
\vspace{0.1cm}

\begin{x}{\small\bf EXAMPLE} \ 
Take $\K = \R$, thus $\K^\ab = \C$ and 
\[
\tN_{\C/\R}(\C^\times) \ = \ \R_{> 0}^\times.
\]
Moreover, 
\[
\Gal(\C/\R) \ = \ \{\id_\C,\sigma\},
\]
where $\sigma$ is the complex conjugation.  
Define now
\[
\rec_\R:\R^\times \lra \Gal(\R^\ab / \R)
\]
by stipulating that
\[
\rec_\R(\R_{> 0}^\times) \ = \ \id_\C, \quad \rec_\R(\R_{< 0}^\times) \ = \ \sigma.
\]
\end{x}

\vspace{0.1cm}

\begin{x}{\small\bf EXAMPLE} \ 
Take $\K = \C$ $-$then $\K^\ab = \C = \K$ and matters in this situation are trivial.
\end{x}
\vspace{0.1cm}

\begin{x}{\small\bf THEOREM} \ 
The arrow 
\[
\LL \lra \tN_{\LL/\K}(\LL^\times)
\]
is a bijection between the finite abelian extensions of $\K$ and the open subgroups of finite index of $\K^\times$.
\end{x}
\vspace{0.1cm}

\begin{x}{\small\bf THEOREM} \ 
The arrow $U \ra \rec_\K^{-1}(U)$ is a bijection between open subgroups of $\Gal(\K^\ab/\K)$ and the open subgroups of finite index of $\K^\times$.
\end{x}

\vspace{0.1cm}

From this point forward, it will be assumed that $\K$ is non-archimedean, hence is a finite extension of $\Q_p$ for some 
$p$ (cf. \S5, \#13).

\vspace{0.2cm}

\begin{x}{\small\bf LEMMA} \ 
$\rec_\K$ is injective and its image is a  proper, dense subgroup of $\Gal(\K^\ab/\K)$.
\end{x}

\vspace{0.1cm}

\begin{x}{\small\bf LEMMA} \ 
\[
(\R^\times,\LL/\K) \ = \ \Gal(\LL/\K_\ur),
\]
where $\K_\ur$ is the largest unramified extension of $\K$ contained in $\LL$ (cf. \S5, \#33).
\end{x}
\vspace{0.1cm}

[Note: \ The image
\[
(1 + p^i,\LL/\K) \ = \ G^i \qquad (i \geq 1),
\]
the 
\un{i$^\text{th}$ ramification group in the upper numbering} 
\index{i$^\text{th}$ ramification group in the upper numbering}
(conventionally, one puts 
\[
G^0 \ = \ \Gal(\LL/\K_\ur)
\]
and refers to it as the 
\index{inertia group}
\un{inertia group}).]
\vspace{0.3cm}

Working within $\K^\sep$, the extension $\K^\ur$ generated by the finite unramified extensions of $\K$ is called the 
\un{maximal unramified extension} 
\index{maximal unramified extension} 
of $\K$.  
This is a Galois extension and
\[
\Gal(\K^\ur / \K) \ \approx \ \Gal(\F_q^\ab / \F_q),
\]
where $\F_q = R/P$ (cf. \S5, \#19).

\vspace{0.25cm}

\begin{x}{\small\bf REMARK} \ 
The finite unramified extensions $\LL$ of $\K$ correspond 1-to-1 with the finite extensions of $R/P = \F_q$ and 
\[
\Gal(\LL / \K)  \ \approx \ \Gal(\F_{q^n} / \F_q) \qquad (n = [\LL:\K]).
\]
\end{x}

\vspace{0.1cm}

\begin{x}{\small\bf LEMMA} \ 
$\K^\ur$ is the field obtained by adjoinging to $\K$ all roots of unity having order prime to $p$.
\end{x}
\vspace{0.1cm}

\begin{x}{\small\bf APPLICATION} \ 
$\K^\ur$ is a subfield of $\K^\ab$.

\vspace{0.1cm}

[Cyclotomic extensions are Galois and abelian.]
\end{x}
\vspace{0.1cm}

\begin{x}{\small\bf THEOREM} \ 
There is a commutative diagram
\[
\begin{tikzcd}[sep=large]
{\K^\times} \ar{d}[swap]{v_\K} \ar{rr}{\rec_\K}     &&{\Gal(\K^\ab / \K)} \ar{d} \\
{\Z}  \ar{rr}[swap]{\rec_q}        &&{\Gal(\F_q^\ab / \F_q)}
\end{tikzcd}
,
\]
the vertical arrow on the right being the composition
\begin{align*}
\Gal(\K^\ab/\K) \ 
&\ra \  \Gal(\K^\ab /\K) / \Gal(\K^\ab/\K^\ur) \\
&\approx \ \Gal(\K^\ur/\K)\\
&\approx \ \Gal(\F_q^\ab/\F_q).
\end{align*}


[Note: \ $\forall \ a \in \K^\times$, 
\[
\mods_\K(a) \ = \ q^{-\ord_\K(a)}.]
\]
\end{x}

\vspace{0.1cm}

\begin{x}{\small\bf \un{N.B.}} \ 
The image of 
\[
\restr{\rec_\K(\pi)}{K^\ur} \in \Gal(\K^\ur / \K)
\]
in $\Gal(\F_q^\ab / \F_q)$ is $\sigma_q$ (cf. \S20, \#7). 

\vspace{0.1cm}

[Note: \ If $\LL$ is a finite unramified extension of $\K$ and if  \ $\widetilde{\sigma}_{q,n}$ is the generator of 
$\Gal(\LL/\K)$ which is the lift of the generator $\sigma_{q,n}$ of $\Gal(\F_{q^n} / \F_q)$ ($n = [\LL:\K]$), then 
\[
(\pi,\LL/\K) \ = \ \widetilde{\sigma}_{q,n}.]
\]
\end{x}
\vspace{0.1cm}

\begin{x}{\small\bf FUNCTORALITY} \ 
Suppose that $\LL / \K$ is a finite extension of $\K$ $-$then the diagram
\[
\begin{tikzcd}[sep=large]
{\LL^\times} \ar{d}[swap]{\tN_{\LL/\K}} \ar{rr}{\rec_\LL}     &&{\Gal(\LL^\ab / \LL)} \ar{d}{\res} \\
{\K^\times}  \ar{rr}[swap]{\rec_\K}        &&{\Gal(\K^\ab / \K)}
\end{tikzcd}
\]
commutes.
\end{x}
\vspace{0.1cm}

\begin{x}{\small\bf DEFINITION} \ 
Given a Hausdorff topological group $G$, let $G^*$ be its commutator subgroup, and put 
$G^\ab = G / \ov{G^*}$ $-$then $\ov{G^*}$ is a closed normal subgroup of $G$ and $G^\ab$ is abelian, the 
\un{topological abelianization}
\index{topological abelianization} 
of $G$.
\end{x}
\vspace{0.1cm}

\begin{x}{\small\bf EXAMPLE} \ 
\[
\Gal(\K^\sep / \K)^\ab \ = \ \Gal(\K^\ab / \K).
\]
\end{x}
\vspace{0.1cm}

\begin{x}{\small\bf CONSTRUCTION} \ 
Let $G$ be a Hausdorff topological group and let $H$ be a closed subgroup of finite index $-$then the 
\un{transfer}
\index{transfer homomorphism} 
homomorphism $\Tee: G^\ab \ra H^\ab$ is defined as follows:  Choose a section 
$s:H \backslash G \ra G$ and for $x \in G$, put
\[
\Tee(x\ov{G^*}) \ = \ \prod\limits_{\alpha \in H \backslash G} h_{x,\alpha} (\modx \ov{H^*}), 
\]
where $h_{x,\alpha} \in H$ is defined by 
\[
s(\alpha)x \ = h_{x,\alpha} s(\alpha x).
\]
\end{x}
\vspace{0.1cm}

\begin{x}{\small\bf EXAMPLE} \ 
Suppose that $\LL / \K$ is a finite extension $-$then $\LL^\sep \ \approx \ \K^\sep$ and 
\[
\Gal(\LL^\sep  / \LL) \ \subset \ \Gal(\K^\sep  / \K)
\]
is a closed subgroup of finite index (viz. $[\LL:\K]$), hence there is a transfer homomorphism
\[
\Tee: \Gal(\K^\ab  / \K) \lra \Gal(\LL^\ab  / \LL).
\]
\end{x}
\vspace{0.1cm}

\begin{x}{\small\bf THEOREM} \ 
The diagram
\[
\begin{tikzcd}[sep=large]
{\LL^\times}  \ar{rr}{\rec_\LL}     &&{\Gal(\LL^\ab / \LL)}  \\
{\K^\times}  \ar{u} \ar{rr}[swap]{\rec_\K}        &&{\Gal(\K^\ab / \K)} \ar{u}[swap]{\Tee}
\end{tikzcd}
\]
commutes.
\end{x}
\vspace{0.1cm}


\chapter{
$\boldsymbol{\S}$\textbf{22}.\quad  WEIL GROUPS: THE ARCHIMEDEAN CASE}
\setlength\parindent{2em}
\setcounter{theoremn}{0}


\begin{x}{\small\bf DEFINITION} \ 
Put $W_\C = \C^\times$, call it the 
\un{Weil group}
\index{Weil group} 
of $\C$, and leave it at that.
\end{x}

\vspace{0.1cm}

\begin{x}{\small\bf DEFINITION} \ 
Put 
\[
W_\R \ = \ \C^\times \ \cup \ \tJ \C^\times \quad \text{(disjoint union) \quad (J a formal symbol)},
\]
where $\tJ^2 = -1$ and $\tJ z \tJ^{-1} = \ov{z}$ (obvious topology on $W_\R$).  Accordingly, there is a nonsplit short exact sequence
\[
1 \lra \C^\times \lra W_\R \lra \Gal(\C / \R) \lra 1, 
\]
the image of J in $\Gal(\C / \R)$ being complex conjugation.

\vspace{0.1cm}

[Note: \ $H^2(\Gal(\C / \R),\C^\times)$ is cyclic of order 2, thus up to equivalence of extensions of $\Gal(\C / \R)$ by 
$\C^\times$ per the canonical action of $\Gal(\C / \R)$ on $\C^\times$, there are two possibilities:

\vspace{0.2cm}

1. \quad A split extension
\[
1 \lra \C^\times \lra E \lra \Gal(\C / \R) \lra 1.
\]

2. \quad A nonsplit extension
\[
1 \lra \C^\times \lra E \lra \Gal(\C / \R) \lra 1.
\]
The Weil group $W$ is a representative of the second situation which is why we took $\tJ^2 = -1$ 
(rather than $\tJ^2 = +1$).]
\end{x}
\vspace{0.1cm}

\begin{x}{\small\bf LEMMA} \ 
The commutator subgroup $W_\R^*$ of $W_\R$ consists of all elements of the form 
$J z \tJ^{-1} z^{-1} = \ds\frac{\ov{z}}{z}$, i.e., $W_\R^* = S$, thus is closed.
\end{x}
\vspace{0.1cm}

Let 
\[
\pr: W_\R \lra \R^\times
\]
be the map sending J to $-1$ and $z$ to $\abs{z}^2$.
\vspace{0.2cm}

\begin{x}{\small\bf LEMMA} \ 
$S$ is the kernel of pr and pr is surjective.
\end{x}
\vspace{0.1cm}

\begin{x}{\small\bf LEMMA} \ 
The arrow 
\[
\pr^\ab: W_\R^\ab \lra \R^\times
\]
induced by pr is an isomorphism.
\end{x}
\vspace{0.1cm}

\begin{x}{\small\bf REMARK} \ 
The inverse $\R^\times \ra W_\R^\ab$ of $\pr^\ab$ is characterized by the conditions 
\[
\begin{cases}
 \ -1 \ra \tJ W_\R^*\\
 \ x \ra \sqrt{x} \ W_\R^* \qquad (x > 0)
\end{cases}
.
\]
\end{x}
\vspace{0.1cm}

\begin{x}{\small\bf NOTATION} \ 
Define
\[
\ncdot : W_\R \lra \R_{>0}^\times
\]
by the prescription 
\[
\norm{z} \ = \ z \ov{z} \quad (z \in \C), \quad \norm{\tJ} = 1.
\]
\end{x}
\vspace{0.1cm}

\begin{x}{\small\bf \un{N.B.}} \ 
$\ncdot$ drops to a continuous homomorphism $W_\R^\ab \ra \R_{>0}^\times$.
\end{x}
\vspace{0.1cm}

\begin{x}{\small\bf DEFINITION} \ 
A 
\un{representation}
\index{representation} 
of $W_\R$ is a continuous homomorphism 
$\rho:W_\R \ra \GL(V)$, where V is a finite dimensional complex vector space.
\end{x}
\vspace{0.1cm}

\begin{x}{\small\bf EXAMPLE} \ 
If $s \in \C$, then the assignment $w \ra \norm{w}^s$ is a 1-dimensional
representation of $W_\R$, i.e., is a character.
\end{x}
\vspace{0.1cm}

\begin{x}{\small\bf \un{N.B.}} \ 
If $\chi$ is a character of $\R^\times$, then $\chi \circ \pr$ is a character of $W_\R$ and all such have this form.

[For any $\rho \in \widetilde{W}_\R$, 
\[
\rho(\ov{z}) \ = \rho(\tJ z \tJ^{-1}) \ = \  \rho(\tJ) \rho(z) \rho(\tJ)^{-1} \ = \  \rho(z).
\]
Therefore
\[
1 \ = \ \rho(-1) \qquad (\text{cf.} \ \S7, \ \#12).
\]
But 
\[
\rho(-1) \ = \ \rho(\tJ^2) \ = \ \rho(\tJ)^2,
\]
so $\rho(\tJ) = \pm1$.  This said, the characters of $\R^\times$ are described in \S7, \#11, thus the 1-dimensional 
representations of $W_\R$ are parameterized by a sign and a complex number $s$:

\qquad \textbullet \quad $(+,s): \rho(z) = \abs{z}^s$, $\rho(\tJ) = + 1$\\

\qquad \textbullet \quad $(-,s): \rho(z) = \abs{z}^s$, $\rho(\tJ) = - 1$.]\\

\vspace{0.1cm}

Let $V$ be a finite dimensional complex vector space.

\vspace{0.1cm}

\end{x}
\vspace{0.1cm}

\begin{x}{\small\bf DEFINITION} \ 
A linear transformation $T:V \ra V$ is 
\un{semisimple} 
\index{semisimple} 
if every $T$-invariant subspace has a complementary $T$-invariant subspace.
\end{x}
\vspace{0.1cm}

\begin{x}{\small\bf FACT} \ 
$T$ is semisimple iff $T$ is diagonalizable, i.e., in some basis $T$ is represented by a diagonal matrix.

\vspace{0.05cm}

[Bear in mind that $\C$ is algebraically closed \ldots \ .]
\end{x}

\vspace{0.1cm}

\begin{x}{\small\bf DEFINITION} \ 
A representation $\rho:W_\R \ra \GL(V)$ is 
\un{semisimple} 
\index{semisimple} 
if $\forall \ w \in W_\R$, 
$\rho(w):V \ra V$ is semisimple.
\end{x}
\vspace{0.1cm}

\begin{x}{\small\bf DEFINITION} \ 
A representation $\rho:W_\R \ra \GL(V)$ is 
\un{irreducible}
\index{irreducible} 
if $V \neq 0$, 
and the only $\rho$-invariant subspaces are 0 and $V$.
\end{x}

\vspace{0.1cm}

The irreducible 1-dimensional representations of $W_\R$ are its characters (which, of course, are automatically semisimple).

\vspace{0.1cm}

\begin{x}{\small\bf LEMMA} \ 
If $\rho:W_\R \ra \GL(V)$ is a semisimple irreducible representation of $W_\R$ of dimension $> 1$, then $\dim V = 2$.

\vspace{0.1cm}

PROOF \ 
There is a nonzero vector $v \in V$ and a charcter $\chi:\C^\times \ra \C^\times$ such that $\forall \ z \in \C^\times$, 
\[
\rho(z) v \ = \ \chi(z) v.
\]
Since the span $S$ of $v$, $\rho(\tJ)v$ is a $\rho$-invariant subspace, the assumption of irreducibility implies that $\dim V = 2$.

[To check the $\rho$-invariance of $S$, note that 
\[
\begin{cases}
 \ \rho(z) \rho(\tJ)v \ = \ \rho(z\tJ)v  \ = \ \rho(\tJ\ov{z})v \ = \ \rho(\tJ) \rho(\ov{z}) v \ = \ \rho(\tJ)\chi(\ov{z}) v\\
 \ \rho(\tJ) \rho(\tJ)v \ = \ \rho(\tJ^2) v \ = \ \rho(-1)v \ = \ \chi(-1) v.]
\end{cases}
.
\]

\vspace{0.1cm}

Given an integer $k$ and a complex number $s$, define a character $\chi_{k,s}:\C^\times \ra \C^\times$ by the prescription
\[
\chi_{k,s}(z) \ = \ \Bigl(\frac{z}{\abs{z}}\Bigr)^k \bigl(\abs{z}^2 \bigr)^s
\]
and let $\rho_{k,s} = \ind \chi_{k,s}$ be the representation of $W_\R$ which it induces.
\end{x}

\vspace{0.1cm}

\begin{x}{\small\bf LEMMA} \ 
$\rho_{k,s}$ is 2-dimensional.
\end{x}

\vspace{0.1cm}

\begin{x}{\small\bf LEMMA} \ 
$\rho_{k,s}$ is semisimple.
\end{x}

\vspace{0.1cm}

\begin{x}{\small\bf LEMMA} \ 
$\rho_{k,s}$ is irreducible iff $k \neq 0$.
\end{x}

\vspace{0.1cm}

\begin{x}{\small\bf DEFINITION} \ 
Let
\[
\begin{cases}
 \ \rho_1:W_\R \ra \GL(V_1)\\
 \ \rho_2:W_\R \ra \GL(V_2)\\
\end{cases}
\]
be representations of $W_\R$ $-$then $(\rho_1,V_1)$ is 
\un{equivalent} 
\index{equivalent representations} 
to $(\rho_2,V_2)$ if there exists an 
isomorphism $f:V_1 \ra V_2$ such that $\forall \ w \in W_\R$, 
\[
f \circ \rho_1(w) \ = \ \rho_2(w) \circ f.
\]
\end{x}

\vspace{0.1cm}

\begin{x}{\small\bf LEMMA} \ 
$\rho_{k_1,s_1}$ is equivalent to $\rho_{k_2,s_2}$ iff $k_1 = k_2$, $s_1 = s_2$ or $k_1 = -k_2$, $s_1 = s_2$.
\end{x}

\vspace{0.1cm}

\begin{x}{\small\bf LEMMA} \ 
Every 2-dimensional semisimple irreducible representation of $W_\R$ is equivalent to a unique $\rho_{k,s}$ $(k > 0)$.
\end{x}

\vspace{0.1cm}

\begin{x}{\small\bf \un{N.B.}} \ 
Therefore the equivalence classes of 2-dimensional semisimple irreducible representations of $W_\R$ are parameterized by 
the points of $\N \times \C$.
\end{x}

\vspace{0.1cm}

\begin{x}{\small\bf DEFINITION} \ 
A representation $\rho:W_\R \ra \GL(V)$ is 
\un{completely reducible}
\index{completely reducible} 
if $V$ is the direct sum of a collection of irreducible 
$\rho$-invariant subspaces.
\end{x}

\vspace{0.1cm}

\begin{x}{\small\bf LEMMA} \ 
Let $\rho:W_\R \ra \GL(V)$ be a semisimple representation $-$then $\rho$ is completely reducible.

\vspace{0.1cm}

PROOF \ 
The characters of $\C^\times$ are of the form $z \ra z^\mu \ov{z}^\nu$ with $\mu$, $\nu \in \C$, $\mu-\nu \in \Z$ and $V$ 
is the direct sum of subspaces $V_{\mu,\nu}$, where 
$\restr{\rho(z)}{V_{\mu,\nu}} =  z^\mu \ov{z}^\nu \ \id_{V_{\mu,\nu}}$.  
Claim: 
\[
\rho(\tJ) V_{\mu,\nu} \ = \ V_{\nu,\mu}.
\]
Proof:
$\forall \ v \in V_{\mu,\nu}$,
\begin{align*}
\rho(z) \rho(\tJ) v \ 
&=\vsy\ \rho(\tJ \ov{z} \tJ^{-1}) \rho(\tJ) v\\
&=\vsy\ \rho(\tJ) \rho(\ov{z}) \rho(\tJ^{-1}) \rho(\tJ) v\\
&=\vsy\ \rho(\tJ) \rho(\ov{z}) v\\
&=\vsy\ \rho(\tJ) \ov{z}^\mu z^\nu v\\
&=\vsy\ \rho(\tJ) z^\nu \ov{z}^\mu v\\
&=\vsy\ z^\nu \ov{z}^\mu \rho(\tJ) v.
\end{align*}

Proceeding:

\vspace{0.1cm}

\qquad \textbullet \quad \un{$\mu = \nu$} \ Choose a basis of eigenvectors for $\rho(\tJ)$ on $V_{\mu,\nu}$ $-$then 
the span of each eigenvector is a 1-dimensional $\rho$-invariant subspace.

\vspace{0.1cm}

\qquad \textbullet \quad \un{$\mu \neq \nu$} Choose a basis $v_1, \ldots v_r$ for $V_{\mu,\nu}$ and put 
$v_i^\prime = \rho(\tJ) v_i$ $(1 \leq i \leq r)$ $-$then $\C v_i \oplus \C v_i^\prime$
is a 2-dimensional $\rho$-invariant subspace and the direct sum 
\[
\bigoplus\limits_{i = 1}^r \ (\C v_i \oplus \C v_i^\prime)
\]
equals 
\[
V_{\mu,\nu} \oplus V_{\nu,\mu}.
\]
\end{x}

\vspace{0.1cm}
\begin{x}{\small\bf REMARK} \ 
Suppose that $\rho:W_\R \ra \GL(V)$ is a representation $-$then 
\begin{align*}
\tJ^2 = -1 
&\implies  (-1) \tJ \cdot \tJ = 1\\
&\implies (-1) \tJ  = \tJ^{-1}
\end{align*}
\qquad\qquad\qquad\qquad $\implies$ 
\begin{align*}
\rho(\tJ)^{-1}\ 
&= \ \rho(\tJ^{-1})\\
&= \ \rho((-1)\tJ)\\
&= \ \rho(-1)\rho(\tJ).
\end{align*}
On the other hand, if $\tJ^2 = 1$ (the split extension situation (cf. \#2)), then 
\begin{align*}
\id_V \ 
&= \  \rho(1) \\
&= \ \rho(\tJ^2)\\
&=\ \rho(\tJ)\rho(\tJ).
\end{align*}
\qquad\qquad\qquad\qquad $\implies$ 
\[
\rho(\tJ)^{-1} \ = \ \rho(\tJ).
\]
\end{x}
\vspace{0.1cm}

\chapter{
$\boldsymbol{\S}$\textbf{23}.\quad  WEIL GROUPS: THE NON-ARCHIMEDEAN CASE}
\setlength\parindent{2em}
\setcounter{theoremn}{0}
\ \indent Let $\K$ be a non-archimedean local field. 

\vspace{0.25cm}

\begin{x}{\small\bf NOTATION} \ 
Put 
\[
\begin{cases}
 \ G_\K \ = \ \Gal(\K^\sep / \K)\\
 \ G_\K^\ab \ = \ \Gal(\K^\ab / \K)
\end{cases}
.
\]
\end{x}

\vspace{0.1cm}

\begin{x}{\small\bf \un{N.B.}} \ 
Every character of $G_\K$ factors through $\ov{G_\K^*}$, hence gives rise to a charcter of $G_\K^\ab$.
\end{x}
\vspace{0.1cm}

To study the characters of $G_\K^\ab$, precompose with the reciprocity map $\rec_\K:\K^\times \ra G_\K^\ab$, thus
\[
\chi_\K \ : \ 
\begin{cases}
 \ (G_\K^\ab)^{\widetilde{\ }} \ra (\K^\times)^{\widetilde{\ }}\\
 \ \chi \ra \chi \circ \rec_\K
\end{cases}
.
\]

\vspace{0.1cm}

\begin{x}{\small\bf LEMMA} \ 
$\chi_\K$ is a homomorphism.
\end{x}
\vspace{0.1cm}

\begin{x}{\small\bf LEMMA} \ 
$\chi_\K$ is injective.

\vspace{0.1cm}

PROOF \ 
Suppose that 
\[
\chi_\K(\chi) \ = \ \chi \circ \rec_\K
\]
is trivial $-$then $\restr{\chi}{\Img \rec_\K} = 1$. But $\Img  \rec_\K$ is dense in $G_\K^\ab$ (cf. \S21, \#13), 
so by continuity, $\chi \equiv 1$.
\end{x}

\vspace{0.1cm}


\begin{x}{\small\bf LEMMA} \ 
$\chi_\K$ is not surjective.

\vspace{0.1cm}

PROOF \ 
$G_\K^\ab$ is compact abelian and totally disconnected.  
Therefore 
$(G_\K^\ab)^{\widetilde{\ }} = (G_\K^\ab)^{\widehat{\ }}$ 
and every $\chi$ is unitary and of finite order 
(cf. \S7, \#7 and \S8, \#2), thus the $\chi_\K(\chi)$ are unitary and of finite order.  
But there are characters of $\K^\times$ for which this is not the case.
\end{x}

\vspace{0.1cm}

\begin{x}{\small\bf \un{N.B.}} \ 
The failure of $\chi_\K$ to be surjective will be remedied below (cf. \#19). 
\end{x}

\vspace{0.1cm}

The kernel of the arrow 
\[
\Gal(\K^\sep / \K) \lra \Gal(\K^\ur / \K)
\]
of restriction is $\Gal(\K^\sep / \K^\ur)$ and there is an exact sequence
\[
1 \lra \Gal(\K^\sep / \K^\ur) \lra \Gal(\K^\sep / \K)  \lra \Gal(\K^\ur / \K) \lra 1.
\]
Identify
\[
\Gal(\K^\ur / \K)
\]
with 
\[
\Gal(\F_q^\ab / \F_q) 
\]
and put 
\[
W(\F_q^\ab / \F_q) \ = \ \langle \sigma_q \rangle \qquad (\text{discrete topology}).
\]

\vspace{0.2cm}

\begin{x}{\small\bf DEFINITION} \ 
The \un{Weil group} $W(\K^\sep / \K)$ is the inverse image of $W(\F_q^\ab / \F_q)$ in $\Gal(\K^\sep / \K)$, 
i.e., the elements of $\Gal(\K^\sep / \K)$ which induce an intgral power of $\sigma_q$.
\end{x}

\vspace{0.1cm}

\begin{x}{\small\bf NOTATION} \ 
Abbreviate $W(\K^\sep / \K)$  to $W_\K$, hence $W_\K \subset G_\K$.
\end{x}

\vspace{0.1cm}

Setting 
\[
I_\K \ = \ \Gal(\K^\sep / \K^\ur) \qquad (\text{the \un{inertia group}}),
\]
there is an exact sequence
\[
\begin{tikzcd}[sep=large]
{1} \ar{r} 
&{I_\K} \ar{r}
&{W_\K} \ar{r}
&{W(\F_q^\ab / \F_q)} \ar{d} \ar{r}
&{1}\\
&&&{\Z} \ar{u}[swap]{\approx}
\end{tikzcd}
.
\]

\vspace{0.1cm}

[Note: \  Fix an element $\widetilde{\sigma}_q \in W_\K$ which maps to $\sigma_q$ $-$then structurally, 
$W_\K$ is the disjoint union
\[
\bigcup\limits_{n \in \Z} (\widetilde{\sigma}_q)^n I_\K.]
\]

Topologize $W_\K$ by taking for a neighborhood basis at the identity the 
\[
\Gal(\K^\sep / \LL) \cap I_\K,
\]
where $\LL$ is a finite Galois extension of $\K$.

\vspace{0.2cm}

\begin{x}{\small\bf REMARK} \ 
$I_\K$ has the relative topology per the inclusion $I_\K \ra G_\K$ and any splitting 
$\Z \ra W_\K$ induces an isomorphism $W_\K \approx I_\K \times \Z$ of topological groups, 
where $\Z$ has the discrete topology.
\end{x}

\vspace{0.1cm}

\begin{x}{\small\bf LEMMA} \ 
$W_\K$ is a totally disconnected locally compact group.

[Note: \ $W_\K$ is not compact \ldots \ .]
\end{x}

\vspace{0.1cm}

\begin{x}{\small\bf LEMMA} \ 
The inclusion  $W_\K \ra G_\K$ is continuous and has a dense image.
\end{x}

\vspace{0.1cm}

\begin{x}{\small\bf LEMMA} \ 
$I_\K$ is open in  $W_\K$.
\end{x}

\vspace{0.1cm}

\begin{x}{\small\bf LEMMA} \ 
$I_\K$ is a maximal compact subgroup of $W_\K$.
\end{x}

\vspace{0.1cm}

Suppose that $\LL / \K$ is a finite extension of $\K$ $-$then $G_\LL \subset G_\K$ is the subgroup of $G_\K$ fixing $\LL$, hence 
\[
W_\LL \subset G_\LL \subset G_\K.
\]

\begin{x}{\small\bf LEMMA} \ 
\[
W_\LL \ = \  G_\LL \cap W_\K \subset W_\K
\]
is open and of finite index in $W_\K$, it being normal in $W_\K$ iff $\LL / \K$ is Galois.
\end{x}

\vspace{0.1cm}

\begin{x}{\small\bf THEOREM} \ 
The arrow 
\[
\LL \ra  W_\LL
\]
is a bijection between the finite extensions of $\K$ and the open subgroups of $W_\K$.

\vspace{0.1cm}

[By contrast, the arrow 
\[
\LL \ra \Gal(\K^\sep / \LL)
\]
is a bijection between the finite extensions of $\K$ and the open subgroups of $G_\K$.]
\end{x}

\vspace{0.1cm}

\begin{x}{\small\bf LEMMA} \ 
\[
\ov{W_\K^*} \ = \ \ov{G_\K^*} .
\]
\end{x}

\vspace{0.1cm}

\begin{x}{\small\bf APPLICATION} \ 
The homomorphism $W_\K^\ab \ra G_\K^\ab$ is 1-to-1.
\end{x}
\vspace{0.1cm}

\begin{x}{\small\bf THEOREM} \ 
The image of 
$\rec_\K:\K^\times \ra G_\K^\ab$ is $W_\K^\ab$ and the induced map $\K^\times \ra W_\K^\ab$ is an isomorphism of topological groups (cf. \S21, \#13).
\end{x}
\vspace{0.1cm}

The characters of $W_\K$ ``are'' the characters of $W_\K^\ab$, so we have:

\vspace{0.2cm}

\begin{x}{\small\bf SCHOLIUM} \ 
There is a bijective correspondence between the characters of $W_\K$ and the characters of $\K^\times$ or still, there is a 
bijective correspondence between the 1-dimensional representations of $W_\K$ and the 
1-dimensional representations of $\GL_1(\K)$.
\end{x}

\vspace{0.1cm}

Suppose that $\LL / \K$ is a finite Galois extension of $\K$ $-$then $G_\LL \subset G_\K$ and 
\[
G_\K / G_\LL \ \approx \ \Gal(\LL / \K)
\]
is finite of cardinality 
$[\LL:\K]$.  Since $W_\K$ is dense in $G_\K$, it follows that the image of the arrow
\[
\begin{cases}
 \ W_\K \lra G_\K / G_\LL\\
 \ w \lra w G_\LL
\end{cases}
\]
is all of $G_\K / G_\LL$, its kernel being those $w \in W_\K$ such that $w \in G_\LL$, i.e., its kernel is 
$G_\LL \cap W_\K$ or still, is $W_\LL$.

\vspace{0.2cm}

\begin{x}{\small\bf LEMMA} \ 
\[
W_\K / W_\LL \ \approx \ G_\K / G_\LL \ \approx \ \Gal(\LL / \K).
\]
\end{x}
\vspace{0.1cm}

\begin{x}{\small\bf LEMMA} \ 
$\ov{W_\LL^*}$ is a normal subgroup of $W_\K$.
\end{x}

\vspace{0.1cm}

[Bearing in mind that $W_\LL$ is a normal subgroup of $W_\K$, if $\alpha, \ \beta \in W_\LL^*$ and if $\gamma \in W_\K$, then 
\[
\gamma \alpha \beta \alpha^{-1} \beta^{-1} \gamma^{-1} \ = \ 
(\gamma \alpha \gamma^{-1}) (\gamma \beta \gamma^{-1}) 
(\gamma \alpha^{-1} \gamma^{-1}) (\gamma \beta^{-1} \gamma^{-1}).]
\]

\vspace{0.1cm}

There is an exact sequence 
\[
1 \lra 
W_\LL /\ov{W_\LL^*} \lra 
W_\K /\ov{W_\LL^*} \lra 
(W_\K /\ov{W_\LL^*})  / (W_\LL /\ov{W_\LL^*}) \lra 
1
\]
or still, there is an exact sequence
\[
1 \lra 
W_\LL /\ov{W_\LL^*} \lra 
W_\K /\ov{W_\LL^*} \lra 
W_\K   / W_\LL \lra 
1.
\]

\vspace{0.2cm}

\begin{x}{\small\bf NOTATION} \ 
Put 
\[
W(\LL,\K) \ = \ W_\K /\ov{W_\LL^*}.
\]
\end{x}

\vspace{0.1cm}

\begin{x}{\small\bf SCHOLIUM} \ 
There is an exact sequence
\[
1 \lra 
W_\LL^\ab \lra
W(\LL,\K) \lra 
W_\K /W_\LL \lra
1
\]
and a diagram
\[
\begin{tikzcd}[sep=large]
&{W_\LL^\ab}   \ar{r}
&{W(\LL,\K)} \ar{r}
&{W_\K / W_\LL} \ar{d}{\approx}\\
{1} \ar{r} &{\LL^\times} \ar{u}{\rec_\LL} &&{\Gal(\LL / \K)} \ar{r} &{1}
\end{tikzcd}
.
\]
\end{x}

\vspace{0.1cm}

\begin{x}{\small\bf NOTATION} \ 
Given $w \in W_\K$, let $\norm{w}$ denote the effect on $w$ of passing
from $W_\K$ to $\R_{>0}^\times$ via the arrows 
\[
\begin{tikzcd}[sep=large]
{W_\K} \ar{r} 
&{W_\K^\ab} \ar{rr}{\rec_\K^{-1}}
&&{\K^\times} \ar{rr}{\mods_\K}
&&{\R_{>0}^\times}
\end{tikzcd}
.
\]
\end{x}

\vspace{0.1cm}

\begin{x}{\small\bf LEMMA} \ 
$\norm{\cdot}:W_\K \ra \R_{>0}^\times$ is a continuous homomorphism and its kernel is $I_\K$.

[Under the arrow 
\[
W_\K \ra W_\K^\ab,
\]
$I_\K$ drops to 
\[
\Gal(\K^\ab / \K^\ur) \ \subset \ W_\K^\ab.
\]
Consider now the arrow
\[
\rec_\K:\K^\times \lra W_\K^\ab.
\]
Then $R^\times$ is sent to $\Gal(\K^\ab / \K^\ur)$ and a prime element $\pi \in R$ is sent to an element 
$\widetilde{\sigma}_q$ in $W_\K^\ab$ whose image in $W(\F_q^\ab / \F_q)$ is $\sigma_q$.  And 
\[
W_\K^\ab \ = \ \bigcup\limits_{n \in \Z} (\widetilde{\sigma}_q)^n \Gal(\K^\ab / \K^\ur).]
\]
\end{x}

\vspace{0.1cm}

\begin{x}{\small\bf DEFINITION} \ 
A \un{representation} of $W_\K$ is a continuous homomorphism 
$\rho:W_\K \ra \GL(V)$, where $V$ is a finite dimensional complex vector space.
\end{x}

\vspace{0.1cm}

\begin{x}{\small\bf LEMMA} \ 
A homomorphism 
$\rho:W_\K \ra \GL(V)$ 
is continuous per the usual topology on $\GL(V)$ iff it is continuous per the discrete topology on $\GL(V)$.

\vspace{0.1cm}

[$\GL(V)$ has no small subgroups.]
\end{x}

\vspace{0.1cm}

\begin{x}{\small\bf SCHOLIUM} \ 
The kernel of every representation of $W_\K$ is trivial on an open subgroup $J$ of $I_\K$.  
Conversely, if $\rho:W_\K \ra \GL(V)$ is a homomorphism which is trivial on an open subgroup $J$ of $I_\K$, then 
the inverse image of any subset of $\GL(V)$ is a union of cosets of $J$, hence is open, hence $\rho$ is continuous, so by definition is a representation of $W_\K$.
\end{x}

\vspace{0.1cm}

\begin{x}{\small\bf EXAMPLE} \ 
Suppose that $\LL / \K$ is a finite Galois extension of $\K$ $-$then 
\begin{align*}
W_\LL \cap I_\K \ 
&=\  G_\LL \cap W_\K \cap I_\K \\
&=\ G_\LL \cap I_\K
\end{align*}
is an open subgroup of $I_\K$.  But
\[
W_\K / W_\LL \ \approx \  \Gal(\LL / \K) \qquad (\text{cf.} \ \#20).
\]
Therefore every homomorphism $\Gal(\LL/\K) \ra \GL(V)$ lifts to a homomorphism $W_\K \ra \GL(V)$ which is trivial on an open subgroup of $I_\K$, hence is a representation of $W_\K$.
\end{x}

\vspace{0.1cm}

\begin{x}{\small\bf \un{N.B.}} \ 
Representations of $W_\K$ arising in this manner are said to be of \un{Galois type}.
\end{x}

\vspace{0.1cm}

\begin{x}{\small\bf LEMMA} \ 
A representation of $W_\K$ is of Galois type iff it has finite image.
\end{x}

\vspace{0.1cm}

\begin{x}{\small\bf EXAMPLE} \ 
$\norm{\cdot}$ is a character of $W_\K$ but as a representation, is not of Galois type.
\end{x}

\vspace{0.1cm}

\begin{x}{\small\bf LEMMA} \ 
Let $\rho:W_\K \ra \GL(V)$ be a representation $-$then the image $\rho(I_\K)$ is finite.

\vspace{0.1cm}

PROOF \ 
Suppose that $J$ is an open subgroup of $I_\K$ on which $\rho$ is trivial.  
Since $I_\K$ is compact and $J$ is open, the quotient $I_\K / J$ is finite, thus 
$\rho(I_\K) = \rho(I_\K / J)$ is finite.
\end{x}

\vspace{0.1cm}

\begin{x}{\small\bf DEFINITION} \ 
A representation $\rho:W_\K \ra \GL(V)$ is \un{irreducible} if $V \neq 0$ 
and the only $\rho$-invariant subspaces are 0 and $V$.
\end{x}

\vspace{0.1cm}

\begin{x}{\small\bf THEOREM} \ 
Given an irreducible representation $\rho$ of $W_\K$, there exists an irreducible representation $\widetilde{\rho}$ of $W_\K$ 
and a complex parameter $s$ such that $\rho \ \approx \ \widetilde{\rho} \ \un{\otimes} \ \norm{\cdot}^s$.
\end{x}

\vspace{0.1cm}

\begin{x}{\small\bf LEMMA} \ 
Let $\rho:W_\K \ra \GL(V)$ be a representation $-$then $V$ is the sum of its irreducible $\rho$-invariant subspaces iff 
every $\rho$-invariant subspace has a $\rho$-invariant complement.
\end{x}

\vspace{0.1cm}

\begin{x}{\small\bf DEFINITION} \ 
Let $\rho:W_\K \ra \GL(V)$ be a representation $-$then $\rho$ is \un{semisimple} if it satisfies either condition of the preceding lemma.
\end{x}

\vspace{0.1cm}

\begin{x}{\small\bf \un{N.B.}} \ 
Irreducible representations are semisimple.
\end{x}

\vspace{0.1cm}

\begin{x}{\small\bf THEOREM} \ 
Let $\rho:W_\K \ra \GL(V)$ be a representation $-$then the following conditions are equivalent

\qquad 1. \quad $\rho$ is semisimple.

\qquad 2. \quad $\rho(\widetilde{\sigma}_q)$ is semisimple.

\qquad 3. \quad $\rho(w)$ is semisimple $\forall \ w \in W_\K$.
\end{x}

\vspace{0.1cm}


\chapter{
$\boldsymbol{\S}$\textbf{24}.\quad  THE WEIL-DELIGNE GROUP}
\setlength\parindent{2em}
\setcounter{theoremn}{0}


\begin{x}{\small\bf DEFINITION} \ 
\  The 
\un{Weil-Deligne}
\index{Weil-Deligne group} 
group \ $WD_\K$ is the semidirect product \  $\C \rtimes W_\K$, the multiplication rule being 
\[
(z_1,w_1) \ (z_2,w_2) \ = \ (z_1  + \norm{w_1} z_2, w_1 w_2).
\]

[Note: \  The identity in $WD_\K$ is $(0,e)$ and the inverse of $(z,w)$ is $(-\norm{w}^{-1}z,w^{-1})$: 
\begin{align*}
(z,w) (-\norm{w}^{-1}z, w^{-1})  \ 
&=\  (z + \norm{w} (-\norm{w}^{-1}z),w w^{-1})\\
&=\ (z - z,e) \\
&=\  (0,e).]
\end{align*}
\end{x}

\vspace{0.1cm}

\begin{x}{\small\bf \un{N.B.}} \ 
The topology on $WD_\K$ is the product topology.
\end{x}

\vspace{0.1cm}

\begin{x}{\small\bf DEFINITION} \ 
A 
\un{Deligne representation}
\index{Deligne representation} 
of $W_\K$ is a triple $(\rho,V,N)$, where $\rho:W_\K \ra \GL(V)$ is a representation of $W_\K$ 
and $N:V \ra V$ is a nilpotent endomorphism of $V$ subject to the relation 
\[
\rho(w) N \rho(w)^{-1} \ = \ \norm{w} N \qquad (w \in W_\K).
\]

[Note: \ $N = 0$ is admissible so every representation of $W_\K$ is a Deligne representation.]
\end{x}

\vspace{0.1cm}

\begin{x}{\small\bf EXAMPLE} \ 
Take $V = \C^n$, hence $\GL(V) = \GL_n(\C)$.  
Let $e_0, e_1, \ldots, e_{n-1}$ be the usual basis of $V$.  Define $\rho$ by the rule
\[
\rho(w) e_i \ = \ \norm{w}^i e_i \qquad (w \in W_\K, \ 0 \leq i \leq n-1)
\]
and define $N$ by the rule 
\[
N e_i \ = \ e_{i+1} \quad (0 \leq i \leq n -2), \quad N e_{n-1} = 0.
\]
Then the triple $(\rho, V, N)$ is a Deligne representation of $W_\K$, the 
\un{$n$-dimensional special} \un{representation}, 
\index{$n$-dimensional special representation}
denoted $\spp(n)$.
\end{x}

\vspace{0.1cm}

\begin{x}{\small\bf DEFINITION} \ 
A 
\un{representation} 
\index{representation} 
of $WD_\K$ is a continuous homomorphism 
$\rho^\prime:WD_\K \ra \GL(V)$ whose restriction to $\C$ is complex analytic, where $V$ is a finite dimensional complex vector space.
\end{x}

\vspace{0.1cm}

\begin{x}{\small\bf LEMMA} \ 
Every Deligne representation $(\rho,V,N)$ of $W_\K$ gives rise to a representation 
$\rho^\prime:WD_\K \ra \GL(V)$ of $WD_\K$.

\vspace{0.1cm}

PROOF \ 
Put
\[
\rho^\prime(z,w) \ = \ \exp(z N) \rho(w).
\]
Then
\begin{align*}
\rho^\prime(z_1,w_1) \rho^\prime(z_2,w_2)   \ 
&=\vsy\  \exp(z_1 N) \rho(w_1) \exp(z_2 N) \rho(w_2) \\
&=\vsy\  \exp(z_1 N) \rho(w_1) \exp(z_2 N) \rho(w_1^{-1}) \rho(w_1)  \rho(w_2) \\
&=\vsy\  \exp(z_1 N) \exp(z_2 \norm{w_1} N) \rho(w_1 w_2)\\
&=\vsy\  \exp(z_1 N + z_2 \norm{w_1} N) \rho(w_1 w_2)\\
&=\vsy\  \exp((z_1 + \norm{w_1} z_2) N) \rho(w_1 w_2)\\
&=\vsy\  \rho^\prime(z_1 + \norm{w_1} z_2, w_1 w_2) \\
&=\vsy\  \rho^\prime((z_1, w_1)(z_2, w_2)).
\end{align*}

[Note: \ 
The continuity of $\rho^\prime$ is manifest as is the complex analyticity of its restriction to $\C$.]
\end{x}
\vspace{0.1cm}

One can also go the other way but this is more involved.

\vspace{0.2cm}

\begin{x}{\small\bf RAPPEL} \ 
If $T:V \ra V$ is unipotent, then 
\[
\log T \ = \ \sum\limits_{n \geq 1} \frac{(-1)^{n+1}}{n} \ (T - I)^n
\]
is nilpotent.
\end{x}
\vspace{0.1cm}

\begin{x}{\small\bf SUBLEMMA} \ 
Let $\rho^\prime:WD_\K \ra \GL(V)$ be a representation of $WD_\K$ $-$then $\forall \ z \neq 0$, $\rho^\prime(z,e)$ is unipotent.
\end{x}

\vspace{0.1cm}

\begin{x}{\small\bf SUBLEMMA} \ 
Let $\rho^\prime:WD_\K \ra \GL(V)$ be a representation of $WD_\K$ $-$then $\forall \ z \neq 0$, 
\[
\log \rho^\prime (z,e)
\]
is nilpotent and 
\[
(\log \rho^\prime (z,e)) / z \qquad (z \neq 0)
\]
is independent of $z$.
\end{x}
\vspace{0.1cm}

\begin{x}{\small\bf LEMMA} \ 
Every representation $\rho^\prime:WD_\K \ra \GL(V)$ of $WD_\K$ gives rise to a Deligne representation 
$(\rho,V,N)$ of $W_\K$.

\vspace{0.1cm}

PROOF \ 
Put 
\[
\rho \ = \ \restr{\rho^\prime}{\{0\}} \times W_{\K}, \  N \ = \ \log \rho^\prime (1,e).
\]
Then $\forall \ w \in W_\K$, 

\allowdisplaybreaks
\begin{align*}
\rho(w) N \rho(w)^{-1}   \ 
&=\vsy\  \rho(w) \log \rho^\prime (1,e) \rho(w)^{-1}\\
&=\vsy\  \rho(w) \bigl( \sum\limits_{n \geq 1} \frac{(-1)^{n+1}}{n} \ (\rho^\prime(1,e) - I)^n\bigr) \rho(w)^{-1}\\
&=\vsy\ \sum\limits_{n \geq 1} \frac{(-1)^{n+1}}{n} (\rho(w) \rho^\prime(1,e) \rho(w)^{-1} - I)^n.
\end{align*}
And 
\allowdisplaybreaks
\begin{align*}
\rho(w) \rho^\prime(1,e) \rho(w)^{-1} \ 
&=\vsy\ \rho^\prime(0,w) \rho^\prime(1,e) \rho^\prime(0,w^{-1})\\
&=\vsy\ \rho^\prime ((0,w)(1,e)(0,w^{-1}))\\
&=\vsy\ \rho^\prime ((\norm{w},w)(0,w^{-1}))\\
&=\vsy\ \rho^\prime(\norm{w},e).
\end{align*}
Therefore
\allowdisplaybreaks
\begin{align*}
\rho(w) N \rho(w)^{-1} \ 
&=\vsy\ \sum\limits_{n \geq 1} \frac{(-1)^{n+1}}{n} (\rho^\prime(\norm{w},e)  - I)^n\\
&=\vsy\ \log \rho^\prime(\norm{w},e)\\
&=\vsy\ \norm{w} \log \rho^\prime (\norm{w}, e)) /  \norm{w}\\
&=\vsy\ \norm{w} \log \rho^\prime (1,e)\\
&=\vsy\ \norm{w} N.
\end{align*}

\end{x}

\vspace{0.1cm}
\begin{x}{\small\bf OPERATIONS} \ 

\vspace{0.3cm}

\index{Deligne representation, direct sum}
\qquad \textbullet \quad \un{Direct Sum}: \ Let $(\rho_1,V_1,N_1)$, $(\rho_2,V_2,N_2)$ be Deligne representations 
$-$then their direct sum is the triple 
\[
(\rho_1 \oplus \rho_2,V_1 \oplus V_2 ,N_1 \oplus N_2).
\]

\index{Deligne representation, tensor product}
\qquad \textbullet \quad \un{Tensor Product}: \ Let $(\rho_1,V_1,N_1)$, $(\rho_2,V_2,N_2)$ be Deligne representations 
$-$then their tensor product is the triple
\[
(\rho_1 \  \otimes \  \rho_2,V_1 \  \otimes \ V_2 ,N_1 \ \otimes \  I_2 + I_1 \  \otimes\  \ N_2).
\]

\index{Deligne representations, contragredient}
\qquad \textbullet \quad \un{Contragredient}: \ Let $(\rho,V,N)$ be a Deligne representation $-$then its contragredient is the triple 
\[
(\rho^\vee, V^\vee, -N^\vee).
\]

[Note: \ 
$V^\vee$ is the dual of $V$ and $N^\vee$ is the transpose of $N$ (thus $\forall \ f \in V^\vee$, 
$N^\vee(f) = f \circ N$).]
\end{x}

\vspace{0.1cm}

\begin{x}{\small\bf REMARK} \ 
The definitions of $\oplus$, $\otimes$, $\vee$ when transcribed to the ``prime picture'' are the usual 
representation-theoretic formalities applied to the group $WD_\K$.  
\end{x}

\vspace{0.1cm}

\begin{x}{\small\bf \un{N.B.}} \ 
Let 
\[
\begin{cases}
\ (\rho_1,V_1,N_1)\\
\ (\rho_2,V_2,N_2)
\end{cases}
\]
be Deligne representations of $W_\K$ $-$then a morphism 
\[
(\rho_1,V_1,N_1) \ra (\rho_2,V_2,N_2)
\]
is a linear map $T:V_1 \ra V_2$ such that 
\[
T \rho_1 (w) \ = \  \rho_2(w) T \qquad (w \in W_\K)
\]
and $T N_1 = N_2 T$.

[Note: \ 
If $T$ is a linear isomorphism, then the Deligne representations 
\[
\begin{cases}
\ (\rho_1,V_1,N_1)\\
\ (\rho_2,V_2,N_2)
\end{cases}
\]
are said to be \un{isomorphic}.]
\index{Deligne representation, isomorphic}
\end{x}

\vspace{0.1cm}

\begin{x}{\small\bf DEFINITION} \ 
Suppose that $(\rho,V,N)$ is a Deligne representation of $W_\K$ $-$then a subspace $V_0 \subset V$ is an 
\un{invariant subspace} 
\index{Deligne representation, invariant subspace}
if it is invariant under $\rho$ and $N$.
\end{x}
\vspace{0.1cm}

\begin{x}{\small\bf LEMMA} \ 
The kernel of $N$ is an invariant subspace. 

\vspace{0.1cm}

PROOF \ 
If $N v = 0$, then $\forall \ w \in W_\K$, 
\[
N \rho (w) v \ = \ \norm{w^{-1}} \rho(w) N v = 0.
\]
\end{x}

\vspace{0.1cm}

\begin{x}{\small\bf DEFINITION} \ 
A Deligne representation $(\rho,V,N)$ of $W_\K$ is 
\un{indecomposable} 
\index{Deligne representation, indecomposable}
if $V$ cannot be written as a direct sum of proper 
invariant subspaces.
\end{x}

\vspace{0.1cm}

\begin{x}{\small\bf EXAMPLE} \ 
Consider $\spp(n)$ $-$then it is indecomposable.

\vspace{0.1cm}

[If $\C^n = S \oplus T$ was a nontrivial decomposition into  proper invariant subspaces, then both
$
\begin{cases}
\ S \cap \ker N\\
\ T \cap \ker N
\end{cases}
$
would be nontrivial.]
\end{x}
\vspace{0.1cm}

\begin{x}{\small\bf DEFINITION} \ 
A Deligne representation $(\rho,V,N)$ of $W_\K$ is 
\un{semisimple} 
\index{Deligne representation, semisimple}
if $\rho$ is semisimple (cf. \S23, \#37).
\end{x}

\vspace{0.1cm}

\begin{x}{\small\bf EXAMPLE} \ 
Consider $\spp(n)$ $-$then it is semisimple.
\end{x}

\vspace{0.1cm}

\begin{x}{\small\bf LEMMA} \ 
Let $\pi$ be an irreducible representation of $W_\K$  $-$then $\spp(n) \ \otimes \  \pi$ is semisimple and indecomposable.  

[Note: \ 
Recall that $\pi$ is identified with $(\pi,0)$.]
\end{x}

\vspace{0.1cm}

\begin{x}{\small\bf THEOREM} \ 
Every semisimple indecomposable Deligne representation of $W_\K$  is equivalent to a Deligne representation of the form 
$\spp(n) \ \otimes \  \pi$, where $\pi$ is an irreducible representation of $W_\K$ and $n$ is a positive integer. 
\end{x}

\vspace{0.1cm}

\begin{x}{\small\bf THEOREM} \ 
Let $(\rho,V,N)$ be a semisimple Deligne representation of $W_\K$ $-$then there is a decomposition 
\[
(\rho,V,N) \ = \ \bigoplus\limits_{i= 1}^{s} \spp(n_i) \ \otimes \ \pi_i,
\]
where $\pi_i$ is an irreducible representation of $W_\K$ and $n_i$ is a positive integer.  
Furthermore, if  
\[
(\rho,V,N) \ = \ \bigoplus\limits_{j = 1}^t \spp(n_j^\prime) \  \otimes \  \pi_j^\prime
\]
is another such decomposition, then $s = t$ and after a renumbering of the summands, 
$\pi_i \approx \pi_i^\prime$ and $n_i = n_i^\prime$.
\end{x}

\vspace{0.1cm}

\[
\textbf{APPENDIX}
\]
\setcounter{theoremn}{0}

Instead of working with 
\[
WD_\K \ = \ \C \rtimes W_\K,
\]
some authorities work with 
\[
SL(2,\C) \times W_\K,
\]
the rationale for this being that the semisimple representations of the two groups are the ``same''.

Given $w \in W_\K$, let 
\[
h_w \ = \ 
\begin{pmatrix}
\norm{w}^{1/2} &0\\
0 &\norm{w}^{-1/2}\\
\end{pmatrix}
\]
and identify $z \in \C$ with 
\[
h_w \ = \ 
\begin{pmatrix}
1 &z\\
0 &1\\
\end{pmatrix}
.
\]

Then
\[
h_w 
\begin{pmatrix}
1 &z\\
0 &1\\
\end{pmatrix}
h_w^{-1}
\ = \ 
\begin{pmatrix}
1 &\norm{w}z\\
0 &1\\
\end{pmatrix}
.
\]
But  conjugation by $h_w$ is an automorphism of $\SL(2,\C)$, thus one can form the semisimple direct product 
$\SL(2,\C) \rtimes W_\K$, the multiplication rule being 
\[
(X_1, w_1) (X_2, w_2) \ = \ (X_1 h_{w_1} X_2 h_{w_1}^{-1}, w_1 w_2).
\]

\vspace{0.1cm}


\begin{x}{\small\bf LEMMA} \ 
The arrow 
\[
(X,w) \lra (Xh_w,w)
\]
from 
\[
\SL(2,\C) \rtimes W_\K \quad \text{to} \quad \SL(2,\C) \times W_\K
\]
is an isomorphism of groups.
\end{x}

\vspace{0.1cm}

\begin{x}{\small\bf DEFINITION} \ 
A 
\un{representation} 
\index{representation of $\SL(2,\C)$} 
of $\SL(2,\C) \times W_\K$ is a continuous homomorphism 
$\rho:\SL(2,\C) \times W_\K \ra \GL(V)$ ($V$ a finite dimensional complex vector space) such that the restriction of $\rho$ to 
$\SL(2,\C)$ is complex analytic.
\end{x}

\vspace{0.1cm}

\begin{x}{\small\bf \un{N.B.}} \ 
$\rho$ is semisimple iff its restriction to $W_\K$ is semisimple.  

\vspace{0.1cm}

[The restriction of $\rho$ to $\SL(2,\C)$ is necessarily semisimple.]
\end{x}

\vspace{0.1cm}

The finite dimensional irreducible representations of $\SL(2,\C)$ are parameterized by the positive integers:
\[
n \longleftrightarrow \sym(n), \qquad \dim \sym(n) = n.
\]

\vspace{0.1cm}

\begin{x}{\small\bf THEOREM} \ 
The isomorphism classes of semisimple Deligne representations of $W_\K$ are in a 1-to-1 correspondence with the isomorphism classes of semisimple representations of $\SL(2,\C) \times W_\K$.
\end{x}

\vspace{0.1cm}

To explicate matters, start with a semisimple indecomposable Deligne representation of $W_\K$, say 
$\spp(n) \ \otimes \  \pi$, and assign to it the external tensor product 
\fboxsep=0cm
$\sym(n)  \boxtimesdmc \hsx \pi$, hence in general 
\[
\bigoplus\limits_{i = 1}^s \spp(n_i)  \ \otimes \  \pi_i \lra 
\bigoplus\limits_{i = 1}^s \sym(n_i) \boxtimesdmc \hsx \pi_i.
\] 

\setcounter{page}{1}
\renewcommand{\thepage}{APPENDIX A-\arabic{page}}
\chapter{
APPENDIX A: TOPICS IN TOPOLOGY}
\setlength\parindent{2em}
\setcounter{theoremn}{0}

\[
\text{NEIGHBORHOODS}
\]

\[
\text{COMPACTNESS}
\]

\[
\text{CONNECTEDNESS}
\]

\[
\text{TOPOLOGICAL GROUPS}
\]
\newpage

\ \indent 

\[
\text{NEIGHBORHOODS}
\]

\begin{x}{\small\bf DEFINITION} \ 
If \mX is a topological space and if $x \in X$, then a 
\un{neighborhood}
\index{neighborhood} 
of $x$ is a set \mU which contains an open set \mV containing $x$, the collection $\sU_x$ of all neighborhoods of $x$ 
being the 
\un{neighborhood system}
\index{neighborhood system} 
at $x$.
\end{x}

\vspace{0.1cm}

Therefore \mU is a neighborhood of $x$ iff $x \in \itr U$.

\vspace{0.2cm}

\begin{x}{\small\bf PROPERTIES} \ 
of $\sU_x$

\vspace{0.1cm}

\un{N-a} \quad If $U \in \sU_x$, then $x \in U$.

\un{N-b} \quad If $U_1, U_2 \in \sU_x$, then $U_1 \cap U_2 \in \sU_x$.

\un{N-c} \quad If $U \in \sU_x$, then there is a $U_0 \in \sU_x$ such that $U \in \sU_{x_0}$ for each $x_0 \in U_0$.

\un{N-d} \quad If $U \in \sU_x$ and $U \subset V$, then $V \in \sU_x$.

\end{x}

\vspace{0.1cm}

\begin{x}{\small\bf FACT} \ 
A subset $G \subset X$ is open iff \mG contains a neighborhood of each of its points.
\end{x}

\vspace{0.1cm}

\begin{x}{\small\bf SCHOLIUM} \ 
If in a set \mX a nonempty collection $\sU_x$ of subsets of \mX is assigned to each $x \in X$ 
so as to satisfy N-a through N-d and if a subset $G \subset X$ is deemed ``open'' provided 
$\forall \ x \in G$, there is a $U \in \sU_x$ such that $U \subset G$, then the result is a topology on \mX 
in which the neighborhood system at each $x \in X$ is $\sU_x$.
\end{x}

\vspace{0.1cm}

\begin{x}{\small\bf DEFINITION} \ 
If \mX is a topological space and if $x \in X$, then a 
\un{neighborhood basis}
\index{neighborhood basis} 
at $x$ is a subcollection $\sB_x$ of $\sU_x$ such that $U \in \sU_x$ contains some $V \in \sB_x$.
\end{x}

\vspace{0.1cm}

\begin{x}{\small\bf EXAMPLE} \ 
Take $X = \R^2$ with the usual topology $-$then the set of all squares with sides parallel to the axes and centered at $x$ 
is a neighborhood basis at $x$.
\end{x}

\vspace{0.1cm}

\begin{x}{\small\bf PROPERTIES} \ 
of $\sB_x$

\vspace{0.1cm}

\un{NB-a} \quad If $V \in \sB_x$, then $x \in V$.

\un{NB-b} \quad If $V_1, V_2 \in \sB_x$, then there is a $V_3 \in \sB_x$ such that $V_3 \subset V_1 \cap V_2$.

\un{NB-c} \quad If $V \in \sB_x$, then there is a $V_0 \in \sB_x$ such that if $x_0 \in V_0$, then there is a 
$W \in \sB_{x_0}$ such that $W \subset V$.
\end{x}

\vspace{0.1cm}

\begin{x}{\small\bf FACT} \ 
A subset $G \subset X$ is open iff \mG contains a basic neighborhood of each of its points.
\end{x}

\vspace{0.1cm}

\begin{x}{\small\bf SCHOLIUM} \ 
If in a set \mX a nonempty collection $\sB_x$ of subsets of \mX is assigned to each $x \in X$ so as to satisfy 
NB-a through NB-c and if a subset $G \subset X$ is deemed ``open'' provided $\forall \ x \in G$, there is a 
$V \in \sB_x$ such that $V \subset G$, then the result is a topology on \mX in which a neighborhood  basis 
at each $x \in X$ is $\sB_x$.

\vspace{0.1cm}

[Put
\[
\sU_x \hsx = \hsx \{U \subset X: V \subset U \ (\exists V \in \sB_x)\}.
\]
Then $\sU_x$ satisfies N-a through N-d above.]

\end{x}

\vspace{0.1cm}

\begin{x}{\small\bf EXAMPLE} \ 
Take $X = \R$ and given $x$, let $\sB_x$ be the $[x,y[$ $(y > x)$ $-$then $\sB_x$ satisfies NB-a through NB-c above, 
from which a topology on the line, the underlying topological space being the 
\un{Sorgenfrey line}.
\index{Sorgenfrey line}.
\end{x}

\vspace{0.1cm}
\begin{x}{\small\bf DEFINITION} \ 
Let \mX be a topological space $-$then a 
\un{basis}
\index{basis} 
for \mX (i.e., for the underlying topology \ldots) is a collection $\sB$ of open sets such that for any open set 
$G \subset X$ and for any point $x \in G$, there is a set $B \in \sB$ such that $x \in B \subset G$.
\end{x}

\vspace{0.1cm}

\begin{x}{\small\bf FACT} \ 
If $\sB$ is a collection of open sets, then $\sB$ is a basis for \mX iff $\forall \ x \in X$, the collection
\[
\sB_x \hsx = \hsx \{B \in \sB : x \in B \}
\]
is a neighborhood basis at $x$.
\end{x}

\vspace{0.1cm}

\begin{x}{\small\bf FACT} \ 
If \mX is a set and if $\sB$ is a collection of subsets of \mX, then $\sB$ is a basis for a topology on \mX iff 
\[
X \hsx = \hsx \bigcup\limits_{B \in \sB} B
\]
and given $B_1, B_2 \in \sB$ and $x \in B_1 \cap B_2$, there exists $B_3 \in \sB$ such that $x \in B_3 \subset B_1 \cap B_2$. 
\end{x}

\vspace{0.1cm}


\newpage
\setcounter{theoremn}{0}


\ \indent 

\[
\text{COMPACTNESS}
\]

\begin{x}{\small\bf DEFINITION} \ 
A topological space \mX is 
\un{compact}
\index{compact} 
if every open cover of \mX has a finite subcover.
\end{x}

\vspace{0.1cm}

\begin{x}{\small\bf EXAMPLE} \ 
The Cantor set is compact.
\end{x}

\vspace{0.1cm}

\begin{x}{\small\bf FACT} \ 
The continous image of a compact space is compact.
\end{x}

\vspace{0.1cm}

\begin{x}{\small\bf FACT} \ 
A one-to-one continuous function from a compact space \mX onto a Hausdorff space \mY is a homeomorphism.
\end{x}

\vspace{0.1cm}

\begin{x}{\small\bf DEFINITION} \ 
A topological space \mX is 
\un{locally compact} 
\index{locally compact} 
if each point in \mX has a neighborhood basis consisting of compact sets.
\end{x}

\vspace{0.1cm}

\begin{x}{\small\bf FACT} \ 
A Hausdorff space \mX is locally compact iff each point in \mX has a compact neighborhood.
\end{x}

\vspace{0.1cm}

\begin{x}{\small\bf APPLICATION} \ 
Every compact Hausdorff space \mX is locally compact.
\end{x}

\vspace{0.1cm}

\begin{x}{\small\bf EXAMPLE} \ 
The Cantor set is a locally compact Hausdorff space.
\end{x}

\vspace{0.1cm}

\begin{x}{\small\bf EXAMPLE} \ 
$\R$ is a locally compact Hausdorff space.
\end{x}

\vspace{0.1cm} 

\begin{x}{\small\bf EXAMPLE} \ 
$\Q$ is a Hausdorff space but it is not locally compact ($\Q$ is first category while a locally compact Hausdorff space is second category).
\end{x}

\vspace{0.1cm}

\begin{x}{\small\bf EXAMPLE} \ 
The Sorgenfrey line is Hausdorff but not locally compact.
\end{x}

\vspace{0.1cm}

\begin{x}{\small\bf FACT} \ 
Suppose that $X_i$ $(i \in I)$ is a nonempty topological space $-$then the product 
$\prod\limits_{i \in I} X_i$ is locally compact iff each $X_i$ is locally compact and all but a finite number of the $X_i$ are compact.
\end{x}


\newpage
\setcounter{theoremn}{0}


\ \indent 

\[
\text{CONNECTEDNESS}
\]

\begin{x}{\small\bf DEFINITION} \ 
A topological space \mX is 
\un{connected} 
\index{connected} 
if it is not the union of two nonempty disjoint open sets.
\end{x}

\vspace{0.1cm}

\begin{x}{\small\bf EXAMPLE} \ 
$\Q$ is not connected (write
\[
\Q \hsx = \hsx \{x: x > \sqrt{2}\} \hsx \cap \hsx \Q \cup \{x: x < \sqrt{2}\} \hsx \cap \hsx \Q). 
\]
\end{x}

\vspace{0.1cm}

\begin{x}{\small\bf EXAMPLE} \ 
$\R$ is connected and the only connected subsets of $\R$ having more than one point are the intervals 
(open, closed, or half-open, half-closed).
\end{x}

\vspace{0.1cm}

\begin{x}{\small\bf FACT} \ 
A topological space \mX is connected iff the only subsets of \mX that are both open and closed are $\emptyset$ and \mX.
\end{x}

\vspace{0.1cm}

\begin{x}{\small\bf FACT} \ 
The continuous image of a connected space is connected.
\end{x}

\vspace{0.1cm}

\begin{x}{\small\bf DEFINITION} \ 
Let \mX be a topological space and let $x \in X$ $-$then the 
\un{component}
\index{component} 
$C(x)$ of $x$ is the union of all connected subsets of $X$ containing $x$.
\end{x}

\vspace{0.1cm}

\begin{x}{\small\bf FACT} \ 
$C(x)$ is a closed subset of \mX.
\end{x}

\vspace{0.1cm}

\begin{x}{\small\bf FACT} \ 
$C(x)$ is a maximal connected subset of \mX.
\end{x}

\vspace{0.1cm}

If $x \neq y$ in \mX, then either $C(x) = C(y)$ or $C(x) \cap C(y) = \emptyset$ (otherwise, 
$C(x) \cup C(y)$ would be a connected set containing $x$ and $y$ and larger than $C(x)$ or $C(y)$, which is 
impossible).  Therefore the set of distinct components of \mX forms a partition of \mX.

\vspace{0.2cm}

\begin{x}{\small\bf EXAMPLE} \ 
Take $X = \Q$ $-$then $\forall \ x \in \Q$, $C(x) = \{x\}$ (under the inclusion $\Q \ra \R$, a connected subset of $\Q$ 
is sent to a connected subset of $\R$).
\end{x}

\begin{x}{\small\bf DEFINITION} \ 
A topological space \mX is 
\un{totally disconnected}
\index{totally disconnected} 
if the components of \mX are singletons, i.e., $\forall \ x \in X$, $C(x) = \{x\}$.

\vspace{0.1cm}

\begin{x}{\small\bf FACT} \ 
A topological space \mX is totally disconnected iff the only nonempty connected subsets of \mX are the one-point 
sets (hence \mX is $T_1$).

\vspace{0.1cm}

[Note: \ In every topological space \mX, the empty set and the one-point sets are connected and in a totally disconnected 
topological space, these are the only connected subsets.]
\end{x}
\end{x}

\vspace{0.1cm}

\begin{x}{\small\bf REMARK} \ 
Let \mE be the equivalence relation defined by writing $x \sim y$ if $x$ and $y$ lie in the same component.  
Equip the set $X/E$ with the identification topology determined by the projection 
$p:X \ra X/E$ $-$then $X/E$ is totally disconnected.
\end{x}

\vspace{0.1cm}

\begin{x}{\small\bf EXAMPLE} \ 
The Cantor set is totally disconnected.
\end{x}

\vspace{0.1cm}

\begin{x}{\small\bf EXAMPLE} \ 
$\Q$ is totally disconnected. 
\end{x}

\vspace{0.1cm}

\begin{x}{\small\bf EXAMPLE} \ 
The Sorgenfrey line is totally disconnected.
\end{x}

\vspace{0.1cm}

\begin{x}{\small\bf FACT} \ 
Every product of totally disconnected topological spaces is totally disconnected.
\end{x}

\vspace{0.1cm}

\begin{x}{\small\bf FACT} \ 
Every subspace of a totally disconnected topological space is totally disconnected.
\end{x}

\vspace{0.1cm}

\begin{x}{\small\bf REMARK} \ 
The continuous image of a totally disconnected space need not be totally disconnected.  
To appreciate the point, recall that evey compact metric space is the continuous image of the Cantor set.
\end{x}


\begin{x}{\small\bf DEFINITION} \ 
A topological space \mX is 
\un{0-dimensional}
\index{0-dimensional} 
if each point of \mX has a neighborhood basis consisting of open-closed sets.
\end{x}

\vspace{0.1cm}

\begin{x}{\small\bf FACT} \ 
A 0-dimensional $T_1$-space is totally disconnected.
\end{x}

\vspace{0.1cm}

\begin{x}{\small\bf EXAMPLE} \ 
The Cantor set is 0-dimensional.
\end{x}

\vspace{0.1cm}

\begin{x}{\small\bf EXAMPLE} \ 
$\Q$ is 0-dimensional.
\end{x}

\vspace{0.1cm}

\begin{x}{\small\bf EXAMPLE} \ 
The Sorgenfrey line is 0-dimensional.
\end{x}

\vspace{0.1cm}

\begin{x}{\small\bf REMARK} \ 
As can be shown by example, a totally disconnected metric space need not be 0-dimensional.
\end{x}

\vspace{0.1cm}

\begin{x}{\small\bf FACT} \ 
A locally compact Hausdorff space is 0-dimensional iff it is totally disconnected.

\vspace{0.1cm}

[Note: \ In such a space, each point has a neighborhood basis consisting of open-compact sets.]
\end{x}

\vspace{0.1cm}

A discrete space is 0-dimensional, hence is totally disconnected, hence a product of discrete spaces is totally disconnected, but an infinite product of nontrivial discrete spaces is never discrete.

\vspace{0.2cm}

\begin{x}{\small\bf DEFINITION} \ 
The 
\un{Cantor space}
\index{Cantor space} 
is the countable product of the two-point discrete space.
\end{x}

\vspace{0.1cm}
\begin{x}{\small\bf FACT} \ 
The Cantor set is homeomorphic to the Cantor space.
\end{x}

\vspace{0.1cm}


\newpage
\setcounter{theoremn}{0}


\ \indent 

\[
\text{TOPOLOGICAL GROUPS}
\]
\begin{x}{\small\bf DEFINITION} \ 
A 
\un{locally compact (compact)} group 
\index{locally compact (compact) group} 
is a topological group \mG that is both locally compact (compact) and Hausdorff.
\end{x}

\vspace{0.1cm}

\begin{x}{\small\bf FACT} \ 
If \mG is a locally compact group and if \mH is a closed subgroup, then $G/H$ is a locally compact Hausdorff space.
\end{x}

\vspace{0.1cm}

\begin{x}{\small\bf FACT} \ 
If \mG is a locally compact group and if \mH is a closed normal subgroup, then $G/H$ is a locally compact group.
\end{x}

\vspace{0.1cm}

\begin{x}{\small\bf FACT} \ 
If \mG is a locally compact group and if \mH is a locally compact subgroup, then \mH is closed in \mG.
\end{x}

\vspace{0.1cm}

\begin{x}{\small\bf FACT} \ 
If \mG is a locally compact 0-dimensional group and if \mH is a closed subgroup of \mG, then $G/H$ is 0-dimensional.
\end{x}

\vspace{0.1cm}

\begin{x}{\small\bf FACT} \ 
If \mG is a totally disconnected locally compact group, then $\{e\}$ has a neighborhood basis consisting of open-compact subgroups.
\end{x}

\vspace{0.1cm}

\begin{x}{\small\bf FACT} \ 
If \mG is a totally disconnected compact group,  then $\{e\}$ has a neighborhood basis consisting of open-compact normal subgroups. 
\end{x}

\vspace{0.1cm}

\begin{x}{\small\bf FACT} \ 
If \mG is a locally compact group, then a subgroup \mH is open iff the quotient $G/H$ is discrete.
\end{x}

\vspace{0.1cm}

\begin{x}{\small\bf FACT} \ 
If \mG is a compact group, then a subgroup \mH is open iff the quotient $G/H$ is finite.
\end{x}

\vspace{0.1cm}

\begin{x}{\small\bf FACT} \ 
If \mG is a locally compact group, then every open subgroup of \mG is closed and every finite index closed subgroup of \mG is open.
\end{x}

\renewcommand{\thepage}{APPENDIX B-\arabic{page}}
\chapter{
APPENDIX B: TOPICS IN ALGEBRA}
\setlength\parindent{2em}
\setcounter{theoremn}{0}
\[
\text{PRINCIPAL IDEAL DOMAINS}
\]

\[
\text{FIELD EXTENSIONS}
\]

\[
\text{ALGEBRAIC CLOSURE}
\]

\[
\text{TRACES AND NORMS}
\]
\newpage

\ \indent 

\[
\text{PRINCIPAL IDEAL DOMAINS}
\]

Let \mA be a commutative ring with unit.

\vspace{0.2cm}

\begin{x}{\small\bf DEFINITION} \ 
An 
\un{ideal}
\index{ideal} 
\mI in \mA is an additive subgroup of \mA such that the relations $a \in A$, $x \in I$ imply that 
$a x$ ($= x a$) belongs to \mI.
\end{x}

\vspace{0.1cm}

\begin{x}{\small\bf DEFINITION} \ 
An ideal \mI in \mA is a 
\un{prime ideal} 
\index{prime ideal} 
if $I \neq A$ and if $a b \in I$ implies that either $a \in I$ or $b \in I$.
\end{x}

\vspace{0.1cm}

\begin{x}{\small\bf DEFINITION} \ 
An ideal \mI in \mA is a 
\un{maximal ideal} 
\index{maximal ideal} 
if $I \neq A$ and there is no larger proper ideal of \mA that contains \mI.
\end{x}

\vspace{0.1cm}

\begin{x}{\small\bf DEFINITION} \ 
\mA is an
\un{integral domain} 
\index{integral domain} 
if $a b = 0$ implies that $a = 0$ or $b = 0$.
\end{x}

\vspace{0.1cm}

\begin{x}{\small\bf \un{N.B.}}  \ 
Every field is an integral domain.
\end{x}

\vspace{0.1cm}

\begin{x}{\small\bf EXAMPLE} \ 
$\Z$ is an integral domain but $Z/n\Z$ is an integral domain iff $n$ is prime.
\end{x}

\vspace{0.1cm}

\begin{x}{\small\bf FACT}  \ 
An ideal $I \neq A$ in \mA is a prime ideal iff $A/I$ is an integral domain.
\end{x}

\vspace{0.1cm}

\begin{x}{\small\bf FACT}  \ 
An ideal $I \neq A$ in \mA is a maximal ideal iff $A/I$ is a field.
\end{x}

\vspace{0.1cm}

\begin{x}{\small\bf EXAMPLE} \ 
Take $A = \Z[X]$ $-$then $\langle X \rangle$ is a prime ideal 
(since $A/\langle X \rangle \approx \Z$ is an integral domain) but $\langle X \rangle$ is not a maximal ideal 
(since $A/\langle X \rangle \approx \Z$ is not a field).
\end{x}

\begin{x}{\small\bf DEFINITION} \ 
An ideal \mI in \mA is a 
\un{principal ideal}
\index{principal ideal} 
if $I = A a_0$ ($\equiv \hsx \langle a_0\rangle$) for some $a_0 \in A$.
\end{x}

\vspace{0.1cm}


\begin{x}{\small\bf DEFINITION} \ 
\mA is a 
\un{principal ideal domain}
\index{principal ideal domain} 
if \mA is an integral domain and if every ideal in \mA is principal.
\end{x}

\vspace{0.1cm}

\begin{x}{\small\bf FACT} \ 
For any field $\K$, the polynomial ring $\K[X]$ is a principal ideal domain.

\vspace{0.1cm}

[If \mI is a nonzero ideal in $\K[X]$, then \mI consists of all the multiples of the monic polynomial in \mI of least degree.]
\end{x}

\vspace{0.1cm}

\begin{x}{\small\bf EXAMPLE} \ 
The polynomial ring $\Z[X]$ is not a principal ideal domain.

\vspace{0.1cm}

[The ideal \mI consisiting of all polynomials with even constant term is not a principal ideal (but it is a maximal ideal).]
\end{x}

\vspace{0.1cm}

\begin{x}{\small\bf FACT} \ 
If \mA is a principal ideal domain, then every nonzero prime ideal is maximal.
\end{x}

\vspace{0.1cm}

\begin{x}{\small\bf FACT} \ 
For any field $\K$, the maximal ideals in $\K[X]$ are the nonzero prime ideals.
\end{x}

\vspace{0.1cm}

\begin{x}{\small\bf DEFINITION} \ 
A 
\un{unit} 
\index{unit} 
in \mA is an element $u \in A$ with a multiplicative inverse, i.e., there is a $v \in A$ such that $u v = 1$.
\end{x}

\vspace{0.1cm}

\begin{x}{\small\bf EXAMPLE} \ 
The units in $\K[X]$ are the nonzero constants.
\end{x}

\vspace{0.1cm}

\begin{x}{\small\bf EXAMPLE} \ 
The units in $\Z$ are 1 and $-1$.
\end{x}

\vspace{0.1cm}

\begin{x}{\small\bf EXAMPLE} \ 
The units in $\Z/n\Z$ are the congruence classes $[a]$ of a mod $n$ such that $(a,n) = 1)$.
\end{x}

\vspace{0.1cm}

\begin{x}{\small\bf DEFINITION} \ 
The elements $a, b \in A$ are said to be 
\un{associates}
\index{associates} 
if there is a unit $u \in A$ such that $a = u b$.
\end{x}

\vspace{0.1cm}


\begin{x}{\small\bf DEFINITION} \ 
A nonzero element $p \in A$ is said to be 
\un{irreducible}
\index{irreducible} 
if $p$ is not a unit and in every factorization $p = a b$, either $a$ or $b$ is a unit.
\end{x}

\vspace{0.1cm}

\begin{x}{\small\bf EXAMPLE} \ 
Take $A = \Z[X]$ $-$then $2X + 2 = 2(X + 1)$ is not irreducible, yet it does not factor into a product of polynomials of lower degree.
\end{x}

\vspace{0.1cm}

\begin{x}{\small\bf SCHOLIUM} \ 
For any field $\K$, a nonzero polynomial $p(X) \in \K[X]$ of degree $\geq 1$ is irreducible iff there is no factorizartion 
$p(X) = f(X) g(X)$ in $\K[X]$ with $\deg f < \deg p$ and $\deg g < \deg p$.
\end{x}

\vspace{0.1cm}

\begin{x}{\small\bf FACT} \ 
If \mA is a principal ideal domain, then the nonzero prime ideals are the ideals $\langle p \rangle$, where $p$ is irreducible.
\end{x}

\vspace{0.1cm}

\begin{x}{\small\bf FACT} \ 
If \mA is a principal ideal domain and if $p \in A$ is irreducible, then $A / \langle p \rangle$ is a field.

\vspace{0.1cm}

[For $\langle p \rangle$ is prime, hence maximal.]
\end{x}

\vspace{0.1cm}

\begin{x}{\small\bf DEFINITION} \ 
\mA is a 
\un{unique factorization domain}
\index{unique factorization domain} 
if \mA is an integral domain subject to:

\vspace{0.1cm}

\un{E}\quad Every nonzero $a \in A$ that is not a unit is a product of irreducible elements.

\un{U}\quad If 
\[
p_1 \cdots p_m \hsx = \hsx q_1 \cdots q_n, 
\]
where the $p$ and $q$ are irreducible, then $m = n$ and there is a one-to-one correspondence between the 
factors such that the corresponding factors are associates.

\end{x}

\vspace{0.1cm}

\begin{x}{\small\bf FACT} \ 
Every principal ideal domain is a unique factorization domain.
\end{x}

\vspace{0.1cm}

\begin{x}{\small\bf APPLICATION} \ 
For any field $\K$, the polynomial ring $\K[X]$ is a unique factorization domain
\end{x}

\vspace{0.1cm}


\begin{x}{\small\bf DEFINITION} \ 
Suppose that \mA is a unique factorization domain $-$then a 
\un{system of representatives of irreducible elements in \mA}
\index{system of representatives of irreducible elements in \mA} 
is a set of irreducible elements having exactly one element in common with the set of all associates of each irreducible element.
\end{x}

\vspace{0.1cm}

\begin{x}{\small\bf SCHOLIUM} \ 
For any field $\K$, the monic irreducible polynomials constitute a system of representatives of irreducible elements in $\K[X]$. 

\vspace{0.1cm}

[Note: \ Let $f$ be a nonconstant polynomial in $\K[X]$ and let $f_1, \ldots, f_n$ be the distinct monic irreducible factors of 
$f$ in $\K[X]$ $-$then 
\[
f \hsx = \hsx C \prod\limits_{k = 1}^n \hsx f_k^{e_k},
\]
where $C$ is the leading coefficient of $f$ and $e_1, \ldots, e_n$ are positive integers.  Moreover, 
this representation of $f$ is unique up to a permutation of $\{1, \ldots, n\}$.]
\end{x}

\vspace{0.1cm}

\begin{x}{\small\bf FACT} \ 
For any field $\K$ and for any irreducible polynomial $p(X)$, the quotient 
$\LL^\prime = \K[X]/\langle p(X) \rangle$ is a field containing an isomorphic copy $\K^\prime$ of $\K$ as a subfield 
and a zero of $p^\prime(X)$.

\vspace{0.1cm}

[Setting $I = \langle p(X) \rangle$, the map $a \ra a + I$ $(a \in \K)$ identifies $\K$ with a subfield $\K^\prime$ 
of $\LL^\prime$.  Write
\[
p(X) \hsx = \hsx a_0 + a_1 X + \cdots + a_n X^n.
\]
Then in $\K^\prime[X]$, 
\[
p^\prime(X) \hsx = \hsx (a_0 + I) + (a_1 + I)X + \cdots + (a_n + I) X^n.
\]
Now put $\theta = X + I$:
\begin{align*}
p^\prime(\theta) \ 
&=\vsx\ (a_0 + I) + (a_1 X + I) + \cdots + (a_n X^n + I)\\
&=\vsx\ a_0 + a_1 X + \cdots + a_n X^n + I\\
&=\vsx\ p(X) + I \\
&=\vsx\ I, 
\end{align*}
the zero element of $\LL^\prime$.]
\end{x}

\vspace{0.1cm}


\setcounter{theoremn}{0}

\newpage

\[
\text{FIELD EXTENSIONS}
\]

\ \indent 

Let $\K$ be a field.

\vspace{0.3cm}

\begin{x}{\small\bf DEFINITION} \ 
A 
\un{field extension} 
\index{field extension} 
of $\K$ is a field $\LL$ having $\K$ as a subfield.  
\end{x}

\vspace{0.1cm}

Given $\LL/\K$ and elements $x_1, \ldots, x_n \in \LL$, write $\K(x_1, \ldots, x_n)$ for the subfield of $\LL$ 
generated by $\K$ and the $x_i$ $(i = 1, \ldots, n)$.  In particular: $\K(x)$ is the subfield generated by $\K$ and $x$.

\begin{x}{\small\bf EXAMPLE} \ 
Take $\K = \Q$, $\LL = \R$, $x = \sqrt{2}$ $-$then $\Q(\sqrt{2})$ consists of all real numbers of the form 
$r + s \sqrt{2}$ $(r, s \in \Q)$.

\vspace{0.1cm}

[Let \mF be the set of all real numbers of the indicated form, thus
\[
\Q \cup \{\sqrt{2}\} \hsx \subset \hsx \F \hsx \subset \hsx \Q(\sqrt{2}),
\]
and, by definition, $\Q(\sqrt{2})$ is the subfield of $\R$  generated by $\Q \cup \{\sqrt{2}\}$.  
Let now 
$x = r + s \sqrt{2}$ $(r, s, \in \Q)$: $r^2 - 2 s^2 \neq 0$ ($\sqrt{2}$ irrational)

\qquad\qquad\qquad\qquad $\implies$
\begin{align*}
\frac{1}{x} \hsx 
&= \hsx \frac{r}{r^2 - 2 s^2} + \frac{-s}{r^2 - 2 s^2} \sqrt{2}\\
&\in \F,
\end{align*}
so $\F$ is a field, so $\F = \Q(\sqrt{2})$.]
\end{x}

\vspace{0.1cm}

\begin{x}{\small\bf EXAMPLE} \ 
Take $\K = \Q$, $\LL = \R$, $x = \sqrt{2}$, $y = \sqrt{3}$ $-$then
\[
\Q(\sqrt{2},\sqrt{3}) \hsx = \hsx \Q(\sqrt{2} + \sqrt{3}).
\]

[Obviously, $\sqrt{2} + \sqrt{3} \in \Q(\sqrt{2},\sqrt{3})$ hence 
$\Q(\sqrt{2} + \sqrt{3}) \subset \Q(\sqrt{2},\sqrt{3})$.  
In the other direction
\[
(\sqrt{2} + \sqrt{3}) (\sqrt{2} - \sqrt{3}) = -1
\]
\qquad\qquad\qquad\qquad $\implies$
\[
\sqrt{3} - \sqrt{2} \hsx = \hsx \ \frac{1}{\sqrt{2} + \sqrt{3}} \in \Q(\sqrt{2} + \sqrt{3})
\]
\qquad\qquad\qquad\qquad $\implies$
\[
\begin{cases}
\ \sqrt{3} \hsx = \hsx ((\sqrt{3} + \sqrt{2}) + (\sqrt{3} - \sqrt{2})) / 2\\
\ \sqrt{2} \hsx = \hsx ((\sqrt{3} + \sqrt{2}) - (\sqrt{3} - \sqrt{2})) / 2
\end{cases}
\in \Q(\sqrt{2} + \sqrt{3}).
\]
Therefore $\Q(\sqrt{2},\sqrt{3}) \subset \Q(\sqrt{2} + \sqrt{3})$.]
\end{x}

\vspace{0.1cm}

Given $\LL \supset \K$, view $\LL$ as a vector space over $\K$ and write $[\LL:\K]$ for its dimension, the 
\un{degree} 
\index{degree} 
of $\LL$ over $\K$.

\vspace{0.1cm}

[Note: \ In this context, the term ``dimension'' refers to the cardinal number of a basis for $\LL$ over $\K$.]

\vspace{0.2cm}

\begin{x}{\small\bf FACT} \ 
Let $\F \subset \K \subset \LL$ be fields $-$then
\[
[\LL:\F] \hsx = \hsx [\LL:\K]  \cdot [\K:\F].  
\]
\end{x}

\vspace{0.1cm}

\begin{x}{\small\bf EXAMPLE} \ 
Take $\F = \Q$, $\K = \Q(\sqrt{2})$, $\LL = \Q(\sqrt{2},\sqrt{3})$ $-$then
\begin{align*}
[\Q(\sqrt{2},\sqrt{3}):\Q] \ 
=\ [\Q(\sqrt{2},\sqrt{3}):\Q(\sqrt{2})] \cdot [\Q(\sqrt{2}):\Q]\\
=\ 2 \cdot 2 \\
=\  4.
\end{align*}
\end{x}

\vspace{0.1cm}

\begin{x}{\small\bf DEFINITION} \ 
$\LL$ is a 
\un{finite extension}
\index{finite extension} 
of $\K$ if $[\LL:\K]$ is finite and $\LL$ is an 
\un{infinite extension}
\index{infinite extension} 
of $\K$ if $[\LL:\K]$ is infinite.
\end{x}

\vspace{0.1cm}


\begin{x}{\small\bf EXAMPLE} \ 
$[\C:\R] = 2$ but $[\C:\Q] = 2^{\aleph_0}$.
\end{x}

\vspace{0.1cm}

Given $\LL/\K$ and $x \in \LL$, the 
\un{ideal $I_x$ of algebraic relations of $x$}
\index{ideal $I_x$ of algebraic relations of $x$} 
is the ideal in $\K[X]$ consisting of all polynomials admitting $x$ as a zero.

\vspace{0.2cm}

\begin{x}{\small\bf DEFINITION} \ 
$x$ is 
\un{algebraic}
\index{algebraic}
over $\K$ 
(\un{transcendental} over $\K$)
\index{transcendental (element over a field $\K$)} 
according to whether $I_x$ is nonzero (zero).  
I.e.: $x$ is algebraic over $\K$ (transcendental over $\K$) according to whether it is (or is not) 
a zero of a nonzero polynomial in $\K[X]$.
\end{x}

\vspace{0.1cm}

\begin{x}{\small\bf EXAMPLE} \ 
Take $\K = \Q$, $\LL = \C$ $-$then $\sqrt{-1}$ is algebraic over $\Q$ but $e$ and $\pi$ are transcendental over $\Q$.
\end{x}

\vspace{0.1cm}

\begin{x}{\small\bf FACT} \ 
Let $x \in \LL$ $-$then $x$ is algebraic over $\K$ iff $I_x$ is a nonzero prime ideal in $\K[X]$ or still, 
is a maximal ideal in $\K[X]$.
\end{x}

\vspace{0.1cm}

\begin{x}{\small\bf FACT} \ 
If $x \in \LL$ is algebraic over $\K$, then $I_x$ has a unique monic polynomial $p_x$ in $\K[X]$ as a generator: 
$I_x = \langle p_x \rangle$, the 
\un{minimal polynomial of $x$ over $\K$}.
\index{minimal polynomial of $x$ over $\K$} 

\vspace{0.1cm}

[Note: \  One can characterize $p_x$ as the monic polynomial in $\K[X]$ that admits $x$ as a zero and divides in $\K[X]$ 
every polynomial admitting $x$ as a zero.]

\end{x}

\vspace{0.1cm}

\begin{x}{\small\bf REMARK} \ 
The minimal polynomial of an element depends on the base field.  
E.g.: If $\K = \Q$ and $\LL = \C$, then $p_{\sqrt{-1}}(X) = X^2 + 1$ but if $\K = \LL = \C$, then 
$p_{\sqrt{-1}}(X) = X - \sqrt{-1}$.
\end{x}

\vspace{0.1cm}

\begin{x}{\small\bf FACT} \ 
If $x \in \LL$ is algebraic over $\K$, then its minimal polynomial $p_x$ is irreducible.
\end{x}

\vspace{0.1cm}


\begin{x}{\small\bf FACT} \ 
If $x \in \LL$ is algebraic over $\K$ and if $n = \deg p_x$, then $p_x$ is the only monic polynomial in $\K[X]$ 
of degree $n$ admitting $x$ as a zero.
\end{x}

\vspace{0.1cm}

\begin{x}{\small\bf FACT} \ 
If $x \in \LL$ is algebraic over $\K$, then the set $\{x^j: 0 \leq j \leq n-1\}$ is a linear basis of $\K(x)$ over $\K$, 
hence $[\K(x):\K] = n$.
\end{x}

\vspace{0.1cm}

\begin{x}{\small\bf EXAMPLE} \ 
Take $\K = \Q$, $\LL = \R$, $x = (2)^{1/3}$ $-$then $\Q((2)^{1/3})$ is a subfield of $\R$ and $(2)^{1/3}$ 
is algebraic over $\Q$, its minimal polynomial being $X^2 - 2$, so $[\Q((2)^{1/3}):\Q] = 3$.
\end{x}

\vspace{0.1cm}

\begin{x}{\small\bf DEFINITION} \ 
$\LL$ is an 
\un{algebraic extension}
\index{algebraic extension} 
of $\K$ if every element of $\LL$ is algebraic over $\K$.
\end{x}

\vspace{0.1cm}

\begin{x}{\small\bf FACT} \ 
If $[\LL:\K] < \infty$, then $\LL$ is an algebraic extension of $\K$.

\vspace{0.1cm}

[If $n = [\LL:\K]$ and if $x \in \LL$, then the sequence $x^j$ $(0 \leq j \leq n)$ is linearly dependent over $\K$, 
so there exists a sequence $a_j$ $(0 \leq j \leq n)$ of elements of $\K$ (not all zero) such that 
$\ds\sum\limits_{j = 0}^n a_j x^j = 0$.]
\end{x}

\vspace{0.1cm}

\begin{x}{\small\bf FACT} \ 
Suppose that $\K$ is infinite and $\LL$ is an algebraic extension of $\K$ $-$then 
\[
\card \K \hsx = \hsx \card \LL.
\]
\end{x}

\vspace{0.1cm}

\begin{x}{\small\bf EXAMPLE} \ 
$\R$ is not an algebraic extension of $\Q$.
\end{x}

\vspace{0.1cm}

\begin{x}{\small\bf DEFINITION} \ 
Let $\K$ be a field and let $\LL_1, \LL_2$ be field extensions of $\K$ $-$then 
a 
\un{$\K$-homomorphism}
\index{$\K$-homomorphism} 
$\phi:\LL_1 \ra \LL_2$ is a ring homomorphism such that 
$\restr{\phi}{\K} = \id_\K$, $\phi$ being called a 
\un{$\K$-isomorphism}
\index{$\K$-isomorphism} 
if it is in addition bijective (injectivity is automatic).

\vspace{0.1cm}

[Note: \ When $\LL_1 = \LL_2$, the term is 
\un{$\K$-automorphism}.]
\index{$\K$-automorphism} 

\end{x}

\vspace{0.1cm}

\begin{x}{\small\bf REMARK} \ 
If $\LL_1 = \LL_2$, call it $\LL$, and if $\LL$ is an algebraic extension of $\K$, then every $\K$-homomorphims 
$\phi:\LL \ra \LL$ is a $\K$-isomorphism.
\end{x}

\vspace{0.1cm}

\begin{x}{\small\bf FACT} \ 
Let $\K$ be a field and let $\LL_1$, $\LL_2$ be field extensions of $\K$.  Suppose that $f$ is an irreducible polynomial in 
$\K[X]$ and suppose that $x_1, x_2$ are, respectively, zeros of $f$ in  $\LL_1$, $\LL_2$  $-$then 
there is a unique $\K$-isomorphism 
$\K(x_1) \ra \K(x)$ such that $x_1 \ra x_2$.
\end{x}

\vspace{0.1cm}

[Note: \ The assumption that $f$ is irreducible cannot be dropped.]

\vspace{0.3cm}

\setcounter{theoremn}{0}
\[
\text{ADDENDUM}
\]

Let $\K$ be a field, $\LL/\K$ a field extension $-$then a sublset $S$ of $\LL$ is a 
\un{transcendence basis}
\index{transcendence basis} 
for $\LL/\K$ if \mS is algebraically independent over $\K$ and if $\LL$ is algebraic over $\K(S)$ 
(the subfield of $\LL$ generated by $\K \cup S$).

\begin{x}{\small\bf FACT} \ 
A transcendence basis for $\LL/\K$ always exists and any two have the same cardinality.
\end{x}

\vspace{0.1cm}

\begin{x}{\small\bf DEFINITION} \ 
The 
\un{transcendence degree}
\index{transcendence degree} 
$\trdeg(\LL/\K)$
\index{$\trdeg(\LL/\K)$} 
is the cardinality of any transcendence basis of $\LL/\K$.
\end{x}

\vspace{0.1cm}


\begin{x}{\small\bf EXAMPLE} \ 
Take $\K = \Q$, $\LL = \C$ $-$then $\trdeg(\C/\Q)$ is infinite (in fact uncountable).
\end{x}

\vspace{0.1cm}

\begin{x}{\small\bf EXAMPLE} \ 
Take $\K = \Q$, $\LL = \Q_p$ $-$then $\trdeg(\Q_p/\Q)$ is infinite (in fact uncountable).
\end{x}

\vspace{0.1cm}


\setcounter{theoremn}{0}

\newpage

\ \indent 

\[
\text{ALGEBRAIC CLOSURE}
\]

Let $\K$ be a field, $\LL/\K$ a field extension.

\vspace{0.2cm}

\begin{x}{\small\bf NOTATION} \ 
$A(\LL/\K)$ is the set of all elements of $\LL$ that are algebraic over $\K$.
\end{x}

\vspace{0.1cm}

\begin{x}{\small\bf DEFINITION} \ 
$A(\LL/\K)$ is the 
\un{algebraic closure}
\index{algebraic closure} 
of $\K$ in $\LL$.
\end{x}

\vspace{0.1cm}

\begin{x}{\small\bf EXAMPLE} \ 
Take $\K = \R$, $\LL = \C$ $-$then $A(\LL/\K) = \C$.

\vspace{0.1cm}

[Given $a + \sqrt{-1} \hsx b$, consider the polynomial
\[
(X - (a + \sqrt{-1} \hsx b)) (X - (a - \sqrt{-1} \hsx b)) \hsx = \hsx 
X^2 - 2 a X + a^2 + b^2.]
\]
\end{x}

\vspace{0.1cm}

\begin{x}{\small\bf FACT} \ 
$\LL$ is an algebraic extension of $\K$ iff $A(\LL/\K) = \LL$.
\end{x}

\vspace{0.1cm}

\begin{x}{\small\bf DEFINITION}  \ 
$\K$ is 
\un{algebraically closed}
\index{algebraically closed} 
in $\LL$ if every element of $\LL$ that is algebraic over $\K$ belongs to $\K$:
\[
A(\LL/\K) \hsx = \hsx \K.
\]
\end{x}

\vspace{0.1cm}

\begin{x}{\small\bf FACT} \ 

\[
\K \hsx \subset\hsx A(\LL/\K) \hsx\subset\hsx \LL.
\]
\end{x}

\vspace{0.1cm}

\begin{x}{\small\bf FACT}  \ 
$A(\LL/\K)$ is a field.
\end{x}

\vspace{0.1cm}

\begin{x}{\small\bf FACT}  \ 
$A(\LL/\K)$ is algebracally closed in $\LL$.

\vspace{0.1cm}

[Spelled out, if $x \in \LL$ is algebraic over $A(\LL/\K)$, then $x \in A(\LL/\K)$.]
\end{x}

\vspace{0.1cm}

\begin{x}{\small\bf SCHOLIUM} \ 
If $\K \subset \E \subset \LL$ and if $\E$ is an algebraic extension of $\K$, then
\[
\E \subset A(\LL/\K).
\]
\end{x}

\vspace{0.1cm}

\begin{x}{\small\bf DEFINITION} \ 
Take $\K = \Q$, $\LL = \C$ $-$then an 
\un{algebraic number}
\index{algebraic number} 
is a complex number which is algebraic over $\Q$, i.e., is an element of $A(\C/\Q)$.
\end{x}

\vspace{0.1cm}

\begin{x}{\small\bf FACT} \ 
$\card A(\C/\Q) \hsx = \hsx \aleph_0$.
\end{x}

\vspace{0.1cm}

\begin{x}{\small\bf FACT} \ 
$[A(\C/\Q):\Q] \hsx = \hsx \aleph_0$.

\vspace{0.1cm}

[Let $n$ be a postive integer $-$then the polynomial $X^n - 2$ is irreducible in $\Q[X]$, thus is the minimal polynomial of 
$(2)^{1/2}$ over $\Q$, so $[Q((2)^{1/2}):\Q] = n$, from which 
\[
[A(\C/\Q):\Q] \hsx \geq \hsx n.
\]
And this implies that 
\[
[A(\C/\Q):\Q] \hsx \geq \hsx \aleph_0.
\]
On the other hand, 
\[
[A(\C/\Q):\Q] \hsx \leq \hsx \card A(\C/\Q) \hsx = \hsx \aleph_0.]
\]
\end{x}

\vspace{0.1cm}

\begin{x}{\small\bf DEFINITION} \ 
A field $\F$ is 
\un{algebraically closed}
\index{algebraically closed} 
if every nonconstant polynomial in $\F[X]$ has a zero in $\F$.

\vspace{0.1cm}

[Note: \ This notion is absolute.]
\end{x}

\vspace{0.1cm}

\begin{x}{\small\bf EXAMPLE} \ 
Neither $\Q$ nor $\R$ is algebraically closed but $\C$ is algebraically closed.
\end{x}

\vspace{0.1cm}

\begin{x}{\small\bf FACT} \ 
$\F$ is algebraically closed iff every irreducible polynomial has degree 1.
\end{x}

\vspace{0.1cm}


\begin{x}{\small\bf FACT} \ 
$\F$ is algebraically closed iff every nonconstant polynomial $f$ in $\F[X]$ splits in $\F[X]$.

\vspace{0.1cm}

[Note: \ I.e.: Given $f$, there exists a postive integer $n$ and elements 
$a, a_1, \ldots, a_n$ (not necessarily distinct) of $\F$ such that 
\[
f(X) \hsx = \hsx a \hsx \prod\limits_{k = 1}^n (X - a_k).]
\]
\end{x}

\vspace{0.1cm}

\begin{x}{\small\bf FACT} \ 
If $\F$ is algebraically closed, then it is its only algebraic extension.
\end{x}

\vspace{0.1cm}

\begin{x}{\small\bf FACT} \ 
If there is an algebraically closed field extension $\F^\prime$ of $\F$ in which $\F$ is algebraically closed, then 
$\F$ is algebraically closed.

\vspace{0.1cm}

[Let $f \in \F[X]$ be a nonconstant polynomial $-$then $f$ has a zero $a^\prime$ in $\F^\prime$, hence 
$a^\prime$ is algebraic over $\F$, hence $a^\prime \in \F$ (since $\F$ is algebraically closed in $\F^\prime$).]
\end{x}

\vspace{0.1cm}

\begin{x}{\small\bf APPLICATION} \ 
Suppose that $\LL/\K$ is an algebraically closed field extension.  
Let $\F = A(\LL/\K)$, $\F^\prime = \LL$ to conclude that $A(\LL/\K)$ is algebraically closed.
\end{x}

\vspace{0.1cm}

\begin{x}{\small\bf EXAMPLE} \ 
Take $\K = \Q$, $\LL = \C$ $-$then $\C$ is algebraically closed, hence $A(\C/\Q)$ is algebraically closed.
\end{x}

\vspace{0.1cm}

\begin{x}{\small\bf FACT} \ 
Let $\K$ be a field, let $\LL$ be an algebraic closure of $\K$, and let $\M$ be an algebraically closed extension of $\K$ $-$then 
there exists a $\K$-monomorphism $\phi:\LL \ra \M$.
\end{x}

\vspace{0.1cm}

\begin{x}{\small\bf EXAMPLE} \ 
Take $\K = \R$, $\LL = \C$, $\M = \C$ $-$then the inclusion $\R \ra \C$ admits two distinct extensions to $\C$, 
viz. the identity and the complex conjugation
(and these are the only $\R$-automorphisms of $\C$).

\vspace{0.1cm}

[Note: \ Therefore uniqueness of the extending $\K$-monomorphism cannot be asserted.]
\end{x}

\vspace{0.1cm}

\begin{x}{\small\bf EXAMPLE} \ 
If $\E \neq \R$ is an algebraic extension of $\R$, then $\E$ is isomorphic to $\C$.

\vspace{0.1cm}

[Take $\K = \R$, $\LL = \E$, $\M = \C$ $-$then there exists an $\R$-monomorphism 
$\phi:\E \ra \C$, hence 
\[
2 \hsx = \hsx [\C:\R] \hsx = \hsx [\C:\phi(\E)] \cdot [\phi(\E):\R],
\]
from which $\C = \phi(\E) \approx \E$.]
\end{x}

\vspace{0.1cm}

\begin{x}{\small\bf DEFINITION} \ 
Given a field $\F$, an 
\un{algebraic closure}
\index{algebraic closure} 
of $\F$ is an algebraicallly closed algebraic extension of $\F$.
\end{x}

\vspace{0.1cm}

\begin{x}{\small\bf EXAMPLE} \ 
$\C$ is an algebraic closure of $\R$ but $\C$ is not an algebraic closure of $\Q$ (since it is not algebraic over $\Q$).
\end{x}

\vspace{0.1cm}

\begin{x}{\small\bf EXAMPLE} \ 
$A(\C/\Q)$ is an algebraic closure of $\Q$.
\end{x}

\vspace{0.1cm}

\begin{x}{\small\bf STEINITZ THEOREM} \ 
Every field $\F$ admits an algebraic closure $\F^\cl$ and any two algebraic closures of $\F$ are 
$\F$-isomorphic.
\end{x}

\vspace{0.1cm}

\begin{x}{\small\bf FACT} \ 
Every automorphism of $\F$ can be extended to an automorphism of $\F^\cl$.

\vspace{0.1cm}

[Note: \ In general, if $\F_1$ and $\F_2$ are fields, then every isomorphism from $\F_1$ to $\F_2$ can be extended to an 
isomorphism from $\F_1^\cl$ to $\F_2^\cl$.]

\end{x}

\vspace{0.1cm}

\begin{x}{\small\bf FACT} \ 
If $\LL / \K$ is an algebraic extension of $\K$, then $\LL$ is $\K$-isomorphic to a subfield of $\K^\cl$.
\end{x}


\setcounter{theoremn}{0}

\newpage

\ \indent 

\[
\text{TRACES AND NORMS}
\]

Let $\K$ be a field, $\LL/\K$ a field extension of $\K$ $-$then each $x \in \LL$ gives rise to a linear transformation 
\[
M_x: \LL \ra \LL
\]
defined by 
\[
M_x(y) \hsx = \hsx x y.
\]

\vspace{0.2cm}

\begin{x}{\small\bf DEFINITION} \ 
The 
\un{trace}
\index{trace}\index{$\tT_{\LL/\K}$} 
of $\LL$ over $\K$ is the function 
\[
\begin{cases}
\ \tT_{\LL/\K}:\LL \ra \K\\
\ \tT_{\LL/\K}(x) \ \hsx = \hsx \tr(M_x).
\end{cases}
\]
\end{x}

\vspace{0.1cm}

\begin{x}{\small\bf DEFINITION} \ 
The 
\un{norm}
\index{norm}\index{$\tN_{\LL/\K}$}
of $\LL$ over $\K$ is the function 
\[
\begin{cases}
\ \tN_{\LL/\K}:\LL \ra \K\\
\ \tN_{\LL/\K}(x) \ \hsx = \hsx \det(M_x).
\end{cases}
\]
\end{x}

\vspace{0.1cm}

\begin{x}{\small\bf PROPERTIES} \ 
$\forall \ x, y \in \LL$, $\forall \ a \in \K$:

\vspace{0.1cm}

1. \ $\tT_{\LL/\K}(x + y) \hsx = \hsx \tT_{\LL/\K}(x) + \tT_{\LL/\K}(y)$.

\vspace{0.1cm}

2. \ $\tT_{\LL/\K}(a) \hsx = \hsx [\LL:\K] a$.

\vspace{0.1cm}

3. \ $\tN_{\LL/\K}(x y) \hsx = \hsx \tN_{\LL/\K}(x) \tN_{\LL/\K}(y)$.

\vspace{0.1cm}

4. \ $\tN_{\LL/\K}(a) \hsx = \hsx a^{[\LL:\K]}$.
\end{x}

\vspace{0.1cm}


\begin{x}{\small\bf FACT} \ 
If $\E$ is a subfield of $\LL$ containing $\K$, then 
\[
\begin{cases}
\ \tT_{\LL/\K}(x) \hsx = \hsx \tT_{\E/\K} ( \tT_{\LL/\E} (x))\\
\ \tN_{\LL/\K}(x) \hsx = \hsx \tN_{\E/\K} ( \tN_{\LL/\E} (x))
\end{cases}
.
\]
\end{x}

\vspace{0.1cm}

\begin{x}{\small\bf EXAMPLE} \ 
Let $\theta \in \K^\times - (\K^\times)^2$ and put $\LL = \K(\sqrt{\theta})$ $-$then $\forall \ a, b \in \K$, 
\[
\begin{cases}
\ \tT_{\LL/\K}(a + b \sqrt{\theta}) \hsx = \hsx 2a\\
\ \tN_{\LL/\K}(x) (a + b \sqrt{\theta}) \hsx = \hsx a^2 - b^2 \theta
\end{cases}
.
\]
\end{x}

\renewcommand{\thepage}{APPENDIX C-\arabic{page}}
\chapter{
TOPICS IN GALOIS THEORY}
\setlength\parindent{2em}
\setcounter{theoremn}{0}

\[
\text{GALOIS CORRESPONDENCES}
\]

\[
\text{FINITE GALOIS THEORY}
\]

\[
\text{INFINITE GALOIS THEORY}
\]

\[
\text{$\K^\sep$ AND $\K^\ab$}
\]
\newpage

\ \indent 
\[
\text{GALOIS CORRESPONDENCES}
\]

Given a field $\F$, $\Aut(\F)$ stands for its associated group of field automorphisms.

\begin{x}{\small\bf EXAMPLE} \ 
Take $\F = \Q$ $-$then $\Aut(\Q)$ is trivial.
\end{x}

\vspace{0.1cm}

\begin{x}{\small\bf EXAMPLE} \ 
Take $\F = \R$ $-$then $\Aut(\R)$ is trivial.

\vspace{0.1cm}

[Let $\phi \in \Aut(\R)$ $-$then $\restr{\phi}{\Q} = \id_\Q$.  Next: 
\begin{align*}
x < y \implies \phi(y) - \phi(x) \ 
&=\ \phi(y - x) \\
&=\ \phi((\sqrt{y - x}\hsx)^2) \\
&=\ \phi(\sqrt{y - x}\hsx)^2 \\
&> \ 0.
\end{align*}
If now $\phi \neq \id_\R$, choose $x$ such that $\phi(x) \neq x$ $-$then there are two possibilities.

\vspace{0.1cm}

\qquad\qquad \textbullet \quad $x < \phi(x)$: \ Choose $q \in \Q$: $x < q < \phi(x)$, so 
$\phi(x) < \phi(q) = q < \phi(x)$.   Contradiction.

\vspace{0.1cm}

\qquad\qquad \textbullet \quad $\phi(x) < x$: \ Choose $q \in \Q$: $\phi(x) < q < x$, so
$\phi(x) < q= \phi(q) < \phi(x)$.  Contradiction.
\end{x}

\vspace{0.1cm}

\begin{x}{\small\bf EXAMPLE} \ 
Take $\F = \C$ $-$then $\Aut(\C)$ is infinite.

\vspace{0.1cm}

[Any automorphism $\phi:\C \ra \C$ will fix $\Q$ and any continuous automorphism $\phi:\C \ra \C$ will fix its closure $\R$, 
there being two such, viz. the identity and the complex conjugation, all others being discontinuous.]

\vspace{0.1cm}

[Note: \ As an illustration, consider the automorphism
\[
a +  b \hsx \sqrt{2} \ra a - b \hsx \sqrt{2} \qquad (a, b \in \Q)
\]
of the field $\Q(\sqrt{2})$ $-$then it can be extended to an automorphism of $\C$ via the following procedure.

\vspace{0.1cm}

\qquad\qquad 1. \ Extend to $\K \equiv \Q(\sqrt{2})^\cl \subset \C$.


\vspace{0.1cm}

\qquad\qquad 2. \ Choose a transcendence basis \mS for $\C/\K$ and extend to $\K(S)$. 

\vspace{0.1cm}

\qquad\qquad 3. \ Extend from $\K(S)$ to $\C$.]

\end{x}

\vspace{0.1cm}

\begin{x}{\small\bf DEFINITION} \ 
Let \mG be a group of automorphisms of $\F$ $-$then the subfield 
\[
\Inv(G) \hsx = \hsx \{x: \sigma x = x\} \qquad (\sigma \in G)
\]
is called the 
\un{invariant field}
\index{invariant field} 
associated with \mG.
\end{x}

\vspace{0.1cm}

\begin{x}{\small\bf DEFINITION} \ 
Given a subfield $\E \subset \F$, the group consisting of all automorphisms of $\F$ leaving every element of $\E$ invariant is 
denoted by $\Gal(\F/\E)$, the 
\un{Galois group}
\index{Galois group} 
of $\F$ over $\E$.
\end{x}

\vspace{0.1cm}

\begin{x}{\small\bf EXAMPLE} \ 
Take $\E = \R$, $\F = \C$ $-$then $\Gal(\C/\R) = \{\id_\C,\sigma\}$, where $\sigma$ is the complex conjugation.
\end{x}

\vspace{0.1cm}

\begin{x}{\small\bf EXAMPLE} \ 
Take $\E = \Q$, $\F = \Q((2)^{1/3})$ $-$then $\Gal(\Q((2)^{1/3}) / \Q)$ is trivial.
\end{x}

\vspace{0.1cm}

\begin{x}{\small\bf EXAMPLE}  \ 
Take $\E = \Q$, $\F = \Q(\omega_n)$ ($\omega_n$ a primitive $n^\nth$ root of unity in $\C$) $-$then 
\[
\Gal(\Q(\omega_n) / \Q) \approx (\Z/n\Z)^\times.
\]
\end{x}

\vspace{0.1cm}

\begin{x}{\small\bf FACT} \ 
We have
\[
G \subset \Gal(\F/\Inv(G)).
\]
\end{x}

\begin{x}{\small\bf FACT} \ 
We have
\[
\E \subset \Inv(\Gal(\F/\E)).
\]
\end{x}

\vspace{0.1cm}

\begin{x}{\small\bf FACT} \ 
\[
G \hsx \subset\hsx \Gal(\F/\E) \hsx\Leftrightarrow\hsx \E \hsx\subset\hsx \Inv(G).
\]
\end{x}

\vspace{0.1cm}


\begin{x}{\small\bf FACT} \ 
\\

\qquad\qquad \textbullet \quad $G_1 \subset G_2 \subset \Aut(\F) \implies \Inv(G_1) \supset \Inv(G_2)$.

\vspace{0.1cm}

\qquad\qquad \textbullet \quad $\E_1 \subset \E_2 \subset \F \implies \Gal(\F/\E_2) \subset \Gal(\F/\E_1)$.

\end{x}

\vspace{0.1cm}

\begin{x}{\small\bf DEFINITION} \ 
Let $\F$ be a field.

\qquad\qquad \textbullet \quad A 
\un{Galois group}
\index{Galois group}
on $\F$ is a group \mG of automorphisms of $\F$ such that 
\[
G \hsx = \hsx \Gal(\F/\Inv(G)).
\]

\vspace{0.1cm}

\qquad\qquad \textbullet \quad An 
\un{invariant field}
\index{invariant field}
in $\F$ is a subfield $\E$ of $\F$ such that 
\[
\E \hsx = \hsx \Inv(\Gal(\F/\E)).
\]

\vspace{0.1cm}
\end{x}

\vspace{0.1cm}

\begin{x}{\small\bf EXAMPLE} \ 
$\Aut(\F)$ is a Galois group on $\F$.

\vspace{0.1cm}

[For 
\begin{align*}
\Aut(\F) \ 
&\subset \Gal(\F/\Inv(\Aut(\F)))\\
&=\ \Aut(\F).]
\end{align*}
\end{x}

\vspace{0.1cm}

\begin{x}{\small\bf EXAMPLE} \ 
$\{\id_\F\}$ is a Galois group on $\F$

\vspace{0.1cm}

[For 
\begin{align*}
\{\id_\F\} \ 
&\subset \Gal(\F/\Inv(\{\id_\F\}))\\
&=\ \Gal(\F/\F) \\
&=\ \{\id_\F\}.]
\end{align*}
\end{x}

\vspace{0.1cm}

\begin{x}{\small\bf EXAMPLE} \ 
$\F$ is an invariant field on $\F$.
\end{x}

\vspace{0.1cm}

\begin{x}{\small\bf REMARK} \ 
Recall that a field is 
\un{prime}
\index{prime} 
if it possesses no proper subfields, these being the fields 
isomorphic to $\Q$ (characteristic 0) or 
isomorphic to $\Z/p\Z$ (characteristic p).  
A prime field admits no automorphism other than the identity.
\end{x}

\vspace{0.1cm}

\begin{x}{\small\bf ABSOLUTE GALOIS CORRESPONDENCE} \ 
Let $\F$ be a field.

\qquad \textbullet \quad If $\E$ is a subfield of $\F$, then $\Gal(\F/\E)$ is a Galois group on $\F$.

\vspace{0.1cm}

\qquad \textbullet \quad If \mG is a group of automorphisms of $\F$, then $\Inv(G)$ is an invariant field in $\F$.

\vspace{0.1cm}

And:  \ The arrow $\E \ra \Gal(\F/\E)$ from the set of all invariant fields in $\F$ to the set of all Galois groups 
on $\F$ and the arrow $G \ra \Inv(G)$ from the set of all Galois groups on $\F$ to the set of all invariant fields in $\F$ 
are mutually inverse inclusion reversing bijections. 
\end{x}

\vspace{0.1cm}

\begin{x}{\small\bf RELATIVE GALOIS CORRESPONDENCE} \ 
Let $\K$ be a field and let $\LL$ be a field extension of $\K$.

\qquad \textbullet \quad If $\K \subset \E \subset \LL$, then $\Gal(\LL/\E)$ is a Galois group on $\LL$ contained in 
$\Gal(\LL/\K)$.

\vspace{0.1cm}

\qquad \textbullet \quad If \mG is a subgroup of $\Gal(\LL/\K)$, then $\Inv(G)$ is an invariant field in $\LL$ containing 
$\K$.

\vspace{0.1cm}

And: \ The arrow $\E \ra \Gal(\LL/\E)$ from the set of all invariant fields in $\LL$ containing $\K$ to the set of all Galois 
groups on $\LL$ contained in $\Gal(\LL/\K)$ and the arrow $G \ra \Inv(G)$ from the set of all Galois groups on $\LL$ 
contained in $\Gal(\LL/\K)$ to the set of all invariant fields in $\LL$ containing $\K$ are mutually inverse inclusion reversing bijections.

\end{x}


\setcounter{theoremn}{0}

\newpage

\ \indent 

\[
\text{FINITE GALOIS THEORY}
\]

\begin{x}{\small\bf DEFINITION} \ 
A field extension $\LL/\K$ is 
\un{Galois over $\K$}
\index{Galois over $\K$} 
(or is a 
\un{Galois extension of $\K$})
\index{Galois extension of $\K$} 
if $\LL$ is algebraic over $\K$ and $\K$ is an invariant field on $\L$ or still, 
\[
\K \hsx = \hsx \Inv(\Gal(\LL/\K)).
\]
\end{x}

\vspace{0.1cm}

\begin{x}{\small\bf FACT} \ 
If $\LL/\K$ is a finite Galois extension and if $\LL \supset \E \supset \K$ is an intermediate field, then $\LL$ is Galois over $\E$.
\end{x}

\vspace{0.1cm}

\begin{x}{\small\bf FACT} \ 
If $\LL/\K$ is a finite Galois extension and if $\LL \supset \E \supset \K$ is an intermediate field, then $\E$ is Galois over $\K$ 
iff $\Gal(\LL/\E)$ is a normal subgroup of $\Gal(\LL/\K)$. 

\vspace{0.1cm}

[Note: \ Under the assumption that $\E$ is Galois over $\K$, there is an arrow of restriction
\[
\Gal(\LL/\K) \ra \Gal(\E/\K).
\]
It is surjective with kernel $\Gal(\LL/\E)$, from which an exact sequence of groups:
\[
1 \ra \Gal(\LL/\E) \ra \Gal(\LL/\K)  \ra \Gal(\E/\K) \ra 1.]
\]
\end{x}

\vspace{0.1cm}

\begin{x}{\small\bf RECOGNITION  PRINCIPLE} \ 
If $\LL/\K$ is a finite extension, then $\LL$ is Galois over $\K$ iff 
\[
\card \Gal(\LL/\K) \hsx = \hsx [\LL:\K].
\]

\vspace{0.1cm}

[Note: \ If $\LL/\K$ is a finite extension, then a priori
\[
\card \Gal(\LL/\K) \hsx \leq \hsx [\LL:\K],
\]
the inequality being strict in general.  Matters break down if it is a question of infinite extensions.  
E.g.: \ If $\Q^\cl$ is an algebraic closure of $\Q$, then 
\[
[\Q^\cl:\Q] \hsx = \hsx \aleph_0
\]
while
\[
\card \Gal(\Q^\cl / \Q) \hsx = \hsx 2^{\aleph_0}.]
\]
\end{x}

\vspace{0.1cm}

\begin{x}{\small\bf EXAMPLE} \ 
Let $\F$ be a field of characteristic 0 and let $a \in \F^\times - (\F^\times)^2$.

Form the quadratic extension $\F(\sqrt{a})$ $-$then $[\F(\sqrt{a}):\F] = 2$, while 
$\Gal(\F(\sqrt{a}) / \F) = \{\id,\sigma\}$ $(\sigma(\sqrt{a}) = - \sqrt{a})$.  
Therefore $\F(\sqrt{a})$ is a Galois extension of $\F$.
\end{x}

\vspace{0.1cm}

\begin{x}{\small\bf EXAMPLE} \ 
Take $\K = \Q$, $\LL = \Q((2)^{1/3})$ $-$then $[\Q((2)^{1/3}):\Q] = 3$ but 
$\Gal(\Q((2)^{1/3}) / \Q)$ is trivial.  Therefore $\Q((2)^{1/3})$ is not a Galois extension of $\Q$.
\end{x}

\vspace{0.1cm}

\begin{x}{\small\bf EXAMPLE} \ 
Take $\K = \Q$, $\LL = \Q((2)^{1/3}, \omega)$, where 
\[
\omega \hsx = \hsx \exp(2 \pi \sqrt{-1} / 3).
\]
Then
\[
[\Q((2)^{1/3}, \omega):\Q] \hsx = \hsx  
[\Q((2)^{1/3}, \omega):\Q((2)^{1/3})]  \cdot [\Q((2)^{1/3}):\Q] \hsx = \hsx
2 \cdot 3 \hsx = \hsx 6.
\]
On the other hand, the six functions 
\begin{align*}
&\vsy(2)^{1/3} \ra (2)^{1/3}, \quad \omega \ra \omega\\
&\vsy(2)^{1/3} \ra \omega(2)^{1/3}, \quad \omega \ra \omega\\
&\vsy(2)^{1/3} \ra (2)^{1/3}, \quad \omega \ra \omega^2\\
&\vsy(2)^{1/3} \ra \omega(2)^{1/3}, \quad \omega \ra \omega^2\\
&\vsy(2)^{1/3} \ra \omega^2 (2)^{1/3}, \quad \omega \ra \omega\\
&\vsy(2)^{1/3} \ra \omega^2 (2)^{1/3}, \quad \omega \ra \omega^2
\end{align*}
extend to distinct automorphisms of $\Q((2)^{1/3}, \omega) /\Q$.
Therefore $\Q((2)^{1/3}, \omega)$ is a Galois extension of $\Q$.
\end{x}

\vspace{0.1cm}

\begin{x}{\small\bf FUNDAMENTAL THEOREM OF FINITE GALOIS THEORY}  \ 
Suppose that $\LL$ is a finite Galois extension of $\K$.

\vspace{0.1cm}

\qquad \textbullet \quad If $\LL \supset \E \supset \K$, then
\[
[\Gal(\LL/\K):\Gal(\LL/\E)] \hsx = \hsx  [\E:\K].
\]
\vspace{0.1cm}

\qquad \textbullet \quad If $G \subset \Gal(\LL/\K)$, then 
\[
[\Inv(G):\K] \hsx = \hsx [\Gal(\LL/\K):G].
\]

\vspace{0.1cm}

And: \ The arow $\E \ra \Gal(\LL/\E)$ from the set of all intermediate fields between $\K$ and $\LL$ 
to the set of all subgroups of $\Gal(\LL/\K)$ and the arrow $G \ra \Inv(G)$ from the 
set of all subgroups of $\Gal(\LL/\K)$ to the set of all intermediate fields between $\K$ and $\LL$ 
are mutually inverse inclusion reversing bijections.
\end{x}

\vspace{0.1cm}

\begin{x}{\small\bf REMARK} \ 
Given a finite Galois extension $\LL/\K$, the problem of determining all intermediate fields 
$\LL \supset \E \supset \K$ amounts to finding all subgroups of $\Gal(\LL/\K)$, a finite problem.

\vspace{0.1cm}

[Note: \ The fact that there are but finitely many intermediate fields cannot be established by a vector space argument alone.]

\end{x}

\vspace{0.1cm}

\begin{x}{\small\bf EXAMPLE} \ 
The field $\Q((2)^{1/3}, \omega)$ is Galois over $\Q$ and its Galois group is a group of order 6, there being two 
possibilities, viz. the cyclic group $\Z/6\Z$ and the symmetric group $S_3$.  
Since $\Q((2)^{1/3})$ is not Galois over $\Q$, the group 
\[
\Gal(\Q((2)^{1/3}, \omega) / \Q((2)^{1/3}))
\]
is not a normal subgroup of $\Gal(\Q((2)^{1/3}, \omega) / \Q$).  
But every subgroup of an abelian group
is normal, so the conclusion is that 
\[
G \hsx \equiv \hsx \Gal(\Q((2)^{1/3}, \omega) / \Q) \hsx \approx \hsx S_3.  
\]
Proceeding, there are $\Q$-automorphisms $\sigma, \tau$ of  $\Q((2)^{1/3}, \omega)$ defined by the specification
\[
\begin{cases}
\ \sigma : (2)^{1/3} \ra \omega(2)^{1/3}, \quad \omega \ra \omega\\
\ \tau : (2)^{1/3} \ra (2)^{1/3}, \quad \omega \ra \omega^2
\end{cases}
.
\]
Then $\sigma$ has order 3, $\tau$ has order 2, and $\sigma \tau \neq \tau \sigma$.  
The subgroups of \mG are 
\[
\langle \id \rangle, \quad 
\langle \sigma \rangle, \quad
\langle \tau \rangle, \quad
\langle \sigma \tau \rangle, \quad
\langle \sigma^2 \tau \rangle, \quad
G
\]
and the corresponding intermediate fields are
\[
\Q((2)^{1/3}, \omega), \quad
\Q(\omega), \quad
\Q((2)^{1/3}), \quad
\Q(\omega^2 (2)^{1/3}), \quad
\Q(\omega(2)^{1/3}), \quad
\Q.
\]
\end{x}

\vspace{0.1cm}

\begin{x}{\small\bf FACT} \ 
Let $\K$ be a finite Galois extension of $\F$ and let $\LL$ be an arbitrary finite extension of $\F$ $-$then 
$\K \vee \LL \supset \LL$ is a Galois extension and 
\[
\Gal(\K \vee \LL / \LL) \hsx \approx \hsx \Gal(\K / \K \cap \LL).
\]
In addition, 
\[
[\K \vee \LL:\LL] \hsx = \hsx [\K:\K \cap \LL].
\]

\vspace{0.1cm}

[Note: \ Tacitly, $\K$ and $\LL$ lie inside some common field $\M$, hence $\K \vee \LL$ is the subfield of $\M$ 
generated by $\K$ and $\LL$.  This said, the arrow 
\[
\Gal(\K \vee \LL /\LL) \ra \Gal(\K / \K \cap \LL)
\]
sends $\sigma$ to its restriction $\restr{\sigma}{\K}$.]
\end{x}

\vspace{0.1cm}

\begin{x}{\small\bf FACT} \ 
Suppose that $\LL$ is a finite Galois extension of $\K$ $-$then

\vspace{0.3cm}

\qquad \textbullet \quad $\tN_{\LL / \K} (x) \hsx = \hsx \ds\prod\limits_{\sigma \in \Gal(\LL / \K)} \sigma x$

\vspace{0.1cm}

\qquad \textbullet \quad $\tT_{\LL / \K} (x) \hsx = \hsx \ds\sum\limits_{\sigma \in \Gal(\LL / \K)} \sigma x$.
\end{x}

\vspace{0.1cm}

\begin{x}{\small\bf NORMAL BASIS THEOREM} \ 
If $\LL / \K$ is finite Galois, then $\exists \ x \in \LL$ such that $\{\sigma x: \sigma \in \Gal(\LL / \K)\}$ is a basis for $\LL / \K$.
\end{x}


\setcounter{theoremn}{0}

\newpage

\ \indent 

\[
\text{INFINITE GALOIS THEORY}
\]

\begin{x}{\small\bf FACT} \ 
If $\K$ is a field and if $\LL$ is an infinite Galois extension of $\K$, then 
\[
\card \Gal(\LL / \K) \hsx \geq \hsx 2^{\aleph_0}.
\]
\end{x}

\vspace{0.1cm}

\begin{x}{\small\bf APPLICATION} \ 
The Galois group of an infinite Galois extension cannot be cyclic.
\end{x}

\vspace{0.1cm}

\begin{x}{\small\bf FACT} \ 
If $\F$ is a field and if $G \subset \Aut(\F)$ is a finite group of automorphisms of $\F$, then \mG is a Galois group on $\F$: 
The a priori containment
\[
G \hsx\subset\hsx \Gal(\F / \Inv(G))
\]
is an equality: 
\[
G \hsx = \hsx \Gal(\F / \Inv(G)).
\]
\end{x}

\vspace{0.1cm}

\begin{x}{\small\bf REMARK} \ 
In general, an infinite group of automorphisms of a field need not be a Galois group.
\end{x}

\vspace{0.1cm}

Given a field $\F$ and an element $a \in \F$, let $D_a$ denote the discrete topological space having $\F$ as its set of points 
$-$then the elements of the product
\[
\prod\limits_{a \in \F} D_a
\]
are just the maps $\F^\F$ from $\F$ to $\F$.

When equipped with the product topology, $\F^\F$  is Hausdorff and totally disconnected 
(but not discrete if $\card \F \geq \aleph_0$).  
Since $\Aut(\F)$ is contained in $\F^\F$, it can be endowed with the relativized product topology, the so-called 
\un{finite topology}.
\index{finite topology}
\vspace{0.2cm}


\begin{x}{\small\bf \un{N.B.}} \ 
Given $\phi \in \Aut(\F)$ and a finite subset \mA of $\F$, let $\Omega_\phi(A)$ be the set of all automorphisms of $\F$ 
that agree with $\phi$ on \mA $-$then $\Omega_\phi(A)$ is open and the collection $\{\Omega_\phi(A)\}$ is a neighborhood basis at $\phi$.
\end{x}

\vspace{0.1cm}

\begin{x}{\small\bf FACT} \ 
In the finite topology, $\Aut(\F)$ is a topological group (as well as being Hausdorff and totally disconnected).
\end{x}

\vspace{0.1cm}

In what follows, if $\Gamma \subset \Aut(\F)$ is a group of automorphisms of $\F$, it will be understood that $\Gamma$ 
carries the relativized finite topology.

\vspace{0.2cm}

\begin{x}{\small\bf FACT} \ 
Suppose that $\Gamma \subset \Aut(\F)$ is compact $-$then $\Gamma$ is a Galois group on $\F$.
\end{x}

\vspace{0.1cm}

\begin{x}{\small\bf REMARK} \ 
A group of automorphisms of $\F$ is compact iff it is closed in $\Aut(\F)$ and has finite orbits.
\end{x}

\vspace{0.1cm}

\begin{x}{\small\bf FACT} \ 
If $\K$ is a field and if $\LL$ is an extension of $\K$, then 
\[
\Gal(\LL / \K) \subset \Aut(\LL)
\]
is closed.
\end{x}

\begin{x}{\small\bf FACT} \ 
If $\K$ is a field and if $\LL$ is an algebraic extension of $\K$, then 
\[
\Gal(\LL / \K) \subset \Aut(\LL)
\]
is compact.

\vspace{0.1cm}

[Note: \ If $\LL$ is finite over $\K$ (hence algebraic), then $\Gal(\LL / \K)$ is discrete.]
\end{x}

\vspace{0.1cm}

\begin{x}{\small\bf REMARK} \ 
The compactness of the Galois group does not characterize algebraic extensions (there exist transcendental extensions with a finite Galois group).


\vspace{0.1cm}

[Note: \ If $\K$ is an infinite field and if $\K(\xi)$ is a simple transcendental extension of $\K$, then $\Gal(\K(\xi) / \K)$ is not compact.]
\end{x}

\vspace{0.1cm}

\begin{x}{\small\bf FUNDAMENTAL THEOREM OF INFINITE GALOIS THEORY} \ 

Suppose that $\LL$ is an infinite Galois extension of $\K$ (hence algebraic, hence $\Gal(\LL / \K)$ compact).

\vspace{0.1cm}

\qquad \textbullet \quad If $\LL \supset \E \supset \K$, then $\Gal(\LL / \E)$ is a closed subgroup of $\Gal(\LL/\K)$ 
(thus is a compact subgroup of $\Gal(\LL/\K)$).  

\vspace{0.1cm}

\qquad \textbullet \quad If \mG is a closed subgroup of $\Gal(\LL/\K)$ (thus is a compact subgroup of $\Gal(\LL/\K))$, then 
$\Inv(G)$ is an intermediate field between $\K$ and $\LL$.

\vspace{0.1cm}

And: \ The arrow $\E \ra \Gal(\LL/\E)$ from the set of all intermediate fields between $\K$ and $\LL$  to the set of all closed  
subgroups of $\Gal(\LL/\K)$ and the arrow $G \ra \Inv(G)$ from the set of all closed  
subgroups of $\Gal(\LL/\K)$ to the set of all intermediate fields between $\K$ and $\LL$
are mutually inverse inclusion reversing bijections. 
\end{x}

\vspace{0.1cm}

\begin{x}{\small\bf REMARK} \ 
Since $\LL / \K$ is an infinite Galois extension, $\Gal(\LL / \K)$ always contains a subgroup that is not closed.

\vspace{0.1cm}

[Any infinite group has a countably infinite subgroup (consider the subgroup generated by a countably infinite subset).  
On the other hand, an infinite compact totally disconnected Hausdorff group has cardinality at least that of the continuum 
(it has a quotient which is homeomorphic to the Cantor set).]
\end{x}

\vspace{0.1cm}

\begin{x}{\small\bf FACT} \ 
$\E/\K$ is finite iff $\Gal(\LL / \E)$ is open.
\end{x}

\vspace{0.1cm}

\begin{x}{\small\bf FACT} \ 
$\E/\K$ is Galois iff $\Gal(\LL / \E)$ is normal.

\vspace{0.1cm}

[Note: \ Canonically, 
\[
\Gal(\E / \K) \hsx \approx \hsx \Gal(\LL / \K) / \Gal(\LL  / \E),
\]
this being a topological identification if $\Gal(\LL / \K)  / \Gal(\LL / \E)$ is given the quotient topology.]
\end{x}

\vspace{0.1cm}

\begin{x}{\small\bf \un{N.B.}} \ 
$\LL$ is Galois over $\E$.
\end{x}

\vspace{0.1cm}

\begin{x}{\small\bf NOTATION} \ 

\vspace{0.2cm}

\qquad \textbullet \quad
$\ds\bigvee\limits_{i \in I} \E_i$ is the subfield generated by the union $\ds\bigcup\limits_{i \in I} \E_i$.

\vspace{0.1cm}

\qquad \textbullet \quad
$\ds\bigvee\limits_{i \in I} G_i$ is the subgroup generated by the union $\ds\bigcup\limits_{i \in I} G_i$.

\vspace{0.1cm}

\end{x}

\vspace{0.1cm}

\begin{x}{\small\bf FACT} \ 
Let $\LL$ be an infinite Galois extension of $\K$.

\vspace{0.2cm}

\qquad \textbullet \quad If $\E_i$ $(i \in I)$ is a nonempty family of intermediate fields between $\K$ and $\LL$, then 
\[
\Gal \biggl(\hsx \LL / \bigcap\limits_{i \in I} \E_i \biggr) \  = \  
\ov{\bigvee\limits_{i \in I} \Gal(\LL / \E_i)}.
\]

\vspace{0.1cm}

\qquad \textbullet \quad If $G_i$ $(i \in I)$ is a nonempty family of closed subgroups of $\Gal(\LL / \K)$, then 
\[
\Inv \biggl(\hsx\bigcap\limits_{i \in I} G_i \biggr) \  = \  
\bigvee\limits_{i \in I} \Inv(G_i).
\]

\vspace{0.1cm}

\end{x}

\vspace{0.1cm}

\begin{x}{\small\bf EXAMPLE} \ 
Take $\K = \Q$, $\LL =\Q(\sqrt{2},\sqrt{3},\sqrt{5}, \ldots)$ (incorporate all primes) $-$then $\LL$ is Galois 
(and infinite) over $\K$ (being the union of 
$\Q$, 
$\Q(\sqrt{2})$, 
$\Q(\sqrt{2},\sqrt{3})$, 
$\Q(\sqrt{2},\sqrt{3},\sqrt{5})$ 
and so on).  Here $\Gal(\LL / \K)$ is a countably infinite direct product of copies of $\Z / 2\Z$.  
Accordingly, every $\K$-automorphism of $\LL$ differet from $\id_\LL$ is an element of order 2.
\end{x}

\vspace{0.1cm}

\begin{x}{\small\bf EXAMPLE} \ 
Take $\K = \Q$, $\LL = A(\C / \Q)$ $-$then $\LL$ is Galois (and infinite) over $\K$.
\end{x}


\setcounter{theoremn}{0}

\newpage

\ \indent 

\[
\text{$\K^\sep$ AND $\K^\ab$}
\]

Let $\K$ be a field, $\LL/\K$ a field extension.

\vspace{0.3cm}

\begin{x}{\small\bf DEFINITION} \ 
An element of $\LL$ is 
\un{separable}
\index{separable (element of a field)} 
if it is algebraic over $\K$ and is a simple zero of its minimal polynomial.
\end{x}

\vspace{0.1cm}

\begin{x}{\small\bf NOTATION} \ 
$S(\LL/\K)$ is the set of all elements of $\LL$ that are separable over $\K$.

\vspace{0.1cm}

[Note: \ Therefore
\[
S(\LL/\K) \subset A(\LL/\K)
\]
and 
\[
S(\LL/\K) \hsx = \hsx A(\LL/\K)
\]
if the characteristic of $\K$ is zero.]
\end{x}

\vspace{0.1cm}

\begin{x}{\small\bf DEFINITION} \ 
$S(\LL/\K)$ is the 
\un{separable closure}
\index{separable closure} 
of $\K$ in $\LL$.
\end{x}

\vspace{0.1cm}

\begin{x}{\small\bf FACT} \ 
$S(\LL/\K)$ is a field.
\end{x}

\vspace{0.1cm}

\begin{x}{\small\bf FACT} \ 
If $\LL \supset \E \supset \K$ and $\E$ is a separable extension of $\K$, then $\E \subset S(\LL/\K)$.
\end{x}

\vspace{0.1cm}

\begin{x}{\small\bf NOTATION} \ 
$\K^\cl$ is the algebraic closure of $\K$.
\end{x}

\vspace{0.1cm}

\begin{x}{\small\bf \un{N.B.}} \ 
If $\K$ is not perfect, then $\K^\cl$ is not Galois over $\K$.
\end{x}

\vspace{0.1cm}

\begin{x}{\small\bf NOTATION}  \ 
$\K^\sep$ is the separable closure of $\K$ in $\K^\cl$:
\[
\K^\sep \hsx = \hsx S(\K^\cl/\K).
\]
\end{x}

\vspace{0.1cm}

\begin{x}{\small\bf FACT} \ 
$\K^\sep$ is the maximal separable extension of $\K$.
\end{x}

\begin{x}{\small\bf FACT} \ 
$\K^\sep$ is a Galois extension of $\K$.
\end{x}

\vspace{0.1cm}

\begin{x}{\small\bf DEFINITION} \ 

\[
\Gal(\K^\sep/\K)
\]
is the 
\un{absolute Galois group}
\index{absolute Galois group} 
of $\K$.
\end{x}

\vspace{0.1cm}

\begin{x}{\small\bf FACT} \ 
If \hsx $\LL/\K$ is Galois, then $\Gal(\LL/\K)$ is a homomorphic image of $\Gal(\K^\sep/\K)$.

\vspace{0.1cm}

[This is because $\Gal(\LL/\K)$ can be identified with the quotient
\[
\Gal(\K^\sep/\K) / \Gal(\K^\sep/\LL).]
\]
\end{x}

\vspace{0.1cm}

\begin{x}{\small\bf EXAMPLE} \ 
Take $\K = \F_p$ $-$then $\Gal(\F_p^\sep /\F_p)$ can be identified with \ 
$\lim\limits_{\substack{\lla\\ n \in \N}} \Z/n\Z$ 
(the set of all (equivalence classes) of sequences $\{a_n\} = \{a_1, a_2, \ldots \}$ 
of natural numbers such that 
\[
a_n \hsx \equiv \hsx a_m \hsx \text{(mod $m$)}
\]
whenever $m|n$).

\vspace{0.1cm}

[Bear in mind that $\forall \ n \in \N$, there is a Galois extension $\K_n/\F_p$ with $[\K_n:\F_p] = n$ and 
$\Gal(\K_n/\F_p) \approx \Z/n\Z$.]

\vspace{0.1cm}

[Note: \ Let $\phi:\F_p^\sep \ra \F_p^\sep$ be the Frobenius automorphism: $\phi(x) = x^p$.  
Let $G = \langle \phi \rangle$ $-$then 
\[
\Inv(G)\hsx = \hsx \F_p, \quad  \Inv(\Gal(\F_p^\sep /\F_p)) \hsx = \hsx \F_p,
\]
yet 
\[
G \hsx \neq \hsx \Gal(\F_p^\sep /\F_p).]
\]
\end{x}

\vspace{0.1cm}


\begin{x}{\small\bf NOTATION} \ 
$\Gal^*(\K^\sep/\K)$ 
\index{$\Gal^*(\K^\sep/\K)$ } 
is the commutator subgroup of $\Gal(\K^\sep/\K)$ .
\end{x}

\vspace{0.1cm}

\begin{x}{\small\bf FACT} \ 
\[
\Inv(\Gal^*(\K^\sep/\K)) \hsx = \hsx \Inv(\ov{\Gal^*(\K^\sep/\K})).
\]

\vspace{0.1cm}

[Put
\[
\Gamma \hsx = \hsx \Gal^*(\K^\sep/\K).
\]
Then
\[
\Gamma \hsx\subset\hsx \ov{\Gamma} \implies \Inv(\ov{\Gamma}) \subset \Inv(\Gamma).
\]
To go the other way, let $x \in \Inv(\Gamma)$, $\ov{\gamma} \in \ov{\Gamma}$ and claim: 
$\ov{\gamma} x = x$ (hence $x \in \Inv(\ov{\Gamma})$).  
If $\ov{\gamma} \in \Gamma$, we are through; otherwise, 
$\ov{\gamma}$ is an accumulation point of $\Gamma$, thus since 
$\Omega_{\ov{\gamma}} (\{x\})$ is a neighborhood of $\ov{\gamma}$, it must contain a 
$\gamma \in \Gamma$ $(\gamma \neq \ov{\gamma}$).  
But
\[
\gamma \in \Gamma \cap \Omega_{\ov{\gamma}} (\{x\}) \implies 
\gamma \in  \Omega_{\ov{\gamma}} (\{x\}) \implies 
\gamma x \hsx = \hsx \ov{\gamma} x.
\]
Meanwhile, 
\[
\gamma \in \Gamma  \ \& \ x \in \Inv(\Gamma) \implies \gamma x = x.
\]
Therefore $\ov{\gamma} x = x$.]
\end{x}

\vspace{0.1cm}

\begin{x}{\small\bf \un{N.B.}} \ 
\[
\ov{\Gal^*(\K^\sep/\K})
\]
is a closed normal subgroup of $\Gal(\K^\sep/\K)$.
\end{x}

\vspace{0.1cm}

\begin{x}{\small\bf DEFINITION} \ 

\[
\Inv(\Gal^*(\K^\sep/\K))
\]
is called the 
\un{maximal abelian extension}
\index{maximal abelian extension} 
of $\K$, denote it by $\K^\ab$.
\index{$\K^\ab$}
\end{x}

\vspace{0.1cm}

\begin{x}{\small\bf FACT} \ 
$\K^\ab$ is a Galois extension of $\K$ and $\Gal(\K^\ab/\K)$ is an abelian group.

\vspace{0.1cm}

[Since
\[
\ov{\Gal^*(\K^\sep/\K)}
\]
is a closed normal subgroup of $\Gal(\K^\sep/\K)$, it follows that
\begin{align*}
\K^\ab \ 
&=\ \Inv ( \Gal^*(\K^\sep/\K) )\\
&=\ \Inv(\ov{\Gal^*(\K^\sep/\K)})
\end{align*}
is a Galois extension of  $\K$ and
\begin{align*}
\Gal(\K^\ab/\K) \ 
&\approx\ \Gal(\K^\sep/\K) / \Gal(\K^\sep/\K^\ab)\\
&=\ \Gal(\K^\sep/\K) / \ov{\Gal^*(\K^\sep/\K)}
\end{align*}
But the group on the RHS is isomorphic to 
\[
\Gal(\K^\sep/\K) / \Gal^*(\K^\sep/\K) / \ov{\Gal^*(\K^\sep/\K)} / \Gal^*(\K^\sep/\K),
\]
thus is a homomorphic image of the abelian group
\[
\Gal(\K^\sep/\K) / \Gal^*(\K^\sep/\K).]
\]
\end{x}

\vspace{0.1cm}

\begin{x}{\small\bf DEFINITION} \ 
A Galois extesnsion $\LL/\K$ is said to be 
\un{abelian}
\index{abelian (field extension)}
if $\Gal(\LL/\K)$ is abelian.
\end{x}

\vspace{0.1cm}

\begin{x}{\small\bf FACT} \ 
The field $\K^\ab$ has no extensions that are abelian Galois extensions of $\K$.

\vspace{0.1cm}

[Let $\LL/\K^\ab$ be an abelian Galois extensions of $\K$:
\[
\LL \hsx = \hsx \Inv(\Gal(\K^\sep/\K)) \supset \K^\ab \hsx = \hsx \Inv(\Gal^*(\K^\sep/\K))
\]
\qquad\qquad $\implies$
\[
\Gal^*(\K^\sep/\K) \hsx\supset\hsx \Gal(\K^\sep/\LL).
\]
On the other hand, $\Gal(\K^\sep/\LL)$ is normal ($\LL/\K$ being Galois) and 
\[
\Gal(\K^\sep/\K) /  \Gal(\K^\sep/\LL) \hsx \approx \hsx \Gal(\LL/\K),
\]
which is abelian by hypothesis, thus
\[
\Gal(\K^\sep/\LL) \supset \Gal^*(\K^\sep/\K).
\]
Therefore
\[
\Gal(\K^\sep/\LL) \ = \  \Gal^*(\K^\sep/\K).
\]
And then
\begin{align*}
\LL \ 
&\vsy=\ \Inv(\Gal(\K^\sep/\LL))\\
&\vsy=\ \Inv(\Gal^*(\K^\sep/\K))\\
&\vsy=\ \K^\ab.]
\end{align*}
\end{x}

\vspace{0.1cm}

\begin{x}{\small\bf FACT} \ 
$\K^\ab$ is generated by the set of finite abelian Galois extensions of $\K$ in $\K^\sep$.

\vspace{0.1cm}

[Every finite Galois extension of $\K$ inside $\K^\ab$ is necessarily abelian.]
\end{x}

\vspace{0.1cm}

\begin{x}{\small\bf DEFINITION} \ 
Take $\K = \Q$ $-$then the splitting field $\Q(n)$ of the polynomial $X^n - 1$ is called the 
\un{cyclotomic field}
\index{cyclotomic field} 
of the $n^\nth$ roots of unity.
\end{x}

\vspace{0.1cm}

\begin{x}{\small\bf FACT} \ 
$\Q(n)$ is a Galois extension of $\Q$ and $\Gal(\Q(n)/\Q)$ is isomorphic to 
$(\Z/n\Z)^\times$, hence $\Gal(\Q(n)/\Q)$ is abelian.
\end{x}
\vspace{0.1cm}

Accordingly, every intermediate field $\E$ between $\Q$ and $\Q(n)$ is abelian Galois (per $\Q$).

\vspace{0.1cm}

[$\Gal(\Q(n) / \Q)$ is abelian, hence every subgroup of $\Gal(\Q(n) / \Q)$ is normal, hence in particular 
$\Gal(\Q(n) / \E)$ is normal, hence $\E/\Q$ is Galois.  And 
\[
\Gal(\E) / \Q) \hsx \approx \hsx \Gal(\Q(n) / \Q) / \Gal(\Q(n) / \E).]
\]

\vspace{0.2cm}

The Kronecker-Weber theorem states that every finite abelian Galois extension of $\Q$ is contained in some $\Q(n)$, 
thus $\Q^\ab$ is the infinite cyclotomic extension $\Q(1, 2, \ldots)$.

\vspace{0.2cm}

\begin{x}{\small\bf SCHOLIUM} \ 
$\Q^\ab$ is generated by the torsion points of the action of $\Z$ on $\C^\times$.

\vspace{0.1cm}

[Note: \ Given $n \in \Z$, $x \in \C^\times$, $(n,x) \ra n \cdot x = x^n$.]
\end{x}

\vspace{0.3cm}

\[
\text{ADDENDUM}
\]
\setcounter{theoremn}{0}

\vspace{0.1cm}

If \mG is a group, then the subgroup $G^*$ generated by the commutators $x y x^{-1}y^{-1}$ is the 
\un{commutator subgroup}
\index{commutator subgroup} 
of \mG.

\vspace{0.2cm}

\qquad\qquad \textbullet \quad $G^*$ is a normal subgroup of \mG.

\vspace{0.1cm}

\qquad\qquad \textbullet \quad $G / G^*$ is abelian.

\vspace{0.1cm}

And if $H \subset G$ is normal and if $G/H$ is abelian, then $H \supset G^*$.

\vspace{0.1cm}


\vspace{0.3cm}

{\small\bf FACT}
If  $\LL/\K$ is an infinite Galois extension and if $N \subset \Gal(\LL/\K)$ is a normal subgroup, then 
$\ov{N} \subset \Gal(\LL/\K)$ is a closed normal subgroup. 

\newpage
\setcounter{page}{1}
\renewcommand{\thepage}{References-\arabic{page}}
\setcounter{page}{1}
\begingroup
\center {\textbf{REFERENCES}}\\
\endgroup
\vspace{0.5cm}

\noindent Arakawa, T. et al.

[1] \quad 
Bernoulli Numbers and Zeta Functions, Springer Verlag, 2004.

\noindent Cassels, J. and Fr\"olich, A 

[2] \quad
Algebraic Number Theory, Academic Press, 1967.

\noindent Edwards, H.

[3] \quad
Galois Theory, Springer Verlag, 1984.

\noindent Gouvea, F.

[4] \quad
$p$-adic Numbers, Springer Verlag, 1991.

\noindent Howe, R. and Tan, E.

[5] \quad
Non-Abelian Harmonic Analysis, Springer Verlag, 1992.

\noindent Iwasawa, K.

[6] \quad
Lectures on $p$-adic $L$-functions, Princeton University Press, 1972.

\noindent Koblitz, N.

[7-(a)] \quad
$p$-adic Analysis: A Short Course on Recent Work, Cambridge University Press, 1980.

[7-(b)] \quad
$p$-adic Numbers, $p$-adic Analysis, and Zeta-Functions, Springer Verlag, 1984.

\noindent K\"orner, T.

[8] \quad
Fourier Analysis, Cambridge University Press, 1988.

\noindent Morandi, P.

[9] \quad
Field and Galois Theory, Springer Verlag, 1996.

\noindent Patterson, S.

[10] \quad
An Introduction to the Theory of the Riemann Zeta-Function, Cambridge University Press, 1988.

\noindent Srivastava, H. and Junesang, C.

[11] \quad
Zeta and $q$-Zeta Functions and Associated Series and Integrals, Elsevier, 2012.

\noindent Weil, A.

[12] \quad
Basic Number Theory, Springer Verlag, 1967.

\setcounter{page}{1}
\renewcommand{\thepage}{Index-\arabic{page}}
\backmatter
\bibliography{}
\printindex
\end{document}